\documentclass[a4paper,10pt,twoside]{report}
\newcommand{\be}{\begin{enumerate}}
\newcommand{\ee}{\end{enumerate}}
\newcommand{\beqn}{\begin{eqnarray*}}
\newcommand{\eeqn}{\end{eqnarray*}}
\newcommand{\disp}{\displaystyle}
\newcommand{\flechen}[3]{%
\xymatrix{
{#1} \ar[r]^-{#2} & {#3} }}
\newcommand{\incl}[1][r]
      {\ar@<-0.2pc>@{^(-}[#1] \ar@<+0.2pc>@{-}[#1]}
\def\C{{\mathbb C}}
\def\D{{\mathbb D}}

\def\N{{\mathbb N}}
\def\P{{\mathbb P}}
\def\Q{{\mathbb Q}}

\def\Z{{\mathbb Z}}
\def\Ac{{\mathcal A}}
\def\Bc{{\mathcal B}}
\def\Cc{{\mathcal C}}
\def\Fc{{\mathcal F}}
\def\Gc{{\mathcal G}}
\def\Hc{{\mathcal H}}
\def\Ic{{\mathcal I}}
\def\Lc{{\mathcal{L}}}
\def\Mc{{\mathcal M}}
\def\Nc{{\mathcal N}}
\def\Oc{{\mathcal O}}
\def\Pc{{\mathcal P}}
\def\Qc{{\mathcal Q}}
\def\Rc{{\mathcal R}}
\def\Sc{{\mathcal S}}
\def\Tc{{\mathcal T}}
\def\Uc{{\mathcal U}}
\def\Vc{{\mathcal{V}}}
\def\Zc{{\mathcal{Z}}}

\def\mgo{{\mathfrak m}}
\def\pgo{{\mathfrak p}}
\def\qgo{{\mathfrak q}}

\def\gm{{\mathbb G}_{\rm{m}}}
\def\mun{{\mu}_n}
\def\bgm{{\rm B}{\mathbb G}_{\rm{m}}}
\def\bmun{{\rm B}\mun}
\def\bmunfppf{({\rm B}\mun)_{{\rm fppf}}}
\def\eps{\varepsilon}
\def\fleche{\longrightarrow}
\def\inj{\hookrightarrow}

\def\ov{\overline}
\newcommand{\cartesien}{\ar@{}[dr]|{\square}}
\def\FlecheNE{\rotatebox[origin=c]{45}{\ensuremath \Rightarrow}}

\def\lind{\lim_{\longrightarrow}}
\def\lpro{\lim_{\longleftarrow}}

\def\pr{{\rm pr}}
\def\Ker{\hbox{\rm Ker}\,}

\def\Mod{\hbox{\rm Mod}}
\def\Ab{\hbox{\rm Ab}}

\def\Aut{\hbox{\rm Aut}}

\def\Auts{\hbox{\rm \underline{Aut}}}

\def\Isom{\hbox{\rm Isom}}
\def\Isoms{\hbox{\rm \underline{Isom}}}
\def\fIsom{{\mathcal I\!som}\,}
\def\fAut{{\mathcal Aut}\,}

\def\ob{\hbox{\rm ob}\,}
\def\cosq{\hbox{\rm cosq}}
\def\sq{\hbox{\rm sq}}
\def\O{\hbox{\rm O}}

\def\Im{\hbox{\rm Im}\,}
\def\Id{{\rm Id}}
\def\id{\hbox{\rm id}}

\def\Der{\hbox{\rm Der}\,}

\def\Hom{\hbox{\rm Hom}}
\def\Ext{\hbox{\rm Ext}\,}
\def\fHom{{\mathcal{H}om}\,}

\def\Spec{{\rm Spec}\,}

\newcommand{\eme}[1]{${#1}^{\textrm{ème}}$}
\newcommand{\iem}[1]{${#1}^{\textrm{ième}}$}
\newcommand{\dd}{\ensuremath{\text{d}}}

\newcommand{\xlist}[2]{\ensuremath{\left({#1}_1,\ldots,{#1}_{#2}\right)}}

\newcommand{\nlist}[1]{\xlist{#1}{n}}

	{\begin{pmatrix}}%
	{\end{pmatrix}}
	{\begin{vmatrix}}%
	{\end{vmatrix}}
	{\ensuremath{\left \{ \hskip -1.5 mm \begin{array}{#1@{\ =\ }l}}}%
	{\end{array}\right.}
	{\ensuremath{\left \{ \hskip -1.5 mm%
	 \begin{array}{#1@{\ }c@{\ }l}}}%
	{\end{array}\right.}
	{\ensuremath{\left \{ \hskip -1.5 mm \begin{array}{#1@{\quad}l}}}%
	{\end{array}\right.}
\newcommand{\fonction}[5]{%
        \ensuremath{#1\colon
        \left\{\hskip -1.5 mm                   %}
        \begin{array}{c@{\ }c@{\ }l}
        \medskip #2 & \longrightarrow & #3 \\
        #4 & \longmapsto & #5 \\
        \end{array}
        \right .
        }}
\usepackage[T1]{fontenc}
\usepackage[latin1]{inputenc}
\usepackage[french]{babel}
\usepackage{verbatim}
\usepackage{latexsym}
\usepackage{amsfonts}
\usepackage{amstext}
\usepackage{amsmath}
\usepackage{graphicx}
\usepackage[all,dvips]{xy}
\input xy
\xyoption{all}
\usepackage{amssymb,euscript,multicol,graphicx}
\usepackage[ec]{aeguill}       %Ca c'est pour que les pdf soient beaux....
\usepackage{delarray, fancyhdr}
\newtheorem{cor}[subsection]{Corollaire}
\newtheorem{souscor}[subsubsection]{Corollaire}

\newtheorem{prop}[subsection]{Proposition}
\newtheorem{sousprop}[subsubsection]{Proposition}

\newtheorem{thm}[subsection]{Th\'eor\`eme}
\newtheorem{sousthm}[subsubsection]{Th\'eor\`eme}

\newtheorem{defi}[subsection]{D\'efinition}
\newtheorem{sousdefi}[subsubsection]{D\'efinition}

\newtheorem{lem}[subsection]{Lemme}
\newtheorem{souslem}[subsubsection]{Lemme}

\newtheorem{pptes}[subsection]{Propri\'et\'es}

\newtheorem{remarque}[subsection]{Remarque}
\newtheorem{sousremarque}[subsubsection]{Remarque}

\newtheorem{exemple}[subsection]{Exemple}
\newtheorem{sousexemple}[subsubsection]{Exemple}

\newenvironment{demo}{{\bf D\'emonstration.}}{$\square$ \vskip .3cm}
\newenvironment{democ}[1]{{\bf D\'emonstration #1}}{$\square$ \vskip .3cm}
\newenvironment{etape}[1]{$\bullet$ \emph{#1}\\}{\vskip .2cm}
\newenvironment{etapefinale}[1]{$\bullet$ \emph{#1}\\}{}
\author{Sylvain Brochard}
\SetMathAlphabet\mathcal{normal}{U}{rsfs}{m}{n} % pour faire des X champêtres...

\renewcommand{\chaptermark}[1]{\markboth{#1}{}}

\fancypagestyle{plain}{ % pages de tetes de chapitre 
\fancyhead{} % supprime l'entete 
}
\title{Champs algébriques et foncteur de Picard}
\date{\empty}

\def\Pic{\hbox{\rm Pic}}
\def\piczar{\Pic_{\X/S \, \rm{(Zar)}}}
\def\picet{\Pic_{\X/S \, \text{(\'Et)}}}
\def\piceth{\Pic_{\X/S \, \emph{(\'Et)}}}
\def\picfppf{\Pic_{\X/S \, \rm{(fppf)}}}
\def\pic{\Pic_{\X/S}}

\def\piczero{\Pic_{\X/S}^0}
\def\champic{\mathcal{P}ic}
\def\liset{\text{{\rm Lis-\'et}}}
\def\lisets{\text{\underline{Lis-\'et}}}
\def\lisetsth{\emph{\underline{Lis-\'et}}}
\def\gll{\text{{\rm Llc}}}
\def\et{\text{\rm {\'Et}}}
\def\fppf{\text{\rm fppf}}
\def\fppfc{\text{\rm fppfc}}
\def\pl{\text{\rm pl}}
\def\O{\mathcal{O}}
\def\A{\mathbb{A}}
\def\X{\mathcal{X}}
\def\Y{\mathcal{Y}}
\def\U{\mathcal{U}}
\def\L{\mathcal{L}}
\def\inv{\hbox{\rm Inv}\,}
\def\aff{\hbox{{\rm (Aff/}S)}}
\def\Xt{\widetilde{\X}}
\def\Yt{\widetilde{\Y}}

\def\Lt{\widetilde{\L}}
\def\Mt{\widetilde{\Mc}}
\def\defm{\hbox{\rm Defm}}
\def\at{\widetilde{\alpha}}
\def\Xzt{\widetilde{X^0}}
\def\Xut{\widetilde{X^1}}
\def\Xdt{\widetilde{X^2}}

\renewcommand{\fleche}{%
\xymatrix@C=1pc{\ar[r] &}}
\def\flechelongue{%
\xymatrix{\ar[r] &}}

\entrymodifiers={!!<0pt,0.7ex>+}
\begin{document}
\thispagestyle{empty}

\begin{center}
\null\vskip 3cm
\LARGE{Champs algébriques et foncteur de Picard}
\vskip 1cm
\large{Sylvain Brochard}
\vskip 2cm
{\bf Comment}\vskip 5mm
\end{center}
{\small This text is my thesis, defended in June 2007, in the status it was at this time. The most important results are contained in the article \emph{Foncteur de Picard d'un champ algébrique} to appear in \emph{Mathematische Annalen} (see the preprint arXiv:0711.4545). In the article, some results have been added, and some previous results have been strengthened. However, the proofs of the results contained in the appendix (concerning the smooth-\'etale cohomology on an algebraic stack) have been removed. The thesis is only put on the ArXiv to provide a more lasting reference than my webpage for these proofs.
}
\vfill
\pagebreak
\tableofcontents \vfill
\chaptermark{Introduction}
\section*{Introduction}
\addcontentsline{toc}{chapter}{Introduction}

Si $X$ est un objet géométrique, son groupe de Picard est le groupe des classes d'isomorphie de faisceaux inversibles (ou, si l'on préfère, de fibrés en droites) sur $X$. Par exemple si $X$ est le spectre d'un anneau de Dedekind $A$, le groupe $\Pic(X)$ est simplement le groupe des classes d'idéaux de $A$. Si $X$ est l'espace projectif $\P^n$ sur un corps, il est connu que le degré des faisceaux inversibles induit un isomorphisme de $\Pic(X)$ vers $\Z$. Dans bien des cas, le groupe de Picard est naturellement muni d'une structure \emph{géométrique}. Par exemple si $X$ est un schéma projectif intègre sur un corps $k$ algébriquement clos, on sait que le groupe de Picard $\Pic(X)$ apparaît en fait comme le groupe sous-jacent à un $k$-schéma en groupes $\Pic_{X/k}$ appelé le \emph{schéma de Picard}  (ou foncteur de Picard) de $X$.

L'objet de cette thèse est de définir et d'étudier le foncteur de Picard d'un \emph{champ algébrique} -- un objet géométrique, apparu vers 1970, qui généralise la notion de schéma. Pour comprendre la nécessité d'étudier un tel objet et pour apprécier à leur juste valeur les travaux de Grothendieck, il n'est sans doute pas inutile de revoir rapidement l'histoire du foncteur de Picard et des champs algébriques. Pour plus de détails sur les aspects historiques évoqués ci-dessous ou pour des références plus précises, nous renvoyons à l'excellente introduction de~\cite{poly_Kleiman} dont nous nous sommes inspiré.

\subsection*{Foncteur de Picard}

Les besoins technologiques du \eme{17} siècle ont poussé les mathématiciens à développer le calcul intégral. Ils ont alors commencé à s'intéresser aux propriétés de certaines fonctions apparaissant comme des intégrales indéfinies. \`A la toute fin de ce siècle, Jacques Bernoulli puis son frère Jean mirent ainsi en évidence des relations algébriques surprenantes entre les arguments des sommes et des différences de logarithmes, de fonctions trigonométriques inverses, ou de certaines fonctions intégrales. Dans la même veine, le travail de Fagnano au début du \eme{18} siècle amena Euler à découvrir en 1757 la \og formule d'addition\fg\ suivante
$$\int_0^{x_1} \frac{\dd x}{\sqrt{1-x^4}}
\pm \int_0^{x_2} \frac{\dd x}{\sqrt{1-x^4}}
= \int_0^{x_3} \frac{\dd x}{\sqrt{1-x^4}}
$$
où $x_1$, $x_2$ et $x_3$ vérifient la relation algébrique :
$$x_1^4x_2^4x_3^4+2x_1^2x_2^2x_3^2(x_1^2+x_2^2+x_3^2)+x_1^4+x_2^4+x_3^4
-2x_1^2x_2^2-2x_1^2x_3^2-2x_2^2x_3^2=0.$$
Il généralisa en 1759 cette formule à d'autres intégrales. Abel donna en 1826 une bien plus large portée à ces considérations en découvrant le fameux \og théorème d'Abel\fg. C'est un théorème d'addition très général pour des intégrales algébriques de ce type, intégrales qui ont vite pris le nom d'\og intégrales abéliennes\fg. Euler comme Abel remarquèrent que ces intégrales étaient en quelque sorte \og multivaluées \fg, un peu comme les fonctions trigonométriques inverses, mais ne poussèrent pas l'analogie plus loin. Abel voyait ce phénomène comme une sorte de \og constante d'intégration\fg\ dont il évitait l'apparition en ne considérant que des chemins d'intégration assez petits.

Riemann éclaircit nettement la situation en 1857. En étudiant les intégrales sur une variété projective lisse complexe $C$ de dimension 1, que l'on appellerait aujourd'hui une \og surface de Riemann\fg, il prouve que certaines de ces intégrales, les intégrales \og de première espèce \fg, forment un $\C$-espace vectoriel de dimension finie $p$ et qu'elles ont $2p$ \og périodes\fg, c'est-à-dire des nombres qui engendrent tous les changements possibles de la valeur de l'intégrale qui peuvent survenir en changeant le chemin d'intégration. Plus précisément, il fixe une base $(\psi_1, \dots, \psi_p)$ de l'espace des intégrales de première espèce. Il choisit aussi $2p$ chemins qui forment une base de l'homologie, puis il forme le réseau $\Lambda$ de $\C^p$ engendré par les $2p$ vecteurs obtenus en intégrant $(\psi_1, \dots, \psi_p)$ le long des $2p$ chemins choisis.
\begin{center}
\includegraphics[scale=0.3]{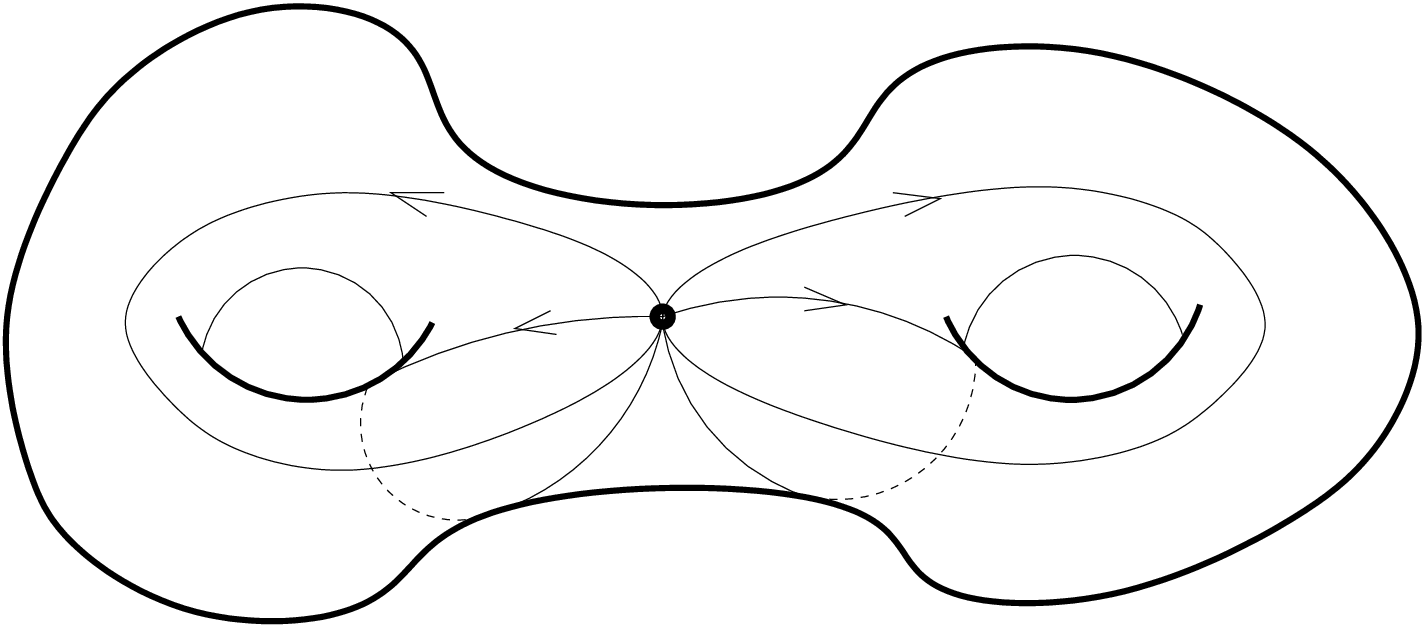}\\
\medskip
\small{Figure 1 : Une surface de Riemann de genre 2 a quatre \og périodes\fg.}
\end{center}
Il appelle alors \og jacobienne\fg\ le quotient 
$$J=\C^p/\Lambda.$$
Il s'intéresse aussi au morphisme dit \og d'Abel-Jacobi\fg\ de la $\mu$-ième puissance symétrique de $C$ vers la jacobienne :
$$\disp \fonction{\Psi_{\mu}}{C^{(\mu)}}{J}{(x_1, \dots, x_{\mu})}{\disp \left(\sum_{i=1}^{\mu} \psi_1 x_i, \dots, \sum_{i=1}^{\mu} \psi_p x_i\right).}$$
Il montre que si les diviseurs $x_1+\dots+x_{\mu}$ et $x'_1+\dots+x'_{\mu}$ sont linéairement équivalents, alors $\psi_{\mu}(x_1,\dots,x_{\mu})$ et $\psi_{\mu}(x'_1,\dots,x'_{\mu})$ sont égaux. La réciproque, vraie également, ne sera prouvée que quelques années plus tard par un autre mathématicien (Clebsch, 1864). Ainsi la jacobienne, dont on commence à voir qu'elle s'identifie au groupe des classes de diviseurs, acquiert avec Riemann une structure de tore complexe de dimension $p$. Elle n'est cependant pour l'instant définie que sur le corps $\C$ des nombres complexes.

Dans les années 1950, Weil fonde la théorie des variétés abéliennes. Il donne en particulier une construction purement algébrique de la jacobienne et prouve ainsi son existence sur un corps quelconque. L'inconvénient de cette démarche est que la structure géométrique n'est pas définie \emph{a priori} mais résulte de la construction algébrique. Il faut par exemple revenir à la construction pour voir que tel objet lié au groupe de Picard est bien un objet algébrique.

Il a fallu attendre Grothendieck pour voir apparaître une définition générale du \og foncteur de Picard\fg\ d'un schéma sur une base quelconque, cette notion se comportant de plus extrêmement bien par changement de base. Il est défini de la manière suivante. Si $X$ est un $S$-schéma, on définit d'abord le foncteur $P$ de la catégorie des $S$-schémas vers celle des groupes qui à un $S$-schéma $T$ associe $P(T)=\Pic(X\times_S T)$. On se demande alors s'il existe un schéma qui représente ce foncteur, c'est-à-dire un $S$-schéma $\underline{\Pic}$ tel que pour tout $T$, $P(T)$ s'identifie à l'ensemble des morphismes de $T$ vers $\underline{\Pic}$. Naturellement, ainsi posé, le problème n'a pas de solution. On sait en effet qu'un foncteur représentable est toujours un faisceau pour la topologie \emph{fpqc}. Or le foncteur $P$ que nous venons de définir n'est jamais un faisceau, même pour la topologie de Zariski : il existe toujours des faisceaux non triviaux qui proviennent de la base. Qu'importe ! Il suffit de considérer, non pas le foncteur $P$ lui-même, mais le faisceau associé pour une topologie suffisamment fine pour que l'on puisse espérer que le faisceau obtenu ait une chance d'être représentable. En pratique on se contente de la topologie étale. On appelle donc foncteur de Picard de $X$, et l'on note $\Pic_{X/S}$, le faisceau associé à $P$ pour la topologie étale. Le schéma de Picard, s'il existe, est alors simplement l'unique schéma qui représente ce foncteur. Avec cette définition, la structure géométrique est intimement liée à la structure algébrique et est en tout cas complètement définie \emph{a priori}. Ce point de vue présente un avantage non négligeable : la construction du foncteur de Picard ne requiert aucune hypothèse sur le schéma $X$. En travaillant directement sur le foncteur de Picard, on peut même étudier certaines propriétés géométriques (lissité, séparation, propreté, dimension\dots) du schéma de Picard sans même savoir s'il existe!! Quant à l'existence du schéma de Picard (c'est-à-dire la représentabilité du foncteur de Picard par un schéma), elle est établie dans de nombreux cas à l'aide de techniques projectives reposant sur une utilisation judicieuse du schéma de Hilbert.

Les travaux ultérieurs s'inscrivent parfaitement dans le prolongement de la démarche \og fonctorielle\fg\ initiée par Grothendieck. Artin invente et développe avec son élève Knutson (voir \cite{Knutson}) la notion d'espace algébrique, qui généralise celle de schéma. Il montre alors avec Raynaud l'existence d'un \og espace algébrique de Picard\fg. Cette fois-ci, la démonstration ne fait aucun usage des techniques projectives encore indispensables quelques années auparavant. Il en résulte un théorème valable dans un cadre bien plus général que celui qui assure l'existence du schéma de Picard, et de surcroît plus facile à obtenir, sans que le produit obtenu (un espace algébrique au lieu d'un schéma) ne soit réellement moins commode.

\subsection*{Champs algébriques}

La thèse de Giraud, vers 1970, donne naissance aux champs, dans la droite lignée des
travaux de Grothendieck. La notion de champ algébrique suit immédiatement avec un article
fondateur de Deligne et Mumford (\cite{Deligne_Mumford}) : ils démontrent, en introduisant les champs qui
portent aujourd'hui leurs noms, que l'espace de module des courbes de genre $g$ est
irréductible. Artin donne peu après une définition un peu plus générale des champs algébriques.

Initialement introduits pour pallier la non-existence de certains espaces de modules ou bien pour effacer
quelques importunes singularités d'iceux, qui n'apparaissaient que par une sorte d'illusion d'optique
due à la mise à mort des automorphismes des objets qu'ils sont censés classifier, les champs algébriques
ont petit à petit pris une place considérable dans l'environnement naturel du géomètre algébriste.
Au point que certains pensent, comme Abramovich et Vistoli (cf.~\cite{Abramovich_Vistoli_note}~1.2), qu'ils sont amenés
à devenir l'objet d'étude de base du géomètre algébriste, au même titre que les schémas l'étaient pour les
\og anciens\fg. De fait, des travaux récents les ont fait apparaître comme tels. Prenons l'exemple des surfaces elliptiques, c'est-à-dire des familles de courbes elliptiques paramétrées par une courbe $C$. Pour compactifier l'espace de modules de ces surfaces, il faut rajouter quelques
points. Or ces points supplémentaires correspondent naturellement à des familles paramétrées par des \og courbes tordues\fg\ introduites
par Abramovich et Vistoli. Ce sont des courbes nodales munies d'une structure supplémentaire d'\og orbispace\fg\ au-dessus de chaque n\oe ud. Elles sont ainsi munies d'une structure de champ de Deligne-Mumford. Nous reparlerons de ces courbes ultérieurement.

\subsection*{Des jacobiennes de courbes champêtres ?}

Toute la force des champs algébriques vient du fait que la plupart des notions géométriques classiques
(lissité, dimension, cohomologie, connexité\dots) s'y étendent sans difficulté.
Pour ce qui concerne le groupe de Picard, le travail a été initié par Mumford dans un article fondateur. Bien que la notion de champ n'ait pas encore été introduite, il y calcule en réalité le groupe de Picard du champ des courbes elliptiques, isomorphe à $\Z/12\Z$. Il semble cependant que personne ne se soit encore réellement attelé à la construction de la jacobienne d'un champ algébrique et à l'étude de ses propriétés géométriques. En l'état actuel des publications, la seule chose que l'on puisse affirmer est que la question de la représentabilité du foncteur de Picard est pour l'essentiel déjà résolue. Masao Aoki montre en effet dans des travaux récents (voir \cite{Aoki_Hom} et \cite{Aoki_erratum}) que si $\X$ et $\Y$ sont des $S$-champs algébriques de présentation finie, et si $\X$ est propre et plat, alors le champ des morphismes $\fHom(\X,\Y)$, dont on trouvera une définition précise dans les articles cités, est un champ (pré)algébrique\footnote{Seule la quasi-séparation n'est pas traitée par Aoki. Nous y reviendrons dans un paragraphe ultérieur.}. Or en prenant pour $\Y$ le champ $\bgm$ qui classifie les $\gm$-torseurs, on obtient le champ $\fHom(\X, \bgm)$, isomorphe au champ des faisceaux inversibles $\champic(\X/S)$ dont la catégorie fibre au-dessus d'un $S$-schéma $U$ est la catégorie des faisceaux inversibles sur $\X\times_S U$. On en déduit facilement que le foncteur de Picard d'un champ algébrique propre, plat et cohomologiquement plat en dimension zéro est représentable par un espace algébrique (cf. infra~(\ref{comparaison_champ_foncteur})). Max Lieblich vient par ailleurs de souligner (cf. \cite{Lieblich_coherent_algebras}) le fait que le champ des faisceaux quasi-cohérents à support propre sur un champ algébrique (avec de bonnes hypothèses) est algébrique. Il en déduit ensuite que le champ des fibrés de rang fini, donc en particulier le champ $\champic(\X/S)$, est algébrique.

Cependant la démonstration d'Aoki ne permet pas réellement de voir ce qui se passe pour le foncteur de Picard. D'une part il n'y est question que de morphismes de champs algébriques, si bien qu'il faut toujours traduire les énoncés et les démonstrations dans le langage des faisceaux inversibles. Mais surtout, l'étude générale des morphismes oblige l'auteur à recourir à des techniques dont on sent bien qu'elles sont par trop élaborées pour le simple cas des faisceaux inversibles.  Max Lieblich utilise implicitement le même genre de techniques pour traiter le cas du champ des faisceaux quasi-cohérents à support propre.

% Enfin le lien entre le groupe de Picard et le groupe des classes de diviseurs n'est toujours pas explicité. Il semble hélas que la notion de groupe des classes de diviseurs ne fonctionne pas très bien sur les champs algébriques. Dans sa thèse \cite{Vistoli_ITAS}, Vistoli étudie les groupes de Chow sur un champ de Deligne-Mumford. Il définit donc en particulier un groupe des classes de diviseurs. Le problème est que ce groupe n'est défini qu'après tensorisation par $\Q$. Ceci provient essentiellement du fait que le degré d'un morphisme de champs algébriques peut être fractionnaire. Il faut donc peut-être renoncer à voir directement le groupe de Picard comme le groupe des classes de diviseurs. On peut néanmoins espérer qu'un tel isomorphisme existe après tensorisation par $\Q$. 

\subsection*{Plan de la thèse et principaux résultats}

\subsubsection*{Généralités}
Le premier chapitre est consacré à la définition et aux premières propriétés \emph{des} foncteurs de Picard d'un champ algébrique sur un schéma $S$. Après avoir revu rapidement ce qu'était un faisceau inversible, nous démontrons que le groupe des classes de faisceaux inversibles est isomorphe au groupe $H^1(\X, \gm)$, retrouvant ainsi le résultat analogue bien connu dans le cas des schémas. Puis nous définissons comme dans le cas des schémas différents foncteurs de Picard, qui sont des faisceaux relativement à des topologies de plus en plus fines sur $S$, et nous en donnons une description cohomologique. Le théorème le plus utile est le théorème de comparaison (\ref{comparaison_des_foncteurs_de_Picard}) entre ces différents foncteurs dans le cas où $f$ est cohomologiquement plat en dimension zéro et quasi-compact. Signalons tout de suite que pour mener à bien cette étude, nous avons besoin d'un certain nombre de résultats élémentaires relatifs à la cohomologie lisse-étale des champs algébriques, qui semblent manquer à la littérature actuelle. Nous avons préféré reléguer ces résultats en annexe afin de ne pas \og polluer\fg\ le texte principal par trop de considérations techniques. Nous reviendrons plus en détail sur cette annexe et son contenu un peu plus loin.

Nous définissons ensuite le champ de Picard. Nous montrons en particulier que sous des hypothèses raisonnables, le champ de Picard est algébrique si et seulement si le foncteur de Picard est représentable.

Cette partie s'achève sur des propriétés élémentaires du foncteur de Picard. Il est tout d'abord localement de présentation finie dès que le morphisme structural de $\X$ est quasi-compact. En outre il commute aux limites projectives. Cette dernière propriété, d'ailleurs facile à obtenir à l'aide du théorème d'existence de Grothendieck généralisé par Olsson aux champs algébriques, trouve sa raison d'être dans le théorème de représentabilité du foncteur de Picard (re)démontré ultérieurement.

\subsubsection*{Propriétés de séparation}
Cette seconde partie est déjà plus technique que la précédente. Nous y étudions les propriétés de séparation du foncteur de Picard et du champ de Picard.
On retiendra principalement deux résultats, et peut-être un lemme qui peut s'avérer utile pour résoudre la question de la quasi-séparation pour d'autres champs que le champ de Picard.

Le premier de ces résultats est la proposition~(\ref{locale_separation_du_foncteur_de_Picard}) : le foncteur de Picard est localement séparé, autrement dit sa diagonale est une immersion quasi-compacte, dès que $\X$ vérifie les hypothèses habituelles. La démonstration repose essentiellement sur les idées d'Artin présentes dans \cite{Global_Analysis_1} et qui apparaissent ici sous la forme du théorème~(\ref{critere_immersionsqc}).

 Pour un champ, les questions de séparation semblent nettement plus délicates en général. Il est déjà hors de question qu'un \og vrai\fg\ champ algébrique, j'entends par là un champ algébrique qui n'est pas un espace algébrique, soit localement séparé au sens habituel. En effet sa diagonale serait alors une immersion donc un monomorphisme, ce qui prouverait que le morphisme structural est représentable. La condition qu'il est raisonnable d'exiger serait donc celle de quasi-séparation : on dit qu'un champ algébrique est quasi-séparé si sa diagonale est quasi-compacte. Cette condition fait d'ailleurs partie de la définition d'un champ algébrique dans le livre \cite{LMB} de Laumon et Moret-Bailly. D'une manière ou d'une autre il est nécessaire d'imposer des conditions de finitude sur la diagonale. On lira à ce sujet la remarque~II~1.9 de Knutson dans \cite{Knutson}. Cependant un certain nombre de champs tout à fait dignes d'intérêt et que l'on aurait envie de qualifier d'\og algébriques\fg\ ne sont pas quasi-séparés. Doit-on considérer que le champ classifiant $B\Z$ est algébrique ? Pour d'autres champs la question n'est tout simplement pas résolue\footnote{Certains auteurs ont même enlevé la condition de quasi-séparation de la définition.}. Par exemple dans \cite{Aoki_Hom}, Masao Aoki ne traite qu'un cas très particulier : il montre que le champ des morphismes $\fHom(\X, \Y)$ est quasi-séparé lorsque $\X$ est un espace algébrique et $\Y$ admet une présentation propre (\cite{Aoki_Hom}~proposition~4.1). Nous montrons dans cette partie (voir le théorème~(\ref{quasi_separation_du_champ_de_Picard})) que dans le cas particulier du champ de Picard, le fait que les groupes d'automorphismes soient quasi-compacts suffit à assurer la quasi-séparation.

 Nous donnons au passage un critère de quasi-séparation pour les champs algébriques fortement inspiré du théorème d'Artin~(\ref{critere_immersionsqc}).
%et peut-être même un ou deux exemples d'application autres que le champ de Picard.

\subsubsection*{Déformation de faisceaux inversibles}
Nous étudions dans cette troisième partie les déformations de faisceaux inversibles, puis nous (re)démontrons que le foncteur de Picard d'un champ algébrique propre, plat et cohomologiquement plat en dimension zéro est représentable et nous donnons enfin des propriétés élémentaires sur sa lissité et sa dimension.

L'étude des déformations est la pierre angulaire des théorèmes de représentabilité d'Artin que l'on peut trouver dans \cite{Global_Analysis_1} ou \cite{Artin_Versal_defm}. C'est souvent le point dont l'étude est la plus délicate. C'est d'ailleurs pour cette raison que Masao Aoki a dû, pour montrer que le champ $\fHom(\X,\Y)$ était algébrique, consacrer un premier article (\cite{Aoki_defm}) à l'étude des déformations de champs algébriques, et s'appuyer sur un article d'Olsson (\cite{Olsson_defm}) concernant les déformations de morphismes \emph{représentables} de champs algébriques et sur les travaux très généraux d'Illusie (\cite{Illusie_CCD}).

Le principal résultat de cette section est le suivant :

\vskip 2mm \noindent
{{\bf Théorème~\ref{thm_defm_fi}} \it
\begin{itemize}
\item[(1)] Il existe un \'el\'ement $\omega\in H^2(\X,I)$
dont l'annulation \'equivaut \`a l'existence d'une d\'eformation de $\L$ \`a
$T$.
\item[(2)] Si $\omega=0$, alors $\ov{\defm_T(\L)}$ est un torseur sous
$H^1(\X,I)$.
\item[(3)] Si $(\Lt,\lambda)$ est une d\'eformation de $\L$, son groupe
d'automorphismes est isomorphe \`a $H^0(\X,I)$.
\end{itemize}
}\vskip 2mm \noindent
Nous l'y étudions sous plusieurs aspects. Nous commençons par montrer que l'on peut le déduire facilement du théorème analogue d'Aoki concernant les déformations de morphismes en prenant comme champ d'arrivée $\Y=\bgm$. Il suffit à peu de choses près de calculer le complexe cotangent du champ classifiant $\bgm$ sur un schéma $T$. Cependant cette démarche nous a paru insatisfaisante car il s'avère que le cas des faisceaux inversibles est nettement plus simple que celui des déformations de morphismes. Nous proposons donc deux autres démonstrations de ce théorème, indépendantes des travaux que nous venons de citer et de nature plus \og géométrique\fg. La première consiste à se ramener au cas des faisceaux inversibles sur un espace algébrique à l'aide d'une présentation du champ $\X$ (elle rejoint en cela les démarches d'Olsson et Aoki). La seconde est plus rapide puisque nous travaillons directement avec des faisceaux inversibles sur $\X$. Elle nécessite cependant un petit travail technique supplémentaire pour relier les groupes de cohomologie lisse-étale d'un champ à ceux d'une extension infinitésimale de ce champ. Nous avons préféré joindre ces résultats de nature purement cohomologique à l'annexe déjà mentionnée plus haut (cf. paragraphe~(\ref{par_coh_et_ext_inf})).

Cette étude directe des déformations de faisceaux inversibles nous permet de donner une démonstration de la représentabilité du foncteur de Picard (sous de bonnes hypothèses) logiquement indépendante des travaux d'Aoki. Bien entendu nous aurions tout aussi bien pu utiliser le critère~\cite{Artin_Versal_defm}~5.3 d'Artin pour qu'un champ soit algébrique au lieu du théorème de représentabilité pour les espaces algébriques~\cite{Global_Analysis_1}~5.3 et obtenir ainsi l'algébricité du champ de Picard. En fait le paragraphe~\ref{champ_de_Picard} montre qu'il est essentiellement équivalent d'étudier le foncteur de Picard ou le champ de Picard (du moins lorsque $\X$ est cohomologiquement plat en dimension zéro).

Ce chapitre s'achève sur un résultat classique reliant la dimension du schéma de Picard d'un champ algébrique sur un corps (lorsqu'il existe) et sa lissité à l'origine. Nous calculons au passage son espace tangent.

\subsubsection{Composante neutre du foncteur de Picard}
Revenons à des exemples plus concrets et considérons le cas d'une courbe elliptique $E$ sur un corps $k$. On sait dans ce cas que la jacobienne $\Pic^0_{E/k}$, c'est-à-dire la composante connexe de l'élément neutre du schéma de Picard $\Pic_{E/k}$, est isomorphe à $E$ et que le schéma de Picard s'identifie à une infinité dénombrable de copies de cette jacobienne indexées par le degré des faisceaux inversibles. Il est bien évident que le schéma de Picard lui-même n'est pas propre. Cependant la jacobienne, elle, est propre. Il en va de même si l'on remplace $E$ par une courbe projective intègre et lisse sur un corps $k$. En fait la jacobienne possède souvent des propriétés de finitude agréables et contient en réalité une grande partie de l'information réellement utile fournie par le schéma de Picard. Sur une base quelconque la situation est naturellement plus complexe car on ne dispose plus d'\emph{un} élément neutre mais d'une \emph{section} neutre. La première idée qui vient à l'esprit est de considérer la composante connexe de l'élément neutre dans chaque fibre (c'est un ouvert de la fibre) et d'appeler jacobienne la réunion de toutes ces composantes. Le problème est qu'il n'est pas évident \emph{a priori} que cette construction fournisse un ouvert du foncteur de Picard. Lorsque ce dernier est représentable par un schéma, on peut appliquer le corollaire~IV~(15.6.5) des EGA (\cite{EGA4_3}) sur la composante connexe des fibres le long d'une section. Mais il ne s'applique pas tel quel lorsque nous avons affaire à un espace algébrique. Notre premier travail a donc été de généraliser ce résultat au cas des espaces algébriques afin de s'assurer que la construction décrite précédemment fournirait bien un sous-espace algébrique ouvert $\piczero$ du foncteur de Picard.

Une fois la jacobienne définie, il semble naturel d'espérer qu'elle soit propre sur $S$ lorsque $X$ est lisse par exemple. Ceci est en effet vrai lorsque $\X$ est un schéma (cf. \cite{poly_Kleiman}~5.20). Nous avons fait un premier pas vers l'obtention de cette propreté en démontrant le résultat lorsque $\X$ est un champ algébrique normal sur un corps.

\vskip 2mm \noindent
{{\bf Théorème~\ref{composante_neutre_propre_sur_k}} \it
On suppose que $\X$ est un champ algébrique propre, géométriquement normal et cohomologiquement plat en
dimension zéro sur $\Spec k$. Alors la composante neutre $\text{Pic}^0_{\X/k}$
du schéma de Picard est propre sur $k$.
}\vskip 0.5mm \noindent

Il suffirait maintenant de relier la propreté de $\piczero$ à la propreté de chacune de ses fibres pour en déduire que $\piczero$ est propre dès que $\X$ est lisse sur $S$. Nous espérons aboutir bientôt à ce résultat.

\subsubsection{Quelques exemples}
Nous avons essayé de voir à travers deux exemples comment l'ajout d'une \og structure champêtre\fg\ à un schéma modifie son foncteur de Picard.

Le premier exemple que nous considérons est celui du champ des \og racines $n^{\text{ièmes}}$ d'un faisceau inversible\fg. On fixe un schéma $X$, un faisceau inversible $\Lc$ sur $X$ et un entier $n$ strictement positif. On peut associer à ces données un champ $[\Lc^{\frac1n}]$ qui classifie les faisceaux inversibles $\Mc$ munis d'un isomorphisme de $\Mc^{\otimes n}$ vers $\Lc$. Alors $[\Lc^{\frac1n}]$ est une $\mun$-gerbe sur $X$, et son foncteur de Picard s'obtient à partir de celui de $X$ en lui adjoignant d'une manière très naturelle une racine \iem{n} de la classe du faisceau $\Lc$.

Plutôt que de mettre une structure champêtre sur le schéma $X$ tout entier, on peut aussi ajouter à une courbe une structure de $\mu_r$-gerbe en quelques points isolés. C'est ce que font Abramovich et Vistoli dans \cite{Abramovich_Vistoli_CSSM}. En essayant de compactifier l'espace des morphismes stables d'une courbe nodale $C$ vers un champ de Deligne-Mumford $\X$ fixé, ils se sont aperçus qu'aux points limites il était nécessaire d'autoriser une telle structure en les n\oe uds de $C$, et même en des points marqués. Une courbe tordue d'Abramovich et Vistoli ressemble à quelque chose comme ceci :
\begin{center}
\includegraphics[scale=0.2]{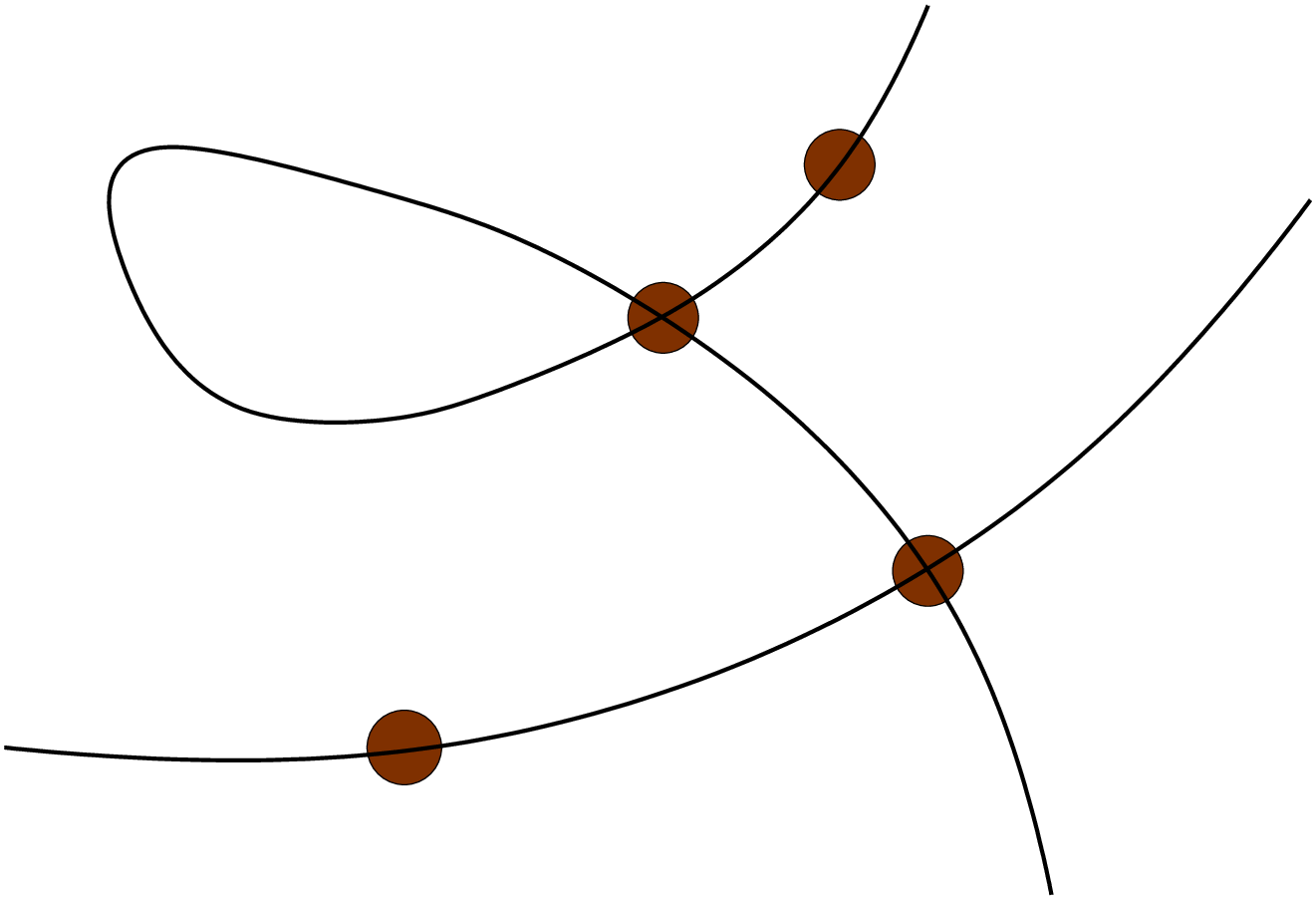}\\
\medskip
\end{center}
où les \og gros \fg\ points signalent la présence d'une action d'un groupe cyclique $\mu_r$. Nous décrivons dans la section~\ref{courbes_tordues} le foncteur de Picard d'une courbe tordue \emph{lisse}.

\subsubsection{Annexe}
Cette dernière partie est un peu différente des autres dans la mesure où elle ne concerne \emph{a priori} pas du tout le foncteur de Picard. C'est ce qui lui a valu l'appellation d'\og annexe\fg\ mais en réalité elle aurait tout aussi bien pu se trouver tout au début et porter le nom de \og résultats préliminaires sur la cohomologie lisse-étale\fg, ce qui aurait d'ailleurs été plus cohérent d'un point de vue purement logique. Nous invitons donc le lecteur soucieux de respecter l'ordre logique des propositions à commencer sa lecture par là.

Cette annexe trouve essentiellement sa raison d'être dans le manque de références sur la cohomologie lisse-étale des champs algébriques. Au fil de nos interrogations sur le foncteur de Picard (et particulièrement sur sa description cohomologique), nous avons naturellement été amené à utiliser des propriétés des groupes de cohomologie qui \og ne pouvaient qu'être vraies\fg\ (sans quoi les définitions auraient pu être qualifiées de mauvaises et définitivement abandonnées) mais pour lesquelles il ne semble pas exister de référence. Il est bon par exemple de vérifier (cf.~(\ref{coh_et_egale_coh_liset})), une fois construits les groupes de cohomologie \emph{lisse-étale} d'un champ algébrique $\X$, qu'ils coïncident avec les groupes de cohomologie \emph{étale} de $\X$ lorsque $\X$ est un schéma (ou plus généralement un champ de Deligne-Mumford). Il convient également de vérifier que si $\Fc$ est un faisceau lisse-étale abélien sur $\X$ et si $f : \X \fleche \Y$ est un morphisme de champs algébriques, alors le faisceau image directe supérieure $R^qf_*\Fc$ est bien, comme on le pense, le faisceau associé au préfaisceau qui à tout ouvert lisse-étale $(U,u)$ de $\Y$ associe le $q^{\textrm{ième}}$ groupe de cohomologie de $\X\times_{\Y} U$ à valeurs dans $\Fc_{\X\times_{\Y} U}$. Contrairement à ce que l'on pourrait croire au premier abord, cette dernière propriété n'est pas une conséquence immédiate de l'étude générale proposée dans SGA4, ceci à cause d'un défaut de fonctorialité du topos lisse-étale des champs algébriques (le foncteur $f^{-1}$ n'est pas toujours exact). Ce défaut de fonctorialité a pour conséquence fâcheuse que la \og machine \fg\ SGA4 ne s'applique pas toujours et qu'il faut par conséquent travailler un peu plus finement pour obtenir un certain nombre de propriétés d'apparence pourtant élémentaire sur la cohomologie lisse-étale des champs algébriques. Ce travail a été fait en grande partie par Olsson et Laszlo pour le cas des faisceaux quasi-cohérents (voir \cite{Olsson_Sheaves_on_Artin_stacks}) ou des coefficients finis (voir \cite{Laszlo_Olsson_Six_operationsI}). Mais les faisceaux abéliens qui ne jouissent pas d'une telle structure ont pour l'instant été laissés de côté.

Or le faisceau abélien que l'on rencontre le plus souvent lorsque l'on s'intéresse au foncteur de Picard n'est pas quasi-cohérent : il s'agit de $\gm$. Nous avons donc été amené à démontrer au fil de nos travaux diverses propriétés. Leur nombre croissant nous a finalement conduit à les regrouper dans un chapitre à part. Cette annexe est donc une sorte de \og fourre-tout \fg\ cohomologique qui, bien loin de prétendre à l'exhaustivité, se contente au contraire de répondre aux strictes exigences des autres chapitres.  Nous espérons tout de même que nous aurons ainsi contribué à combler une lacune de la littérature existante.
Nous renvoyons à l'introduction de ladite annexe pour un exposé plus détaillé des propriétés que le lecteur intéressé pourra y trouver. Signalons juste l'introduction d'un site un peu plus gros que le site lisse-étale, le site lisse-lisse champêtre, qui induit le même topos mais se comporte de manière un peu plus agréable à certains égards (notamment vis-à-vis des images directes). Il nous a rendu de fiers services et pourra sans doute encore se montrer utile.

\subsection*{Conventions}
\vskip -1.5mm
Suivant \cite{LMB}, sauf mention expresse du contraire, tous les champs algébriques (\emph{a fortiori} tous les schémas et tous les espaces algébriques) seront quasi-séparés. Un champ \og algébrique \fg\ non-quasi-séparé, \emph{i.e.} un champ dont la diagonale est représentable et localement de type fini, et qui admet une présentation lisse, sera dit préalgébrique.

\chapter{Foncteurs de Picard}

En 1965, c'est-à-dire quatre ans avant la parution des premiers articles fondant la théorie des champs, Mumford publie l'article~\cite{Mumford_Picard_groups}. Il définit en particulier une notion de faisceau inversible sur un \og problème de modules\fg\ $\Mc$. Un tel faisceau inversible est la donnée d'un faisceau inversible sur $S$ pour toute famille de courbes $X\fleche S$, et d'isomorphismes de transition entre ces faisceaux inversibles vérifiant une condition de compatibilité naturelle. Bien que le champ considéré par Mumford se trouve être algébrique, on voit bien que cette condition n'est en rien nécessaire  et que la notion de faisceau inversible prend un sens sur un champ quelconque. Mumford définit ensuite le groupe de Picard de $\Mc$ comme étant le groupe des classes d'isomorphie de faisceaux inversibles sur $\Mc$, montre que ce groupe est isomorphe au groupe $H^1(\Mc, \gm)$, puis le calcule à titre d'exemple dans le cas où $\Mc$ est le champ $\Mc_{1,1}$ des courbes elliptiques.

La notion de faisceau inversible sur un champ est devenue aujourd'hui commune. Nous présentons dans une première section différentes réalisations de la catégorie des faisceaux inversibles sur un champ (parfois algébrique) et nous montrons l'équivalence entre ces différents points de vue. Puis nous définissons le groupe de Picard d'un champ et nous en donnons une interprétation cohomologique. Nous suivons ensuite Grothendieck pour définir les foncteurs de Picard relatifs d'un champ et nous comparons (cf. thm.~(\ref{comparaison_des_foncteurs_de_Picard})) les différents foncteurs obtenus. Nous introduisons enfin le champ de Picard, c'est-à-dire le champ des faisceaux inversibles, et tentons de démêler les liens étroits qui l'unissent au foncteur de Picard. Cette partie s'achève sur des propriétés de commutation aux limites inductives ou projectives.

\section{Groupe de Picard d'un champ}

\subsection{Différents visages des faisceaux inversibles}

Suivant Mumford (cf. \cite{Mumford_Picard_groups}~p.~64), nous adoptons la définition ci-dessous.

\begin{sousdefi}
\label{definition_fi}
Soit $\X$ un $S$-champ. On note $p$ son morphisme structural. Un faisceau inversible $L$ sur $\X$ est la donn\'ee de :
\begin{itemize}
\item[(i)] pour tout $U\in\ob\aff$ et tout $x\in \ob \X_U$, un faisceau inversible $L(x)$ sur $U$ ;
\item[(ii)] pour toute fl\`eche 
$\varphi : y \fleche x$ dans $\X$, un isomorphisme
$$\xymatrix{L(\varphi) : L(y) \ar[r]& p(\varphi)^*L(x)} ;$$
\end{itemize}
v\'erifiant la condition de compatibilit\'e suivante : si $\varphi : y\fleche x$ et $\psi : z\fleche y$ sont deux fl\`eches composables, on a
$$p(\psi)^*L(\varphi) \circ L(\psi) = L(\varphi \circ \psi)$$
(modulo un isomorphisme canonique que nous nous dispenserons d'\'ecrire).
\end{sousdefi}

\begin{sousexemple}\rm
On d\'efinit le faisceau trivial $\O_{\X}$ en posant $\O_{\X}(x)=\O_{p(x)}=\O_U$ pour tout objet $x$ de $\X$ et $\O_{\X}(\psi)=id_{\O_{p(y)}}$ pour toute fl\`eche $\psi :y\fleche x$.
\end{sousexemple}

Un morphisme entre deux faisceaux inversibles est simplement une collection de morphismes compatibles avec les isomorphismes de changement de base.

\begin{sousdefi}
Soient $L$ et $M$ deux faisceaux inversibles sur $\X$. Un morphisme $\Phi:
L\fleche M$ est la donn\'ee pour tout $x\in \ob \X$ d'un morphisme
$$\xymatrix{\Phi(x) : L(x) \ar[r]& M(x)}$$
tel que pour tout $\varphi : y \fleche x$ dans $\X$, le carr\'e suivant
commute (o\`u $\tilde{\varphi}=p(\varphi)$) :
$$\xymatrix{L(y) \ar[r]^{\sim}_{L(\varphi)} \ar[d]_{\Phi(y)} &
   \tilde{\varphi}^*L(x) \ar[d]^{\tilde{\varphi}^*\Phi(x)} \\
   M(y) \ar[r]^{\sim}_{M(\varphi)} &\tilde{\varphi}^*M(x)}.$$
\end{sousdefi}

On obtient de cette mani\`ere une cat\'egorie $\inv(\X)$ des faisceaux
inversibles sur $\X$. On notera $\inv'(\X)$ la sous-catégorie de $\inv(\X)$ dans
laquelle on ne garde que les isomorphismes.
Dans le cas où $\X$ est algébrique, on peut aussi voir un faisceau inversible, de mani\`ere peut-\^etre plus naturelle, comme \'etant un $\O_{\X}$-module localement libre
de rang 1 sur le site
lisse-\'etale\footnote{La définition du site lisse-étale est rappelée en annexe.} de $\X$ d\'efini dans \cite{LMB} en (12.1), où $\O_{\X}$ désigne
cette fois le faisceau structural de $\X$ défini dans loc.cit. en (12.7.1).
En fait, la proposition suivante montre que ces deux constructions sont \'equivalentes. On rappelle qu'un $\O_{\X}$-module $\Mc$ est dit localement libre de rang~1 si pour tout ouvert lisse-étale $(U,u)$ de $\X$, le $\O_{U_{\textrm{ét}}}$-module $\Mc_{U,u}$ est localement libre de rang~1. Cela équivaut à dire qu'il existe une présentation $P : X \fleche \X$ de $\X$ telle que le $\O_{X_{\textrm{ét}}}$-module $\Mc_{X,P}$ soit 
isomorphe à $\O_{X_{\textrm{ét}}}$. Un tel $\O_{\X}$-module est en particulier cartésien au sens de \cite{LMB}~(12.7.3) d'après la proposition~(13.2.1) du même ouvrage.

\begin{sousprop} \label{equiv cat de fi}
Soit $\X$ un $S$-champ algébrique. On a une \'equivalence naturelle de cat\'egories entre la cat\'egorie $\inv(\X)$ des faisceaux inversibles sur $\X$ et la cat\'egorie des $\O_{\X}$-modules localement libres de rang~1 sur le champ alg\'ebrique annel\'e $(\X, \O_{\X})$.
\end{sousprop}
\begin{demo}
Soit $\Mc$ un $\O_{\X}$-module localement libre de rang~1 pour la topologie
lisse-\'etale sur $\X$.
% Par définition il existe une présentation $P : X \fleche
% \X$ telle que le module lisse-étale $P^*\Mc$ soit isomorphe à
% $\O_{X_{\textrm{lis-ét}}}$. En particulier (cf. proposition~(13.2.1) de
% \cite{LMB}), $\Mc$ est un faisceau cart\'esien et le
% module étale
% $\Mc_{X,P}$ est isomorphe \`a $\O_{X_{\textrm{\'et}}}$.
Pour tout $U\in \ob \aff$ et
tout morphisme $u : U \fleche \X$, on pose $$M(u)=(u^*\Mc)_{\textrm{zar}}.$$
Les isomorphismes de changement de base $M(v)\fleche
\varphi^*M(u)$ se d\'efinissent de mani\`ere \'evidente. On v\'erifie facilement
que $M$ ainsi construit est bien un faisceau inversible, et que cette
construction se prolonge de mani\`ere naturelle en un foncteur $F$ de la
cat\'egorie des $\O_{\X}$-modules localement libres de rang~1, vers la cat\'egorie $\inv(\X)$.

R\'eciproquement, \'etant donn\'e un syst\`eme $(M(u),M(\varphi))$, montrons
comment cons\-truire un $\O_{\X}$-module
localement libre de rang~1. D'apr\`es
le lemme~(12.1.2) de~\cite{LMB}, il suffit de le construire sur le site
\underline{Lis-\'et}$(\X)$, dont les objets sont les $(U,u)$ de Lis-\'et$(\X)$
avec $U\in\ob\aff$, et o\`u les familles couvrantes sont les familles finies de
morphismes qui forment une famille couvrante dans Lis-\'et$(\X)$. Pour un tel
$(U,u)$, on pose
$$\Mc(U,u)=\Gamma(U,M(u)).$$
Les flèches de restriction sont définies de manière évidente à partir des isomorphismes $M(\varphi)$. Remarquons que pour tout ouvert lisse-étale $(U,u)$ de $\X$, le faisceau étale $\Mc_{U,u}$ n'est autre que $\varepsilon^*M(u)$, où $\varepsilon$ est le morphisme naturel de topos de $U_{\textrm{ét}}$ dans $U_{\textrm{zar}}$.\footnote{Le foncteur $\varepsilon^*$ induit une équivalence de catégories entre les faisceaux inversibles sur $U$ au sens étale et les faisceaux inversibles au sens de Zariski. \`A partir de maintenant nous identifierons ces deux catégories.}

Le foncteur $G$ ainsi d\'efini est clairement un inverse \`a gauche de $F$. Il
reste \`a montrer que c'est un inverse \`a droite. On part d'un faisceau
inversible $M$ donn\'e par un syst\`eme $(M(u),M(\varphi))$, et on lui associe
$\Mc=G(M)$. Il faut montrer que pour tout $U\in \ob\aff$ et tout $u\in\ob\X_U$,
$M(u)$ est canoniquement isomorphe \`a $u^*\Mc$. Pour ce faire, il suffit de
montrer que pour tout ouvert $(V,v)$ de \underline{Lis-\'et}$(U)$,
$(u^*\Mc)_{V,v}$ est canoniquement isomorphe \`a $M(u)_{V,v}=v^*M(u)$. Or
$$(u^*\Mc)_{V,v}=\lind (f^*\Mc_{U',u'})$$ o\`u la limite
inductive est prise sur l'ensemble des diagrammes 2-commutatifs
$$\xymatrix{V \ar[r]^v \ar[d]_f & U\ar[d]^u \\
U' \ar[r]^{u'} & \X
}$$
o\`u $(U',u')\in\ob$\,\underline{Lis-\'et}$(\X)$.
Pour un tel diagramme, la d\'efinition de $\Mc$ et le fait que $(U',u')$ soit un
ouvert lisse-\'etale montrent que $\Mc_{U',u'}\simeq M(u')$. Alors gr\^ace aux
isomorphismes de changement de base sous-jacents \`a $M$, on obtient un
isomorphisme entre $f^*\Mc_{U',u'}$ et $v^*M(u)$, donc $(u^*\Mc)_{V,v}\simeq
v^*M(u)$, ce qui montre le r\'esultat attendu.
\end{demo}

\begin{sousremarque}\rm
\label{rem_fi_cas_des_schemas}
Si $\X=X$ est un $S$-espace algébrique, l'équivalence de catégories
$$\flechen{\Mod_{\textrm{qcoh}}(\O_{X_{\textrm{ét}}})}{\sim}
{\Mod_{\textrm{qcoh}}(\O_{X_{\textrm{lis-ét}}})}$$
mentionnée dans \cite{LMB} (13.2.3) induit une équivalence entre la catégorie
$\inv(\X)$ définie ci-dessus et la catégorie des faisceaux inversibles
sur $X$ au sens étale.

En particulier si $X$ est un $S$-schéma on a une
équivalence de catégories naturelle entre la catégorie $\inv(\X)$ et la catégorie des
faisceaux inversibles sur $X$ au sens de Zariski.
\end{sousremarque}

Le résultat ci-dessous découle de l'étude de la descente fidèlement plate des
modules quasi-cohérents proposée dans \cite{LMB}. Nous ne
l'énonçons que pour le cas particulier des faisceaux inversibles car c'est là
l'objet principal de
notre étude, mais en réalité le paragraphe (13.5) de \cite{LMB} montre qu'il est
encore valable tel quel en remplaçant \og faisceau
inversible \fg\ par \og module quasi-cohérent\fg. Il nous sera surtout utile
lorsque le morphisme $f$ est une présentation de $\X$ par un espace algébrique
ou un schéma.

\begin{sousprop}[\cite{LMB}~(13.5)]
\label{descente_fidelement_plate}
Soit $f : \Y\fleche \X$ un morphisme fidèlement plat de champs algébriques, que l'on
suppose de plus quasi-compact ou localement de présentation finie. Alors la
catégorie des faisceaux inversibles sur $\X$ est équivalente à
la catégorie suivante. Un objet est un couple $(\L,\alpha)$ où $\L$ est un
faisceau inversible sur $\Y$ et
$\alpha : p_1^*\L \fleche p_2^*\L$
est un isomorphisme tel que, à des isomorphismes canoniques près,
$(p_{23}^*\alpha) \circ (p_{12}^*\alpha) = p_{13}^*\alpha$. Un morphisme
de $(\L,\alpha)$ vers $(\Mc, \beta)$ est un morphisme $\gamma :
\L \fleche \Mc$ tel que $(p_2^*\gamma)\circ \alpha= \beta\circ
(p_1^*\gamma)$.
\end{sousprop}
\begin{demo}
Avec les notations de \cite{LMB}~(13.5), la catégorie décrite ci-dessus est
équivalente à la catégorie des $\O_{\Y_{\bullet}}$-modules quasi-cohérents
cartésiens et localement libres de rang 1. Par~(13.5.5), le morphisme
$f$ est de descente cohomologique effective, et donc par~(13.5.4), il induit une
équivalence de catégories entre la catégorie des $\O_{\X}$-modules
quasi-cohérents et la catégorie des $\O_{\Y_{\bullet}}$-modules
quasi-cohérents cartésiens. Il est clair que cette équivalence de catégories
préserve les faisceaux inversibles.
\end{demo}

\noindent {\sc Les faisceaux inversibles comme morphismes de $\X$ dans $\bgm$.}
\vskip 2mm

Si $\X$ et $\Y$ sont deux $S$-champs alg\'ebriques, $\Hom(\X,\Y)$ est la
cat\'egorie dont les objets sont les 1-morphismes de $\X$ vers $\Y$ et dont les
fl\`eches sont les 2-isomorphismes entre les 1-fl\`eches. On rappelle aussi que
l'on note $\bgm$ le $S$-champ dont la fibre
en $U$ est la cat\'egorie des $\gm$-torseurs sur $U$.

\begin{sousprop} \label{fi_equiv_morphisme_de_X_dans_BGm}
La cat\'egorie $\inv'(\X)$ des faisceaux inversibles sur $\X$ est \'equivalente
\`a la cat\'egorie $\Hom(\X,\bgm)$.
\end{sousprop}
\begin{demo}
Un morphisme $F :\X\fleche \bgm$ est la donn\'ee pour tout $U\in\ob\aff$ et pour tout $x\in\ob\X_U$ d'un
$\gm$-torseur $F(x)$ sur $U$, ces donn\'ees \'etant de plus compatibles au
changement de base. Autrement dit pour toute fl\`eche
$\varphi : V \fleche U$
dans $\aff$ et pour tout $x\in\ob\X_U$, on se donne un isomorphisme $F(\varphi,
x) : F(\varphi^*x)=F(x\circ \varphi) \fleche \varphi^*F(x)$. Le r\'esultat
d\'ecoule donc simplement du fait que la donn\'ee d'un $\gm$-torseur sur $U$ est
\'equivalente \`a la donn\'ee d'un faisceau inversible sur $U$.
\end{demo}

\subsection{Description du groupe de Picard}
On d\'efinit de mani\`ere naturelle le produit tensoriel de deux faisceaux
inversibles en posant\footnote{Il est préférable de procéder ainsi plutôt que de considérer le produit tensoriel des faisceaux localement libres de rang 1 sur le site lisse-étale d'un champ algébrique, car ce dernier n'aurait plus de sens sur un champ quelconque.} :
$$(L\otimes M)(x)=L(x)\otimes M(x),$$
$$(L\otimes M)(\varphi)=L(\varphi)\otimes M(\varphi).$$
On v\'erifie alors que le produit tensoriel pr\'eserve les classes d'isomorphie,
que le faisceau structural $\O_{\X}$ est neutre, et que tout faisceau inversible
a un inverse donn\'e par :
$$L^{-1}(x)=L(x)^{-1}=\fHom(L(x),\O_U)$$
et $L^{-1}(\varphi)$ est induit par $L(\varphi)$. Il est clair que ce produit
tensoriel coïncide avec le produit tensoriel des faisceaux quasi-cohérents via
les équivalences de catégories ci-dessus. On peut alors poser la
d\'efinition suivante.

\begin{sousdefi}
On appelle groupe de Picard de $\X$ l'ensemble $\Pic \X$ des classes
d'isomorphie de faisceaux inversibles muni de la loi de groupe induite par le
produit tensoriel.
\end{sousdefi}

\begin{sousremarque}\rm
Si $\X$ est un sch\'ema ou un espace algébrique, la remarque
(\ref{rem_fi_cas_des_schemas}) montre que l'on retrouve ainsi le groupe de
Picard usuel.
\end{sousremarque}

\noindent {\bf Fonctorialit\'e}
Soit $F:\X \fleche \Y$ un morphisme de champs alg\'ebriques, et soit $L$ un
faisceau inversible sur $\Y$. On d\'efinit l'image inverse $F^*L$ de $L$ en
posant pour tout $x\in \ob\X$, $F^*L(x)=L(F(x))$, et pour tout
$\xymatrix@C=1pc{\varphi : y \ar[r]&x}$ dans $\X$, $F^*L(\varphi)=L(F(\varphi))$. On v\'erifie
que $F^*L$ ainsi d\'efini est bien un faisceau inversible sur $\X$. Le foncteur
$F^*$ se d\'efinit de mani\`ere tout aussi \'evidente sur les fl\`eches. De plus
il pr\'eserve le produit tensoriel, donc il induit un morphisme de groupes $\Pic
F : \Pic \Y \fleche \Pic \X$.

\begin{sousremarque}\rm
Via l'\'equivalence de cat\'egories \ref{equiv cat de fi}, le foncteur $F^*$
d\'efini ci-dessus correspond au foncteur image inverse d\'efini au chapitre 12
de \cite{LMB}.
\end{sousremarque}

La description cohomologique traditionnelle du groupe de Picard est encore
valable.

\begin{sousprop}
\label{pic_egal_h1}
Soit $\X$ un $S$-champ alg\'ebrique. Alors :
$$\Pic \X \simeq H^1(\X,\gm),$$
o\`u $H^1(\X, \gm)$ est le premier groupe de cohomologie du faisceau $\gm$ sur
$\X$ muni de la topologie lisse-\'etale calcul\'e au sens des foncteurs
d\'eriv\'es.
\end{sousprop}
\begin{demo}
Soit $T\fleche \X$ une pr\'esentation de $\X$. On reprend ici les notations
de l'annexe (\ref{annexe_desc_coh}) qui pr\'ec\`edent le th\'eor\`eme
(\ref{suite_spectrale_de_descente}). On consid\`ere en particulier le premier
groupe de cohomologie \og \`a la \v{C}ech \fg\ associ\'e \`a cette
pr\'esentation :
$$\check{H}^1(H^0(T^{\bullet},\gm))=
\frac{\Ker(p_{23}^*-p_{13}^*+p_{12}^*)}{\Im(p_1^*-p_2^*)}.$$ En observant que
pour tout $i$, $\gm(T^i)=\Aut(\O_{T^i})$, on voit que se donner un 1-cocycle de
\v Cech \`a valeurs dans $\gm$ revient \`a se donner une donn\'ee de descente
sur $\O_T$, et que deux telles donn\'ees de descente $g_1, g_2$ d\'efinissent le
m\^eme \'el\'ement dans $\check{H}^1(H^0(T^{\bullet},\gm))$ si et seulement si
$(\O_T,g_1)$ et $(\O_T,g_2)$ sont
isomorphes dans la cat\'egorie des faisceaux inversibles sur $T$ munis d'une
donn\'ee de descente relativement \`a $T\fleche \X$. Compte tenu de
l'\'equivalence (\ref{descente_fidelement_plate}) entre cette cat\'egorie et la
cat\'egorie des faisceaux inversibles sur $\X$, on voit que le groupe
$\check{H}^1(H^0(T^{\bullet},\gm))$ s'identifie \`a l'ensemble des classes
d'isomorphie de faisceaux inversibles sur $\X$ dont l'image inverse sur $T$ est
isomorphe \`a $\O_T$, donc au groupe $\Ker(\Pic(\X) \fleche \Pic(T))$. De plus
cet isomorphisme est fonctoriel en $T$. En passant \`a la limite inductive, on
obtient un isomorphisme
$$\lind \check{H}^1(H^0(T^{\bullet},\gm)) \simeq \Pic \X,$$
o\`u la limite inductive est prise sur l'ensemble des pr\'esentations lisses
$T\fleche \X$ de $\X$. On dispose par ailleurs de la suite exacte en bas
degr\'es associ\'ee \`a la suite spectrale de descente
(\ref{suite_spectrale_de_descente}) :
$$0\fleche \check{H}^1(H^0(T^{\bullet},\gm)) \fleche H^1(\X,\gm) \fleche
\check{H}^0(H^1(T^{\bullet},\gm)).$$
Cette suite exacte \'etant fonctorielle en $T$, on voit que les injections
des $\check{H}^1(H^0(T^{\bullet},\gm))$ dans $H^1(\X,\gm)$ induisent
un morphisme
injectif $$\lind \check{H}^1(H^0(T^{\bullet},\gm)) \fleche H^1(\X,\gm).$$
Pour conclure, il ne reste plus qu'\`a montrer, vu la suite exacte en bas
degr\'es ci-dessus, que pour tout $x\in H^1(\X,\gm)$, il existe une
pr\'esentation $T\fleche \X$ telle que l'image de $x$ dans $H^1(T,\gm)$ soit
nulle. Soit $0\fleche \gm \fleche \Ic \fleche \Qc \fleche 0$ une suite exacte de
faisceaux ab\'eliens, avec $\Ic$ injectif. On a alors un diagramme commutatif 
$$\xymatrix{H^0(\X,\Ic) \ar[r] \ar[d]& H^0(\X,\Qc) \ar[r] \ar[d]&
  H^1(\X,\gm) \ar[r] \ar[d]& H^1(\X,\Ic)=0\\
  H^0(T,\Ic_{|_T}) \ar[r] & H^0(T,\Qc_{|_T}) \ar[r] &H^1(T,\gm)}$$
dans lequel les lignes sont exactes. Soit $s\in H^0(\X,\Qc)$ un ant\'ec\'edent
de $x$. Comme $\Qc$ est le faisceau quotient $\Ic/\gm$ on sait par d\'efinition
du faisceau associ\'e \`a un pr\'efaisceau que $s$ provient localement de $\Ic$.
Autrement dit il existe une pr\'esentation de $\X$ telle que $s_{|_T}\in
H^0(T,\Qc_{|_T})$ provienne de $H^0(T,\Ic_{|_T})$, auquel cas l'image de $x$
dans $H^1(T,\gm)$ est nulle, conform\'ement \`a nos exigences.
\end{demo}

\section{Les foncteurs de Picard relatifs}

\subsection{Définitions et premières propriétés}
\begin{sousdefi}
Soit $\X$ un $S$-champ. Le foncteur de Picard na\"if de $\X/S$ est le foncteur :
$$\fonction{\Pic_{\X}}{(Sch/S)^{\circ}}{(Gr)}{T}{\Pic (\X\times_S T)}.$$
\end{sousdefi}

Le foncteur de Picard na\"if (parfois aussi appel\'e foncteur de Picard absolu,
m\^eme si cette terminologie peut sembler \'etrange dans la mesure o\`u il
d\'epend
\'evidemment de $S$) ayant aussi peu de chances d'\^etre repr\'esentable
que dans le cas des sch\'emas, on d\'efinit encore des foncteurs de Picard
relatifs obtenus par faisceautisation relativement \`a des topologies de plus en
plus fines.

\begin{sousdefi}
Nous noterons $P_{\X/S}$ le foncteur suivant :
$$\fonction{P_{\X/S}}{(Sch/S)^{\circ}}{(Gr)}{T}
{\Pic (\X\times_S T)/\Pic(T)}.$$
Nous noterons $\piczar$ (resp. $\piceth, \picfppf$) le faisceau associ\'e \`a
$P_{\X/S}$ pour la topologie de Zariski (resp. \'etale, fppf).
\end{sousdefi}

\begin{sousremarque}\rm
Nous avons donc cinq foncteurs de Picard avec des morphismes naturels 
$\Pic_{\X} \fleche P_{\X/S} \fleche \piczar \fleche \picet \fleche \picfppf$,
et il est facile de voir que les trois derniers sont les faisceaux associ\'es,
pour la topologie indiqu\'ee, \`a chacun des foncteurs pr\'ec\'edents, y compris
le foncteur de Picard absolu $\Pic_{\X}$.
\end{sousremarque}

\begin{sousremarque}\rm
Il est tout aussi formel de voir que si $\pic$ d\'esigne l'un des cinq
foncteurs d\'efinis ci-dessus, alors
$$\pic(T) \simeq \Pic_{(\X\times_S T)/T}(T).$$
En particulier la formation de ces foncteurs commute au changement de base.
\end{sousremarque}

\begin{sousremarque}\rm
\label{P_isom_Pet_si_k_alg_clos}
On peut d'ores et déjà comparer les groupes de sections de ces différents
foncteurs dans quelques cas particuliers grâce à la remarque de topologie
suivante. Supposons que $U$ soit un $S$-schéma affine tel que, pour toute
famille couvrante $(U_i\fleche U)_{i\in I}$ pour une certaine topologie
\emph{(Top)}, il existe un indice $i\in I$ tel que le morphisme $U_i\fleche U$
ait une section. Alors le morphisme $\Pic_{\X}(U) \fleche
\Pic_{\X/S \, \rm{(Top)}}(U)$ est un isomorphisme. Ceci est vrai en particulier
pour la topologie de Zariski lorsque $U=\Spec A$ est le spectre d'un anneau
local. On obtient donc dans ce cas un isomorphisme :
$$\xymatrix{\Pic_{\X}(A) \ar[r]^-{\sim}& \piczar(A)}.$$
Si $A$ est un anneau local strictement hensélien, la condition précédente est vérifiée pour
$U=\Spec A$ avec la topologie étale, d'où un isomorphisme :
$$\xymatrix{\Pic_{\X}(A) \ar[r]^-{\sim}& \picet(A)}.$$
Il en va encore de même lorsque $U$ est le spectre d'un corps algébriquement
clos pour la topologie \emph{(fppf)}.
\end{sousremarque}

\begin{sousprop}
\label{description_coh_de_picet_et_picfppf}
Soit $\X$ un $S$-champ alg\'ebrique. Alors pour tout $S$-sch\'ema $T$:
$$\piceth(T)\simeq H^0(T,R^1{f_T}_* \gm)$$
$$\picfppf(T)\simeq H_{\text{\rm pl}}^0(T,R^1f_{T*}^{\text{\rm pl}} \gm).$$
\end{sousprop}
\begin{demo}
D'après la proposition~(\ref{prop_image_directe_sup}) appliquée au morphisme
$f_T : \X_T=\X\times_S T \fleche T$
la restriction au site étale de $T$ du faisceau $R^1f_{T*}\gm$ est le
faisceau \'etale associ\'e \`a
$$\xymatrix{U \ar@{|->}[r] & H^1(\X_T\times_T U, \gm)=\Pic_{\X_T}(U)}.$$
Donc $(R^1f_{T*}\gm)_{\text{ét}}=\Pic_{\X_T/T \, \text{(\'Et)}}$ et en particulier
$$H^0(T, R^1f_{T*}\gm)=\Pic_{\X_T/T \, \text{(\'Et)}}(T)=\picet(T).$$
On obtient la seconde assertion en appliquant exactement le même raisonnement pour la topologie \emph{fppf}. On utilise~(\ref{images_directes_superieures_fppf})
au lieu de~(\ref{prop_image_directe_sup}).
\end{demo}

\subsection{Morphismes cohomologiquement plats en dimension zéro}

Afin de comparer ces diff\'erents foncteurs, comme dans
le cas des sch\'emas, la notion de morphisme cohomologiquement plat en
dimension z\'ero nous sera utile. Commen\c cons par remarquer que si
$F : \X\fleche \Y$ est un morphisme de champs alg\'ebriques, on a un
morphisme naturel
\begin{equation} \label{morph}
F^{\sharp} : \O_{\Y} \fleche F_* \O_{\X}.
\end{equation}
En effet, compte tenu des d\'efinitions, un tel morphisme correspond \`a la
donn\'ee, pour tout ouvert lisse-\'etale $(U,u)$ de $\Y$, d'un morphisme
$\disp \O_U \fleche \lpro \varphi_*\O_V$
o\`u la limite projective est prise sur
l'ensemble des carr\'es 2-commutatifs
$$\xymatrix{
V \ar[d]_{\varphi} \ar[r] & \X \ar[d]^F\\
U\ar[r]^u &\Y},$$
avec $(V,v)\in \ob \liset(\X)$.
Pour chaque $(U,u)$, on prend le morphisme
de $\O_U$ vers $\disp \lpro \varphi_*\O_V$
induit par le syst\`eme compatible des $\O_U \fleche \varphi_*\O_V$.

\begin{sousdefi}
On dit qu'un morphisme de champs alg\'ebriques $F : \X \fleche \Y$ est
cohomologiquement plat en dimension z\'ero s'il est plat et si le
morphisme~\eqref{morph}
ci-dessus est un isomorphisme universellement.
\end{sousdefi}

\begin{sousremarque} \rm
On rappelle que si $\X$ est un champ alg\'ebrique, et si $\Fc$ est un faisceau
sur $\X$, on d\'efinit l'ensemble des sections globales de $\Fc$ par
$$\Gamma(\X, \Fc)=\Gamma(S,(A_*\Fc)_{S,\Id}),$$
o\`u $A:\X \fleche S$ est le
morphisme structural de $\X$. Alors si $F : \X \fleche \Y$ est
cohomologiquement plat en dimension z\'ero, il r\'esulte des d\'efinitions que
$$\Gamma(\X,\O_{\X})=\Gamma(\Y, \O_{\Y}).$$
On en d\'eduit $$\Aut(\O_{\X})=\Aut(\O_{\Y}).$$
\end{sousremarque}

\begin{sousremarque} \rm
Un morphisme cohomologiquement plat en dimension z\'ero est surjectif. En effet,
montrons que pour tout point $s : \Spec k \fleche \Y$ de $\Y$, la fibre $\X_s$
de $\X$ au-dessus de $s$ est non vide. Le morphisme $\X_s \fleche \Spec k$ est
encore cohomologiquement plat en dimension zéro et en particulier on a
$\Gamma(\X_s,\O_{\X_s})=\Gamma(\Spec k, \O_{\Spec k})=k$.
L'assertion r\'esulte alors du fait qu'un
champ alg\'ebrique est vide si et seulement si l'anneau de ses sections globales
est r\'eduit \`a z\'ero.
\end{sousremarque}

\begin{sousremarque} \rm
\label{remarque_section}
On déduit de la remarque précédente que si $S$ est un schéma et
$f~:~\X \fleche S$ un morphisme localement de présentation finie et
cohomologiquement plat en dimension zéro de $S$-champs algébriques
alors $f$ a une section localement pour la topologie \emph{(fppf)} sur $S$. On
vient en effet de voir que $f$ est surjectif,
donc il est fidèlement plat. Soit $S' \fleche \X$ une présentation de $\X$ par
un schéma. Alors la famille $(S'\fleche S)$ est une famille couvrante pour la
topologie \emph{(fppf)}, et il est clair que le morphisme induit par $f$ sur
$S'$ a une section.
\end{sousremarque}

\begin{sousremarque}\rm
La notion de platitude cohomologique est intimement liée à la conne\-xité des fibres géométriques. En effet on peut montrer que si $\X$ est cohomologiquement plat en dimension zéro sur $\Y$, ses fibres sont géométriquement connexes. Il suffit pour cela de montrer que $\X$ est connexe dans le cas où $\Y$ est le spectre d'un corps $L$. C'est immédiat puisque l'anneau $\Gamma(\X,\O_{\X})$ est alors isomorphe à $L$, donc intègre.
\end{sousremarque}

\begin{sousremarque}\rm
Réciproquement, dans le cas où le morphisme $f$ est propre et plat, il suffit que les fibres de $\X$ soient géométriquement connexes et géométriquement réduites pour que $\X$ soit cohomologiquement plat en dimension zéro sur $\Y$. En effet, commençons par montrer que pour tout point géométrique $y : \Spec \Omega \fleche \Y$ de $\Y$, l'anneau $H^0(\X_y)$ des fonctions globales sur la fibre géométrique $\X_{y}$ est égal à $\Omega$. D'après le théorème de finitude pour les morphismes propres, c'est une algèbre de dimension finie sur $\Omega$. Comme $\X_y$ est connexe et réduit, $\Spec H^0(\X_y)$ l'est aussi, donc c'est le spectre d'un corps (puisqu'il est par ailleurs fini sur $\Omega$). Comme $\Omega$ est algébriquement clos on a bien $H^0(\X_y)=\Omega$. Ceci prouve en particulier que la fonction qui à $y$ associe $\dim_\kappa(y) H^0(\X_y)$ est constante sur $\Y$. On en conclut avec les résultats classiques sur la cohomologie et les changements de base que le morphisme $f^{\sharp} : \Oc_{\Y} \fleche f_*\Oc_{\X}$ est un isomorphisme (au moins lorsque $\Y$ est noethérien et réduit).
\end{sousremarque}

Le lemme suivant donne un critère de \og trivialité relative \fg\  pour un faisceau
inversible sur un $\Y$-champ algébrique $\X$ cohomologiquement plat en dimension
zéro.

\begin{souslem}
\label{critere_relative_trivialite}
Soit $f : \X \fleche \Y$ un morphisme quasi-compact et cohomologiquement plat en
dimension zéro de $S$-champs algébriques. Soit $\Lc$ un faisceau inversible sur
$\X$. Les propositions suivantes sont équivalentes :
\begin{itemize}
\item[(1)] Le morphisme naturel d'adjonction
$$a_{\Lc} : f^*f_*\Lc \fleche \Lc$$
est un isomorphisme.
\item[(2)] Il existe un faisceau inversible $\Mc$ sur $\Y$ tel que $f^*\Mc$ soit
isomorphe à $\Lc$.
\end{itemize}
\end{souslem}
\begin{demo}
Pour montrer l'implication $(1) \Rightarrow (2)$, il suffit de montrer que le
faisceau $f_*\Lc$ est inversible. Par \cite{LMB}~(13.2.6)~(iii), on sait déjà
que c'est un faisceau quasi-cohérent, et vu l'hypothèse ci-dessus, son
image par $f^*$ est un faisceau inversible. De plus $f$ est un morphisme
fidèlement plat et quasi-compact, donc par descente fidèlement plate (cf.
(\ref{descente_fidelement_plate})) on en déduit que $f_*\Lc$ est inversible.

Pour la réciproque, considérons tout d'abord le cas où $\Lc$ est le faisceau
structural $\O_{\X}$. Le morphisme $f^{\sharp} : \O_{\Y} \fleche f_*\O_{\X}$ est
un isomorphisme par hypothèse, et son adjoint $f^*\O_{\Y}\fleche \O_{\X}$ en est
évidemment un aussi. De plus le diagramme
$$\shorthandoff{!;:?}
\xymatrix@!0 @R=3pc @C=3pc{f^*\O_{\Y} \ar[rd]\ar[rr]^{f^*(f^{\sharp})}&&
f^*f_*\O_{\X} \ar[ld]^{a_{\O_{\X}}}\\
& \O_{\X}}$$
commute, ce qui prouve que $a_{\O_{\X}}$ est un isomorphisme.

\begin{comment}
Passons maintenant au cas où $\Lc$ est trivial. Soit $\varphi : \Lc \fleche
\O_{\X}$ un isomorphisme. Alors les propriétés formelles de l'adjonction
montrent que le diagramme
$$\xymatrix{f^*f_*\Lc \ar[r]^{a_{\Lc}}\ar[d]_{f^*f_*\varphi}&
\Lc \ar[d]^{\varphi} \\
f^*f_*\O_{\X} \ar[r]^{a_{\O_{\X}}}&\O_{\X}}$$
commute, ce qui, vu le premier cas, prouve que $a_{\Lc}$ est un isomorphisme.
\end{comment}

Le cas où $\Lc$ est trivial s'en déduit immédiatement. %en utilisant les propriétés formelles de l'adjonction.
Dans le cas général, soit $\Mc$ un faisceau inversible sur $\Y$ tel que $\Lc$
soit isomorphe à $f^*\Mc$, et soit $\pi : Y \fleche \Y$ une présentation de $\Y$
telle que $\pi^* \Mc$ soit trivial. On note $\X'$ le champ algébrique obtenu par
changement de base.
$$\xymatrix{\X'\ar[d]_{f'} \ar[r]^{\pi'}\cartesien&\X \ar[d]^f\\
Y \ar[r]^{\pi}& \Y}$$
Notons $\psi$ le morphisme naturel de $\pi^*f_*\Lc$ vers $f'_*\pi'^*\Lc$. On
a de manière purement formelle un diagramme commutatif :
$$\shorthandoff{!;:?}
\xymatrix@!0 @R=2pc @C=8pc{&\pi'^*f^*f_*\Lc \ar[rd]^{\pi'^*a_{\Lc}}\\
f'^*\pi^*f_*\Lc\ar[ru]^{\sim}_{\text{can.}} \ar[rd]_{f'^*\psi} && \pi'^*\Lc.\\
&f'^*f'_*(\pi'^*\Lc)\ar[ru]_{a_{\pi'^*\Lc}}}$$
Or $a_{\pi'^*\Lc}$ est un isomorphisme d'après le cas précédent puisque
$\pi'^*\Lc$ est trivial. Par ailleurs, comme $\pi$ est plat, le morphisme $\psi$
est un isomorphisme en vertu de la
proposition~(\ref{images_directes_commutent_chgt_base_plat}), ce qui montre que
$\pi'^*a_{\Lc}$ est un isomorphisme. Enfin, en utilisant la théorie de la descente
fidèlement plate des modules quasi-cohérents (cf. \cite{LMB}~(13.5)),
on en déduit que $a_{\Lc}$ est un
isomorphisme.
\end{demo}

\subsection{Le théorème de comparaison}

Le théorème suivant permet de comparer les différents foncteurs de Picard d'un champ algébrique. Il généralise \cite{poly_Kleiman} 2.5.
\begin{sousthm}
\label{comparaison_des_foncteurs_de_Picard}
Soit $f : \X\fleche S$ un $S$-champ alg\'ebrique cohomologiquement plat en
dimension z\'ero. Alors les morphismes naturels
$$\xymatrix{P_{\X/S} \ar[r]^-{i_1} &\piczar \ar[r]^{i_2} &
\piceth}$$
sont injectifs, et le morphisme naturel
$$\flechen{\piceth}{}{\picfppf}$$
est un isomorphisme. Si de plus $f$ a une section localement pour la topologie
de Zariski, alors $i_2$ est un isomorphisme. Enfin si $f$ a une section, $i_1$
est lui aussi un isomorphisme.
\end{sousthm}
\begin{demo}
Soit $T$ un $S$-sch\'ema.
La suite spectrale de Leray
$$H^p(T, R^qf_{T*}\gm) \Rightarrow H^{p+q}(\X_T,\gm)$$
induit la suite exacte longue en bas degr\'es suivante :
\begin{multline} 
0\fleche H^1(T,f_{T*}\gm) \fleche H^1(\X_T,\gm)\fleche H^0(T, R^1f_{T*}\gm)
\fleche  \\
\fleche H^2(T,f_{T*}\gm) \fleche H^2(\X_T,\gm) \nonumber
\end{multline}
 (o\`u tous les calculs sont effectu\'es pour la topologie lisse-\'etale sur le
champ alg\'ebrique consid\'er\'e. En particulier lorsqu'il s'agit d'un sch\'ema,
cela revient d'après~(\ref{coh_et_egale_coh_liset}) \`a calculer sa cohomologie
pour la topologie \'etale.) Le morphisme
$f$ \'etant cohomologiquement plat en dimension z\'ero, on a $f_{T*}\gm=\gm$.
D'apr\`es les propri\'et\'es pr\'ec\'edentes, le d\'ebut de la suite exacte
ci-dessus fournit la suite exacte :
$$0\fleche \Pic(T)\fleche \Pic(\X_T)\fleche \picet(T)$$
ce qui montre que le morphisme naturel $P_{\X/S} \fleche \picet$
est injectif. Il en r\'esulte, d'une part, que le morphisme naturel
$P_{\X/S} \fleche \piczar$
est injectif, et d'autre part, en appliquant le foncteur \og faisceau associ\'e
pour la topologie de Zariski \fg, que le morphisme naturel
$\piczar \fleche \picet$
est injectif. En effet, ce foncteur est exact (c'est un fait g\'en\'eral dans le
cadre des topologies de Grothendieck) et laisse $\picet$ invariant
puisque c'est d\'ej\`a un faisceau pour la topologie de Zariski.

Dans le cas o\`u $f$ a une section, le morphisme induit par $f$ de
$H^2(T,f_{T*}\gm)$ vers $H^2(\X_T,\gm)$ a une r\'etraction, donc il est
injectif. Mais alors la fl\`eche de
$H^0(T, R^1f_{T*}\gm)$ vers $H^2(T,f_{T*}\gm)$
est nulle, et donc $H^1(\X_T,\gm)\fleche H^0(T, R^1f_{T*}\gm)$
est surjectif, d'o\`u une suite exacte courte :
$$0\fleche \Pic(T)\fleche \Pic(\X_T)\fleche \picet(T) \fleche 0$$
ce qui montre que $i_2\circ i_1 : P_{\X/S}\fleche \picet$ est un
isomorphisme, et donc $i_2$ et $i_1$ sont aussi des isomorphismes.

Si $f$ a une section localement pour la topologie de Zariski, montrons que
$$\xymatrix{\piczar \ar@{^{(}->}[r]^{i_2} &\picet}$$
est un isomorphisme. Comme il
s'agit d'un morphisme entre deux faisceaux pour la topologie de Zariski, la
question est locale pour cette topologie et l'on peut supposer que $f$ a une
section, ce qui nous ram\`ene au cas pr\'ec\'edent.

Il nous reste \`a montrer que $\picet \fleche \picfppf$ est un
isomorphisme.
On a un diagramme commutatif dans lequel les lignes sont les suites exactes de bas degré associées aux suites
spectrales de Leray pour la cohomologie lisse-étale et pour la cohomologie \emph{fppf} :
$$\hskip -5pt \shorthandoff{!;:?}
\xymatrix@!0 @R=3.5pc @C=7pc{ 
H^1(T,\gm) \ar[r] \ar[d] & H^1(\X_T,\gm)\ar[r]\ar[d] & H^0(T, R^1f_{T*}\gm)
\ar[r]\ar[d] & H^2(T,\gm) \ar[r]\ar[d] & H^2(\X_T,\gm)\ar[d] \\
H^1_{\text{pl}}(T,\gm) \ar[r] & H^1_{\text{pl}}(\X_T,\gm)\ar[r] & H^0_{\text{pl}}(T, R^1f_{T*}^{\text{pl}}\gm)
\ar[r] & H^2_{\text{pl}}(T,\gm) \ar[r] & H^2_{\text{pl}}(\X_T,\gm).
}$$
Le diagramme ci-dessus tient compte du fait que, par platitude cohomologique, $f_{T*}\gm=\gm$ et
$f_{T*}^{\text{pl}}\gm=\gm$. Comme $\gm$ est un groupe lisse
sur la base, le théorème (\ref{coh_fppf_groupe_lisse}) nous dit que les deux flèches verticales de gauche
et les deux flèches verticales de droite sont des isomorphismes. D'après le lemme des 5, celle du milieu en est un aussi.
Or d'après (\ref{description_coh_de_picet_et_picfppf}) cette flèche est précisément le morphisme
$\picet(T)\fleche \picfppf(T)$.
\end{demo}

\begin{sousremarque}\rm
On peut démontrer que le morphisme $\picet\fleche \picfppf$ est injectif d'une manière plus directe, sans utiliser
la cohomologie \emph{fppf}. Voici comment.
On note $P$ le foncteur $\picet$. Il faut montrer que si
$(U_i\fleche U)$ est une famille couvrante pour la topologie \emph{(fppf)}, alors le
morphisme $P(U) \fleche \disp \prod_iP(U_i)$ est injectif, autrement dit que $P$ est un
préfaisceau séparé pour la topologie \emph{(fppf)}.
Vu que $P$ est un faisceau pour la topologie de Zariski, il suffit
clairement de traiter le cas d'une famille couvrante à un élément $(U' \fleche
U)$.

Soit $s$ un élément de $P(U)$ dont l'image dans $P(U')$ est nulle.
Considérons dans un premier temps le cas où $s$ provient de $P_{\X/S}(U)$.
Il existe alors un faisceau inversible $\Lc$ sur $\X\times_S U$ dont la classe
dans $P_{\X/S}(U)$ est égale à $s$. Avec les notations du diagramme
ci-dessous,
$$\xymatrix{\X\times_S U' \ar[d]_{f_{U'}} \ar[r]^{v}\cartesien&
\X\times_S U \ar[d]^{f_U}\\
U' \ar[r]^{u}& U}$$
dire que l'image de $s$ dans $P(U')$ est nulle signifie que $v^*(\Lc)$ provient
de la base, autrement dit d'après le critère
(\ref{critere_relative_trivialite}),
que le morphisme d'adjonction $a_{\Lc_{U'}} : f^*_{U'} {f_{U'}}_*\Lc_{U'}
\fleche \Lc_{U'}$ est un isomorphisme. Or on a montré au cours de la
démonstration de ce critère que la formation du morphisme d'adjonction
$a_{\Lc}$ commute au changement de base plat. On en déduit donc que le morphisme
$$v^*a_{\Lc_{U}} : v^*f^*_{U}{f_{U}}_*\Lc_{U}\fleche v^*\Lc_{U}$$
est lui aussi un isomorphisme. Maintenant, comme $v$ est fidèlement plat et
localement de
présentation finie, il résulte de la théorie de la descente fidèlement plate des
modules quasi-cohérents (cf.
\cite{LMB} (13.5)) que $a_{\Lc_{U}}$ est un isomorphisme, ce qui signifie en
vertu du critère (\ref{critere_relative_trivialite}) que $s=0$. Dans le cas
général, il existe une famille couvrante étale $(V\fleche U)$ telle que
$s_{|_V}$ provienne de $P_{\X/S}(V)$. Alors, en notant $V'=V\times_U U'$ on
voit que $(s_{|_V})_{|_{V'}}=(s_{|_{U'}})_{|_{V'}}=0$ ce qui implique la
nullité de $s_{|_V}$ d'après le premier cas puisque $(V'\fleche V)$ est une
famille couvrante \emph{(fppf)}, et donc celle de $s$ puisque $V\fleche U$ est
une famille couvrante étale.
\end{sousremarque}

Par la suite, lorsque nous parlerons \emph{du} foncteur de Picard, il s'agira de
$\picfppf$. Nous le noterons simplement $\pic$.

\section{Le champ de Picard d'un champ algébrique}
\label{champ_de_Picard}

Avec l'apparition des champs, un nouvel objet \og classifiant \fg\ pour les faisceaux inversibles vient s'ajouter aux cinq foncteurs précédents : le champ de Picard. La proposition~(\ref{comparaison_champ_foncteur}) et le corollaire~(\ref{foncteur_pic_repres_implique_champic_alg}) montrent que dans le cas cohomologiquement plat, les questions de l'algébricité du champ de Picard et de la représentabilité du foncteur de Picard sont essentiellement équivalentes.

\begin{defi}
Soit $\X$ un $S$-champ. On appelle \emph{champ de Picard de $\X$},
et on note $\champic(\X/S)$ le champ des faisceaux inversibles. Pour tout $U \in \ob \aff$, la cat\'egorie $\champic(\X/S)_U$ est donc
la cat\'egorie dont les objets sont les faisceaux inversibles sur $\X\times_S U$
et dont les fl\`eches sont les isomorphismes de faisceaux inversibles. Les
changements de base sont d\'efinis de mani\`ere \'evidente.
\end{defi}

\begin{remarque}\rm
Il r\'esulte ais\'ement de (\ref{descente_fidelement_plate}) que le
$S$-groupoïde $\champic(\X/S)$ défini ci-dessus est bien un
$S$-champ.
\end{remarque}

\begin{remarque}\rm
La terminologie adoptée ici, qui semble s'imposer naturellement, entre presque en conflit avec celle de \cite{LMB}~(14.4.2), où l'on appelle \og champ de Picard \fg\ un champ muni d'un morphisme d'addition qui en fait une sorte de \og champ en groupes\fg. Heureusement, si $\X$ est un $S$-champ, \emph{le} champ de Picard de $\X$ défini ci-dessus est bien \emph{un} champ de Picard au
sens de~\cite{LMB}~(14.4.2), le 1-morphisme d'addition
$$+ : \champic(\X/S)\times_S \champic(\X/S) \fleche \champic(\X/S)$$
étant donné par le produit tensoriel de faisceaux inversibles.
\end{remarque}

\begin{remarque}\rm
Il est clair que la formation de $\champic(\X/S)$ commute au changement de base.
Autrement dit, si $S'\fleche S$ est un morphisme de sch\'emas, le champ
$\champic(\X\times_S S'/S')$ sur $({\rm Aff}/S')$ est canoniquement isomorphe au
$S'$-champ $\champic(\X/S)\times_S S'$.
\end{remarque}

\begin{remarque}\rm
En notant $\fHom(\X,\Y)$ le $S$-champ dont la fibre en $U$ est la cat\'egorie
$\Hom(\X\times_SU,\Y\times_SU)$, la proposition
(\ref{fi_equiv_morphisme_de_X_dans_BGm}) montre que
$\champic(\X/S)$ est isomorphe à $\fHom(\X,\bgm)$.
\end{remarque}

Le champ $\champic(\X/S)$ est muni d'un morphisme naturel vers le
$S$-groupo\"\i de
associ\'e au pr\'efaisceau $\Pic_{\X}$, d\'efini sur $\champic(\X/S)_U$
par le foncteur qui envoie un faisceau inversible sur sa classe d'isomorphie
dans $\Pic_{\X}(U)$. On en d\'eduit par composition un
morphisme naturel vers le $S$-espace $P$ (où $P$ désigne $\picet$ ou $\picfppf$
selon l'envie) :
$$\flechen{\champic(\X/S)}{\pi}{P}.$$

\begin{prop}
Le morphisme de $S$-champs $\pi$ ci-dessus est une gerbe (resp. une gerbe
\emph{(fppf)}), autrement dit $\pi$
et le morphisme diagonal $$\Delta :\champic(\X/S) \fleche
\champic(\X/S)\times_P \champic(\X/S)$$ sont des \'epimorphismes.
\end{prop}
\begin{demo}
Soient $U\in \ob\aff$ et $x\in P(U)$. Le $S$-espace $P$ \'etant le faisceau
\'etale (resp. \emph{(fppf)}) associ\'e
\`a $\Pic_{\X}$, il existe une famille couvrante dans $\aff$, que l'on peut
supposer r\'eduite \`a un \'el\'ement $(U'\fleche U)$, et un \'el\'ement $l\in
\Pic_{\X}(U')$, tels que l'image de $l$ dans $P(U')$ soit \'egale \`a $x_{|U'}$.
En d'autres termes il existe un faisceau inversible $\Lc$ dans
$\champic(\X/S)_{U'}$ dont l'image dans $P(U')$ est $x_{|U'}$, ce qui montre que
$\pi$ est un \'epimorphisme.

Par ailleurs, si $U$ est un objet de $\aff$ et si
$\Lc_1, \Lc_2 \in \champic(\X/S)_U$ sont
tels que $\pi(\Lc_1)=\pi(\Lc_2)$, alors il existe une famille couvrante \`a un
\'el\'ement telle que les images de $\Lc_{1|U'}$ et $\Lc_{2|U'}$ dans
$\Pic_{\X}(U')$ soient \'egales, c'est-\`a-dire telle que $\Lc_{1|U'}$ et
$\Lc_{2|U'}$ soient isomorphes, ce qui montre que $\Delta$ est un
\'epimorphisme.
\end{demo}

\begin{prop}
\label{comparaison_champ_foncteur}
Si
$\champic(\X/S)$ est un champ alg\'ebrique (resp. préalgébrique), et si f est cohomologiquement plat
en dimension z\'ero, alors $\pic$ est repr\'esentable par un $S$-espace
alg\'ebrique (resp. préalgébrique), et le 1-morphisme $\pi : \champic(\X/S)\fleche \pic$ est
fid\`element plat et localement de pr\'esentation finie.
\end{prop}
\begin{demo}
D'apr\`es le corollaire (10.8) de \cite{LMB}, il suffit de v\'erifier que pour
tout $U\in\ob\aff$ et pour tout faisceau inversible $\Lc$ sur $\X\times_S U$, le
$U$-espace
alg\'ebrique en groupes $\fIsom(\Lc,\Lc)$ est plat et localement de
pr\'esentation finie. Or $\fIsom(\Lc,\Lc)$ est le faisceau qui \`a tout
$V\in\ob {\rm (Aff/}U)$ associe $\Aut(\Lc_V)$, et $\Aut(\Lc_V)$ s'identifie
canoniquement \`a $\Gamma(\X\times_S V, \O_{\X\times_S V})^{\times}$. Mais par
platitude cohomologique en dimension z\'ero, on a $$\Gamma(\X\times_S V,
\O_{\X\times_S V})^{\times} = \Gamma(V, \O_V)^{\times}=\gm(V).$$ Finalement
$\fIsom(\Lc,\Lc)$ est isomorphe \`a $\gm$, donc il est lisse et de
pr\'esentation finie.
\end{demo}

Réciproquement, on peut essayer de déduire l'algébricité du champ de Picard à partir de celle du foncteur de Picard. Voyons comment.
Le champ des faisceaux inversibles sur $S$, c'est-à-dire le champ dont la catégorie fibre au-dessus d'un $S$-schéma $U$ est la catégorie des faisceaux inversibles sur $U$, s'identifie naturellement au champ $\bgm$ sur $S$. Le foncteur \og image inverse par $f$\fg\ induit donc un morphisme de $S$-champs algébriques
$$\flechen{\bgm}{f^*}{\champic(\X/S).}$$
De plus ce morphisme est un monomorphisme (i.e. les foncteurs $f^*_U$ sont pleinement fidèles) dès que $f$ est cohomologiquement plat en dimension zéro. On a donc dans ce cas une \og suite exacte \fg\ de champs de Picard :
$$\xymatrix{0 \ar[r]& \bgm \ar[r]^-{f^*} &\champic(\X/S) \ar[r]^-{\pi}& \pic \ar[r]&0}.$$
Autrement dit $f^*$ est un monomorphisme, $\pi$ est un épimorphisme, un objet de $\champic(\X/S)$ est envoyé sur $0$ par $\pi$ si et seulement s'il provient de $\bgm$, et tous les morphismes sont compatibles aux morphismes d'addition des champs de Picard considérés\footnote{Si $f$ n'est pas cohomologiquement plat en dimension zéro, $f^*$ n'est plus un monomorphisme, mais le reste est presque vrai : il faut juste remplacer \og provient \fg\ par \og provient localement pour la topologie étale\fg.}. La proposition ci-dessous nous permet de déterminer quand cette suite exacte est scindée.

\begin{prop}
\label{comparaison_champ_foncteur2}
Soit $\X$ un $S$-champ algébrique dont le morphisme structural $f$ est cohomologiquement plat en dimension zéro.
\begin{itemize}
\item[(i)] Les propositions suivantes sont équivalentes :
\begin{itemize}
\item[a)] La suite exacte ci-dessus est scindée, autrement dit il existe un morphisme $s$ de $\pic$ dans $\champic(\X/S)$ tel que $\pi\circ s$ soit égal à l'identité.
\item[b)] Il existe un isomorphisme
$$\flechen{\bgm\times_S \pic}{\sim}{\champic(\X/S)}$$
compatible avec les projections $\pr_2$ et $\pi$ sur $\pic$.
\item[c)] Il existe sur le $S$-champ (non nécessairement algébrique) $\X\times_S\pic$ un faisceau inversible \og universel\fg\ $\Pc$  qui représente le foncteur $\pic$, i.e. tel que pour tout $U\in\ob\aff$ et tout élément $l$ de $\pic(U)$ on ait $l=[\Pc_{|_{\X\times_S U}}]=\pi(\Pc_{|_{\X\times_S U}})$.
$$\xymatrix{\X\times_S U \ar[r] \ar[d] \cartesien &
\X\times_S \pic \ar[r] \ar[d] \cartesien &\X \ar[d] \\
U \ar[r]_l & \pic \ar[r]& S}$$
\end{itemize}
\item[(ii)] Les propositions suivantes sont équivalentes :
\begin{itemize}
\item[d)] Le morphisme naturel de $P_{\X/S}$ dans $\pic$ est un isomorphisme.
\item[e)] Il existe un élément de $P_{\X/S}(\pic)$ (c'est-à-dire $\Hom(\pic,P_{\X/S})$ lorsque $\pic$ n'est pas représentable) qui a pour image l'identité dans $\pic(\pic)$ (ou, ce qui revient au même, $P_{\X/S}\fleche \pic$ a une section).
\end{itemize}
\item[(iii)] Les conditions de (i) impliquent celles de (ii), et la réciproque est vraie si $\pic$ est représentable par un schéma. Toutes ces conditions sont vérifiées si le morphisme structural $\X \fleche S$ a une section.
\end{itemize}
\end{prop}
\begin{demo}
Montrons que a) implique b). On suppose que $\pi$ a une section $s$. On vérifie alors facilement que le morphisme composé 
$$\xymatrix{\bgm \times_S \pic \ar[r]^-{(f^*,s)}&
\champic(\X/S)\times_S \champic(\X/S) \ar[r]^-{\otimes} & \champic(\X/S)}$$
est un isomorphisme. (On utilise ici la platitude cohomologique en dimension zéro.) Un quasi-inverse est donné par le foncteur de $\champic(\X/S)$ vers $\bgm \times_S \pic$ qui à un faisceau inversible $\Mc$ associe le couple
$(\Mc\otimes s(\pi(\Mc))^{-1},\pi(\Mc))$, en identifiant $\bgm$ à son image essentielle dans $\champic(\X/S)$. 

La réciproque est claire : un isomorphisme comme dans b) induit une section de $\pi$ par composition avec la section évidente de $\pr_2$.

Montrons maintenant que a) implique c)\footnote{C'est évident si $\pic$ est représentable par un schéma ! Ce que l'on ne suppose pas ici\dots}. On se donne une section $s$ de $\pi$ et l'on va construire un faisceau inversible $\Pc$ (au sens de (\ref{definition_fi})) sur le champ $\X\times_S \pic$. Il faut donc construire pour tout $S$-schéma affine $U$ et tout objet $x$ de $(\X\times_S \pic)_U$ un faisceau inversible $\Pc(x)$ sur $U$, et des isomorphismes de transition entre les $\Pc(x)$. Soit $x$ un objet de $(\X\times_S \pic)_U$. On note $l=f_P\circ x$ et $t_x$ la section de $f_U$ induite par $x$, comme dans le diagramme ci-dessous.
$$\xymatrix{\X\times_S U \ar[r] \ar[d]^{f_U} &
\X\times_S \pic \ar[r] \ar[d]^{f_P} \cartesien &\X \ar[d]^f \\
U \ar[r]_l \ar[ru]_x \ar@/^6mm/[u]^{t_x}& \pic \ar[r]& S}$$
On a un faisceau inversible $s(l)$ sur $\X\times_S U$. On pose $\Pc(x)=t_x^*s(l)$. Les isomorphismes de transition sont définis de manière évidente et l'on obtient ainsi un faisceau inversible $\Pc$ sur $\X\times_S \pic$ dont il ne reste plus qu'à montrer qu'il représente $\pic$. Soit $l$ un élément de $\pic(U)$. Il faut montrer que
$l=\pi(\Pc_{|_{\X\times_S U}})$. Il est évident, vu la construction de $\Pc$, que $\Pc_{|_{\X\times_S U}}$ est isomorphe à $s(l)$. Comme $s$ est une section de $\pi$, on a bien le résultat attendu.

Réciproquement on suppose donné un faisceau inversible universel $\Pc$. Alors le morphisme qui à $l$ associe $\Pc_{|_{\X\times_S U}}$ définit une section de $\pi$.

a) $\Rightarrow$ e) est évident puisque $\pi$ se factorise par $P_{\X/S}$.

L'équivalence entre d) et e) est évidente puisque le morphisme en question est de toute manière injectif par (\ref{comparaison_des_foncteurs_de_Picard}).

Supposons que $\pic$ soit représentable par un schéma et montrons que e) implique c). Sous l'hypothèse e), il existe par définition du foncteur $P_{\X/S}$ un faisceau inversible $\Pc$ sur $\X\times_S \pic$ qui induit l'élément identité de $\pic(\pic)$. Il est clair que $\Pc$ est universel.

Enfin dans le cas où $f$ a une section $\sigma$, montrons que la condition b) est vérifiée. Le morphisme $f^*$ admet dans ce cas une rétraction $\sigma^*$.
$$\xymatrix{0 \ar[r]& \bgm \ar[r]^-{f^*} &\champic(\X/S) \ar@/^5mm/[l]^{\sigma^*}\ar[r]^-{\pi}& \pic \ar[r]&0}.$$
On vérifie alors facilement que le 1-morphisme
$$\fonction{(\sigma^*,\pi)}{\champic(\X/S)}%
{\bgm\times_S \pic}{\Mc}{(\sigma^*\Mc,\pi(\Mc))}$$
est un isomorphisme.
\end{demo}

\begin{remarque}\rm
On se donne à la fois une section $\sigma$ de $f$ et un faisceau inversible universel $\Pc$ sur $\X\times_S \pic$ (on ne suppose pas $\pic$ représentable). Alors l'isomorphisme $\bgm\times_S \pic \fleche \champic(\X/S)$ induit par $\Pc$ est un quasi-inverse de $(\sigma^*, \pi)$ si et seulement si $\sigma^*\Pc$ est trivial. En particulier on obtient toujours un quasi-inverse en \og rigidifiant \fg\ $\Pc$ le long de $\sigma$, c'est-à-dire en remplaçant $\Pc$ par $\Pc\otimes (f^*\sigma^*\Pc)^{-1}$.
\end{remarque}

\begin{cor}
\label{foncteur_pic_repres_implique_champic_alg}
Soit $\X$ un $S$-champ algébrique dont le morphisme structural $f$ est cohomologiquement plat en dimension zéro. On suppose que le foncteur de Picard $\pic$ est représentable par un espace algébrique et que $f$ a une section localement pour la topologie \emph{fppf} sur $S$. Alors le champ de Picard $\champic(\X/S)$ est algébrique.
\end{cor}
\begin{demo}
Pour un $S$-champ, être algébrique est une condition locale pour la topologie \emph{fppf}. On peut donc supposer que $f$ a une section et il suffit d'appliquer la proposition précédente.
\end{demo}

% \begin{exemple}\rm
% Attention cependant, le champ de Picard n'est pas toujours isomorphe à $\bgm\times_S \pic$. Par exemple si $\X$ est le $\R$-schéma ... alors $\pic$ est isomorphe au schéma constant $\Z$ (voir ...). Il ne peut pas y avoir de faisceau inversible universel sur $\X \times_S \Z$ car ...
% \end{exemple}

\section{Propriétés de finitude relative}

Nous terminons cette première partie en montrant que le foncteur de Picard commute aux limites inductives (i.e. il est de présentation finie) et aux limites projectives. Ces deux propriétés nous seront utiles pour étudier sa représentabilité.

\subsection{Pr\'esentation finie}

\begin{sousprop}
\label{lpf}
Soient $S$ un sch\'ema affine, et $\X$ un $S$-champ alg\'ebrique quasi-compact.
Alors le foncteur de Picard relatif $P_{\X/S}$ d\'efini ci-dessus est
localement de pr\'esentation finie. 
Si de plus le morphisme structural $f : \X \fleche S$ est cohomologiquement plat
en dimension zéro et admet une section
localement pour la topologie (fppf), alors le foncteur de Picard relatif
$\pic$ est localement de pr\'esentation finie. 
\end{sousprop}
\begin{demo}
Il s'agit de montrer que $P_{\X/S}$ commute aux limites inductives filtrantes
d'anneaux. Autrement dit, on veut montrer que pour tout ensemble inductif
filtrant $I$, et pour tout syst\`eme inductif de $S$-anneaux $(A_i)_{i\in I}$
index\'e par $I$, le morphisme canonique
$$\lind P_{\X/S}(A_i) \flechelongue P_{\X/S}(A)$$
o\`u $\disp A=\lind A_i$ est la limite inductive des $A_i$, est un isomorphisme.
Rappelons que pour tout $S$-sch\'ema $S'$, on a
$P_{\X/S}(S')=\frac{\Pic(\X \times_S S')}{\Pic(S')}.$ Le foncteur
$\disp \lind$ \'etant exact, il suffit donc de montrer que $\Pic_{\X}$
lui-m\^eme est localement de pr\'esentation finie. On peut d\'ecrire la
cat\'egorie 
$\disp \lind\inv(\X\otimes_S A_i)$ de la mani\`ere suivante. Les objets sont les
couples $(i,\Lc_i)$ o\`u $i\in I$ et $\Lc_i$ est un faisceau inversible sur
$\X\otimes_S A_i$. Pour deux tels objets $(i,\Lc_i)$ et $(j,\Lc_j)$,
$\Hom((i,\Lc_i),(j,\Lc_j))$ est
$\disp
\lind \Hom((\Lc_i)_{|_{(\X\otimes_S A_l)}},(\Lc_j)_{|_{(\X\otimes_S A_l)}})$
o\`u la limite inductive est prise sur l'ensemble (non vide) des indices $l$ qui
majorent $i$ et $j$. Pour montrer que $\Pic_{\X}$ est localement de
pr\'esentation finie, il suffit de montrer que le foncteur naturel 
$$\lind\inv(\X\otimes_S A_i) \flechelongue \inv(\X\otimes_S A)$$ est une \'equivalence
de cat\'egories. Ceci r\'esulte de la proposition (4.18) de \cite{LMB},
compte tenu de (\ref{fi_equiv_morphisme_de_X_dans_BGm}) et du fait que le champ
$\bgm$ est localement de pr\'esentation finie.

Supposons maintenant que $f : \X \fleche S$ est cohomologiquement plat
en dimension zéro et a une section localement pour la
topologie \emph{(fppf)}. Dans le cas o\`u $f$ a une section, on sait
d'apr\`es (\ref{comparaison_des_foncteurs_de_Picard}) que $\pic$ s'identifie \`a 
$P_{\X/S}$, dont on vient de montrer qu'il est localement de pr\'esentation
finie. Soit $S'\fleche S$ un morphisme fid\`element plat et localement de
pr\'esentation finie, tel que le morphisme $f' : \X\times_S S' \fleche S'$
obtenu par changement de base ait une section.
Comme la formation de $\pic$ commute au changement de base, on en
d\'eduit que le foncteur $\pic\times_S S'$ est localement de pr\'esentation finie,
ce qui montre, en utilisant le fait que $\pic$ est un faisceau pour la topologie
\emph{(fppf)} et l'exactitude du foncteur $\disp \lind$, que $\pic$ lui-m\^eme
est localement de pr\'esentation finie.
\end{demo}

\begin{sousremarque}\rm
L'hypoth\`ese de quasi-compacit\'e sur $\X$ est r\'eellement n\'ecessaire,
m\^eme lorsque $\X$ est un sch\'ema, comme
le montre l'exemple suivant. Soit $E$ une courbe elliptique sur $\Q$. On prend
$\X=\coprod_{n\in \N} E_n$, o\`u pour tout $n$, $E_n=E$. Alors on a
$\underline{\Pic}^0_{E/\Q}\simeq E$, donc $\underline{\Pic}^0_{\X/\Q}\simeq
E^{\N}$. En particulier pour toute extension $L$ de $\Q$ on a
$\underline{\Pic}^0_{\X/\Q}(L)=(E(L))^{\N}$.
On peut voir que le sch\'ema de Picard ainsi obtenu n'est pas localement de
pr\'esentation finie en regardant par exemple les points \`a valeurs dans
$\ov{\Q}=\disp \lind K$, o\`u la limite inductive est prise sur l'ensemble des
corps de nombres $K$ sur $\Q$. La fl\`eche canonique
$\disp \lind \underline{\Pic}^0_{\X/\Q}(K)\fleche
\underline{\Pic}^0_{\X/\Q}(\ov{\Q})$ s'identifie \`a :
$$\disp \lind E(K)^{\N} \flechelongue E(\ov{\Q})^{\N}.$$
Cette application n'est clairement pas bijective puisque $E$ a des $K$-points
pour des corps de nombres $K$ de degr\'e aussi grand que l'on veut.
\end{sousremarque}

\subsection{Commutation aux limites projectives}

\begin{sousprop}
\label{lim_proj}
Soit $\X$ un $S$-champ alg\'ebrique propre et
cohomologiquement plat en dimension z\'ero.
Soient $\ov{A}$ un anneau noeth\'erien local complet sur
$S$ et $\mgo$ son id\'eal maximal. 
Alors le morphisme canonique
$$P_{\X/S}(\ov{A}) \flechelongue \lpro P_{\X/S}(\ov{A}/\mgo^n)$$
est un isomorphisme.
\end{sousprop}
\begin{demo}
On peut supposer $S=\Spec \ov{A}$. On note
$S_n= \Spec(\ov{A}/\mgo^{n+1})$, et $\X_n=\X\times_S S_n$.
$$\xymatrix{
\X_0 \ar[r] \ar[d] & \dots \ar[r] & \X_n \ar[r] \ar[d] &
	\X_{n+1} \ar[r] \ar[d] & \dots \ar[r] & \X \ar[d] \\
S_0 \ar[r] & \dots \ar[r] & S_n \ar[r] &
	S_{n+1} \ar[r] & \dots \ar[r] & S
}$$
Comme $S$ et tous les
$S_n$ ont des groupes de Picard triviaux, nous avons \`a montrer que 
$$\Pic(\X) \flechelongue \lpro \Pic(\X_n)$$
est un isomorphisme.\\
\emph{Injectivit\'e :}
Soient $\Lc$ un faisceau inversible sur $\X$, et $\Lc_n$ les faisceaux induits
par $\Lc$ sur les $\X_n$. On suppose que $\Lc$ est trivial sur
chacun des $\X_n$ et on veut montrer qu'il est trivial. D'apr\`es le
théorème~(1.4) de \cite{Olsson_lemme_chow},
on a une \'equivalence de cat\'egories
$$\xymatrix{(\textrm{faisceaux coh\'erents sur }\X) \ar[d]\\
(\textrm{familles compatibles de faisceaux coh\'erents sur les }\X_n).}$$
Cette \'equivalence de cat\'egories induit un foncteur
pleinement fid\`ele :
$$\xymatrix{(\textrm{faisceaux inversibles sur }\X) \ar[d]\\
(\textrm{familles compatibles de faisceaux inversibles sur les }\X_n).}$$
Il suffit donc de montrer que la famille compatible de faisceaux inversibles
$(\Lc_n)_{n\geq 0}$ est isomorphe \`a la famille triviale, autrement dit que
l'on a une famille compatible d'isomorphismes
$\xymatrix@C=1pc{\O_{\X_n} \ar[r]^{\sim} & \Lc_n}$. Par
hypoth\`ese on sait que pour tout $n$ il existe un tel isomorphisme, mais la
famille ainsi obtenue n'est \emph{a priori} pas compatible au changement de base. Deux
isomorphismes $\xymatrix@C=1pc{\O_{\X_n} \ar[r]^{\sim}& \Lc_n}$ diff\`erent par un
automorphisme de
$\O_{\X_n}$, donc pour conclure il suffit de montrer que
$$\Aut(\O_{\X_n}) \flechelongue \Aut(\O_{\X_{n-1}})$$
est surjective. Or par platitude cohomologique,
$\Aut(\O_{\X_n})=\Gamma(\X_n, \O_{\X_n})^{\times}$ n'est autre que
$(\ov{A}/\mgo^{n+1})^{\times}$, et le morphisme
$(\ov{A}/\mgo^{n+1})^{\times} \fleche (\ov{A}/\mgo^{n})^{\times}$ est clairement
surjectif.\\
\emph{Surjectivit\'e :}
Soit $(\Lc_n)_{n\geq 0}$ une famille compatible de faisceaux inversibles sur les
$\X_n$. Il nous faut montrer qu'il existe un faisceau inversible sur $\X$, qui
induit, \`a isomorphisme pr\`es, la famille $(\Lc_n)_{n\geq 0}$. Il existe en
vertu de l'\'equivalence de cat\'egories mentionn\'ee ci-dessus un faisceau
quasi-coh\'erent $\L$ sur $\X$ ayant cette propri\'et\'e, et il suffit de
montrer que $\L$ est un faisceau inversible. On construit un inverse en
consid\'erant un faisceau cohérent $\Mc$ qui induit la famille compatible
$(\Lc_n^{-1})_{n\geq 0}$. Pour montrer que
$\Mc\otimes \Lc$ est isomorphe à $\O_{\X}$ il suffit de
montrer que la famille induite par $\Mc\otimes \Lc$ est isomorphe \`a la famille
triviale. Mais ceci est \'evident puisque que l'on obtient un syst\`eme
compatible d'isomorphismes en prenant pour tout $n$ l'isomorphisme :
$$(\Mc\otimes \Lc)_{|_{\X_n}}\simeq (\Mc_{|_{\X_n}})\otimes(\Lc_{|_{\X_n}})
\simeq \L_n^{-1}\otimes \Lc_n \simeq \O_{\X_n}.$$
\vskip-0.4cm\noindent
\end{demo}

\chapter{Propri\'et\'es de s\'eparation}

\section{Pr\'eliminaires}

\subsection{Critère valuatif de locale séparation d'Artin}

Artin a caché dans les méandres techniques de son théorème de représentabilité pour les espaces algébriques un très joli critère valuatif permettant de caractériser les immersions quasi-compactes de schémas\footnote{ou plutôt, ce qui revient essentiellement au même, un critère de locale séparation pour les morphismes de schémas} (cf. \cite{Global_Analysis_1}~pp.~58-59). Nous nous proposons d'en donner ci-dessous un énoncé précis assorti d'une démonstration autonome (thm.~\ref{critere_immersionsqc}). Au passage, nous le généralisons aux morphismes représentables de champs algébriques.

Si un morphisme $f : Y \fleche X$ est une immersion, il doit vérifier la propriété suivante :
\begin{gather}
\forall y_0, y_1 \in Y, \ y_0\in \overline{\{y_1\}} \Longleftrightarrow
f(y_0)\in \overline{\{f(y_1)\}} \tag{$*$}
\end{gather}
En particulier, pour tout point $y$ de $Y$ d'image $x$ dans $X$ et tout plongement d'un petit morceau de courbe $S$ passant par $x$ dans $X$, le point $y$ doit être dans l'adhérence de tout relèvement (unique, s'il existe !) de $S\setminus\{x\}$ à $Y$.

\begin{center}
\includegraphics[scale=0.3]{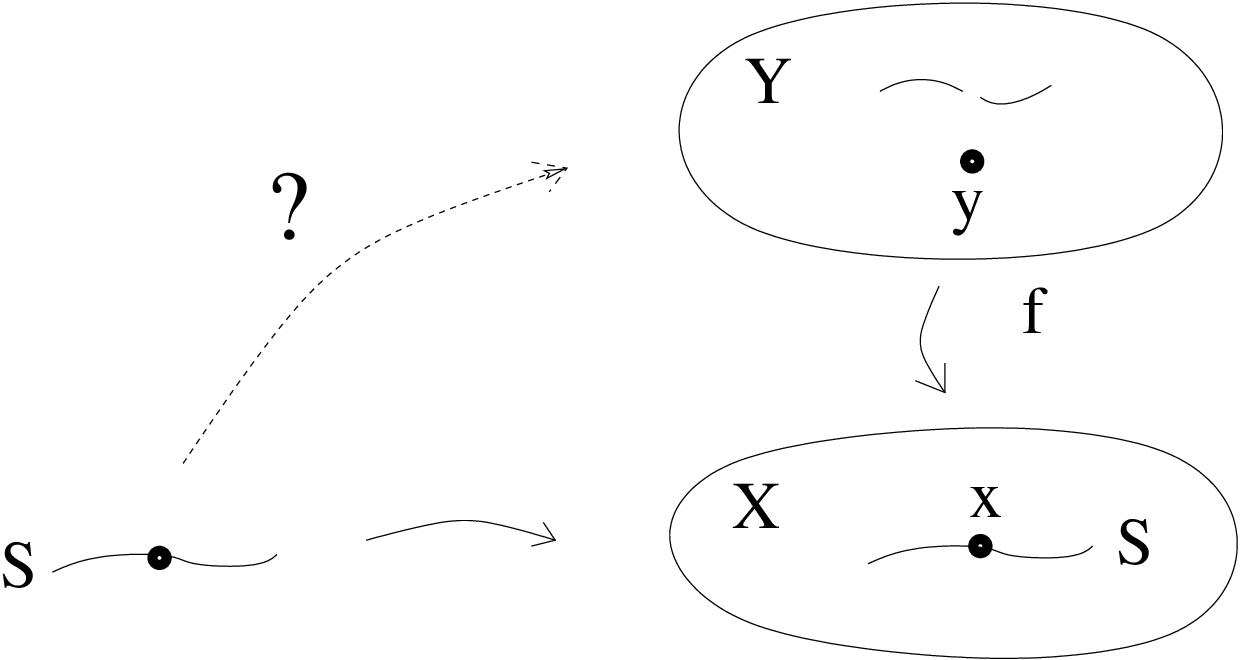}\\
\medskip
\end{center}

Le lemme suivant dit essentiellement que cette dernière propriété est équivalente à la propriété ($*$).

%Le but de cette section est de donner un crit\`ere facile \`a v\'erifier pour
%caract\'eriser les immersions quasi-compactes. Nous appliquerons 
%ult\'erieurement ces r\'esultats \`a la diagonale du champ de Picard (ou du
%foncteur de Picard) pour en d\'eduire qu'ils sont localement s\'epar\'es (donc
%quasi-s\'epar\'es).
\begin{souslem}
Soit $f : Y\fleche X$ un monomorphisme de sch\'emas.
Les propri\'et\'es
suivantes sont \'equivalentes.
\begin{itemize}
\item[(1)] Le morphisme $f$ v\'erifie ($*$) universellement.
\item[(1 bis)] Le morphisme $f$ v\'erifie ($*$) apr\`es tout changement de base
localement de type fini.
\item[(2)] Soient $A$ un anneau de valuation, $k$ son corps r\'esiduel, $K$ son corps
des fractions, $S=\Spec A$, $s_g=\Spec K$, $s_f=\Spec k$, et $S\fleche X$ un
morphisme, tels que $s_g$ et $s_f$ se rel\`event \`a $Y$. Alors $S$ se rel\`eve
\`a $Y$. En d'autres termes, tout diagramme commutatif en traits pleins :
$$\xymatrix{&&&Y \ar[d] \\ s_g\ar[urrr] \ar@/_1.3pc/[rr] & s_f\ar[r] \ar[urr] &
     S \ar[r] \ar@{.>}[ur]&X,}$$ \vskip0.2cm
\noindent se prolonge par un morphisme de $S$ vers $Y$ en un diagramme
commutatif.
\end{itemize}

Si de plus X est localement noeth\'erien, alors ces conditions sont
encore \'equivalentes \`a :
\begin{itemize}
\item[(2 bis)] m\^eme condition que (2) mais en se restreignant \`a des anneaux de
valuation discr\`ete. 
\end{itemize}
\end{souslem}
\begin{demo}

$\underline{(\text{1 bis}) \Longrightarrow (1) :}$
Soit $\varphi : X' \fleche X$ un morphisme de changement de base. On consid\`ere
le carr\'e cart\'esien :
$$\xymatrix{Y'\ar[d]_{f'} \ar[r]^{\varphi'} \cartesien&Y \ar[d]^f \\
     X' \ar[r]^{\varphi} &X.}$$
Soient $y_0', y_1' \in Y'$. On note $x_0', x_1'$ (resp. $y_0, y_1, x_0, x_1$)
leurs images dans $X'$ (resp. $Y$, $X$). On suppose que $x_0'\in \ov{\{x_1'\}}$
et il faut montrer que $y_0'\in \ov{\{y_1'\}}$. Soit $U$ un ouvert affine
de $X$ contenant $x_0$. Il contient aussi $x_1$. En notant $U'=\varphi^{-1}(U)$
et $V'=f'^{-1}(U')$, on voit qu'il suffit de montrer que
$f'_{|V'} : V'\fleche U'$ v\'erifie $(*)$. Or $V'=Y\times_X U'=f^{-1}(U)\times_U
U'$. La condition (1 bis) \'etant clairement stable par changement
de base par une immersion ouverte, on peut supposer $X$ affine. De
m\^eme, le fait qu'une immersion ouverte v\'erifie (1) (donc (1 bis)) montre que
l'on peut supposer $Y$ affine. Il est clair que l'on peut aussi supposer $X'$
affine. On \'ecrit $Y=\Spec B$, $X=\Spec A$, $X'=\Spec A'$ et $\lambda :
A'\fleche A'\otimes_A B$. Soient $\pgo,\qgo \in \Spec(A'\otimes_A B)$ tels que
$\lambda^{-1}(\qgo)\subset \lambda^{-1}(\pgo)$. Il s'agit de montrer que $\qgo
\subset \pgo$. Soit $x=\sum_i a_i'\otimes b_i \in \qgo$. On note $A''$ la
sous-$A$-alg\`ebre de $A'$ engendr\'ee par les $a_i'$.
$$\xymatrix{A'\otimes_A B & A''\otimes_A B \ar[l]_{\gamma} & B\ar[l]\\
A' \ar[u]^{\lambda}& A'' \ar[u]_{\psi} \ar[l] & A \ar[l] \ar[u]}$$
Alors dans $A''\otimes_A B$, l'\'el\'ement $y$ d\'efini par $\sum_i a_i'\otimes
b_i$ appartient \`a $\gamma^{-1}(\qgo)$. Par ailleurs, comme $A''$ est de type
fini sur $A$, l'inclusion $\psi^{-1}(\gamma^{-1}(\qgo))\subset
\psi^{-1}(\gamma^{-1}(\pgo))$ implique
$$\gamma^{-1}(\qgo)\subset\gamma^{-1}(\pgo)~.$$ Donc $x=\gamma(y)\in \pgo$.

$\underline{(1) \Longrightarrow (2) :}$ Avec les notations ci-dessus, les
donn\'ees de la condition (2) induisent un diagramme
$$\xymatrix{&& Y_S\ar[d]^{f_S} \ar[r] &Y \ar[d]^f \\
     s_g\ar[urr] \ar@/_1.3pc/[rr] & s_f\ar[r] \ar[ur] &
     S \ar[r] &X,}$$\vskip0.2cm
\noindent 
o\`u $Y_S=Y\times_S X$. On note $y_f$, $y_g$ les images de $s_f$, $s_g$ dans
$Y_S$, et $Z=\overline{\{y_g\}}$ muni de sa structure de sous-sch\'ema ferm\'e
r\'eduit. On a $f_S(y_f)=s_f$ et $f_S(y_g)=s_g$, donc la condition (1) implique
que $y_f\in Z$. On note alors $\O=\O_{Z,y_f}$ l'anneau local de $Z$ en $y_f$.
Le morphisme $\Spec \O \fleche Z \fleche S$
correspond \`a un morphisme local d'anneaux locaux $R\fleche \O$, qui fait
commuter le diagramme suivant :
$$\xymatrix{R \ar[r] \incl[d] & \O\incl[d] \\ K \ar[r] &\kappa(y_g).}$$
Par ailleurs comme $s_g$ se rel\`eve on a $\kappa(y_g)=\kappa(s_g)=K$. Donc $\O$
est un anneau local qui domine $R$ et par maximalit\'e de $R$ pour la relation
de domination, $R=\O$, d'o\`u un morphisme $S\fleche \Spec \O \fleche Z \subset
Y_S \fleche Y$ qui a les propri\'et\'es escompt\'ees.

$\underline{(2) \Longrightarrow (1) :}$
La condition (2) est stable par changement de base donc on peut oublier
l'adverbe de la condition (1). Soient $y_0, y_1\in Y$ tels que
$f(y_0)\in \overline{\{f(y_1)\}}$. On
note $x_0=f(y_0)$, $x_1=f(y_1)$, $Z=\overline{\{x_1\}}$ muni de sa structure de
sous-sch\'ema ferm\'e r\'eduit, et $\O=\O_{Z,x_0}$ l'anneau local de $Z$ en
$x_0$. $\O\subset \kappa(x_1) \subset \kappa(y_1)=K$ est un anneau local. Il
existe un anneau de valuation $R$ de $K$ qui domine $\O$. Autrement dit on a un
diagramme commutatif :
$$\xymatrix{\O \ar@{^{(}->}[r]^i \incl[d] & R\incl[d] \\
\kappa(x_1) \ar@{^{(}->}[r] & K=\kappa(y_1)}$$
et $i$ est un morphisme local d'anneaux locaux. D'o\`u un morphisme de $S=\Spec
R$ vers $Z$ qui envoie le point g\'en\'erique $s_g$ de $S$ sur $x_1$ et le point
ferm\'e $s_f$ sur $x_0$. Comme $K=\kappa(y_1)$, $s_g$ se rel\`eve clairement en
un morphisme $s_g \fleche Y$ d'image $y_1$. Par ailleurs comme $Y\fleche X$ est
un monomorphisme, l'extension r\'esiduelle $\kappa(x_0) \fleche \kappa(y_0)$ est
triviale de sorte que le point ferm\'e se rel\`eve aussi. On obtient donc un
diagramme commutatif :
$$\xymatrix{&&&&Y \ar[d] \\
s_g\ar[urrrr] \ar@/_1.3pc/[rr] & s_f\ar[r] \ar[urrr] &
     S \ar[r] &Z \ar@{^{(}->}[r] &X,}$$
     
\noindent
et la condition (2) nous donne un morphisme $S\fleche Y$ qui envoie $s_g$ sur
$y_1$ et $s_f$ sur $y_0$, ce qui montre que $y_0\in \overline{\{y_1\}}$.

\medskip
$\underline{(\text{2 bis}) \Longrightarrow (\text{1 bis}) :}$ Ici on suppose de
plus que $X$ est localement noeth\'erien. Alors dans le raisonnement
pr\'ec\'edent $\O$ est un anneau local int\`egre \emph{noeth\'erien} de corps
des fractions $\kappa(x_1)$, et l'extension $\kappa(x_1) \fleche K=\kappa(y_1)$
est triviale puisque $f$ est un monomorphisme, donc l'anneau $R$ peut \^etre
choisi de valuation discr\`ete, et la suite du raisonnement est encore valable.
\end{demo}

Le théorème ci-dessous reprend les idées d'Artin présentes pp.~58-59 de~\cite{Global_Analysis_1}.

\begin{sousthm}
\label{critere_immersionsqc}
Soit $f : \Y \fleche \X$ un morphisme
repr\'esentable et localement de pr\'esentation finie de champs
alg\'ebriques. On suppose de plus que $\X$ est localement noeth\'erien.
Alors $f$ est une immersion quasi-compacte si et seulement si les
conditions suivantes sont v\'erifi\'ees :

(a) $f$ est un monomorphisme.

(b) Soient $A$ un anneau de valuation (resp. un anneau de valuation discr\`ete),
$k$ son corps r\'esiduel, $K$ son corps
des fractions, $S=\Spec A$, $s_g=\Spec K$, $s_f=\Spec k$, et $S\fleche \X$ un
morphisme, tels que $s_g$ et $s_f$ se rel\`event \`a $\Y$. Alors $S$ se rel\`eve
\`a $\Y$. Autrement dit, tout diagramme 2-commutatif en traits pleins :
$$\xymatrix{&&&\Y \ar[d] \\ s_g\ar[urrr] \ar@/_1.3pc/[rr] & s_f\ar[r] \ar[urr] &
     S \ar[r] \ar@{.>}[ur]&\X.}$$
\noindent se prolonge en un diagramme 2-commutatif par la fl\`eche en
pointill\'es.

(c) Soit $A$ un anneau int\`egre. Pour tout morphisme $\Spec A
\fleche \X$ tel qu'il existe un ensemble dense $\Sc$ de points de $\Spec A$ qui
se rel\`event \`a $\Y$, il existe un ouvert non vide $U$ de $\Spec A$ qui se
rel\`eve \`a $\Y$ :
$$\xymatrix{&&&\Y \ar[d] \\ s\in \Sc\ar[urrr] \ar@/_1.5pc/[rr] &  U\ar[r]
	\ar@{.>}[urr] & \Spec A \ar[r] &\X.}$$
\end{sousthm}
\begin{sousremarque}\rm
La d\'emonstration ci-dessous montre aussi que si les conditions (a) et (c) sont
v\'erifi\'ees, alors $f$ est un morphisme quasi-compact. Cependant, la condition
(c) seule ne suffit pas pour assurer la quasi-compacité de $f$. En effet, elle
est vérifiée par exemple dès que $f$ a une section.
\end{sousremarque}
\begin{demo}
Supposons que $f$ est une immersion quasi-compacte. Alors $f$ est un
monomorphisme, et il est de plus sch\'ematique. Avec les notations de (b), en
posant $Y_S=\Y\times_{\X} S$, on a un diagramme commutatif
$$\xymatrix{&& Y_S\ar[d]^{f_S} \\
     s_g\ar[urr] \ar@/_1.3pc/[rr] & s_f\ar[r] \ar[ur] &S,}$$
o\`u $f_S : Y_S \fleche S$ est une immersion quasi-compacte de sch\'emas, et il
s'agit de montrer que $f_S$ a une section. D'apr\`es le lemme
pr\'ec\'edent, il suffit de montrer qu'une immersion quasi-compacte $\varphi
:U\fleche V$ v\'erifie la condition 
$$\forall u_0, u_1 \in U \ u_0\in \overline{\{u_1\}} \Longleftrightarrow
\varphi(u_0)\in \overline{\{\varphi(u_1)\}},$$
ce qui r\'esulte trivialement du fait qu'un tel morphisme induit un
hom\'eomorphisme sur une partie localement ferm\'ee. Montrons la
condition (c). On se donne un schéma affine int\`egre $S=\Spec A$ et un
morphisme $S\fleche \X$ comme dans (c).
En notant $Y_S=\Y\times_{\X} S$, le morphisme $f_S : Y_S \fleche S$ est une
immersion quasi-compacte de sch\'emas, et par hypoth\`ese son image contient un
ensemble dense de $S$. Son image est donc un ouvert $U$ de $S$ et il se
factorise en $Y_S \fleche U \fleche S$ o\`u $U\fleche S$ est une immersion
ouverte, et $Y_S \fleche U$ une immersion ferm\'ee surjective. Mais comme $U$
est r\'eduit, cette derni\`ere fl\`eche est un isomorphisme. On obtient donc
un morphisme $U\fleche Y_S$ qui par composition avec $Y_S \fleche \Y$ donne le
rel\`evement souhait\'e.

R\'eciproquement, supposons les conditions (a), (b), et (c) satisfaites et
montrons que $\Y\fleche \X$ est une immersion quasi-compacte. Comme $f$ est un
monomorphisme localement de pr\'esentation finie, il est s\'epar\'e et
localement quasi-fini, donc il est sch\'ematique d'apr\`es \cite{LMB}, (A.2).
En vertu de \cite{EGA4_2} (2.7.1), il suffit de montrer le r\'esultat apr\`es
changement de base par une pr\'esentation lisse $X\fleche \X$ de $\X$. On
choisit une pr\'esentation par un sch\'ema $X$ (localement noeth\'erien). Alors
$Y=\Y\times_{\X} X$ est lui aussi un sch\'ema, et le morphisme $Y\fleche X$
obtenu par changement de base v\'erifie les m\^emes hypoth\`eses que $f$ (il
sera encore not\'e $f$ par la suite). \^Etre une immersion quasi-compacte est
une propri\'et\'e locale \`a la base, donc on peut supposer que $X$ est un
sch\'ema affine (noeth\'erien).

\begin{etape}{Quasi-compacit\'e de $Y$.}
\indent Pour montrer que $Y$ est quasi-compact, nous allons montrer par r\'ecurrence
noeth\'erienne que pour tout ferm\'e $F$ de $X$, $f^{-1}(F)$ est un espace
topologique noeth\'erien. Il faut montrer que si $A$ est un ferm\'e de $X$ tel
que pour tout ferm\'e strict $B$ de $A$, $f^{-1}(B)$ soit noeth\'erien, alors
$f^{-1}(A)$ est noeth\'erien. Or le morphisme $f^{-1}(A) \fleche A$ obtenu par
changement de base v\'erifie les m\^emes hypoth\`eses que $Y\fleche X$ (les
conditions (a), (b), et (c) sont stables par changement de base). Donc quitte
\`a remplacer $X$ par $A$ on peut supposer que l'image r\'eciproque de tout
ferm\'e strict de $X$ est noeth\'erienne et il s'agit de montrer que $Y$ est
quasi-compact. C'est clair si $X$ est r\'eductible, ou si $f(Y)$ n'est pas
dense, donc on peut supposer $X$ irr\'eductible et $f$ dominant. On peut
supposer $Y$ r\'eduit, puisque la question est purement topologique, et les
sch\'emas qui entrent en jeu dans les conditions (b) et (c) sont tous r\'eduits.
Alors $f$ se factorise par $X_{\text{r\'ed}}$, et le morphisme $Y\fleche
X_{\text{r\'ed}}$ ainsi obtenu v\'erifie encore les conditions (a), (b), (c), si
bien que l'on peut supposer $X$ int\`egre. Soit $U$ un ouvert affine
non vide de $X$. Soit $\Sc$ l'ensemble des
points de $U$ qui sont
dans l'image de $Y$. $\Sc$ est dense dans $U$, et de plus, comme $Y\fleche X$
est un monomorphisme, l'extension r\'esiduelle $\kappa(f(y)) \inj \kappa(y)$
est triviale pour tout $y\in Y$, de sorte que chaque point $s\in \Sc$ se
rel\`eve \`a $Y$. En appliquant la condition (c) \`a $U$, on obtient, quitte
\`a r\'eduire $U$, un rel\`evement de $U$ \`a $Y$, rendant le diagramme suivant
commutatif :
$$\xymatrix{&U \incl[d] \ar@{.>}[dl]_g\\ Y \ar[r] &X.}$$
Alors $Y=g(U)\cup f^{-1}(X\setminus U)$ est quasi-compact.
\end{etape}

\begin{etapefinale}{$f$ est une immersion.}
%Montrons \`a pr\'esent que $Y\fleche X$ est une immersion.
\indent C'est une question
locale, donc il suffit de le faire au voisinage de chaque point $x_0$ de $X$. On
note encore $Z=\overline{f(Y)} \subset X$. Il y a trois cas \`a traiter :
\begin{itemize}
\item $x_0\in X\setminus Z$,
\item $x_0\in f(Y)$,
\item $x_0\in Z\setminus f(Y)$.
\end{itemize}
Le premier cas est trivial. Montrons que le troisi\`eme r\'esulte du second.
Il suffit pour cela de montrer que $Z\setminus f(Y)$ est ferm\'e. En effet,
$f(Y)$ sera
alors un ouvert de $Z$, et s'\'ecrira donc $f(Y)=Z\cap U$ o\`u $U$ est un ouvert
de $X$. Le morphisme $f$ se factorisera par $U$ et pour v\'erifier que c'est une
immersion, il suffira de v\'erifier qu'il en est ainsi au voisinage de chaque
point $x_0$ de $U$, ce qui nous ram\`enera \`a l'un des deux premiers cas
suivant que $x_0$ appartient \`a $f(Y)$ ou non.
Maintenant, si $Z\setminus f(Y)$ n'\'etait pas ferm\'e, il existerait un point
$x_0\in \overline{Z\setminus f(Y)}\cap f(Y)$. D'apr\`es le second cas, $f$ est
une immersion au voisinage de $x_0$, donc il existe un voisinage $U$ de $x_0$
tel que $f(Y)\cap U$ soit ferm\'e dans $U$, ce qui contredit $x_0\in
\overline{Z\setminus f(Y)}$ (remarquer que $Z\cap U=f(Y)\cap U$).

Il nous reste donc \`a montrer que $f$ est une immersion au voisinage de tout
point de $f(Y)$. Soient $y_0\in Y$ et $x_0=f(y_0)$. Il suffit de le faire pour
un voisinage \'etale de $x_0$, car \^etre une immersion ferm\'ee est une
propri\'et\'e de nature locale pour la topologie \'etale sur $X$ (\cite{EGA4_2},
(2.7.1)). 

\emph{Premier cas : $X$ est un schéma local hens\'elien et
noeth\'erien dont $x_0$ est le point ferm\'e.}

Alors $f : Y\fleche X$ est s\'epar\'e et de type fini, et $y_0$ est
isol\'e dans sa fibre (car $f$ est un monomorphisme) donc d'apr\`es
\cite{Artin_Montreal}, th\'eor\`eme (1.10) p.26, il existe un splittage
$Y=Y_0\sqcup Y_1$ avec $Y_0\fleche X$ fini et $y_0\in Y_0$. Comme $f$
est un monomorphisme, il est injectif et $f^{-1}(x_0)\cap Y_1$ est vide. De
plus, $Y_0\fleche X$ est un monomorphisme propre, donc d'apr\`es \cite{EGA4_4}
(18.12.6) c'est une immersion ferm\'ee. Il suffit donc de montrer que
$Y_1$ est vide. Dans le cas contraire, $x_0\in \overline{f(Y_1)}$.
Le morphisme $f : Y\fleche X$ est de pr\'esentation
finie donc d'apr\`es le th\'eor\`eme de Chevalley (\cite{EGA1}, (7.1.4) p.329)
$f(Y_1)$ est constructible. Il s'\'ecrit alors
$f(Y_1)=(U_1\cap F_1) \cup \dots \cup (U_k\cap F_k)$, o\`u pour tout $i$, $F_i$
est un ferm\'e irr\'eductible, $U_i$ un ouvert, et $(U_i\cap F_i)$ est non vide.
On a $x_0\in \overline{f(Y_1)} \subset F_1\cup \dots \cup F_k$, par exemple
$x_0\in F_1$. On note $x_1$ le point g\'en\'erique de $F_1$. Alors
$x_0\in\overline{\{x_1\}}$. Par ailleurs, $U_1\cap F_1$ est un ouvert non vide
de $F_1$, donc il contient $x_1$ et il existe $y_1\in Y_1$ tel que $x_1=f(y_1)$.
Le lemme pr\'ec\'edent et la condition (b) impliquent alors que
$y_0\in\overline{\{y_1\}}$, ce qui est absurde car $Y=Y_0\sqcup Y_1$.

\emph{Cas g\'en\'eral.}

On note $X=\Spec A$, $X^h=\Spec A^h$ o\`u $A^h$ est l'hens\'elis\'e de $A$ en
$x_0$, et $f^h : Y^h\fleche X^h$ le morphisme obtenu par changement de base.
Ce dernier v\'erifie encore les conditions (a), (b), (c), il est de type fini et
$A^h$ est local hens\'elien et noeth\'erien. Donc d'apr\`es ce qui pr\'ec\`ede
$f^h$ est une immersion ferm\'ee. Maintenant, comme $Y\fleche X$ est de
pr\'esentation finie, et $X$ quasi-compact, le th\'eor\`eme (8.10.5) de
\cite{EGA4_3} nous assure qu'il existe un voisinage \'etale $X'$ de $x_0$ tel
que le morphisme induit $Y'\fleche X'$ soit une immersion ferm\'ee. Ceci montre
que $f$ est une immersion ferm\'ee au voisinage de $x_0$ (\cite{EGA4_2},
(2.7.1)) et ach\`eve la d\'emonstration.
\end{etapefinale}
\end{demo}

On rappelle qu'un morphisme d'espaces alg\'ebriques est dit localement
s\'epar\'e si sa diagonale est une immersion quasi-compacte. Le lecteur érudit reconnaîtra dans le corollaire ci-dessous les conditions~[3']~(a) et~(b) du théorème~(5.3) (p.48) de représentabilité pour les espaces algébriques de~\cite{Global_Analysis_1}.

\begin{souscor}
\label{critere_locale_separation}
Soient $S$ un sch\'ema noeth\'erien, et $X$, $Y$ des espaces alg\'ebriques
localement de pr\'esentation finie sur $S$. Soit $f : Y\fleche X$ un morphisme
localement de pr\'esentation finie. Alors $f$ est localement s\'epar\'e si et
seulement si :
\begin{itemize}
\item[(a)] Pour tout anneau $A$ de valuation (resp. de valuation discr\`ete), de corps
des fractions $K$ et de corps r\'esiduel $k$, pour
tout couple $(\xi, \eta)\in (Y\times_X Y)(A)$ tel que $\xi_K = \eta_K$ et
$\xi_k = \eta_k$ on a $\xi = \eta$.
\item[(b)] Soient $A$ un anneau int\`egre et $(\xi, \eta)\in (Y\times_X Y)(A)$. On
suppose qu'il existe un ensemble dense $\Sc$ de points de $\Spec A$ tel que pour
tout $s\in \Sc$, on ait $\xi_s = \eta_s$. Alors il existe un ouvert non vide $U$
de $\Spec A$ tel que $\xi_U = \eta_U$.
\end{itemize}
\end{souscor}
\begin{demo}
La diagonale $\Delta$ de $f$ est un monomorphisme localement de pr\'esentation
finie, et
$Y\times_X Y$ est localement noeth\'erien. De plus les conditions (a) et (b)
ci-dessus signifient pr\'ecis\'ement que $\Delta$ v\'erifie les conditions (b)
et (c) du th\'eor\`eme (\ref{critere_immersionsqc}), de sorte que notre
assertion en r\'esulte.
\end{demo}
\begin{sousremarque}\rm
Pour les champs algébriques, on ne peut pas obtenir de critère de locale séparation, ni même de quasi-séparation, à partir du théorème (\ref{critere_immersionsqc}), pour la
simple raison que la diagonale d'un morphisme de champs algébriques n'est
généralement pas un monomorphisme. En fait, si $\X$ est un $S$-champ algébrique,
sa diagonale est un monomorphisme si et seulement si $\X$ est un $S$-espace
algébrique. En conséquence nous devrons développer d'autres techniques pour
aboutir à la quasi-séparation du champ de Picard.
\end{sousremarque}

\subsection{Sections globales et trivialit\'e des faisceaux inversibles}
\label{critere_trivialite}

Soient $\X$ un champ alg\'ebrique, et $\L$ un faisceau inversible sur $\X$. Le
but de ce paragraphe est de donner des crit\`eres portant sur $\Gamma(\X,\L)$ et
\'eventuellement sur $\Gamma(\X,\L^{-1})$
permettant de voir facilement si $\L$ est trivial. Le premier de ces crit\`eres
montrera plus pr\'ecis\'ement que $\L$ est trivial si et seulement s'il a une
section globale \og partout non nulle\fg. On rappelle que $\Gamma(\X,\L)$
est l'ensemble des collections $(s_{(V,v)})$ de sections de $\L$ sur les
$(V,v)\in \ob\lisets(\X)$ telles que pour toute fl\`eche
$\xymatrix@C=1pc{\varphi : (U,u) \ar[r] & (V,v)}$ dans $\lisets(\X)$ on ait
$\text{r\'es}_{\varphi}\, s_{(V,v)}=s_{(U,u)}$.

% Commen\c cons par remarquer que $\Gamma(\X,\L)$ s'identifie naturellement
% \`a l'ensemble des morphismes de $\O_{\X}$ vers $\L$. En effet, $\Gamma(\X,\L)$
% est l'ensemble des collections $(s_{(V,v)})$ de sections de $\L$ sur les
% $(V,v)\in \ob\lisets(\X)$ telles que pour toute fl\`eche
% $\flechen{(U,u)}{\varphi}{(V,v)}$ dans $\lisets(\X)$ on ait
% $\text{r\'es}_{\varphi}\, s_{(V,v)}=s_{(U,u)}$.
% Par ailleurs, un morphisme de $\O_{\X}$ vers $\L$
% correspond \`a la donn\'ee, pour tout $(U,u)\in \ob \lisets(\X)$, d'un morphisme
% $\O_U \fleche \L_{U,u}$, avec une condition de compatibilit\'e avec les
% morphismes du topos $\lisets(\X)$ (cf. \cite{LMB}, (12.1.2) et (12.2.1)).
% Or un morphisme $\O_U \fleche \L_{U,u}$ correspond
% \`a une section globale $s_{(U,u)}$ de $\L_{U,u}$, et la condition de
% compatibilit\'e se traduit exactement par la condition $\text{r\'es}_{\varphi}\,
% s_{(V,v)}=s_{(U,u)}$.

% Notons $\tilde{s}$ le morphisme associ\'e \`a $s=(s_{(V,v)})\in \Gamma(\X,\L)$.
% Il est clair que $\tilde{s}$ est un isomorphisme si et seulement si pour tout
% $(V,v)\in \ob\lisets(\X)$, le morphisme induit $\O_V \fleche
% \L_{V,v}$ par $s_{(V,v)}$ en est un, c'est-\`a-dire si et seulement si pour tout
% $(V,v)$, la section $s_{(V,v)}$ est partout non nulle sur $V$ (i.e. en tout
% point $x$, son germe engendre la fibre en $x$ de $\L_{V,v}$, ou encore est un
% \'el\'ement inversible de $(\L_{V,v})_x \simeq \O_{V,x}$).

\begin{sousdefi}
Soit $s$ une section globale de $\Lc$.
\begin{itemize}
\item[(i)] On dit que $s$ est partout non nulle sur $\X$ si
pour tout ouvert lisse-étale $(V,v)$ de $\X$, la section globale $s_{(V,v)}$ de
$\L_{V,v}$ est partout non nulle sur $V$.
\item[(ii)] Soit $x\in |\X|$ un point de $\X$, repr\'esent\'e par un morphisme
$\varphi$ de $\Spec K$ vers $\X$ (cf. \cite{LMB}, chapitre 5). La section
$s$ induit une section $\varphi^*s$ de
$\Gamma(\Spec K,\varphi^*\L)$. On dit que $s$ est nulle en $x$ si $\varphi^*s=0$.
(Il est clair que ceci ne d\'epend pas du repr\'esentant $\varphi$ choisi.)
\end{itemize}
\end{sousdefi}

% La discussion pr\'ec\'edente montre que $\tilde{s}$ est un isomorphisme si et
% seulement si $s$ est partout non nulle sur $\X$. Nous allons maintenant relier
% cette notion \`a l'ensemble $|\X|$ sous-jacent \`a $\X$.

\begin{sousprop}
\label{prop_sections_partout_non_nulles}
Soit $s$ une section globale de $\Lc$. Les propositions suivantes sont équivalentes :
\begin{itemize}
\item[(i)] $s$ est partout non nulle sur $\X$ ;
\item[(ii)] $s$ est non nulle en tout point de $\X$ ;
\item[(iii)] le morphisme de $\Oc_{\X}$ dans $\Lc$ correspondant à $s$ est un isomorphisme.
\end{itemize}
De plus l'ensemble des points de $\X$ o\`u $s$ est
non nulle est un ouvert de $|\X|$.
\end{sousprop}
\begin{demo}
L'équivalence de~(i) et~(iii) résulte immédiatement du fait qu'une section globale $s_{(V,v)}$ de $\Lc_{V,v}$ est partout non nulle si et seulement si le morphisme correspondant de $\Oc_V$ dans $\Lc_{V,v}$ est un isomorphisme.

Le fait que (i) implique (ii) est clair puisque tout point de $|\X|$ admet un représentant qui se factorise par un ouvert lisse-étale de $\X$. La réciproque est évidente aussi : plus généralement si $s$ est non nulle en tout point de $\X$ alors pour tout morphisme $F : \Y \fleche \X$ la section $F^*s$ est non nulle en tout point de $\Y$.

Pour la dernière assertion, soit $x\in |\X|$ un point de $\X$ en lequel $s$ est non nulle et soit $U$ un ouvert lisse-étale de $\X$ le contenant. Il existe alors un voisinage $V$ de $x$ dans $U$ tel que $s$ soit non nulle en tout point de $V$, et son image est un voisinage ouvert de $x$ dans $\X$.
\end{demo}

\begin{sousremarque}\rm
La section $s$ définit en fait un idéal quasi-cohérent $\Ic$ de $\Oc_{\X}$, donc un sous-champ fermé de $\X$. Le fermé de $|\X|$ correspondant est l'ensemble des points où $s$ est nulle.
\end{sousremarque}

Il est bien connu qu'un faisceau inversible sur un sch\'ema propre est trivial si et seulement
s'il est engendr\'e par ses sections globales, de m\^eme que son inverse. Nous aurons besoin d'un résultat analogue pour les champs. L'hypothèse de platitude cohomologique sous laquelle nous nous plaçons souvent rend le résultat presque trivial et nous dispense même de l'hypothèse de propreté (cf.~\ref{critere_trivialite_2}).

\begin{sousdefi}
Soient $\X$ un $S$-champ alg\'ebrique, et $\Fc$ un faisceau de
$\O_{\X}$-modules. On dit que $\Fc$
est engendr\'e par ses sections globales s'il existe une famille de sections
globales $\{s^i\}_{i\in I}$, $s^i\in\Gamma(\X,\Fc)$, telle que pour tout
$(U,u)\in\ob\lisetsth(\X)$, $\Fc_{U,u}$ soit engendr\'e par les sections
$\{s^i_{U,u}\}_{i\in I}$.
\end{sousdefi}

\begin{sousremarque} \rm
Si $\Fc$ est un faisceau inversible et si $\{s_i\}_{i\in I}$ est une famille de sections globales de $\Fc$, alors la famille $\{s^i_{U,u}\}_{i\in I}$
engendre $(\Fc_{U,u})_x$ en un point $x\in U$ o\`u $(U,u)\in\ob \lisets(\X)$ si
et seulement s'il existe un indice $i$ tel que $s^i$ soit non nulle au point de $|\X|$
repr\'esent\'e par $x$. Donc $\Fc$ est engendr\'e par ses sections globales si
et seulement si pour tout point $x$ de $|\X|$, $\Fc$ a une section globale
non nulle en $x$.
\end{sousremarque}

\begin{sousprop}
\label{critere_trivialite_2}
Soit $\X$ un champ alg\'ebrique cohomologiquement plat en
dimension z\'ero sur le spectre d'un corps $k$. Soit $\Fc$ un faisceau inversible sur $\X$.
Alors $\Fc$ est trivial si et seulement si $\Fc$ et $\Fc^{-1}$ sont engendr\'es
par leurs sections globales.
\end{sousprop}
\begin{demo}
La partie \og seulement si\fg\ est \'evidente puisque le faisceau $\O_{\X}$ est
engendr\'e par sa section globale constante $1\in \Gamma(\X,\O_{\X})$. Pour la
r\'eciproque, supposons que $\Fc$ et $\Fc^{-1}$ soient engendr\'es par leurs
sections globales. Soit $x$ un point de $\X$. Il existe des sections
$s$ de $\Gamma(\X,\Fc)$ et $s'$ de $\Gamma(\X,\Fc^{-1})$ non nulles en $x$. Alors
$s\otimes s'\in \Gamma(\X,\O_{\X})$ est non nulle en $x$. Mais
$\Gamma(\X,\O_{\X})=k$ par platitude cohomologique, donc $s\otimes s'$ est
partout non nulle sur $\X$, et donc $s$ aussi, ce qui
montre que $\Fc$ est trivial.
\end{demo}

\begin{souslem}
Soient $\X$ un champ alg\'ebrique propre sur $S=\Spec k$, $L/k$ une extension de
$k$, et $\Fc$ un faisceau inversible sur $\X$. Alors $\Fc$ est engendr\'e par
ses sections globales si et seulement si $\Fc_L$ l'est.
\end{souslem}
\begin{demo}
D'apr\`es \cite{Varietes_abeliennes}, cor. 5 p.53 et \cite{Aoki_Hom}, annexe A,
on a un isomorphisme naturel $\xymatrix@C=1pc{H^0(\X,\Fc)\otimes_k
L \ar[r]^{\sim} & H^0(\X\otimes_k L, \Fc_L)}.$ Supposons que $\Fc$ soit engendr\'e par
ses sections globales. Soient $y\in|\X\otimes_k L|$ un point de $\X\otimes_k L$,
et $x\in |\X|$ son image. Il existe par hypoth\`ese une section globale $s\in
H^0(\X,\Fc)$ non nulle en $x$. Alors $s\otimes 1\in H^0(\X\otimes_k L, \Fc_L)$
est non nulle en $y$. R\'eciproquement, si $\Fc$ n'est pas engendr\'e par ses
sections globales, soit $x\in |\X|$ un point de $\X$ en lequel toute section
globale de $\Fc$ est nulle. Soit $y\in |\X\otimes_k L|$ un ant\'ec\'edent de $x$
($\X\otimes_k L\fleche \X$ est surjectif). Alors il est clair que toute section
globale $t\in H^0(\X,\Fc)\otimes_k L$ de $\Fc_L$ est nulle en $y$, ce qui montre
que $\Fc_L$ n'est pas engendr\'e par ses sections globales.
\end{demo}

\begin{souscor}
\label{trivialite_apres_extension_de_corps}
Soit $\X$ un champ alg\'ebrique propre et cohomologiquement plat en
dimension z\'ero sur un corps $k$. Soit $L/k$ une
extension de $k$. Alors $\Fc$ est isomorphe \`a $\O_{\X}$ si et seulement si
$\Fc_L$ est isomorphe \`a $\O_{\X\times_{k} L}$. $\square$
\end{souscor}

\subsection{Un peu de descente fid\`element plate}

\begin{sousprop}
\label{cond_sep_locales}
Soit
% $$\shorthandoff{!;:?}
% \xymatrix@!0 @R=1.5pc @C=3pc{\X'\ar[rr]\ar[rd]^{F'}\ar[dd] &&
% \X\ar[rd]^F\ar[dd]|!{[ld];[rd]}\hole &\\
% & \Y'\ar[rr]\ar[ld] && \Y\ar[ld] \\ \Zc'\ar[rr]^q && \Zc&}$$
$$\xymatrix{\X'\ar[r]\ar[d]_{F'}\cartesien& \X\ar[d]^F\\
\Y'\ar[r]_q&\Y}$$
un diagramme 2-cart\'e\-sien de $S$-champs alg\'ebriques.
\begin{itemize}
\item[(i)] On suppose que le morphisme de changement de base 
$q:\Y'\fleche \Y$ est surjectif et universellement ouvert (c'est en particulier le cas d'une
pr\'esentation de $\Y$, ou encore d'un morphisme fid\`element plat et
localement de pr\'esentation finie). Alors $F$ est quasi-compact (resp.
quasi-s\'epar\'e) si et seulement si $F'$ l'est.
\item[(ii)] On suppose $q$ fid\`element plat et quasi-compact, et $F$ repr\'esentable.
Alors $F$ est une immersion quasi-compacte si et seulement si $F'$ en est une.
\item[(iii)] On suppose $q$ fid\`element plat et quasi-compact. Alors $F$ est
localement s\'epar\'e si et seulement si $F'$ est localement s\'epar\'e.
\end{itemize}
\end{sousprop}
\begin{demo}
(i) Pour la premi\`ere assertion, la d\'emonstration de \cite{LMB}, corollaire
(5.6.3), est encore valable dans notre cas sans aucune modification. La seconde
s'en d\'eduit imm\'ediatement en constatant que le morphisme diagonal
$\Delta_{F'}$ s'identifie au morphisme obtenu \`a partir de $\Delta_F$ par le
changement de base $q$.

(ii) $F$ \'etant repr\'esentable, on se ram\`ene formellement au cas o\`u $\Y$
est un sch\'ema affine. Quitte \`a consid\'erer une pr\'esentation de $\Y'$, on
peut aussi supposer que $\Y'$ est un sch\'ema. Alors $\X$ est un $S$-espace
alg\'ebrique, et $\X'$ est un sch\'ema puisqu'une immersion est toujours
sch\'ematique. Dans ce cas le r\'esultat d\'ecoule de l'effectivit\'e de la
descente des morphismes quasi-affines par morphismes fid\`element plats
quasi-compacts.

(iii) se d\'eduit formellement de (ii) en tenant compte du fait qu'un morphisme
diagonal est toujours repr\'esentable.
\end{demo}

\section{Locale s\'eparation du foncteur de Picard}

On vérifie dans le lemme ci-dessous les conditions~(a) et~(b) du corollaire~\ref{critere_locale_separation} pour le foncteur de Picard.

\begin{lem}
\label{lemme_locale_separation_du_foncteur_de_Picard}
Soient $A$ un anneau int\`egre, $T=\Spec A$, et $\X \fleche T$ un $T$-champ
alg\'ebrique propre, plat, de pr\'esentation finie, et cohomologiquement
plat en dimension z\'ero. Soit $\L$ un faisceau inversible sur $\X$.
\begin{itemize}
\item[(i)] Si $A$ est un anneau local (en particulier si $A$ est un anneau de
valuation discr\`ete), de corps
r\'esiduel $k$ et de corps des fractions $K$, alors $\L$
est trivial si et seulement si $\L_k$ et $\L_K$ le sont.
\item[(ii)] S'il existe un ensemble dense $\Sc$ de points de $T$ tel que pour tout
$t\in \Sc$, $\L_t$ soit trivial, alors il existe un ouvert non vide $U$ de $T$
tel que $\L_U$ soit trivial.
\end{itemize}
\end{lem}
\begin{demo}
D'apr\`es l'annexe A de~\cite{Aoki_Hom} les corollaires~1 et 2 du paragraphe 5
de~\cite{Varietes_abeliennes} sont encore valables. Ils nous apprennent que la
fonction $$t\longmapsto d(t):=\dim_{\kappa(t)} H^0(\X_t, \L_t)$$ est semi-continue
sup\'erieurement, et que si elle est constante sur
$T$ (qui est bien r\'eduit et connexe), alors le morphisme naturel
$$(f_*\L)\otimes_{\O_T} \kappa(t) \flechelongue H^0(\X_t, \L_t)$$ est un isomorphisme
pour tout $t\in T$. De plus en un point o\`u $\L_t$ est trivial, on a $H^0(\X_t,
\L_t)=H^0(\X_t, \O_{\X_t})\simeq \kappa(t)$ par platitude cohomologique, donc
$d(t)=1$.

\medskip
(i) Ici comme $\L_K$ et $\L_k$ sont triviaux, la fonction $d$ vaut 1 au point
g\'en\'erique et au point ferm\'e, donc elle est constante
\'egale \`a 1 sur $T$. D'apr\`es ce qui pr\'ec\`ede, on a en particulier un
isomorphisme $H^0(T, f_*\L)\otimes_A k=H^0(\X,\L)\otimes_A k\fleche
H^0(\X_k,\L_k)$. Le faisceau $\L_k$ \'etant trivial, il a une section globale
partout non nulle sur $\X_k$. D'apr\`es l'isomorphisme pr\'ec\'edent,
cette section provient d'une
section globale $s$ de $\L$ (comme $k$ est un quotient de $A$, les \'el\'ements
de $H^0(\X,\L)\otimes_A k$ peuvent tous s'\'ecrire sous la forme $s\otimes 1$).
Alors $s$ est non nulle en tout point $x$ de $\X$ d'image $\Spec k$ dans $T$.
D'apr\`es (\ref{prop_sections_partout_non_nulles}), l'ensemble $C$ des points de
$\X$ o\`u $s$ est nulle est un ferm\'e de $|\X|$. Comme $f$ est propre, l'image
de $C$ est un ferm\'e de $T$, qui de plus ne contient pas le point ferm\'e.
N\'ecessairement $C$ est vide, donc $s$ est partout non nulle et d'apr\`es
(\ref{critere_trivialite}) $\L$ est trivial.

(ii) La fonction $d$ \'etant semi-continue sup\'erieurement, il existe un ouvert
non vide sur lequel elle est constante et quitte \`a remplacer $T$ par cet
ouvert, on peut supposer qu'elle est constante. D\`es lors pour tout $t\in T$ le
morphisme naturel
$$H^0(\X,\L)\otimes_A \kappa(t) \flechelongue H^0(\X_t, \L_t)$$
est un isomorphisme. Soit $t\in \Sc$. Alors $\L_t$ est trivial donc il a une
section globale $s_t\in H^0(\X_t, \L_t)$ partout non nulle sur $\X_t$. D'apr\`es
l'isomorphisme pr\'ec\'edent, elle provient d'un \'el\'ement de
$H^0(\X,\L)\otimes_A \kappa(t)$ que l'on peut \'ecrire sous la forme
$s\otimes(\frac1f)$ o\`u $f\in A$ n'est pas dans l'id\'eal premier correspondant
\`a $t$. Quitte \`a remplacer $T$ par $D(f)$, on peut supposer $f=1$ et on a
donc trouv\'e une section $s\in H^0(\X,\L)$ dont la r\'eduction $s_t$ \`a $\X_t$
est partout non nulle.

Maintenant l'ensemble $C$ des points de $\X$ o\`u $s$ est nulle est un ferm\'e
de $|\X|$ d'apr\`es (\ref{prop_sections_partout_non_nulles}), et son image est
un ferm\'e $F$ de $T$, qui ne contient pas $t$. Son compl\'ementaire
$U=T\setminus F$ est alors un ouvert non vide de $T$, et la section $s_U$
induite par $s$ sur $\X\times_T U$ est partout non nulle sur $\X\times_T U$, de
sorte que $\L_U$ est trivial.
\end{demo}

\begin{prop}
\label{locale_separation_du_foncteur_de_Picard}
Soient $S$ un sch\'ema noeth\'erien et $\X$ un $S$-champ alg\'ebrique propre,
plat, et cohomologiquement plat en dimension z\'ero. Alors le foncteur
$\pic$ est un $S$-espace alg\'ebrique localement s\'epar\'e.
\end{prop}
\begin{demo}
On sait d\'ej\`a d'apr\`es \cite{Aoki_Hom} et la proposition
(\ref{comparaison_champ_foncteur}) que $\pic$ est repr\'esentable par un $S$-espace
pr\'ealg\'ebrique. Il s'agit de montrer que le morphisme $$\pic\flechelongue S$$ est
localement
s\'epar\'e. Par~(\ref{cond_sep_locales})~(iii) et~(\ref{remarque_section}), on peut supposer que
le morphisme $\X\fleche S$ a une section, ce qui nous ram\`ene par
(\ref{comparaison_des_foncteurs_de_Picard}) au cas o\`u $\pic$ est donn\'e par
$$\pic(T)=\frac{\Pic(\X\times_S T)}{\Pic(T)}.$$
On peut appliquer le corollaire~(\ref{critere_locale_separation}) au morphisme
$\pic\fleche S$, qui est bien localement de pr\'esentation finie par~(\ref{lpf}). Les conditions (a) et (b) sont v\'erifi\'ees gr\^ace au lemme
pr\'ec\'edent.
\end{demo}

\section{Quasi-s\'eparation du champ de Picard}

Nous allons montrer dans ce paragraphe que sous de bonnes hypoth\`eses, le
champ de Picard $\champic(\X/S)$ est quasi-s\'epar\'e. En fait ce r\'esultat se
d\'eduit sans trop d'efforts du r\'esultat de locale s\'eparation du foncteur de
Picard (\ref{locale_separation_du_foncteur_de_Picard}).
Cependant, cette d\'emarche peut sembler insatisfaisante, et il convient
de donner de ce r\'esultat une d\'emonstration directe. Nous avons pour cela
d\'evelopp\'e un crit\`ere de quasi-s\'eparation pour les champs alg\'ebriques
(\ref{critere_quasi_separation}),
inspir\'e du crit\`ere de locale s\'eparation (\ref{critere_locale_separation}),
qu'il nous a sembl\'e pertinent d'\'enoncer s\'epar\'ement.
Rappelons tout d'abord quelques notations.

Si $\X$ est un $S$-champ, $U\in\ob\aff$, et si $x,y$ sont deux
objets de $\X_U$, on notera respectivement $\Auts(x)$ et $\Isoms(x,y)$ les
foncteurs des automorphismes et des isomorphismes dont les d\'efinitions sont
rappel\'ees ci-dessous :
$$\fonction{\Auts(x)}{({\rm Aff}/U)}{({\rm Gr})}{(V\fleche U)}{\Aut(x_V)}$$
$$\fonction{\Isoms(x,y)}{({\rm Aff}/U)}{({\rm Ens})}{(V\fleche U)}{
   \Isom(x_V,y_V)}.$$

Le r\'esultat principal concernant le champ de Picard est le suivant.

\begin{thm}
\label{quasi_separation_du_champ_de_Picard}
Soient $S$ un sch\'ema noeth\'erien, et $\X$ un $S$-champ alg\'ebrique propre,
plat, et cohomologiquement plat en dimension z\'ero. Alors le champ de Picard
$\champic(\X/S)$ est un $S$-champ \emph{alg\'ebrique} (c'est-\`a-dire qu'il est
quasi-s\'epar\'e).
\end{thm}
{\bf Premi\`ere d\'emonstration (utilisant
(\ref{locale_separation_du_foncteur_de_Picard}))}

\noindent
Notons $\Pc=\champic(\X/S)$, et $P=\pic$. D'apr\`es
(\ref{locale_separation_du_foncteur_de_Picard}), le morphisme diagonal
$\Delta_P : P\fleche P\times_S P$ est une immersion quasi-compacte. On veut
montrer que le morphisme diagonal $\Delta_{\Pc}$ de $\Pc$ est quasi-compact,
i.e. que pour tout $U\in\ob \aff$ et tout couple
$(\Lc, \Mc)$ de faisceaux inversibles sur $\X\times_S U$, le morphisme de $U$-espaces
alg\'ebriques $\Isoms(\Lc, \Mc) \fleche U$ est quasi-compact. Les faisceaux
$\Lc$ et $\Mc$ (plus pr\'ecis\'ement leurs classes d'isomorphie) d\'efinissent
un morphisme $([\Lc], [\Mc])$ de $U$ dans $P\times_S P$. Notons $U_0$ le
sous-sch\'ema quasi-compact de $U$ sur lequel $\Lc$ et $\Mc$ sont isomorphes,
i.e. le sous-sch\'ema d\'efini par le carr\'e cart\'esien
$$\xymatrix{U_0 \ar@{^{(}->}[r] \ar[d] \cartesien & U \ar[d]^{([\Lc], [\Mc])}\\
P \ar@{^{(}->}[r]^{\Delta_P} &P\times_S P.}$$
Alors on a une factorisation :
$$\shorthandoff{!;:?}
\xymatrix@!0 @R=1.5pc @C=3pc{\Isoms(\Lc,\Mc) \ar[rr] \ar@{.>}[rd]&&U\\
    &U_0\ar@{^{(}->}[ur]&}$$
qui provient de la commutativit\'e du diagramme suivant :
$$\xymatrix{\Isoms(\Lc,\Mc) \cartesien \ar[r]\ar[d]
    &\Pc\ar[r]^{\pi}\ar[d]^{\Delta_{\Pc}} &P\ar[d]^{\Delta_P}\\
    U\ar[r]^{(\Lc,\Mc)}& \Pc\times_S\Pc \ar[r]^{\pi\times\pi} & P\times_S P}$$
(o\`u $\pi$ d\'esigne le morphisme naturel $\Pc \fleche P$).
Il suffit donc de montrer que
$$\Isoms(\Lc,\Mc) \flechelongue U_0$$
est quasi-compact.
Or c'est un \'epimorphisme (on peut le voir soit de mani\`ere directe, soit en
constatant qu'on a un diagramme 2-cart\'esien
$$\xymatrix{\Isoms(\Lc,\Mc) \cartesien\ar[r]\ar[d] &U_0 \ar[d]\\
    \Pc \ar[r]^{\Delta_{\pi}}& \Pc\times_P\Pc}$$
et en se rappelant que $\pi$ est une gerbe). $\Isoms(\Lc,\Mc) \fleche U_0$ est
donc un torseur sous $\Auts(\Mc)$, c'est-\`a-dire un $\gm$-torseur puisque par
platitude cohomologique, $\Auts(\Mc)$ est isomorphe à $\gm$. Alors $\Isoms(\Lc,\Mc) \fleche
U_0$ est quasi-compact puisque $\gm$ est quasi-compact. $\square$ \vskip .3cm

La proposition ci-dessous donne un crit\`ere g\'en\'eral de quasi-s\'eparation,
gr\^ace auquel nous pourrons donner une seconde d\'emonstration, directe, du
fait que le champ de Picard est quasi-s\'epar\'e.

\begin{prop}
\label{critere_quasi_separation}
Soient $S$ un sch\'ema localement noeth\'erien et $\X$ un $S$-champ
pr\'ealg\'ebrique
localement de pr\'esentation finie. On suppose que les deux conditions suivantes
sont remplies.

(i) Pour tout $U\in\ob\aff$ et tout $x\in\ob \X_U$, le morphisme
$\Auts(x)\fleche U$ est quasi-compact.

(ii) Soit $U\in\ob\aff$ int\`egre, et soient $x, y\in \ob \X_U$. On suppose
qu'il existe un ensemble dense de points $t$ de $U$, tels qu'il existe une
extension $L(t)$ de $\kappa(t)$ telle que $x_{L(t)} \simeq y_{L(t)}$. Alors $x$
et $y$ sont isomorphes sur un ouvert dense de $U$.

Alors $\X$ est quasi-s\'epar\'e. (Autrement dit c'est un $S$-champ
\emph{alg\'ebrique}.)
\end{prop}
\begin{demo}
Le champ $\X\times_S\X$ est localement noeth\'erien (car localement de pr\'esentation
finie sur $S$). Soit $U\fleche \X\times_S\X$ une pr\'esentation de
$\X\times_S\X$, o\`u $U$ est un sch\'ema localement noeth\'erien. Alors par
(\ref{cond_sep_locales}) (i), pour montrer que $\Delta: \X \fleche \X\times_S
\X$ est quasi-compact, il suffit de montrer que le morphisme obtenu par
changement de base
$$\xymatrix{I\ar[r] \ar[d] \cartesien &U\ar[d]^{(x,y)}\\
  \X\ar[r]& \X\times_S \X}$$
est quasi-compact. Comme c'est une question locale sur $U$, on peut supposer $U$
affine. Il s'agit donc de montrer que si $U$ est un sch\'ema affine
noeth\'erien, et si $x,y$ sont deux objets de $\X_U$, alors
$\xymatrix@C=1pc{I=\Isoms(x,y)\ar[r]^{f}& U}$ est un morphisme quasi-compact.

Soit $G=\Auts(y)$. L'action \`a gauche naturelle de $G$ sur $I$ fait clairement
de $I\fleche U$ un pseudo-torseur\footnote{Autrement dit la fl\`eche
$I\fleche U$ est invariante sous $G$, et le morphisme naturel $G\times I \fleche
I\times_U I$ est un isomorphisme.} sous $G$. Par r\'ecurrence noeth\'erienne sur
les ferm\'es de $U$, on peut supposer que pour tout ferm\'e strict $F$ de $U$,
$I\times_U F$ est quasi-compact. Il suffit maintenant de trouver un ouvert non
vide $V$ de $U$ tel que $I\times_U V$ soit quasi-compact. En effet, si $V$ est
un tel ouvert, alors en notant $F=U\setminus V$, l'espace alg\'ebrique
$I'=(I\times_U V )\coprod (I\times_U F)$ est quasi-compact, et le morphisme
$I'\fleche I$ obtenu par changement de base \`a partir de $V\coprod F \fleche U$
est surjectif, donc $I$ est quasi-compact.

Posons $Z=\overline{f(I)}$. Soit $W$ un ouvert affine int\`egre de
$Z_{\text{r\'ed}}$. Alors le morphisme $f_W: I\times_U W \fleche W$ est
dominant, et la condition (ii) signifie pr\'ecis\'ement que dans ce cas, quitte
\`a choisir $W$ plus petit, $f_W$ a une section. Ceci prouve que $I\times_U W$
est un torseur sous $G$, donc est quasi-compact. Ceci \'etant, on peut \'ecrire
$W=Z_{\text{r\'ed}}\cap V=Z_{\text{r\'ed}}\times_U V$ o\`u $V$ est un ouvert
(non vide) de $U$. On v\'erifie sans peine que le morphisme induit $I\times_U W
\fleche I\times_U V$ est surjectif, ce qui montre que $I\times_U V$ est
quasi-compact, et ach\`eve la d\'emonstration.
\end{demo}

\noindent {\bf Seconde d\'emonstration du th\'eor\`eme
(\ref{quasi_separation_du_champ_de_Picard})}
On sait d\'ej\`a d'apr\`es~\cite{Aoki_Hom} que le $S$-champ $\champic(\X/S)$ est pr\'ealg\'ebrique. Il est localement de pr\'esentation finie d'apr\`es~\ref{lpf}. La condition (i) r\'esulte du fait que, par platitude
cohomologique, pour tout faisceau inversible $\Lc$ sur $\X\times_S U$, le
$U$-espace en groupes $\Auts(\Lc)$ est isomorphe \`a ${\mathbb G}_{{\rm m},U}$,
donc il est quasi-compact. La condition (ii) r\'esulte du lemme
(\ref{lemme_locale_separation_du_foncteur_de_Picard}) (ii) et de
(\ref{trivialite_apres_extension_de_corps}). $\square$

\chapter{Th\'eorie des d\'eformations}

\def\modtf{\rm{Mod}^{\rm{pl}}}
\section{Conditions de Schlessinger}

Dans toute cette section, on suppose que $S=\Spec R$ est un sch\'ema affine, et
on consid\`ere un carr\'e cart\'esien de $R$-alg\`ebres :
$$\xymatrix{B' \ar[r] \ar[d] \cartesien &A' \ar[d]^p \\ B \ar[r] &A.}$$
$B'$ est donc isomorphe au produit fibr\'e $A'\times_A B$. On suppose de plus
que $p$ est surjectif.

On rappelle le th\'eor\`eme suivant, d\^u \`a Ferrand, o\`u pour tout anneau $R'$, $\modtf(R')$ d\'esigne la
cat\'egorie des $R'$-modules plats.

\begin{thm}[\cite{Ferrand}, th. 2.2]
\label{thm_ferrand}
Le foncteur $$\modtf(B') \flechelongue \modtf(A')\times_{\modtf(A)} \modtf(B)$$ qui
\`a un $B'$-module $M$ associe le triplet $(M\otimes_{B'}A', M\otimes_{B'}B,
\varphi)$, o\`u $\varphi$ est l'isomorphisme canonique
$\xymatrix@C=1pc{(M\otimes_{B'}A')\otimes_{A'}A \ar[r]^{\sim} &
(M\otimes_{B'}B)\otimes_{B}A}$, est une \'equivalence de cat\'egories. Un
quasi-inverse est donn\'e par le foncteur qui \`a un triplet $(M',N,\varphi)$,
associe le $A'\times_A B$-module
$$S(M',N,\varphi)=
        \{(m',n)\in M'\times N\ |\ \varphi(1\otimes m')=1\otimes n\}.$$
L'\'enonc\'e est encore valable en rempla\c cant la condition \og plats\fg\ par
\og plats et de type fini\fg, \og projectifs et de type fini\fg, ou encore
\og localement libres de rang 1\fg.
\end{thm}

Nous allons en d\'eduire un th\'eor\`eme de comparaison entre la
cat\'egorie des faisceaux inversibles sur $\X_{B'}$ et le produit fibr\'e des
cat\'egories des faisceaux inversibles sur $\X_{A'}$ et sur $\X_B$. (Ici un 
indice comme $\X_A$ d\'esigne le changement de base $\X\times_S \Spec A$.)

\begin{thm}
\label{equiv_cat_fi_pour_defo}
Avec les notations pr\'ec\'edentes, si $\X$ est un champ alg\'ebrique
\emph{plat} sur $S$, on a une \'equivalence de cat\'egories
$$\inv(\X_{B'}) \flechelongue \inv(\X_{A'}) \times_{\inv(\X_{A})} \inv(\X_{B}).$$
\end{thm}
\begin{demo}

Commen\c cons par l'observation suivante. Soient $\U$ un $S$-champ
alg\'ebrique, $U^1$ un
$S$-espace alg\'ebrique, et soit $U^1 \fleche \U$ un \'epimorphisme. On note $U^2=U^1\times_{\U} U^1$ et
$U^3=U^1\times_{\U} U^1\times_{\U} U^1$.
Alors d'après~\ref{descente_fidelement_plate} la cat\'egorie $\inv(\U)$ est \'equivalente \`a la
cat\'egorie ${\rm Desc}(U^1,\U)$ d\'efinie comme suit :\\
\indent - un objet de ${\rm Desc}(U^1,\U)$ est la donn\'ee d'un faisceau
inversible $\L$
sur $U^1$ et d'un isomorphisme $\varphi: pr_1^*\L \fleche pr_2^*\L$ dans
$\inv(U^1\times_{\U} U^1)=\inv(U^2)$ v\'erifiant la condition de cocycle
$$pr_{1,3}^*\varphi=pr_{2,3}^*\varphi \circ pr_{1,2}^*\varphi$$
dans $\inv(U^3)$ ;\\
\indent - un morphisme de $(\L,\varphi)$ vers $(\L',\varphi')$ est la donn\'ee
d'un isomorphisme $\alpha$ de $\L$ vers $\L'$ faisant commuter le diagramme
$$\xymatrix{pr_1^*\L \ar[r]^{\varphi} \ar[d]_{pr_1^*\alpha} &pr_2^*\L
\ar[d]^{pr_2^*\alpha} \\
pr_1^*\L' \ar[r]_{\varphi'} &pr_2^*\L'.}$$
De plus la construction de ${\rm Desc}(U^1,\U)$ commute aux produits fibr\'es de
cat\'egories, de sorte que si les trois foncteurs naturels
$$\inv(U^i_{B'}) \flechelongue \inv(U^i_{A'}) \times_{\inv(U^i_{A})} \inv(U^i_{B}),
\quad i=1, 2, 3,$$
sont des \'equivalences de cat\'egories, il en est de m\^eme du foncteur
$$\inv(\U_{B'}) \flechelongue \inv(\U_{A'}) \times_{\inv(\U_{A})} \inv(\U_{B}).$$

Revenons \`a la d\'emonstration proprement dite de notre th\'eor\`eme.
\vskip .2cm
\emph{Premier cas : $\X$ est un $S$-sch\'ema affine $\Spec R'$.}\\
Alors l'assertion r\'esulte du th\'eor\`eme (\ref{thm_ferrand}) et du lemme
suivant.

\begin{lem}
Si $R'$ est une $R$-alg\`ebre plate, alors le carr\'e obtenu par changement de
base
$$\xymatrix{B'\otimes_R R'\ar[r]\ar[d] &
A'\otimes_R R'\ar[d]^{p_{R'}} \\
B\otimes_R R' \ar[r] &A\otimes_R R'}$$
est encore cart\'esien, et $p_{R'}$ est surjectif.
\end{lem}
\begin{demo}
Un carr\'e de $R$-alg\`ebres est cart\'esien si et seulement si le carr\'e
sous-jacent de $R$-modules l'est. Or un carr\'e de $R$-modules
$$\xymatrix{E \ar[r] \ar[d] &F \ar[d]^f \\ G \ar[r]_g &H}$$
est cart\'esien si et seulement si la suite induite
$$\xymatrix{0 \ar[r]& E \ar[r] &F\oplus G \ar[r]^{f-g} &H}$$
est exacte. Cette condition est \'evidemment stable par changement de base plat.
La surjectivit\'e de $p_{R'}$ est claire.
\end{demo}

\emph{$2^{\textrm{\`eme}}$ cas : $\X$ est un $S$-sch\'ema qui est union
disjointe de sch\'emas affines.}\\
On a $\X=\disp \coprod_i X_i$, o\`u chaque $X_i$ est affine. Alors
$\inv(\X)=\disp \prod_i \inv(X_i)$ universellement, et l'assertion est vraie
pour chacun des $S$-sch\'emas $X_i$, donc il suffit de montrer que les produits
fibr\'es de cat\'egories commutent aux produits arbitraires, ce qui est purement
formel et \'evident.

\emph{$3^{\textrm{\`eme}}$ cas : $\X$ est un sch\'ema s\'epar\'e.}\\
Soit $X$ l'union disjointe d'une famille d'ouverts affines recouvrant $\X$.
Alors $X$, $X\times_{\X} X$, et $X\times_{\X} X\times_{\X} X$ sont unions
disjointes de sch\'emas affines, donc d'apr\`es la remarque qui amor\c cait
notre d\'emonstration, l'assertion r\'esulte du second cas.

\emph{$4^{\textrm{\`eme}}$ cas : $\X$ est un $S$-sch\'ema quelconque.}\\
De m\^eme que pr\'ec\'edemment soit $X$ la somme des ouverts d'un recouvrement
affine de $\X$. Alors $X$ est s\'epar\'e, et, $\X$ \'etant localement s\'epar\'e
(comme tout sch\'ema !), les sch\'emas $X\times_{\X} X$ et $X\times_{\X}
X\times_{\X} X$ sont s\'epar\'es, donc il suffit d'appliquer le troisi\`eme cas.

\emph{$5^{\textrm{\`eme}}$ cas : $\X$ est un $S$-espace alg\'ebrique.}\\
On choisit un morphisme \'etale surjectif $X\fleche \X$ o\`u $X$ est un
sch\'ema. La diagonale d'un espace alg\'ebrique \'etant sch\'ematique (par
d\'efinition), on voit que $X\times_{\X} X$ et $X\times_{\X}
X\times_{\X} X$ sont encore des sch\'emas, de sorte que l'assertion r\'esulte du
cas pr\'ec\'edent.

\emph{$6^{\textrm{\`eme}}$ cas : $\X$ est un $S$-champ alg\'ebrique.}\\
On choisit une pr\'esentation $P:X\fleche \X$ par un espace alg\'ebrique. Alors 
$X\times_{\X} X$ et $X\times_{\X}
X\times_{\X} X$ sont encore des espaces alg\'ebriques, donc d'apr\`es la
discussion pr\'eliminaire le r\'esultat se d\'eduit du cinqui\`eme cas.
\end{demo}

\begin{cor}
\label{commutation_au_produit_fibre}
On garde les hypoth\`eses et notations du th\'eor\`eme pr\'ec\'edent, mais on
suppose de plus que $p$ est une extension infinit\'esimale (i.e. un morphisme
surjectif \`a noyau nilpotent), et non seulement un morphisme surjectif. On
suppose \'egalement que le morphisme $f : \X \fleche S$ est cohomologiquement
plat en dimension z\'ero. Alors les morphismes naturels
$$\Pic(\X_{B'}) \flechelongue \Pic(\X_{A'})\times_{\Pic(\X_{A})}\Pic(\X_{B})$$
et
$$P_{\X/S}(B') \flechelongue P_{\X/S}(A')\times_{P_{\X/S}(A)}P_{\X/S}(B)$$
sont des isomorphismes.
\end{cor}
\begin{demo}
On a un carr\'e commutatif
$$\xymatrix{\X_{B'} &\X_{A'}\ar[l]\\
\X_{B} \ar[u] & \X_{A}\ar[u] \ar[l]}$$
qui induit le premier morphisme \'evoqu\'e dans l'\'enonc\'e apr\`es
application du foncteur $\Pic$.
La surjectivit\'e r\'esulte tautologiquement du th\'eor\`eme pr\'ec\'edent.
Montrons que c'est un morphisme injectif. Soit $\L_{B'}$ un faisceau inversible
sur $\X_{B'}$ tel que ${\L_{B'}}_{|\X_{A'}}$ et
${\L_{B'}}_{|\X_{B}}$ soient triviaux. On veut montrer que $\L_{B'}$
est isomorphe à $\O_{\X_{B'}}$. Vu l'\'equivalence de cat\'egories
(\ref{equiv_cat_fi_pour_defo}), il suffit de trouver des isomorphismes
$$\flechen{{\L_{B'}}_{|\X_{A'}}}{\varphi_{A'}}{\O_{\X_{A'}}}$$
$$\flechen{{\L_{B'}}_{|\X_B}}{\varphi_B}{\O_{\X_B}}$$
tels que, aux isomorphismes canoniques pr\`es,
${\varphi_{A'}}_{|\X_{A}}$ et ${\varphi_B}_{|\X_{A}}$ soient égaux. On sait par hypoth\`ese
qu'il existe un $\varphi_{A'}$ et un $\varphi_B$, mais \emph{a priori}
${\varphi_{A'}}_{|\X_{A}}\neq{\varphi_B}_{|\X_{A}}$. Soit
$\psi=({\varphi_{A'}}_{|\X_{A}})({\varphi_B}_{|\X_{A}})^{-1}$. $\psi$ est un
isomorphisme de $\O_{\X_A}$, donc correspond \`a un \'el\'ement de
$\Gamma(\X_A,\O_{\X_A})^{\times}$. On voudrait que $\psi$ soit l'identit\'e de
$\O_{\X_A}$. Il suffirait pour cela de corriger $\varphi_{A'}$ (ou $\varphi_B$)
par un automorphisme de $\O_{\X_{A'}}$ (ou $\O_{\X_B}$) qui induise $\psi$ sur
$\O_{\X_A}$. Il suffit donc pour conclure de montrer que l'une des deux
fl\`eches
$$\Gamma(\X_{A'},\O_{\X_{A'}})^{\times}\flechelongue\Gamma(\X_A,\O_{\X_A})^{\times},$$
$$\Gamma(\X_B,\O_{\X_B})^{\times} \flechelongue \Gamma(\X_A,\O_{\X_A})^{\times},$$
est surjective. Or par platitude cohomologique en dimension z\'ero, on a
$\Gamma(\X_T,\O_{\X_T})=T$ pour toute $R$-alg\`ebre $T$, donc il suffit que
$A'^{\times} \fleche A^{\times}$ ou $B^{\times} \fleche A^{\times}$
soit surjective. Comme $A' \fleche A$ est une extension infinit\'esimale,
$A\simeq A'/I$ avec $I$ nilpotent et on voit facilement que $A'^{\times} \fleche
A^{\times}$ est surjective.

Passons maintenant au foncteur de Picard relatif. Il s'agit de montrer que le
morphisme naturel
$$\xymatrix{\frac{\text{Pic}(\X_{B'})}{\text{Pic}(B')} \ar[r]&
\frac{\text{Pic}(\X_{A'})}{\text{Pic}(A')}\times_{\frac{\text{Pic}(\X_{A})}{\text{Pic}(A)}}
\frac{\text{Pic}(\X_{B})}{\text{Pic}(B)}}$$
est un isomorphisme. On sait d\'ej\`a gr\^ace \`a la platitude cohomologique de
$f$ que pour tout $S$-sch\'ema $T$, $\Pic(T)$ s'injecte dans
$\Pic(\X\times_S T)$ (cf.
th. (\ref{comparaison_des_foncteurs_de_Picard})). On identifie dans la suite
$\Pic(T)$ \`a un sous-groupe de $\Pic(\X\times_S T)$.

\emph{Injectivit\'e :} Soit $b'\in\Pic(\X_{B'})$. On note $a'$, $b$, et $a$ ses
images dans $\Pic(\X_{A'})$, $\Pic(\X_B)$, et $\Pic(\X_{A})$, et on suppose que
$a'\in\Pic(A')$ et $b\in\Pic(B)$. Il s'agit de montrer que $b'\in\Pic(B')$. Les éléments $a'$ et $b$ ont m\^eme image dans $\Pic(A)$ donc gr\^ace \`a l'isomorphisme
$$\Pic(B')\simeq \Pic(A')\times_{\Pic(A)} \Pic(B)$$
il existe un \'el\'ement $\beta'$ de $\Pic(B')$ qui induit le couple $(a',b)$.
Maintenant $b'$ et $\beta'$ ont la m\^eme image dans
$\Pic(\X_{A'})\times_{\Pic(\X_{A})}\Pic(\X_{B})$ donc d'apr\`es ce qui
pr\'ec\`ede $b'=\beta'$.

\emph{Surjectivit\'e :} Soient $a'\in\Pic(\X_{A'})$ et $b\in\Pic(\X_B)$
induisant le m\^eme \'el\'ement dans $\Pic(\X_A)$ modulo $\Pic(A)$. Comme
$\Pic(A') \fleche \Pic(A)$ est surjectif (cf. lemme
(\ref{pic_extension_inf_est_isom})) on peut supposer, quitte \`a corriger
$a'$ par un \'el\'ement de $\Pic(A')$ bien choisi, que
$(a',b)\in\Pic(\X_{A'})\times_{\Pic(\X_A)}\Pic(\X_B)$. Alors il existe un
\'el\'ement $b'$ de $\Pic(\X_{B'})$ qui induit $(a',b)$, ce qui montre le
r\'esultat attendu.
\end{demo}

\begin{lem}
\label{pic_extension_inf_est_isom}
Si $A'\fleche A$ est une extension infinit\'esimale, alors le morphisme
$$\xymatrix{\Pic(A')\ar[r]& \Pic(A)}$$ est un isomorphisme.
\end{lem}
\begin{demo}
Signalons tout d'abord que l'on peut supposer que le noyau de $A'\fleche A$ est de carré nul, puisqu'une extension infinit\'esimale peut s'\'ecrire comme la
compos\'ee d'extensions infinit\'esimales \`a noyaux de carrés nuls.
On a une suite exacte de $A'$-modules :
$$0 \flechelongue I \flechelongue A' \flechelongue A \flechelongue 0$$
o\`u $I$ est le noyau de $A'\flechelongue A$. On en d\'eduit une suite exacte de
$\O_{\Spec A'}$-modules
$$0 \flechelongue \tilde{I} \flechelongue \O_{\Spec A'} \flechelongue \O_{\Spec A} \flechelongue 0,$$
d'o\`u une suite exacte de faisceaux de groupes ab\'eliens sur $\Spec A'$ (ou
sur $\Spec A$, n'oublions pas qu'ils ont le m\^eme espace topologique
sous-jacent) via l'application exponentielle $x\mapsto 1+x$ :
$$0 \flechelongue \tilde{I} \flechelongue \O_{\Spec A'}^{\times} \flechelongue
\O_{\Spec A}^{\times} \flechelongue 0.$$
La suite exacte longue de cohomologie associ\'ee, et la nullit\'e des groupes
$H^1(\Spec A', \tilde{I})$ et $H^2(\Spec A', \tilde{I})$ nous donnent le
r\'esultat.
\end{demo}

\section{D\'eformations de faisceaux inversibles}

Commen\c cons par une petite remarque d'alg\`ebre que nous utiliserons
abondamment par la suite et que, par commodit\'e, nous \'enon\c cons sous la
forme d'un lemme.
%\footnote{auquel d'ailleurs, vu sa trivialit\'e, adjoindre une
%d\'emonstration serait faire, d'une certaine mani\`ere, offense au lecteur.}

\begin{lem}
\label{defm_lemme_alg}
Soient $A$ un anneau et $I$ un id\'eal de carr\'e nul de $A$. On note $\pi$ la
projection canonique $\pi : A \fleche A/I$.\\
1) Le morphisme de groupes ab\'eliens $\pi^{\times} : A^{\times} \fleche
(A/I)^{\times}$ induit par $\pi$ est surjectif.\\
2) L'application $x\mapsto 1+x$ induit un isomorphisme de groupes ab\'eliens de
$I$ sur $\Ker \pi^{\times}$. $\square$
\end{lem}

Soit $\X$ un champ alg\'ebrique sur un sch\'ema $T$ et soit $\L$ un faisceau
inversible sur $\X$. On consid\`ere une immersion ferm\'ee
$$\xymatrix{i : \X \ar[r] & \Xt}$$
d\'efinie par un id\'eal quasi-coh\'erent $I$ de $\Xt$ de carr\'e nul.

\begin{remarque} \rm
Si $\X$ et $\Xt$ sont des champs de Deligne-Mumford, le morphisme $i$ induit une équivalence de sites entre les sites étales de $\X$ et de $\Xt$, ce qui permet d'identifier les topos étales de $\X$ et de $\Xt$. Il est alors évident que la catégorie des $\Oc_{\X}$-modules quasi-cohérents est équivalente à la catégorie des $\Oc_{\Xt}$-modules quasi-cohérents annulés par $I$. Lorsque l'on travaille avec des champs d'Artin (donc avec leurs sites lisses-étales) il faut être plus prudent. En effet le foncteur naturel du site lisse-étale de $\Xt$ vers celui de $\X$ n'est même plus fidèle si bien que ces derniers ne sont \emph{a priori} pas équivalents. Cependant, la descente fidèlement plate des modules quasi-cohérents (\emph{cf.} par exemple \ref{descente_fidelement_plate}) nous permet encore d'identifier la catégorie des $\Oc_{\X}$-modules à la catégorie des $\Oc_{\Xt}$-modules annulés par $I$. En particulier, l'idéal $I$ peut être vu comme un $\Oc_{\X}$-module. 
\end{remarque}

On note $\defm(\L)$ la cat\'egorie des d\'eformations de $\L$ \`a $\Xt$
d\'efinie de la mani\`ere suivante. Un objet de $\defm(\L)$ est un couple
$(\Lt,\lambda)$ o\`u $\Lt$ est un faisceau inversible sur $\Xt$ et $\lambda$ est
un isomorphisme $\xymatrix@C=1pc{\lambda : i^*\Lt\ar[r]^-{\sim}&\L}$. Un morphisme de
$(\Lt,\lambda)$ vers $(\widetilde{\Mc},\mu)$ est un isomorphisme $\alpha :
\raisebox{.7ex}{\xymatrix@C=1pc{\Lt\ar[r]^-{\sim}&\widetilde{\Mc}}}$ tel que
$\mu \circ i^*\alpha=\lambda$. On
note $\ov{\defm(\L)}$ l'ensemble des classes d'isomorphie de $\defm(\L)$.

Dans \cite{Aoki_Hom}, Aoki d\'efinit la
cat\'egorie $\defm_T(f)$ des d\'eformations d'un morphisme
de $T$-champs alg\'ebriques $f : \X\fleche \Y$ et en donne une description
compl\`ete au th\'eor\`eme 2.1.1, que nous rappelons ci-dessous. Nous renvoyons
\`a \emph{loc. cit.} pour les notations pr\'ecises et la d\'efinition de $\defm_T(f)$. 

\begin{thm}[Aoki]
\label{thm_defm_aoki}
\begin{enumerate}
\item[(1)] Il existe un \'el\'ement $\omega\in \Ext^1({\rm L}f^*{\rm L}_{\Y/T},I)$
dont l'annulation \'equivaut \`a l'existence d'une d\'eformation de $f$.
\item[(2)] Si $\omega=0$, alors $\ov{\defm_T(f)}$ est un torseur sous
$\Ext^0({\rm L}f^*{\rm L}_{\Y/T},I)$.
\item[(3)] Si $(\widetilde{f},\lambda)$ est une d\'eformation de $f$, son groupe
d'automorphismes est isomorphe \`a $\Ext^{-1}({\rm L}f^*{\rm L}_{\Y/T},I)$.
\end{enumerate}
\end{thm}

En appliquant ce th\'eor\`eme avec
$\Y=\bgm$, et en
tenant compte du fait que la cat\'egorie des morphismes de $\X$ vers $\bgm$ est
\'equivalente \`a la cat\'egorie des faisceaux inversibles sur $\X$, on en
d\'eduit une description de la cat\'egorie $\defm(\L)$ en termes des groupes
$\Ext^{i}({\rm L}f^*{\rm L}_{\bgm/T},I)$ (o\`u $f$ d\'esigne le morphisme
$\X \fleche \bgm$ correspondant au faisceau inversible $\L$ sur $\X$), que nous
allons promptement remplacer par les groupes de cohomologie de $\X$ \`a valeurs
dans $I$ gr\^ace au lemme suivant.

\begin{lem}
Pour tout $i\in \Z$, on a $\Ext^{i}(Lf^*L_{\bgm/T},I)=
{\rm H}^{i+1}(\X,I)$.
\end{lem}
\begin{demo}
Commençons par calculer le complexe cotangent de $\bgm$ sur $T$. On va montrer qu'il est représenté
par le complexe $\O_{\bgm}[-1]$, qui a pour seul terme non nul $\O_{\bgm}$ situé en degré 1.
Vu que $\bgm$ est lisse, on sait d'après \cite{LMB}~(17.8) que $L_{\bgm/T}$ est canoniquement isomorphe
à $L_{\Delta}[-1]$, où $L_{\Delta}$ est le complexe cotangent du morphisme diagonal
$$\Delta : \bgm \flechelongue \bgm \times_T \bgm.$$
Mais comme $\Delta$ est représentable et lisse, la proposition~(17.5.8) de \cite{LMB} montre que le système projectif $L_{\Delta}$
est essentiellement constant, et peut être représenté par le $\O_{\bgm}$-module quasi-cohérent
$\Omega_{\Delta}=\Omega_{\bgm/\bgm\times_T\bgm}$ placé en degré 0 (rappelons que par définition
$\Omega_{\Delta}$ est la limite dans $\Mod_{\text{\rm qcoh}}(\O_{\X})$ du ind-objet $\Hc^0(L_{\Delta})$).
Plus précisément, si $\Y \fleche \Zc$ est un morphisme
représentable et lisse, son complexe cotangent peut être calculé de la manière suivante. On choisit une présentation
$Z\fleche \Zc$ de $\Zc$, et on note $Y=\Y\times_{\Zc} Z$ l'espace algébrique obtenu par changement de base. Le
morphisme $Y\fleche \Y$ est automatiquement une présentation de $\Y$. On note $Y^{\bullet}$ (resp. $Z^{\bullet}$)
l'espace algébrique simplicial obtenu en prenant le 0-cosquelette du morphisme $Y\fleche \Y$ (resp. $Z\fleche \Zc$).
$$\xymatrix{Y^{\bullet} \ar[r] \ar[d] & Z^{\bullet}\ar[d]\\ \Y\ar[r]& \Zc}$$
\`A l'étage $i$, le morphisme $Y^i\fleche Z^i$ est lisse donc son complexe cotangent est représenté par le $\O_{Y^i}$-module
quasi-cohérent $\Omega_{Y^i/Z^i}$ placé en degré 0. Le $\O_{Y^{\bullet}}$-module simplicial
$\Omega_{Y^{\bullet}/Z^{\bullet}}$ se descend en un $\O_{\Y}$-module quasi-cohérent unique à unique isomorphisme
près (\cite{LMB}~(13.5.4) et~(13.5.5)), qui n'est autre que $\Omega_{\Y/\Zc}$. Pour le cas qui nous occupe,
on a $\Y=\bgm$ et $\Zc=\bgm\times_T\bgm$. Le morphisme $T\fleche \bgm$ correspondant au torseur trivial fournit
par composition avec le morphisme diagonal une présentation de $\bgm\times_T\bgm$. De plus on a un carré 2-cartésien
$$\xymatrix{Y=\gm \ar[r] \ar[d] \cartesien & Z=T\ar[d]\\ \bgm\ar[r]& \bgm\times_T\bgm}$$
dont on peut facilement calculer le 0-cosquelette. En effet, on a
\begin{eqnarray*}
   Y^1 &=& \gm\times_{\bgm}\gm \\
   &\simeq& \gm\times_T T\times_{\bgm}T\times_T\gm\\
   &\simeq& \gm\times_T\gm\times_T\gm = (\gm/T)^3
\end{eqnarray*}
puis par récurrence on en déduit $Y^i\simeq (\gm/T)^{2i+1}$. Par ailleurs,
\begin{eqnarray*}
   Z^1 &=& T\times_{(\bgm\times_T\bgm)}T \\
   &\simeq& (T\times_{(\bgm\times_T\bgm)}\bgm)\times_{\bgm}T\\
   &\simeq& \gm\times_{\bgm}T  \\
   &\simeq& \gm\times_T (T\times_{\bgm}T) \\
   &\simeq& \gm\times_T\gm = (\gm/T)^2
\end{eqnarray*}
puis par récurrence on en déduit $Z^i\simeq (\gm/T)^{2i}$.
\`A l'étage $i$, le morphisme $Y^i \fleche Z^i$ s'identifie donc à $(\gm/T)^{2i+1}\fleche (\gm/T)^{2i}$,
autrement dit c'est le morphisme structural ${\gm}_{,S_i} \fleche S_i$ du groupe $\gm$ sur la base
$S_i=(\gm/T)^{2i}$, et son module des différentielles relatives est canoniquement isomorphe à $\O_{(\gm/T)^{2i+1}}$.
En conséquence, le $\O_{Y^{\bullet}}$-module simplicial $\Omega_{Y^{\bullet}/Z^{\bullet}}$ est isomorphe à
$\O_{Y^{\bullet}}$, si bien que $\Omega_{\Y/\Zc}$ est isomorphe à $\O_{\Y}$, ce qui achève le calcul du
complexe cotangent de $\bgm$.

Maintenant, si $f$ est un morphisme de $\X$ dans $\bgm$,
le système projectif $Lf^*L_{\bgm/T}$ est représenté
par le complexe $\O_{\X}[-1]$ qui a pour seul terme non nul $\O_{\X}$ situé en degré~1.
On en déduit, par définition,
\begin{eqnarray*}
   \Ext^{i}(Lf^*L_{\bgm/T},I)&=& \Ext^i(\O_{\X}[-1],I)\\
   &=& \Hom_{D_{\text{\rm qcoh}}(\O_{\X})}(\O_{\X}[-1],I[i]) \\
   &=& \Hom_{D_{\text{\rm qcoh}}(\O_{\X})}(\O_{\X},I[i+1]) \\
   &=& \Ext^{i+1}(\O_{\X},I)\\
   &=& H^{i+1}(\X,I)
\end{eqnarray*}
\end{demo}

Le th\'eor\`eme ci-dessus devient alors :

\begin{thm}
\label{thm_defm_fi}
\begin{enumerate}
\item[(1)] Il existe un \'el\'ement $\omega\in H^2(\X,I)$
dont l'annulation \'equivaut \`a l'existence d'une d\'eformation de $\L$ \`a
$\Xt$.
\item[(2)] Si $\omega=0$, alors $\ov{\defm(\L)}$ est un torseur sous
$H^1(\X,I)$.
\item[(3)] Si $(\Lt,\lambda)$ est une d\'eformation de $\L$, son groupe
d'automorphismes est isomorphe \`a $H^0(\X,I)$.
\end{enumerate}
\end{thm}

Comme on peut s'en douter, le th\'eor\`eme (\ref{thm_defm_fi}), qui n'est qu'un
cas particulier de (\ref{thm_defm_aoki}), est en r\'ealit\'e beaucoup plus
facile \`a obtenir que ce dernier. De plus la d\'emonstration du th\'eor\`eme
(\ref{thm_defm_aoki}) fait appel \`a un certain nombre de r\'esultats non
triviaux concernant les d\'eformations de morphismes d'espaces alg\'ebriques
(\cite{Illusie_CCD}), les d\'eformations de morphismes repr\'esentables de
champs alg\'ebriques (\cite{Olsson_defm}), les d\'eformations de champs
alg\'ebriques (\cite{Aoki_defm}) et la th\'eorie du complexe
cotangent (\cite{Illusie_CCD}, \cite{LMB} et
\cite{Olsson_Sheaves_on_Artin_stacks}). Il nous a donc paru n\'ecessaire de
donner une d\'emonstration directe, et plus \'el\'ementaire que celle qui
consiste \`a faire appel \`a \cite{Aoki_Hom}, du th\'eor\`eme
(\ref{thm_defm_fi}). Commen\c cons par en d\'emontrer le troisi\`eme point (pour
mesurer l'ampleur des simplifications, le lecteur pourra comparer avec la
d\'emonstration du r\'esultat analogue pour les d\'eformations de morphismes,
qui est essentiellement la proposition 2.8 de \cite{Aoki_Hom}).

\vskip .3cm \noindent {\bf D\'emonstration de (\ref{thm_defm_fi}), (3)}
Par d\'efinition un automorphisme de $(\Lt,\lambda)$ est un automorphisme
$\varphi$ de $\Lt$ tel que $i^*\varphi=\Id_{i^*\Lt}$. Autrement dit,
$\Aut_{\defm(\L)}(\Lt,\lambda)$ est le noyau du morphisme $\Aut(\Lt) \fleche
\Aut(\L)$ induit par $i^*$. Or $\Aut(\Lt)=\Gamma(\Xt,\O_{\Xt})^{\times}$ et
$\Aut(\L)=\Gamma(\X,\O_{\X})^{\times}$. Donc
$\Aut_{\defm(\L)}(\Lt,\lambda)$ est le noyau du morphisme
$$H^0(\Xt,\O_{\Xt})^{\times} \flechelongue
\left(\frac{H^0(\Xt,\O_{\Xt})}{H^0(\Xt,I)}\right)^{\times}.$$
On en d\'eduit d'apr\`es le lemme (\ref{defm_lemme_alg}) 2), que
le groupe $\Aut_{\defm(\L)}(\Lt,\lambda)$ est isomorphe \`a $H^0(\Xt,I)$ via
l'application $x\mapsto 1+x$.
$\square$ \vskip .3cm

Nous allons en fait donner deux d\'emonstrations des parties (1) et (2) du
th\'eor\`eme (\ref{thm_defm_fi}). La premi\`ere d\'emonstration propos\'ee
ci-dessous est une traduction dans le langage des faisceaux inversibles
de la d\'emonstration du th\'eor\`eme 2.1.1 de \cite{Aoki_Hom}. Un certain
nombre d'arguments se simplifient, voire deviennent inutiles, et les outils
invoqu\'es diff\`erent sensiblement, mais le canevas
sous-jacent est le m\^eme :
\begin{itemize}
\item prendre une pr\'esentation $X_0$ de
$\X$ et l'espace alg\'ebrique simplicial associ\'e \`a cette pr\'esentation ;
\item relier les d\'eformations dans la cat\'egorie des champs alg\'ebriques
\`a celles dans la cat\'egorie des espaces alg\'ebriques simpliciaux ;
\item comparer les groupes Ext (dans le cas des morphismes) ou les groupes de
cohomologie (dans le cas des faisceaux inversibles).
\end{itemize}

La seconde d\'emonstration, plus rapide, g\'en\'eralise au cas des champs
alg\'ebriques l'argumentaire classique valable dans le cadre des espaces
alg\'ebriques.

\vskip .3cm \noindent {\bf Premi\`ere d\'emonstration de (\ref{thm_defm_fi}),
(1) et (2) :}
Commençons par une remarque générale. Soit $X^0\fleche \X$ une présentation de
$\X$. Alors la catégorie des faisceaux inversibles sur $\X$ est équivalente à
la catégorie suivante. Un objet est un couple $(\L^0,\alpha)$ où $\L^0$ est un
faisceau inversible sur $X^0$ et
$\xymatrix@C=1pc{\alpha : p_1^*\L^0\ar[r]^-{\sim}& p_2^*\L^0}$
est un isomorphisme tel que, à des isomorphismes canoniques près,
$(p_{23}^*\alpha) \circ (p_{12}^*\alpha) = p_{13}^*\alpha$. Un (iso)morphisme
$(\L^0,\alpha) \fleche (\Mc^0, \beta)$ est un isomorphisme $\xymatrix@C=1pc{\gamma :
\L^0\ar[r]^-{\sim}&\Mc^0}$ tel que $(p_2^*\gamma)\circ \alpha= \beta\circ
(p_1^*\gamma)$.

On choisit une présentation $P : X^0 \fleche \X$ de $\X$ telle que
$\L^0:=P^*\L$ soit trivial et telle que $X^0$ soit union disjointe de schémas
affines. D'après le lemme (\ref{defm_pres}) ci-dessous,
il existe alors un morphisme lisse et surjectif
$\widetilde{P} : \Xzt \fleche \Xt$ et un diagramme 2-cartésien :
$$\xymatrix{X^0\ar[r]^{i_0} \ar[d]_{P} \cartesien &
   \Xzt \ar[d]^{\widetilde{P}} \\
   \X \ar[r]^i & \Xt.}$$
Quitte à prendre une présentation de $\Xzt$, on peut supposer que $\Xzt$ est
union disjointe de sch\'emas affines, et donc que $i_0$ est une somme disjointe
d'immersions ferm\'ees entre des sch\'emas affines.
Soit $\xymatrix@C=1pc{\alpha : p_1^*\L^0\ar[r]^-{\sim}& p_2^*\L^0}$ l'isomorphisme
correspondant à $\L$ via l'équivalence de catégories ci-dessus (disons qu'on
fixe une fois pour toute une telle équivalence de catégories ainsi qu'un
quasi-inverse). Afin d'alléger un peu les notations, on supposera dans la suite
que le faisceau $\L^0$ est \emph{égal}, et non seulement isomorphe, à
$\O_{X^0}$. (Ceci est un abus clairement inoffensif. En toute rigueur il
faudrait choisir une fois pour toute un isomorphisme entre $\O_{X^0}$ et $\L^0$
et le transporter tout au long de la démonstration.)
On note alors $\defm(\alpha)$ l'ensemble des isomorphismes
$\xymatrix@C=1pc{\at : p_1^*\O_{\Xzt}\ar[r]^-{\sim}&p_2^*\O_{\Xzt}}$
tels que
$$\left\{ \begin{array}{l}
(p_{23}^*\at) \circ (p_{12}^*\at) = p_{13}^*\at \\
\rm{et} \quad \at_{|_{X^1}}=i_1^*\at=\alpha
\end{array}\right.$$
(à des isomorphismes canoniques près, que nous omettrons pour rester lisible),
où $X^1=X^0\times_{\X}X^0$, $\Xut=\Xzt\times_{\Xt}\Xzt$, et $i_1$ est le
morphisme induit $X^1 \fleche \Xut$.

Il existe alors une application $A : \defm(\alpha) \fleche \ov{\defm(\L)}$
qui à $\at$ associe la classe du couple $(\Lt,\lambda)$ constitué du faisceau
inversible $\Lt$ sur $\Xt$ associé à $(\O_{\Xzt},\at)$ via l'équivalence de
catégories mentionnée plus haut, et de l'isomorphisme $\lambda : i^*\Lt \fleche
\L$ induit par l'isomorphisme canonique $i_0^*\O_{\Xzt}\fleche \O_{X^0}$.

\begin{etape}{L'application $A$ ainsi construite est surjective.}
Soit $(\Lt,\lambda)$ une d\'eformation de $\L$. Soit
$\widetilde{\L^0}=\widetilde{P}^*\Lt$. On a un isomorphisme $P^*\lambda :
i_0^*\widetilde{\L^0} \fleche \L^0=\O_{X^0}$, donc en vertu du lemme
(\ref{pic_extension_inf_est_isom}) le faisceau
inversible $\widetilde{\L^0}$ est trivial. On fixe un isomorphisme $\varphi :
\widetilde{\L^0} \fleche \O_{\Xzt}$. Soit $\at : p_1^*\widetilde{\L^0}
\fleche p_2^*\widetilde{\L^0}$ l'isomorphisme correspondant \`a $\Lt$.
On pose $\widehat{\alpha}=(p_2^*\varphi)\circ \at \circ (p_1^*\varphi)^{-1} :
p_1^*\O_{\Xzt} \fleche p_2^*\O_{\Xzt}$. Il est
clair que $\widehat{\alpha}$ appartient \`a $\defm(\alpha)$. Soit
$A(\widehat{\alpha})=(\widehat{\L},\widehat{\lambda})$ la d\'eformation
correspondante.
On va montrer que, pour un choix de $\varphi$ convenable, $(\Lt,\lambda)$ et
$(\widehat{\L},\widehat{\lambda})$ sont isomorphes.

L'isomorphisme $\varphi$, par d\'efinition de $\widehat{\alpha}$, fait commuter
le diagramme
$$\xymatrix{p_1^*\widetilde{\L^0} \ar[r]^{p_1^*\varphi} \ar[d]_{\at} &
p_1^*\O_{\Xzt} \ar[d]^{\widehat{\alpha}}\\
p_2^*\widetilde{\L^0} \ar[r]_{p_2^*\varphi} & p_2^*\O_{\Xzt}}$$
donc il induit un isomorphisme
$\xymatrix@C=1pc{\Lt\ar[r]^-{\sim}&\widehat{\L}}$.
C'est un isomorphisme de d\'eformations si et seulement si
$\widehat{\lambda}\circ i^*\varphi=\lambda$. Vu la construction de
l'isomorphisme $\xymatrix@C=1pc{\Lt\ar[r]^-{\sim}&\widehat{\L}}$, cela revient pr\'ecis\'ement
\`a exiger que $i_0^*\varphi=P^*\lambda$. Or le morphisme de groupes
$\Aut(\O_{\Xzt})\fleche \Aut(\O_{X^0})$ est surjectif en vertu du lemme
(\ref{defm_lemme_alg}) 1),
donc on peut toujours choisir un tel $\varphi$, ce qui prouve notre assertion.
\end{etape}

On note maintenant $K$ le noyau
$$K=\Ker(\xymatrix{H^0(X^1,I^1) \ar[rr]^{p_{23}^*-p_{13}^*+p_{12}^*} &&
H^0(X^2,I^2)}),$$
o\`u $X^2=X^0\times_{\X}X^0\times_{\X}X^0$, et $I^1$, $I^2$ sont les images
inverses de $I$ sur $X^1$, $X^2$.
Ce groupe agit sur $\defm(\alpha)$ de la manière suivante : si $x\in K$ et
$\at\in\defm(\alpha)$, on pose
$$x.\at = \mu_{1+x}\circ \at,$$
où $\mu_{1+x}$ désigne l'automorphisme de $p_2^*\O_{\Xzt}$ de multiplication par
$1+x$. 

\begin{etape}{Muni de cette action, $\defm(\alpha)$
est soit vide, soit un torseur sous $K$.}
On suppose que $\defm(\alpha)$ est non vide. C'est alors clairement un torseur
sous le groupe $B$ des automorphismes $\beta\in\Aut(\O_{\Xut})$ tels que
$\beta_{|_{X^1}}=\Id_{\O_{X^1}}$ et $(p_{23}^*\beta)\circ (p_{12}^*\beta)=
p_{13}^*\beta$. Le lemme (\ref{defm_lemme_alg}), 2) montre que l'application
$x\mapsto 1+x$ induit un
isomorphisme de $H^0(X^1,I^1)$ sur le groupe des $\beta\in\Aut(\O_{\Xut})$ tels
que $\beta_{|_{X^1}}=\Id_{\O_{X^1}}$. En consid\'erant sa restriction \`a $K$ on
obtient un isomorphisme de $K$ sur $B$, d'o\`u le r\'esultat.
\end{etape}

On a par ailleurs un morphisme $\flechen{D : H^0(X^0,I^0)}{p_1^*-p_2^*}{K}$, qui
induit donc une action de $H^0(X^0,I^0)$ sur $\defm(\alpha)$.

\begin{etape}{Les fibres de $A$ sont les orbites de $\defm(\alpha)$ sous
$H^0(X^0,I^0)$.}
Soient $\at_1, \at_2\in\defm(\alpha)$. Il s'agit de montrer que $\at_1$ et
$\at_2$ sont dans la m\^eme orbite sous $H^0(X^0,I^0)$ si et seulement si
$A(\at_1)=A(\at_2)$. On note $(\Lt_1,\lambda_1),(\Lt_2,\lambda_2)$ les
d\'eformations de $\L$ associ\'ees. Un isomorphisme
$(\Lt_1,\lambda_1)\fleche(\Lt_2,\lambda_2)$ est un isomorphisme
$\xymatrix@C=1pc{\Lt_1 \ar[r]^{\gamma} & \Lt_2}$ tel que $\lambda_2\circ i^*\gamma=\lambda_1$.
Un tel isomorphisme correspond \`a un isomorphisme $\gamma : \O_{\Xzt} \fleche
\O_{\Xzt}$ tel que les diagrammes
$$\xymatrix{p_1^*\O_{\Xzt} \ar[r]^{p_1^*\gamma} \ar[d]_{\at_1} &
p_1^*\O_{\Xzt} \ar[d]^{\at_2} \\
p_2^*\O_{\Xzt} \ar[r]^{p_2^*\gamma} & p_2^*\O_{\Xzt}}
\qquad \rm{et} \qquad
\xymatrix{i_0^*\O_{\Xzt} \ar[rr]^{i_0^*\gamma} \ar[rd]_{\rm{can.}} &&
i_0^*\O_{\Xzt} \ar[ld]^{\rm{can.}} \\
& \L^0=\O_{X^0}}$$
commutent. Le second diagramme signifie simplement que
$i_0^*\gamma=\Id_{i_0^*\O_{\Xzt}}$. Cela équivaut encore \`a dire que
$\gamma$ est de la forme $\mu_{1+x}$ o\`u $x\in H^0(X^0,I^0)$. Donc
$A(\at_1)=A(\at_2)$ si et seulement s'il existe $x\in H^0(X^0,I^0)$ tel que
$\at_2=(p_2^*\mu_{1+x})\circ\at_1 \circ (p_1^*\mu_{1+x})^{-1}
=(p_2^*x-p_1^*x).\at_1$, i.e. si et seulement si $\at_1$ et
$\at_2$ sont dans la m\^eme orbite sous $H^0(X^0,I^0)$.
\end{etape}

Le fait que $A$ soit surjective et que ses fibres s'identifient aux orbites de
$\defm(\alpha)$ (qui, ne l'oublions pas, est soit vide soit un torseur sous
$K$) sous $H^0(X^0,I^0)$ montre que, lorsqu'il est non vide,
$\ov{\defm(\L)}$ est naturellement un torseur sous le conoyau du morphisme
$D : H^0(X^0,I^0) \fleche K$. Or, comme $P : X^0 \fleche \X$ est une
présentation de $\X$, la proposition (\ref{suite_spectrale_de_descente}) de
l'annexe montre que l'on a une suite spectrale (cf. \emph{loc. cit.} pour les
notations) :
$$E_2^{p,q}=\check{H}^p
(H^q(X^{\bullet},I^{\bullet}))
 \Rightarrow H^{p+q}(\X,I).$$
La suite exacte des termes de bas degré associée à cette suite spectrale est la
suivante (\cite{Cartan_Eilenberg}, chapitre XV, paragraphe 5, p.329):
\begin{equation}
0 \fleche E^{1,0}_2\fleche H^1 \fleche E^{0,1}_2 \fleche E^{2,0}_2 \fleche H^2.
\end{equation}
Or, comme $X^0$ est union disjointe de schémas affines, le groupe
$H^1(X^0,I^0)$ est trivial, ce qui montre que le terme $E^{0,1}_2$ de la suite
exacte ci-dessus est nul. On obtient donc d'une part une injection de
$\check{H}^2(H^0(X^{\bullet},I^{\bullet}))$ dans $H^2(\X,I)$, dont nous nous
servirons bientôt, et
d'autre part un isomorphisme entre $H^1(\X,I)$ et
$\check{H}^1(H^0(X^{\bullet},I^{\bullet}))$. Ce dernier groupe n'étant autre que
le conoyau de $D$, on en déduit la deuxième assertion du théorème.

\begin{etapefinale}{Classe d'obstruction.}
Vu que $A$ est surjective, $\ov{\defm(\L)}$ est vide si et seulement si
$\defm(\alpha)$ l'est. Il ne nous reste donc plus qu'à construire 
un élément $\omega$ de $H^2(\X,I)$ dont la nullité caractérise
l'existence d'une déformation $\at\in\defm(\alpha)$. Or, via l'identification
$\Aut(\O_{\Xzt}) \simeq \Gamma(\Xzt,\O_{\Xzt})^{\times}$, chercher un tel $\at$
revient à chercher $\at\in \Gamma(\Xut,\O_{\Xut})^{\times}$ tel que
$$\left\{ \begin{array}{l}
i_1^*\at=\alpha \in \Gamma(X^1,\O_{X^1})^{\times} \\
(p_{23}^*\at)(p_{12}^*\at)(p_{13}^*\at)^{-1}=1 \in
\Gamma(\Xdt,\O_{\Xdt})^{\times}
\end{array}\right.$$
On rappelle que
$\Gamma(X^1,\O_{X^1})=\frac{\Gamma(\Xut,\O_{\Xut})}{H^0(X^1,I^1)}$.
Choisissons un antécédent $\beta\in\Gamma(\Xut,\O_{\Xut})$ de $\alpha$. Alors
$\beta$ est automatiquement inversible, et les équations ci-dessus deviennent :
$$\left\{ \begin{array}{l}
i_1^*(\at\beta^{-1})=1 \in \Gamma(X^1,\O_{X^1})^{\times} \\
(p_{23}^*(\at\beta^{-1}))(p_{12}^*(\at\beta^{-1}))(p_{13}^*(\at\beta^{-1}))^{-1}
=(p_{23}^*\beta)(p_{12}^*\beta)(p_{13}^*\beta)^{-1}=:\xi
\end{array}\right.$$
Le fait que $\alpha$ vérifie la condition de cocycle
$(p_{23}^*\alpha)(p_{12}^*\alpha)(p_{13}^*\alpha)^{-1}=1$
montre que $i_2^*\xi=1$,
ce qui signifie d'après le lemme (\ref{defm_lemme_alg}) que $\xi$ est de la
forme $1+u$ avec $u\in H^0(X^2,I^2)$.
Posant $\gamma=\at\beta^{-1}$ les équations deviennent
$$\left\{ \begin{array}{l}
i_1^*\gamma=1 \in \Gamma(X^1,\O_{X^1})^{\times} \\
(p_{23}^*\gamma)(p_{12}^*\gamma)(p_{13}^*\gamma)^{-1}=1+u \in
\Gamma(\Xdt,\O_{\Xdt})^{\times}
\end{array}\right.$$
Toujours d'après (\ref{defm_lemme_alg}) l'ensemble des éléments qui vérifient la
première équation
s'identifie à $H^0(X^1,I^1)$ via $x\mapsto \gamma=1+x$. Donc chercher un
élément $\at\in\defm(\alpha)$ revient à chercher un $x\in H^0(X^1,I^1)$ tel
que $(p_{23}^*-p_{13}^*+p_{12}^*)x=u$. Par ailleurs, vu l'expression de
$\xi=(p_{23}^*\beta)(p_{12}^*\beta)(p_{13}^*\beta)^{-1}$, on a clairement
$(p_{234}^*-p_{134}^*+p_{124}^*-p_{123}^*)u=0$. On note alors $\omega$ la classe
de $u$ dans
$$\check{H}^2(H^0(X^{\bullet},I^{\bullet}))=
\frac{\Ker(p_{234}^*-p_{134}^*+p_{124}^*-p_{123}^*)}
{\Im(p_{23}^*-p_{13}^*+p_{12}^*)}.$$
L'injection $\check{H}^2(H^0(X^{\bullet},I^{\bullet})) \inj H^2(\X,I)$ obtenue
ci-dessus à partir de la suite spectrale de descente cohomologique permet de
voir $\omega$ comme un élément de $H^2(\X,I)$, et la discussion que nous venons
de mener montre que $\omega=0$ si et seulement si $\L$ admet une déformation à
$\Xt$.
\end{etapefinale}
$\square$ \vskip .3cm

\begin{lem}
\label{defm_pres}
Soit $i : \X \fleche \Xt$ une immersion ferm\'ee de champs alg\'ebriques
d\'efinie par un id\'eal quasi-coh\'erent $I$ de $\Xt$ de carr\'e nul (ou
nilpotent). Soit $P : X^0\fleche \X$ un morphisme lisse, o\`u $X^0$ est union
disjointe de sch\'emas affines. Alors il
existe un morphisme lisse de champs alg\'ebriques
$\widetilde{P} : \Xzt\fleche \Xt$ et un diagramme 2-cart\'esien :
$$\xymatrix{X^0\ar[r]^{i_0} \ar[d]_{P} \cartesien &
    \Xzt \ar[d]^{\widetilde{P}} \\
    \X \ar[r]^i & \Xt.}$$
Si de plus $P$ est surjectif (resp. \'etale), il en est de m\^eme
de $\widetilde{P}$.
\end{lem}
\begin{demo}
D'apr\`es \cite{Olsson_defm}, thm. 1.4, il existe $o\in \Ext^2(L_{X^0/\X},P^*I)$
tel que $o=0$ si et seulement s'il existe un diagramme 2-cart\'esien comme
ci-dessus avec $\widetilde{P}$ plat. Comme $P$ est lisse, d'apr\`es \cite{LMB},
(17.5.8), le complexe cotangent $L_{X^0/\X}$ est quasi-isomorphe au
$\O_{X^0}$-module quasi-coh\'erent $\Omega^1_{X^0/\X}$, qui est localement libre
de rang fini. On a donc :
\begin{eqnarray*}
\Ext^i(L_{X^0/\X},P^*I) &=& \Ext^i(\Omega^1_{X^0/\X},P^*I)\\
   &=& \Ext^i(\O_{X^0},(\Omega^1_{X^0/\X})^{\vee}\otimes P^*I)\\
   &=& H^i(X^0,(\Omega^1_{X^0/\X})^{\vee}\otimes P^*I)\\
   &=& 0
\end{eqnarray*}
pour $i>0$ puisque $X^0$ est union disjointe de sch\'emas affines. Il reste juste à montrer que $\widetilde{P}$ est lisse, et qu'il est m\^eme \'etale
(resp. surjectif) d\`es que $P$ l'est. Or $\widetilde{P}$ est repr\'esentable
d'apr\`es \cite{Olsson_defm}, lemme (2.1).
Par descente fid\`element plate, on peut supposer, quitte \`a faire un
changement de base par une pr\'esentation de $\Xt$, que $\Xt$ est un espace
alg\'ebrique. Le r\'esultat d\'ecoule alors de \cite{Aoki_defm}, lemme (2.2.4).
\end{demo}

\begin{remarque}\rm
En réalité, dans le lemme qui précède, $\Xzt$ est automatiquement une union disjointe de schémas affines. Nous n'aurons pas besoin de ce fait par la suite.
\end{remarque}

%Afin de faciliter la comparaison entre cette preuve et celle de
%\cite{Aoki_Hom}, nous avons
%dress\'e un dictionnaire permettant de passer de l'une \`a l'autre, et de voir
%ce que sont devenus les objets apparaissant dans la d\'emonstration originale.
%
%TABLEAU

La seconde d\'emonstration que nous proposons repose sur une g\'en\'eralisation
au cas des champs alg\'ebriques des arguments avanc\'es par Artin dans
\cite{Global_Analysis_1}. Elle est plus directe que la pr\'ec\'edente dans la
mesure o\`u elle ne repose pas sur une r\'eduction au cas des
espaces alg\'ebriques. Elle n\'ecessite cependant un petit travail technique
pour relier la cohomologie des faisceaux ab\'eliens sur $\X$ et sur $\Xt$ (ce
qui \'etait trivial lorsque c'\'etaient des sch\'emas puisqu'ils avaient le
m\^eme espace topologique sous-jacent).

\vskip .3cm \noindent {\bf Seconde d\'emonstration de (\ref{thm_defm_fi}),
(1) et (2) :}

Montrons d'abord que $\ov{\defm(\L)}$ est isomorphe \`a l'ensemble 
$\Pic_{[\L]}(\Xt)$ des \'el\'ements de $\Pic(\Xt)$ qui
sont envoy\'es sur $[\L]\in \Pic(\X)$. On a une application naturelle
$\phi : \ov{\defm(\L)} \fleche \Pic_{[\L]}(\Xt)$ qui
\`a une d\'eformation $(\Lt,\lambda)$ associe la classe de $\Lt$ dans
$\Pic(\Xt)$. Elle est clairement surjective par d\'efinition m\^eme de
$\defm(\L)$. Montrons qu'elle est injective. Soient $(\Lt,\lambda)$ et
$(\Mt,\mu)$ tels que $[\Lt]=[\Mt]$. Autrement dit, il existe un isomorphisme
$\alpha : \Lt \fleche \Mt$. On veut montrer que l'on peut choisir $\alpha$ de
telle sorte que $i^*\alpha=\mu^{-1} \circ \lambda$. Il suffit pour cela de voir
que le morphisme de groupes $\Aut(\Lt) \fleche \Aut(i^*\Lt)$ est surjectif : il
n'y aura plus alors qu'\`a corriger $\alpha$ par un automorphisme convenablement
choisi de $\Lt$. Or ce morphisme de groupes s'identifie \`a
$\Aut(\O_{\Xt})\fleche\Aut(\O_{\X})$, qui est surjectif en vertu de
(\ref{defm_lemme_alg}), 1).

L'ensemble $\ov{\defm(\L)}$, qui s'identifie \`a $\Pic_{[\L]}(\Xt)$, est donc
naturellement un torseur sous le noyau du morphisme $\Pic(\Xt)\fleche \Pic(\X)$.
Calculons ce noyau. On a une suite exacte de faisceaux quasi-coh\'erents sur
$\Xt$ :
$$0 \flechelongue I \flechelongue \O_{\Xt} \flechelongue i_*\O_{\X} \flechelongue 0.$$
Le morphisme $\O_{\Xt} \fleche i_*\O_{\X}$ est donn\'e sur un ouvert
lisse-\'etale $(U,u)$ de $\Xt$ par
$$\xymatrix{\Gamma((U,u),\O_{\Xt}) \ar[r] & \frac{\Gamma((U,u),\O_{\Xt})}{\Gamma((U,u),I)}}.$$
Il induit donc d'apr\`es
(\ref{defm_lemme_alg}) un morphisme surjectif
$$(\Gamma((U,u),\O_{\Xt}))^{\times}
\flechelongue \left(\frac{\Gamma((U,u),\O_{\Xt})}{\Gamma((U,u),I)}\right)^{\times}$$
dont le noyau s'identifie via l'application exponentielle \`a $\Gamma((U,u),I)$.
En d'autres termes, on obtient une suite exacte de faisceaux de groupes
ab\'eliens sur $\Xt$:
$$0 \flechelongue I \flechelongue \O_{\Xt}^{\times} \flechelongue i_*\O_{\X}^{\times} \flechelongue 0.$$
La suite exacte longue de cohomologie associ\'ee nous donne :
$$H^1(\Xt, I) \flechelongue H^1(\Xt,\O_{\Xt}^{\times}) \flechelongue
H^1(\Xt, i_*\O_{\X}^{\times}) \flechelongue H^2(\Xt, I).$$
Or d'apr\`es \ref{pic_egal_h1}, le groupe $H^1(\Xt,\O_{\Xt}^{\times})$ est
isomorphe \`a $\Pic(\Xt)$. De plus le lemme \ref{coh_et_ext_inf}
fourni en annexe
montre que $H^1(\Xt, i_*\O_{\X}^{\times})\simeq H^1(\X,\O_{\X}^{\times})\simeq
\Pic(\X)$. Montrons enfin que la premi\`ere fl\`eche de la suite exacte longue
ci-dessus est injective. Il suffit pour cela de voir que le morphisme
$H^0(\Xt,\O_{\Xt}^{\times}) \fleche
H^0(\Xt, i_*\O_{\X}^{\times})$ est surjectif, ce qui r\'esulte encore du lemme
(\ref{defm_lemme_alg}). On obtient donc la suite exacte 
$$0\flechelongue H^1(\Xt, I) \flechelongue \Pic(\Xt) \flechelongue
\Pic(\X) \flechelongue H^2(\Xt, I),$$
ce qui ach\`eve notre d\'emonstration.
$\square$ \vskip .3cm

\begin{remarque} \rm
\label{rem_noyau}
Nous avons montr\'e au passage que le noyau de $\Pic(\Xt)\fleche \Pic(\X)$
s'identifie \`a $H^1(\Xt, I)$.
\end{remarque}

\section{Lissit\'e et dimension}

Soient $k$ un corps, $P$ un $k$-espace alg\'ebrique, et $x : \Spec k \fleche P$
un $k$-point de $P$. Soit $k[\eps]=k[X]/(X^2)$. On note $i : \Spec k \fleche
\Spec k[\eps]$ l'immersion ferm\'ee d\'efinie par le morphisme qui envoie $\eps$
sur 0. L'espace tangent \`a $P$ en $x$, not\'e $T_xP$ est par d\'efinition un
$k$-espace vectoriel dont l'ensemble sous-jacent est l'ensemble des 
morphismes de $\Spec k[\eps]$ vers $P$ qui induisent $x$ par composition avec
$i$. Autrement dit c'est l'ensemble des \'el\'ements de $P(k[\eps])$ qui
sont envoy\'es sur $x$ par la fl\`eche $P(k[\eps])\fleche P(k)$ induite par
$i$, muni d'une structure de $k$-espace vectoriel naturelle, que nous
n'expliciterons pas ici (pour plus de d\'etails, voir~\cite{poly_Kleiman}~(5.11)). 

Si $f : X\fleche P$ est un morphisme de $k$-espaces alg\'ebriques, et si $x$ est
un $k$-point de $X$, on note $T_xf : T_xX \fleche T_{f(x)}P$ le 
morphisme naturel :
$$\fonction{T_xf}{T_xX}{T_{f(x)}P=\{\psi : \Spec k[\eps] \fleche P\
|\ \psi \circ i
=f\circ x\}}{\varphi}{f\circ \varphi.}$$
C'est un morphisme de $k$-espaces vectoriels. Nous rappelons ci-dessous quelques
r\'esultats bien utiles dans l'\'etude des espaces tangents.

\begin{lem}
\label{lemme_espace_tangent}
\begin{itemize}
\item[1)] Soit $f : X \fleche P$ un morphisme formellement lisse (resp.
formellement non ramifi\'e, formellement \'etale) de $k$-espaces alg\'ebriques
et soit $x : \Spec k \fleche X$ un point $k$-rationnel de $X$. Alors
$T_xf : T_xX \fleche T_{f(x)}P$ est surjective (resp. injective, resp. un
isomorphisme).
\item[2)] Soient $X$ un $k$-espace alg\'ebrique localement de type fini, $x$ un
$k$-point de $X$, et $L$ une
extension de $k$. On note $X_L$ le $L$-espace alg\'ebrique obtenu par changement
de base, et $x_L$ le point de $X_L$ induit par $x$.
Alors le morphisme naturel
$$(T_xX)\otimes_k L \flechelongue T_{x_L}(X_L)$$
est un isomorphisme.
\end{itemize}
\end{lem}
\begin{demo}
1) 
Dire que $T_xf$ est surjective (resp. injective, resp. un
isomorphisme) revient \`a dire qu'\'etant donn\'e un diagramme
commutatif en traits pleins :
$$\xymatrix{\Spec k \ar@{^(->}[r]^i \ar[d]_x & \Spec k[\eps]\ar[d]^{\psi}
\ar@{.>}[ld]_{\varphi}\\
X \ar[r]^f &P,}$$
il existe une (resp. au plus une, resp. une unique) fl\`eche $\varphi$ en
pointill\'es qui le rend commutatif.
Ceci r\'esulte clairement du fait que $f$ est formellement lisse (resp.
formellement non ramifi\'e, formellement \'etale).

2) D'apr\`es \cite{Knutson} II (6.4), il existe un diagramme commutatif
$$\xymatrix{& X' \ar[d]^f \\
\Spec k \ar[ur]^{x'} \ar[r]^x & X}$$
o\`u $X'$ est un sch\'ema affine \'etale au-dessus de $X$. Compte tenu du point
pr\'ec\'edent, on peut donc supposer que $X$ est un sch\'ema affine d'anneau
$A$. Via des identifications bien connues, tout
revient \`a montrer que le morphisme $\Der_k(A,k)\otimes_k L \fleche
\Der_k(A,L)$
est un isomorphisme. On note $\pi : A \fleche k$ le morphisme correspondant au
point $x$. La $k$-alg\`ebre $A$ est de type fini, donc s'\'ecrit
$A=\frac{k[X_1, \dots, X_n]}{(P_1, \dots, P_r)}$. Si on note
$M : k^n \fleche k^r$ l'application lin\'eaire de matrice
$$M=\left(\begin{matrix}
\pi\left(\frac{\partial P_1}{\partial X_1}\right) &\dots &
\pi\left(\frac{\partial P_1}{\partial X_n}\right)\\
\vdots && \vdots \\
\pi\left(\frac{\partial P_r}{\partial X_1}\right) &\dots &
\pi\left(\frac{\partial P_r}{\partial X_n}\right)
\end{matrix}\right),$$
on voit facilement que $T_{x_L}(X_L)=\Der_k(A,L)=\Hom_A(\Omega_{A/k},L)$
s'identifie (fonctoriellement en $L$) au noyau de $M\otimes_k L$. La suite
exacte
$$\xymatrix{0 \ar[r] & T_xX \ar[r] & k^n \ar[r]^{M} & k^r}$$
induit une suite exacte
$$\xymatrix{0 \ar[r] & (T_xX)\otimes_k L \ar[r] & L^n \ar[r]^{M\otimes_k L}
& L^r}$$
ce qui prouve notre assertion.
\end{demo}

Nous rappelons également le r\'esultat suivant, valable aussi pour un espace algébrique puisqu'un espace algébrique en groupes sur un corps est toujours un schéma.  

\begin{prop}[\cite{poly_Kleiman} 5.13 et 5.14]
\label{prop_espaces_alg_en_groupes}
Soient $P$ un $k$-schéma en groupes localement de type fini, et $e$
le $k$-point neutre. Alors
$P$ a la m\^eme dimension en tout point. De plus cette dimension est
inf\'erieure \`a $\dim_k(T_eP)$, et les propri\'et\'es suivantes sont
\'equivalentes :
\begin{itemize}
\item[(i)] $\dim P =\dim_k(T_eP)$;
\item[(ii)] $P$ est lisse en 0;
\item[(iii)] $P$ est lisse.
\end{itemize}
Elles sont v\'erifi\'ees lorsque $k$ est de caract\'eristique nulle.
\end{prop}

\begin{thm}
Soient $k$ un corps, $S=\Spec k$, et $\X$ un $S$-champ alg\'ebrique.
On note $\Pic_{\X/k}$ le foncteur de Picard relatif $\piceth$ et
on suppose qu'il est repr\'esentable par un $S$-espace
alg\'ebrique localement de type fini.

a) Alors l'espace tangent \`a l'origine est $$T_0\Pic_{\X/k}
=H^1(\X,\O_{\X}).$$

b) L'espace alg\'ebrique $\Pic_{\X/k}$ a la m\^eme dimension en tout point. De
plus cette dimension est inf\'erieure \`a $\dim_k
H^1(\X,\O_{\X})$, et il y a \'egalit\'e si et seulement si $\Pic_{\X/k}$ est
lisse \`a l'origine. Dans ce cas, $\Pic_{\X/k}$ est lisse de dimension $\dim_k
H^1(\X,\O_{\X})$ partout. Il en est toujours ainsi lorsque $k$ est de caract\'eristique nulle.
\end{thm}
\begin{demo}
a) On note encore $P=\Pic_{\X/k}$. On rappelle que l'espace tangent \`a
$P$ en 0 est d\'efini par 
$T_0P:=\Ker\left(P(k[\eps]) \fleche P(k)\right)$, o\`u
$k[\eps]=\frac{k[X]}{(X^2)}$. Il est naturellement muni d'une structure de
$k$-espace vectoriel. Si l'on pose $\Xt=\X\times_{\Spec k}\Spec k[\eps]$, alors
$\Xt$ est une d\'eformation de $\X$. On note
$\xymatrix{i:\X \ar@{^{(}->}[r] &\Xt}$
l'injection canonique. On voit facilement que l'id\'eal $I$ sur $\Xt$ de carr\'e
nul qui d\'efinit $\X$ comme sous-champ ferm\'e de $\Xt$ est isomorphe \`a
$i_*\O_{\X}$. La remarque~(\ref{rem_noyau}) montre alors que l'on a un
isomorphisme naturel
$$\Ker(\Pic(\Xt)\rightarrow \Pic(\X)) \fleche H^1(\Xt,i_*\O_{\X}).$$
Or d'apr\`es le lemme~(\ref{coh_et_ext_inf}), $H^1(\Xt,i_*\O_{\X})$ est
isomorphe \`a $H^1(\X,\O_{\X})$. On a par ailleurs un carr\'e commutatif :
$$\xymatrix{ \Pic(\Xt)\ar[r] \ar[d] & \Pic(\X)\ar[d] \\
P(k[\eps]) \ar[r]& P(k),}$$
qui induit un morphisme entre les noyaux des fl\`eches horizontales, et donc,
d'apr\`es ce qui pr\'ec\`ede, un morphisme $v: H^1(\X,\O_{\X}) \fleche
T_0P$. Montrons que $v$ est un isomorphisme. Remarquons tout d'abord que
si $k$ est alg\'ebriquement clos, les fl\`eches verticales du carr\'e ci-dessus
sont des isomorphismes en vertu de (\ref{P_isom_Pet_si_k_alg_clos}), de sorte
que $v$ en est un aussi.
Passons maintenant au cas g\'en\'eral. Soit $\ov{k}$ une cl\^oture alg\'ebrique
de $k$. Le carr\'e ci-dessus s'envoie alors sur le carr\'e correspondant obtenu
apr\`es extension du corps de base \`a $\ov{k}$. On en d\'eduit en regardant les
noyaux un carr\'e commutatif de $k$-espaces vectoriels :
$$\xymatrix{H^1(\X,\O_{\X})\ar[r] \ar[d]_v &
H^1(\X_{\ov{k}},\O_{\X_{\ov{k}}})\ar[d]^{\ov{v}} \\
T_0P \ar[r] &T_0(P_{\ov{k}}).}$$
Par adjonction ce carr\'e induit le diagramme commutatif suivant :
$$\xymatrix{H^1(\X,\O_{\X})\otimes_k \ov{k}\ar[r] \ar[d]_{v\otimes_k \ov{k}} &
H^1(\X_{\ov{k}},\O_{\X_{\ov{k}}})\ar[d]^{\ov{v}} \\
(T_0P)\otimes_k \ov{k} \ar[r] &T_0(P_{\ov{k}}).}$$
Clairement pour montrer que $v$ est un isomorphisme il suffit de montrer que
$v\otimes_k \ov{k}$ en est un. Or nous avons vu ci-dessus que $\ov{v}$ en est
un. La fl\`eche du haut est un isomorphisme d'apr\`es
(\ref{lemme_coh_et_chgt_base}), et celle du bas en est un d'apr\`es le lemme
(\ref{lemme_espace_tangent}) ci-dessus. Pour b) il suffit d'appliquer
(\ref{prop_espaces_alg_en_groupes}).
\end{demo}

\section{Repr\'esentabilit\'e}

\begin{thm}
\label{thm_representabilite}
Soient $S$ un sch\'ema et $\X$ un $S$-champ alg\'ebrique propre, plat,
de pr\'esentation finie et cohomologiquement plat en dimension z\'ero. Alors le
foncteur de Picard relatif $P=\pic$ est repr\'esentable par un
$S$-espace alg\'ebrique localement de pr\'esentation finie et localement séparé.
\end{thm}
\begin{demo}
On peut évidemment supposer que $S$ est un schéma affine d'anneau $R$. On peut même supposer que $R$ est de type fini sur $\Z$ grâce à la proposition~(4.18)~(ii) de~\cite{LMB}. 
La question de la représentabilité est locale pour la topologie \emph{(fppf)}
sur $S$ (\cite{Artin_implicit_fct_thm} corollaire~(7.2)). On se ramène donc
avec la remarque~(\ref{remarque_section}) au cas où $f$ a une
section. Dès lors, il résulte du
théorème (\ref{comparaison_des_foncteurs_de_Picard}) que $P$ est isomorphe au
foncteur de Picard relatif $P_{\X/S}$, i.e. que
$$P(S')=\frac{\Pic(\X\times_S S')}{\Pic(S')}.$$
Nous allons maintenant utiliser le th\'eor\`eme (5.3) de
\cite{Global_Analysis_1},
et il nous faut pour cela en v\'erifier les conditions [0'] \`a [5']. On sait
déjà que les conditions [0'] (\og $P$ est un faisceau \emph{fppf}\fg), [1']
(\og $P$ est localement de présentation finie \fg), et [2'] (\og $P$ commute
aux limites projectives \fg) sont vérifiées, respectivement
d'après (\ref{comparaison_des_foncteurs_de_Picard}), (\ref{lpf}), et
(\ref{lim_proj}). Les conditions de séparation [3'] a) et b) résultent
immédiatement du lemme (\ref{lemme_locale_separation_du_foncteur_de_Picard}). Il
nous reste donc à construire une \og théorie des déformations \fg\ au sens
d'Artin \cite{Global_Analysis_1} (5.2) et à vérifier les conditions du théorème
(5.3) s'y rapportant. Or d'après le corollaire
(\ref{commutation_au_produit_fibre}), si $A\fleche A'$ est une extension
infinitésimale et $B\fleche A$ un morphisme quelconque, le morphisme canonique
\begin{equation}
\label{morph_produit_fibre}
P(A'\times_A B) \flechelongue P(A')\times_{P(A)} P(B)
\end{equation}
est un isomorphisme. On en déduit en utilisant la remarque au bas de la page~47
de \cite{Global_Analysis_1} que l'on obtient une théorie des déformations pour
$P$ en posant pour tout triplet $(A_0,M,\xi_0)$, où $A_0$ est une $R$-algèbre
intègre, $M$ un $A_0$-module de type fini, et $\xi_0\in P(A_0)$,
$$D(A_0,M,\xi_0)=\Ker(P(A_0[M])\flechelongue P(A_0)).$$
Ici $A_0[M]$ désigne la $A_0$-algèbre $A_0\oplus M$ munie de la multiplication
donnée par $(a_0,m)\cdot(a_0',m')=(a_0a_0',a_0m'+a_0'm)$. (Le $A_0$-module $M$
s'identifie donc à un idéal de carré nul de $A_0[M]$.) En utilisant le fait que
le morphisme $\Pic(A_0[M]) \fleche \Pic(A_0)$ est un isomorphisme
(lemme~(\ref{pic_extension_inf_est_isom})), on voit facilement que le noyau de
$P(A_0[M])\fleche P(A_0)$ s'identifie au noyau du morphisme
$\Pic(\X_{A_0[M]})\fleche \Pic(\X_{A_0})$. D'après la
remarque~(\ref{rem_noyau}), on obtient donc :
$$D(A_0,M,\xi_0)=H^1(\X_{A_0[M]},I)$$ où $I$ est l'idéal de carré nul
définissant l'immersion fermée $\X_{A_0}\fleche \X_{A_0[M]}$. Il est clair que
$I=f_{A_0[M]}^*j_*M\simeq i_*f_{A_0}^*M$, où les notations sont celles du
diagramme :
$$\xymatrix{\X_{A_0} \ar[r]^i \ar[d]_{f_{A_0}} &\X_{A_0[M]}\ar[d]^{f_{A_0[M]}}\\
\Spec A_0 \ar[r]^j &\Spec A_0[M].}$$
Le lemme~(\ref{coh_et_ext_inf}) nous donne alors :
$$D(A_0,M,\xi_0)=H^1(\X_{A_0},f_{A_0}^*M).$$

Vérifions maintenant les conditions [4'] a) à [5'] c) portant sur la théorie des
déformations ainsi construite. Notons que les conditions [4'] b) et [5'] a)
sont automatiquement vérifiées grâce à la bijectivité de
\eqref{morph_produit_fibre}, comme le fait remarquer Artin
(\cite{Global_Analysis_1}, page~48).

\emph{[4'] a)} Il s'agit de vérifier que la formation de $D(A_0,M,\xi_0)$ \og
commute avec la localisation en $A_0$ \fg, et que si $M$ est un $A_0$-module
libre de rang 1, alors $D(A_0,M,\xi_0)$ est un $A_0$-module de type fini. La
première assertion
résulte de la proposition~(\ref{lemme_coh_et_chgt_base}), et la seconde de la
finitude de la cohomologie des faisceaux cohérents sur un champ algébrique
propre et localement noethérien (\cite{Olsson_lemme_chow} théorème~(1.2), ou
\cite{Faltings_finiteness_coco}).

\emph{[4'] c)} Nous devons montrer que si $A_0$ est une $R$-algèbre intègre de
type fini, il existe un ouvert non vide $U$ de $\Spec A_0$ tel que pour tout
point fermé $s$ dans $U$ le $\kappa(s)$-espace vectoriel
$D(k,M\otimes_{A_0}\kappa(s),{\xi_0}_s)$ soit isomorphe à
$D(A_0,M,\xi_0)\otimes_{A_0}\kappa(s)$. En d'autres termes on demande que pour
un tel point $s$, le morphisme canonique
\begin{equation}
\label{morph_condition_4c}
H^1(\X_{A_0},f_{A_0}^*M)\otimes_{A_0}\kappa(s) \fleche
H^1(\X_{s},f_{s}^*(M\otimes_{A_0}\kappa(s)))
\end{equation}
soit un isomorphisme. Montrons d'abord qu'il existe un ouvert non vide de $\Spec
A_0$ sur lequel $M$ est libre de rang fini. On note $K$ le corps des fractions
de $A_0$. Comme $M$ est de type fini, $M_K$ est un $K$-espace vectoriel de
dimension finie. Soient $\nlist{x}$ des éléments de $M$ dont les images dans
$M_K$ engendrent $M_K$. Soit $(y_1, \dots, y_n)$ une famille génératrice de $M$.
Il existe $f\in A_0$ non nul tel que pour tout $j$ on puisse écrire $y_j$ comme
combinaison linéaire des $\frac{x_i}{f}$ dans $M_K$. Donc quitte à localiser par
$f$ on peut supposer que la famille $\nlist{x}$ engendre $M$. C'est alors
aussitôt une base de $M$. En effet, on a un morphisme surjectif $\varphi : A^n
\fleche M$ qui induit un isomorphisme $\varphi_K : A^n_K\fleche M_K$, et le
morphisme de localisation $A^n \fleche K^n$ est injectif puisque $A$ est
intègre, donc $\varphi$ est aussi injectif. Maintenant, $M$ est libre de rang
fini sur $\Spec A_0$, donc $f_{A_0}^*M\simeq (\O_{\X\times_S \Spec A_0})^n$ est
cohérent et plat sur $\Spec A_0$ (car $f$ est plat). Donc le résultat de Mumford
(\cite{Varietes_abeliennes} paragraphe~5) généralisé par Aoki aux champs
algébriques (\cite{Aoki_Hom}, théorème~(A.1)) s'applique et il existe un ouvert
non vide $U$ de $\Spec A_0$ tel que, pour tout point $s$ de $U$, le morphisme
\eqref{morph_condition_4c} ci-dessus soit un isomorphisme. 

Pour les conditions [5'] b) et [5'] c), on peut reprendre telles quelles les
démonstrations proposées dans \cite{Global_Analysis_1}, page~70 dans le cas des
espaces algébriques, en utilisant la classe d'obstruction que nous avons
construite au théorème~(\ref{thm_defm_fi}).
\end{demo}

\begin{comment}
\begin{lem}
\label{condition_4c}
Soit $\X$ un $S$-champ algébrique (de Deligne-Mumford ?) propre, où
$S=\Spec A_0$ est le spectre d'un
anneau noethérien intègre, et soit $\Fc$ un faisceau cohérent sur $\X$. Soit $q$
un entier. Alors il existe un ouvert non vide $U$ de $S$ tel que pour tout $s\in
U$ le morphisme canonique
$$H^q(\X,\Fc)\otimes_{A_0}\kappa(s) \fleche
H^q(\X_{s},\Fc_s)$$
soit un isomorphisme.
\end{lem}
\begin{demo}
\end{demo}
\end{comment}

\chapter{Composante neutre du foncteur de Picard}

\section{Préliminaires}

\subsection{Composante des fibres le long d'une section}

\begin{souslem}
\label{CC_des_fibres_lemme_projection}
Soient $k$ un corps, $X$ un $k$-schéma connexe localement de type fini, $L$ une extension de $k$ et $p$ la projection de $X_L$ dans $X$.
Alors les fibres de $p$ rencontrent toutes les composantes connexes de $X_L$. Autrement dit, pour toute composante
connexe $U$ de $X_L$, le morphisme induit $p_{|_U} : U \fleche X$ est surjectif.
\end{souslem}
\begin{demo} En dévissant l'extension, il suffit clairement de traiter le cas d'une extension algébrique et le cas d'une extension
transcendante pure. Dans le premier cas, le morphisme $\Spec L \fleche \Spec k$ est universellement ouvert (\cite{EGA4_2}~2.4.9)
et universellement fermé (\cite{EGA2}~6.1.10). Si $U$ est une composante connexe de $X_L$, son image $p(U)$ dans $X$ est ouverte, fermée
et non vide, donc c'est $X$ tout entier.

Supposons maintenant l'extension $L/k$ transcendante pure. Si $\Omega$ est une autre extension de $k$, alors l'anneau
$L\otimes_k \Omega$ est intègre. En effet, en écrivant $L=k(\underline{T})$, où $\underline{T}$ est une famille
d'indéterminées, on voit que $L\otimes_k \Omega $ est isomorphe à $S^{-1}\Omega[\underline{T}]$ avec
$S=k[\underline{T}]$. Il en résulte que les fibres de $p$ sont géométriquement intègres.
Notons $(X_i)_{i\in I}$ la famille des composantes connexes de $X_L$. Vu que les fibres de $p$ sont connexes, chaque
$X_i$ est une réunion de fibres. Donc les $p(X_i)$ forment une partition de $X$, et ils sont ouverts puisque $p$
est ouvert. Par connexité un seul d'entre eux est non vide, donc $X_L$ est connexe.
\end{demo}

\begin{souslem}
\label{CC_des_fibres_lemmeA}
Soit
$$\xymatrix{& X\ar[d]^f \\ S' \ar[r]_g \ar[ur]^e & S}$$
un diagramme commutatif de schémas. On suppose que $f$ est lisse et quasi-compact, et que $g$ est
universellement ouvert. Alors il existe un unique ouvert $U$ de $X$ tel que pour
tout $s\in S$, $U_s$ soit la réunion des composantes connexes de $X_s$ qui rencontrent $e(S')$.
\end{souslem}
\begin{demo}
On note $X'$ le produit fibré $X\times_S S'$ et $e' : S' \fleche X'$ la section induite par $e$.
D'après \cite{EGA4_3}~(15.6.5), on a un ouvert $U'$ de $X'$ tel que pour tout $s'\in S'$, la fibre de $U'$
au-dessus de $s'$ soit la composante connexe de $X'_{s'}$ contenant $e'(s')$. Soit $U$ l'image de $U'$
dans $X$. C'est un ouvert de $X$ puisque $g$ est universellement ouvert.

Si $s$ est un point de $S$, il est clair que la fibre $U_s$ est l'image de $U'_s$ par le morphisme
$X'_s \fleche X_s$ obtenu par changement de base. Donc pour montrer que $U$ vérifie la propriété annoncée,
on peut supposer que $S$ est le spectre d'un corps $k$. Pour tout $s'\in S'$, la fibre
$U'_{s'}$ est connexe et contient $e'(s')$, de sorte que son image dans $X$ est incluse dans la
composante connexe de $e(s')$. Comme $U$ est la réunion des images des $U'_{s'}$,
on en déduit que $U$ est inclus dans la réunion des composantes connexes de $X$ qui rencontrent
$e(S')$. Réciproquement, étant donné un point $s'$ de $S'$ et un point $x$ de $X$ qui est dans la composante
connexe de $e(s')$, montrons que $x$ appartient à $U$. 
Notons $C'$ la composante connexe de $X'_{s'}$ qui contient $e'(s')$, et $C$ la composante connexe de
$e(s')$ dans $X$. D'après le lemme~(\ref{CC_des_fibres_lemme_projection}), le morphisme induit de $C'$ vers $C$
est surjectif. En particulier $x$ appartient à l'image de $C'$, donc à $U$.

L'unicité de $U$ est claire, puisque l'ensemble sous-jacent à $U$ est déterminé de manière unique
par ses fibres. 
\end{demo}

\begin{sousremarque}\rm
\label{CC_des_fibres_lemmeA_rem_adherence}
Dans le lemme précédent, pour tout $s\in S$, $X_s$ est lisse sur $\kappa(s)$ donc ses composantes connexes sont irréductibles.
On en déduit que si $S'$ est un ouvert de $X$, alors $U_s$ est l'adhérence de $S'_s$ dans $X_s$.
\end{sousremarque}

\begin{sousremarque}\rm
L'hypothèse de quasi-compacité sur $f$ n'est pas nécessaire. Elle est cependant présente dans l'énoncé~\cite{EGA4_3}~(15.6.5),
et c'est le seul endroit où nous l'utilisons. Nous allons nous en affranchir dans le lemme~(\ref{CC_des_fibres_lemmeB})
ci-dessous. Ceci prouve en particulier, en reprenant la démonstration du lemme~(\ref{CC_des_fibres_lemmeA}) et en 
y remplaçant \og D'après \cite{EGA4_3}~(15.6.5) \fg\ par \og D'après le lemme~(\ref{CC_des_fibres_lemmeB})\fg, que
l'énoncé~(\ref{CC_des_fibres_lemmeA}) est encore valable même lorsque $f$ n'est plus supposé quasi-compact.
\end{sousremarque}

\begin{souslem}
\label{CC_des_fibres_lemmeB}
Soient $S$ un schéma et $X$ un espace algébrique lisse sur $S$ muni d'une section $e : S \fleche X$. 
Alors il existe un sous-espace algébrique ouvert de $X$, noté $U(X/S)$, tel que pour
tout $s\in S$, la fibre de $U(X/S)$ au-dessus de $s$ soit la composante connexe de $e(s)$ dans $X_s$.
\end{souslem}
\begin{demo}
Il suffit de traiter le cas où $S$ est quasi-compact. On note $f: X \fleche S$ le morphisme structural de $X$.
Soit $\pi$ un morphisme étale d'un schéma quasi-compact $X_1$ vers $X$ tel que le morphisme composé $f\circ \pi$ soit surjectif. On note $S_1$ le produit fibré
$S \times_{e,X,\pi} X_1$. Le diagramme commutatif
$$\xymatrix{& X_1 \ar[d]^{f\circ \pi} \\ S_1 \ar[ru]^{e'} \ar[r] & S}$$
vérifie bien les hypothèses du lemme~(\ref{CC_des_fibres_lemmeA}), donc il existe un ouvert $U_1$ de $X_1$ tel que pour tout
$s\in S$, $U_{1s}$ soit la réunion des composantes connexes de $X_{1s}$ qui rencontrent $e'(S_1)$. On note $U$ l'image de $U_1$
par $\pi$ (c'est un ouvert de $X$), puis $V_1$ l'image réciproque de $U$ par $\pi$ (c'est un ouvert de $X_1$).
Autrement dit $V_1$ est le saturé de $U_1$ pour la relation d'équivalence définie par $\pi$.
On applique maintenant le lemme~(\ref{CC_des_fibres_lemmeA}) au diagramme commutatif
$$\xymatrix{& X_1 \ar[d]^{f\circ \pi} \\ V_1 \ar@{^{(}->}[ru] \ar[r] & S}$$
et on obtient ainsi un ouvert $W_1$ de $X_1$ tel que, pour tout $s\in S$, $W_{1s}$ soit la réunion des composantes connexes
de $X_{1s}$ qui rencontrent $V_{1s}$, c'est-à-dire l'adhérence de $V_{1s}$ dans $X_{1s}$ (cf.
remarque~(\ref{CC_des_fibres_lemmeA_rem_adherence})). On note maintenant $W$ l'image de $W_1$ par $\pi$. C'est un ouvert de $X$.

Notons $C_s$ la composante connexe de $X_s$ qui contient $e(s)$. L'ouvert $U_{1s}$ est la réunion des composantes
connexes de $X_{1s}$ qui rencontrent $e'(S_{1s})$. Donc son image $U_s$ est une réunion de parties connexes
qui contiennent $e(s)$. En particulier $U_s$ est inclus dans $C_s$. De plus $U_s$ est non vide car tout point de $S$ est dans l'image de $X_1$. Comme $C_s$ est irréductible (puisque $X_s$
est lisse), on a $C_s=\ov{U_s}=\ov{\pi_s(V_{1s})}$.

Par ailleurs, on a la suite d'inclusions et d'égalités ensemblistes suivante :
\begin{eqnarray*}
\pi_s^{-1}(W_s) &=& \pi_s^{-1}(\pi_s(W_{1s}))\\
&=& \pi_s^{-1}(\pi_s(\ov{V_{1s}})) \quad \textrm{par construction de } W_{1s}\\
&\subset& \pi_s^{-1}(\ov{\pi_s(V_{1s})}) \quad \textrm{par continuité de } \pi_s\\
&\subset& \ov{\pi_s^{-1}(\pi_s(V_{1s}))} \quad \textrm{parce que } \pi_s \textrm{ est un morphisme ouvert}\\
&=& \ov{V_{1s}} \\
&=& W_{1s} \\
&\subset& \pi_s^{-1}(W_s).
\end{eqnarray*}
D'où l'égalité entre $\pi_s^{-1}(W_s)$ et $\pi_s^{-1}(\ov{\pi_s(V_{1s})})$, qui peut encore s'écrire :
$$W_s= \pi_s(X_{1s}) \cap C_s.$$

Maintenant, soit $(\pi_i : X_i \fleche X)_{i\in I}$ une famille couvrante étale, où les $X_i$ sont des schémas affines et où chacun des morphismes composés $f\circ \pi_i$ est surjectif (c'est possible : il suffit de veiller à ce que chacun des $X_i$ recouvre l'image, quasi-compacte, de la section $e$ dans $X$).
Pour chaque $i\in I$, on note $W_i$ l'ouvert de $X$ obtenu par la construction précédente. Soit $U(X/S)$ la réunion
des $W_i$. Il est clair que pour tout $s\in S$, la fibre de $U(X/S)$ au-dessus de $s$ est égale à $C_s$.
\end{demo}

\begin{sousremarque}\rm
\label{rem_caractérisation_fonctorielle_CC_neutre}
Comme sous-foncteur de $X$, l'ouvert $U=U(X/S)$ est caractérisé par la propriété suivante.
Pour tout $S$-schéma $T$ et pour tout $\xi \in X(T)$, $\xi$ appartient à $U(T)$ si et seulement si 
pour tout $s\in S$, le point $\xi_s\in X_s(T_s)$ obtenu par changement de base appartient à $C_s(T_s)$, autrement
dit le morphisme correspondant $T_s \fleche X_s$ se factorise par $C_s$.
La vérification est immédiate et laissée au lecteur. En particulier on voit que la propriété vérifiée par les
fibres détermine $U(X/S)$ de manière unique.
\end{sousremarque}

\begin{sousremarque} \rm La formation de $U(X/S)$ commute à tout changement de base. En effet, soit 
$S' \fleche S$ un morphisme de changement de base. Notons $U=U(X/S)$, $X'=X\times_S S'$, $U'=U\times_S S'$,
et $e' : S' \fleche X'$ la section obtenue par changement de base.
D'après la remarque précédente il suffit de montrer que pour tout point $s'$ de $S'$, la fibre
$U'_{s'}$ est la composante connexe de $X'_{s'}$ contenant $e'(s')$. Soient $s'$ un point de $S'$ et $s$ son image dans $S$.
On a un diagramme commutatif à carrés cartésiens :
$$\xymatrix{U'_{s'} \ar[r]\ar[d] \cartesien & X'_{s'} \ar[r]\ar[d]^p \cartesien & \Spec \kappa(s')\ar[d] \\
U_s \ar[r] & X_s \ar[r] & \Spec \kappa(s).}$$
L'ouvert $U_s$ est un $\kappa(s)$-schéma connexe qui a un $\kappa(s)$-point, donc d'après~\cite{EGA4_2}~4.5.13
il est géométriquement connexe. En particulier $U'_{s'}$ est connexe, et comme il contient $e'(s')$ il est inclus
dans la composante connexe $C'_{s'}$ de $e'(s')$. D'autre part, $C'_{s'}$ est connexe donc son image aussi et elle contient $e(s)$,
donc $p(C'_{s'})$ est inclus dans $C_s$, ou, ce qui revient au même, $C'_{s'}$ est inclus dans $p^{-1}(U_s)=U'_{s'}$.
\end{sousremarque}

\subsection{Cohomologie à coefficients dans $\Z$ des champs normaux}

Le but du présent paragraphe est de démontrer le résultat suivant, qui nous sera utile pour
ramener, lorsque $S$ est le spectre d'un corps, l'étude de la propreté de $\piczero$ au cas où
$\X$ est un schéma.

\begin{sousthm}
\label{H1_normal}
Soit $\X$ un champ algébrique localement noethérien et normal. Alors $H^1(\X,\Z)=0$.
\end{sousthm}

Nous allons pour cela montrer que tout $\Z$-torseur sur $\X$ est trivial. Notre démarche est
fortement inspirée de l'étude des préschémas constants tordus quasi-isotriviaux proposée dans SGA3
(\cite{SGA3} exposé~X, paragraphe~5).

\begin{souslem} Soit $f : \X \fleche \Y$ un morphisme représentable de $S$-champs algébriques. Soit $y$ un point
de l'espace topologique $|\Y|$ sous-jacent à $\Y$. Les propositions suivantes sont équivalentes :
\begin{itemize}
\item[(i)] Pour un représentant $\Spec K \fleche \Y$ de $y$, le morphisme $f_K$ induit par changement de base
de $\X_K=\X\times_{\Y} \Spec K$ vers $\Spec K$ est fini.
\item[(ii)] Pour tout représentant $\Spec K \fleche \Y$ de $y$, le morphisme $f_K$ est fini.
\end{itemize}
Lorsqu'elles sont vérifiées, on dit que $f_y : \X_y \fleche y$ est fini, ou encore que $\X_y$ est fini.
\end{souslem}
\begin{demo}
Il suffit clairement de montrer que si $L$ est une extension de $K$ et si $f_L$ est fini, alors $f_K$ est fini. Or si
$f_L$ est fini, $\X_L$ est un schéma affine, et $\X_K$ aussi par descente \emph{fpqc} pour les morphismes affines.
Maintenant la descente fidèlement plate pour les morphismes finis de schémas (\cite{EGA4_2}~(2.7.1)) assure que $f_K$ est fini.
\end{demo}

Le lemme suivant généralise le lemme~5.13 de~\cite{SGA3}, exposé~X.

\begin{souslem}
\label{morphisme_fini_condition_ouverte_et_fermee}
Soit $p :\Pc \fleche \X$ un morphisme représentable de $S$-champs algébriques, avec $\X$ localement noethérien.
On suppose qu'il existe une présentation $X\fleche \X$ de $\X$ telle que $P=\Pc \times_{\X} X$ soit une union
disjointe de copies de $X$. Soit $\Cc$ un sous-champ ouvert et fermé de $\Pc$. Alors l'ensemble
des points $x$ de $|\X|$ tels que $\Cc_x$ soit fini est ouvert et fermé dans $|\X|$. Si on note $\Uc$ le sous-champ
ouvert et fermé que cet ensemble définit, le champ $\Cc_{\Uc}=\Cc\times_{\X}\Uc$ est fini sur $\Uc$.
\end{souslem}
\begin{demo}
Notons $C$ le sous-espace algébrique ouvert et fermé de $P$ obtenu par changement de base à partir de $\Cc$.
Les propriétés que l'on veut montrer sont clairement de nature locale pour la topologie lisse sur $\X$, 
donc il suffit de montrer qu'elles sont vérifiées par $C \fleche P \fleche X$. Comme $X$ est lui aussi localement
noethérien, ses composantes connexes sont ouvertes et fermées donc on peut supposer que $X$ est un schéma connexe.
Dans ce cas, vu que $P$ est une union disjointe de copies de $X$, le sous-schéma ouvert et fermé $C$ est
lui-même l'union disjointe de certaines de ces copies. Si elles sont en nombre fini alors $C$ est fini sur $X$, sinon
l'ensemble des points $x$ de $X$ où $C_x$ est fini est vide.
\end{demo}
\begin{democ}{du théorème \ref{H1_normal}}
Nous allons utiliser la description du premier groupe de cohomologie en termes de torseurs
(cf. paragraphe~\ref{Cohomologie_et_torseurs}).
Il suffit en vertu de~(\ref{thm_coh_torseur}) de montrer que tout $\Z$-torseur sur $\X$ au sens de~(\ref{def_torseur}) est trivial.
Le cas où $\X$ est le spectre d'un corps est bien connu et nous nous en servirons par la suite.
Remarquons tout d'abord que dans tous les champs algébriques (\emph{a fortiori} tous les schémas ou espaces algébriques) qui vont
intervenir au cours de la démonstration, les composantes connexes sont irréductibles.
En effet, ils seront tous normaux et localement noethériens car ce sont là des propriétés de nature locale pour la topologie lisse.
Notre affirmation résulte alors de la proposition~4.13 de~\cite{LMB}.

Donnons-nous donc un $\Z$-torseur sur $\X$, c'est-à-dire un morphisme $p : \Pc \fleche \X$ représentable et lisse muni d'une
action de $\Z$ qui en fait un torseur. Pour montrer que $\Pc$ est trivial, on peut supposer $\X$ connexe donc irréductible.
Soit $\Cc$ une composante connexe de $\Pc$. Notons $\eta$ le point générique de
l'espace $|\X|$ sous-jacent à $\X$.

\begin{etape}{La fibre générique $\Cc_{\eta}$ est finie.}
Soit $u : U\fleche \X$ un morphisme lisse, où $U$ est un schéma affine irréductible (donc intègre, puisque $U$ est normal).
On note $s$ le point générique de $U$ et on adopte encore les notations du diagramme suivant :
$$\xymatrix{\Cc_s \ar[r] \ar[d] \cartesien & \Cc_U \ar[r] \ar[d] \cartesien & \Cc \ar[d]\\
\Pc_s \ar[r] \ar[d] \cartesien & \Pc_U \ar[r] \ar[d] \cartesien & \Pc \ar[d]^p \\
\Spec \kappa(s) \ar[r] & U\ar[r]^u & \X.}$$
Notons $(\Cc_i)_{i\in I}$ les composantes connexes de $\Cc_U$ et pour chaque $i\in I$ notons $\eta_i$ le point
générique de $|\Cc_i|$. Notons enfin $\xi$ le point générique de $|\Cc|$ et $\Spec K \fleche \Cc$ l'un de ses
représentants. Pour tout $i$ le morphisme de $\Cc_i$ vers $\Cc$ est lisse donc générisant (\cite{LMB}~(5.8)) si bien qu'il envoie
le point $\eta_i$ sur $\xi$. On en déduit (cf. par exemple \cite{LMB}~(5.4)~(iv)) que le champ $\Cc_{i,K}=\Cc_i \times_{\Cc} \Spec K$
est non vide. Maintenant, les $|\Cc_{i,K}|$ forment une partition ouverte de $|\Cc_{U,K}|$, où $\Cc_{U,K}$ est le produit
fibré $\Cc_U\times_{\Cc} \Spec K$. Or le champ
$$\Cc_{U,K} = U\times_{\X} \Spec K$$
est quasi-compact car $\X$ est quasi-séparé donc les $|\Cc_{i,K}|$ sont en nombre fini et finalement $\Cc_U$ n'a qu'un
nombre fini de composantes connexes.

Par ailleurs, pour chaque $i$ le morphisme naturel de $\Cc_i$ vers $U$ est lui aussi générisant donc il envoie $\eta_i$
sur $s$. En particulier sa fibre générique est irréductible donc $\Cc_s$ est une union finie d'irréductibles. Comme
$\Pc_s$ est un $\Z$-torseur sur un corps, il est nécessairement trivial, donc $\Cc_s$ est une union disjointe finie
de copies de $\Spec \kappa(s)$. Or le morphisme de $\Spec \kappa(s)$  vers $\X$ est un représentant
de $\eta$ donc $\Cc_{\eta}$ est fini.
\end{etape}

\begin{etape}{Montrons que $\Cc$ est fini.}
D'après le lemme~(\ref{morphisme_fini_condition_ouverte_et_fermee}), l'ensemble des points $x$ de $|\X|$ où
$\Cc_x$ est fini est ouvert et fermé dans $|\X|$. Or il est non vide puisqu'il contient $\eta$, donc par
connexité c'est $|\X|$ tout entier. Le même lemme prouve alors que le morphisme de $\Cc$ vers $\X$ est fini.
\end{etape}

\begin{etape}{Montrons que $\Cc \fleche \X$ est étale.}
Il s'agit d'un morphisme fini. En particulier il est schématique. Notre assertion résulte alors du
fait qu'un morphisme fini et lisse de schémas est étale (\cite{SGA1}~II~1.4).
\end{etape}

\begin{etape}{Montrons que $\Cc \fleche \X$ est radiciel.}
Supposons qu'il existe un corps $K$ et un morphisme $\Spec K \fleche \X$ tel que le schéma $\Cc_K$ obtenu par changement
de base contienne au moins deux points $x_1$ et $x_2$. Notons $c_1$ et $c_2$ leurs images
dans $|\Cc| \subset |\Pc|$. Soit $n$ l'unique élément de $\Z\setminus\{0\}$ tel que l'automorphisme $\tau_n$ correspondant
envoie $x_1$ sur $x_2$. On note encore $\tau_n$ l'automorphisme de $|\Pc|$ qui correspond à $n$. Il est clair qu'il envoie
$c_1$ sur $c_2$. Or $\tau_n(|\Cc|)$ est un connexe qui contient $c_2$ donc il est inclus dans la composante connexe de
$c_2$, à savoir $|\Cc|$. On montre de même que $\tau_n^{-1}(|\Cc|)$ est inclus dans $|\Cc|$ donc $\tau_n$ induit
un automorphisme de $|\Cc|$ et tous les $\tau_n^k(c_1)$, $k\in \Z$, sont dans $|\Cc|$. Donc tous les $\tau_n^k(x_1)$
sont dans l'ouvert $\Cc_{K}$ de $\Pc_K$, ce qui contredit le fait qu'il est de type fini.
\end{etape}

\begin{etapefinale}{Conclusion.}
Le morphisme de $\Cc$ vers $\X$ est schématique, étale, radiciel et de type fini donc par~\cite{SGA1}~I~5.1 c'est une
immersion ouverte. De plus il est aussi fermé puisqu'il est fini donc par connexité de $\X$ il est surjectif.
Ceci prouve que c'est un isomorphisme, donc le torseur $\Pc \fleche \X$ est trivial.
\end{etapefinale}
\end{democ}

\section{La composante connexe de l'identité}

On rappelle le résultat suivant concernant la composante neutre d'un $k$-schéma
en groupes localement de type fini.

\begin{thm}[\cite{poly_Kleiman}~5.1]
Soient $k$ un corps et $G$ un $k$-schéma en groupes localement de type fini. Alors $G$ est séparé. Soit
$G^0$ la composante connexe de l'élément neutre de $G$. Alors $G^0$ est un sous-schéma en groupes ouvert et
fermé de $G$. Il est de type fini et géométriquement irréductible. De plus, la formation de $G^0$ commute
aux extensions du corps $k$.
\end{thm}

Naturellement, tout ceci est valable pour le foncteur de Picard $\Pic_{\X/k}$ dès
qu'il est représentable par un espace algébrique localement de type fini sur
$k$. En effet, c'est alors automatiquement un schéma en vertu d'un lemme
d'Artin (\cite{Global_Analysis_1} lemme~4.2 p.~43). En particulier, c'est le cas
dès que $\X$ est un champ algébrique propre et cohomologiquement plat en dimension zéro sur $k$ (théorème~(\ref{thm_representabilite})).

On sait que si $\X$ est un schéma propre et géométriquement normal sur un corps, alors la composante neutre du groupe de Picard est propre (voir par exemple \cite{poly_Kleiman}~théorème~5.4).
Le théorème suivant généralise ce résultat au cas où $\X$ est un champ algébrique.

\begin{thm}
\label{composante_neutre_propre_sur_k}
On suppose que $\X$ est un champ algébrique propre, géométriquement normal et cohomologiquement plat en
dimension zéro sur $\Spec k$. Alors la composante neutre $\Pic^0_{\X/k}$
du schéma de Picard est propre sur $k$.
\end{thm}
\begin{demo}
Si $\ov{k}$ est une clôture algébrique de $k$ alors le champ $\X_{\ov{k}}$ obtenu par changement de base est normal d'après~\cite{EGA4_2}~(6.7.7). Il vérifie clairement les autres hypothèses du théorème donc par descente fidèlement plate on peut supposer le corps $k$ algébriquement clos.

Il suffit (cf. argumentaire de Kleiman au cours de la démonstration du théorème~5.4 de~\cite{poly_Kleiman}) de montrer que tout morphisme de schémas de $G$ vers $\Pic_{\X/k}$ est constant, avec $G=\mathbb{G}_a$ ou $G=\gm$. (En fait il suffirait même de le faire pour $\gm$ puisqu'alors c'est aussi vrai pour $G=\mathbb{G}_a$, mais
cette restriction n'apporte pas grand-chose.) En effet, il suffit de montrer que le
réduit $(\Pic^0_{\X/k})_{\textrm{réd}}$ est propre, or ce dernier est lisse donc on peut
lui appliquer le théorème de structure de Chevalley et Rosenlicht (cf. par
exemple~\cite{Conrad_thm_Chevalley}, théorème~1.1). On en déduit que
$(\Pic^0_{\X/k})_{\textrm{réd}}$ a un sous-groupe algébrique linéaire $H$, fermé et
distingué dans $(\Pic^0_{\X/k})_{\textrm{réd}}$, tel que le quotient
$(\Pic^0_{\X/k})_{\textrm{réd}}/H$ soit une variété abélienne. Il suffit de montrer
que $H$ est trivial. Le groupe $H$ est commutatif donc résoluble. Il est dès lors triangulable d'après le théorème de Lie-Kolchin. On en déduit que s'il était non trivial, il contiendrait un sous-groupe isomorphe à $\gm$ ou $\mathbb{G}_a$ (voir par exemple le livre de Springer~\cite{Springer_Linear_Algebraic_Groups}, lemme~6.3.4).

Vu que $f$ a une section, et vu que les groupes de Picard de $\Spec k$ et
de $G$ sont triviaux, dire que tout morphisme de schémas de $G$ vers $\Pic_{\X/k}$ est constant revient à dire, grâce au théorème~(\ref{comparaison_des_foncteurs_de_Picard}), que
le morphisme naturel
$$\xymatrix{\Pic(\X) \ar[r]& \Pic(\X\times G)}$$
est un isomorphisme.
La suite exacte des termes de bas degré associée à la suite spectrale de Leray du morphisme $p : \X\times G \fleche \X$
s'écrit :
$$\xymatrix{0 \ar[r]& H^1(\X,p_*\gm) \ar[r] & H^1(\X\times G, \gm) \ar[r] & H^0(\X, R^1p_*\gm).}$$
Commençons par montrer que le faisceau $R^1p_* \gm$ est nul. On sait d'après le calcul des images directes supérieures
effectué en annexe, que c'est le faisceau associé au préfaisceau qui à tout ouvert lisse-étale $(U,u)$ de $\X$
associe $H^1(U\times G,\gm)$. Donc d'après le lemme~(\ref{lemme_nullite_faisceau_associe}) il suffit de montrer que
pour tout schéma affine  $U$ lisse sur $\X$ et pour tout $\xi\in H^1(U\times G,\gm)$,
il existe une famille couvrante étale $V\fleche U$ telle
que l'élément $\xi_{|_V}$ de $H^1(V\times G,\gm)$ soit nul. Mais pour démontrer ceci il suffit clairement de savoir que 
le morphisme $\Pic(U) \fleche \Pic(U\times G)$ est surjectif. Nous sommes donc ramenés à montrer que si $U$ est un schéma affine normal
(que l'on peut aussi supposer intègre, vu que $U$ est de toute manière somme disjointe finie
de schémas intègres) sur $\Spec k$, alors $\Pic(U\times G)$ s'identifie à $\Pic(U)$.
Ce fait est démontré par Kleiman dans~\cite{poly_Kleiman}, au cours de la démonstration du théorème~(5.4).
% Nous rappelons ici brièvement la démonstration, pour la commodité du lecteur.
%{\bf (...)}

Calculons maintenant le faisceau $p_*\gm$. Considérons d'abord le cas où $G=\gm$. Nous allons montrer que
$p_*\gm$ s'identifie à $\gm\times \Z$. Il suffit bien évidemment de vérifier que ces deux faisceaux coïncident
sur le site $\lisets(\X)$. Si $U=\Spec A$ est un schéma affine lisse sur $\X$, il est en particulier normal, donc somme disjointe
finie de schémas affines intègres, si bien que l'on peut supposer $U$ intègre. Alors $\gm\times \Z (U)=A^{\times}\times \Z$.
Par ailleurs, vu que $p$ est représentable, on a
$$p_*\gm(U)=\gm((\X\times G) \times_{\X} U)=\gm(U\times G)=A[X,X^{-1}]^{\times}.$$
Or lorsque $A$ est intègre, il est clair que les éléments inversibles de l'anneau $A[X,X^{-1}]$ sont les
éléments de la forme $aX^n$ avec $a\in A^{\times}$ et $n\in \Z$.
Dans le cas où $G=\mathbb{G}_a$, le lecteur vérifiera facilement que l'on trouve $p_*\gm=\gm$.

Or le groupe $H^1(\X,\Z)$ est réduit à zéro d'après le théorème~(\ref{H1_normal}). On a alors, que $G$ soit égal à $\mathbb{G}_a$ ou $\gm$,
$$H^1(\X,p_*\gm)=H^1(\X,\gm)=\Pic(\X),$$
ce qui, vu la suite exacte évoquée ci-dessus, fournit l'isomorphisme désiré.
\end{demo}

La notion de faisceaux inversibles algébriquement équivalents se généralise très bien aux champs algébriques, et comme
dans le cas des schémas elle permet de caractériser les $k$-points du foncteur de Picard qui sont dans la composante neutre.

\begin{defi}[\cite{poly_Kleiman}~5.9]
Soit $\X$ un champ algébrique sur un corps $k$. Soient $\Lc$ et $\Nc$ deux faisceaux inversibles sur $\X$. On dit que $\Lc$ et $\Nc$ sont
\emph{algébriquement équivalents} s'il existe une suite de $k$-schémas connexes de type fini $T_1, \dots, T_n$, des points
géométriques $s_i, t_i$ de $T_i$ (pour tout $i$) ayant tous le même corps, et un faisceau inversible $\Mc_i$ sur
$\X\times_k T_i$ tels que
$$\Lc_{s_1} \simeq \Mc_{1,s_1}, \quad\Mc_{1,t_1}\simeq \Mc_{2,s_2},\quad \dots, \quad\Mc_{n-1,t_{n-1}}\simeq \Mc_{n,s_n},
\quad\Mc_{n,t_n}\simeq \Nc_{t_n}$$
\end{defi}

\begin{thm}[\cite{poly_Kleiman}~5.10]
Soit $\X$ un champ algébrique sur un corps $k$. On suppose que $\Pic_{\X/k}$ est représentable par un schéma localement
de type fini. Soit $\Lc$ un faisceau inversible sur $\X$ et soit $\lambda$ le point correspondant de $\Pic_{\X/k}$. Alors
$\Lc$ est algébriquement équivalent à $\O_{\X}$ si et seulement si $\lambda$ est dans la composante neutre $\Pic^0_{\X/k}$.
\end{thm}
\begin{demo} On peut recopier telle quelle la démonstration qui accompagne l'énoncé référencé ci-dessus.
\end{demo}

Comme toujours, la définition de la composante neutre est plus délicate lorsque la base n'est plus un corps, mais un schéma quelconque.
La remarque~(\ref{rem_caractérisation_fonctorielle_CC_neutre}) motive la définition suivante.

\begin{defi} Soient $S$ un schéma et $\X$ un $S$-champ algébrique. On suppose que le foncteur de Picard
$\pic$ est représentable
par un espace algébrique localement de type fini. On désigne alors par $\piczero$ le sous-foncteur de $\pic$ défini de la
manière suivante. Pour tout $S$-schéma $S'$ et pour tout $\xi \in \pic(S')$, on dit que $\xi$ appartient à
$\piczero(S')$ si et seulement si pour tout point $s'$ de $S'$, l'élément $\xi_{|_{s'}}$
de $\Pic_{\X_{s'}/\kappa(s')}(\kappa(s'))$ est dans la composante neutre $\Pic^0_{\X_{s'}/\kappa(s')}(\kappa(s'))$.
\end{defi}

\begin{remarque}\rm
Il est clair que $\piczero$ est bien un sous-foncteur en groupes de $\pic$. De plus, la formation de $\piczero$ commute
au changement de base. Par ailleurs, si $S$ est le spectre d'un corps $k$ et si $\Pic_{\X/k}$ est représentable par un schéma
en groupes localement de type fini, alors le sous-foncteur $\Pic^0_{\X/k}$ ainsi défini coïncide avec le sous-foncteur
ouvert défini par la composante connexe de l'élément neutre dans $\Pic_{\X/k}$. Il en résulte dans le cas d'une base $S$ quelconque,
que pour tout point $s$ de $S$ la fibre de $\piczero$ au-dessus de $s$ coïncide avec $\Pic^0_{\X_s/\kappa(s)}$.
En particulier, si l'on suppose que pour tout $s\in S$ la composante neutre $\Pic^0_{\X_s/\kappa(s)}$
du schéma en groupes $\Pic_{\X_s/\kappa(s)}$ est lisse sur $\kappa(s)$, alors 
la définition que nous avons adoptée
coïncide avec celle donnée par Kleiman dans la proposition~5.20 de~\cite{poly_Kleiman}.
\end{remarque}

\begin{prop} Soient $S$ un schéma et $\X$ un $S$-champ algébrique. On suppose que le foncteur de Picard $\pic$ 
est représentable par un espace algébrique lisse sur $S$.  Alors le morphisme naturel 
$$\xymatrix{\piczero  \ar[r]& \pic}$$
est une immersion ouverte. De plus, si $S$ est localement noethérien, alors $\piczero$ est de type fini sur $S$.
\end{prop}
\begin{demo} Le fait que $\piczero$ soit un ouvert de $\pic$ résulte immédiatement du lemme~(\ref{CC_des_fibres_lemmeB}).
Comme il est par hypothèse localement de type fini sur $S$, il ne reste plus qu'à montrer qu'il est quasi-compact
lorsque $S$ est localement noethérien. Nous reprenons pour cela, en l'adaptant à notre cas, la démonstration
présentée par Kleiman dans~\cite{poly_Kleiman} (proposition~5.20). C'est une question locale sur $S$, donc on peut supposer que
$S$ est affine noethérien. Notons $\sigma$ le morphisme naturel de $\piczero$ vers $S$. Par récurrence noethérienne
sur les fermés de $S$, on peut supposer que pour tout fermé strict $Z$ de $S$, $\sigma^{-1}(Z)=\piczero\times_S Z$ est 
quasi-compact. Il suffit alors, par le même raisonnement que celui que nous avions suivi lors de la démonstration
du théorème~(\ref{critere_immersionsqc}), de construire un ouvert non vide $U$ de $S$ tel que $\sigma^{-1}(U)$ soit quasi-compact.

Soit $\pi : V \fleche \piczero$ un morphisme étale, où $V$ est un schéma affine non vide, et soit $U=\sigma(\pi(V))$.
Comme $\sigma \circ \pi$ est un morphisme
lisse, $U$ est un ouvert de $S$. Pour montrer que cet ouvert $U$ convient, on va construire un morphisme
surjectif de $V\times_S V$ vers $\piczero \times_S U$. Ce sera suffisant, puisque $V\times_S V$ est quasi-compact.
On a un diagramme commutatif de $S$-espaces algébriques :
$$\xymatrix{V\times_S V \ar[r]^-{\pi\times \pi} \ar[d] & \piczero\times_S \piczero \ar[r]^-{\alpha} & \piczero \ar[d]^{\sigma}\\
U \ar[rr] && S}$$
où $\alpha$ est défini fonctoriellement par $\alpha(g,h)=g.h^{-1}$.
On en déduit un morphisme $\alpha' : V\times_S V \fleche \piczero\times_S U$, dont il ne nous reste plus qu'à montrer
la surjectivité. En vertu du lemme~(\ref{lemme_surjectivite}) ci-dessous,
on peut supposer pour cela que $S$ est le spectre d'un corps algébriquement clos.
Mais dans ce cas on a $U=S$ et il faut montrer que le morphisme $\alpha \circ (\pi \times \pi)$ est surjectif.
De plus, vu que $S$ est le spectre d'un corps, $\piczero$ est un schéma en groupes et on peut supposer que $V$ est un
ouvert de $\piczero$. Alors le résultat découle du lemme~(\ref{lemme_alpha_surjectif}).
\end{demo}

\begin{lem}
\label{lemme_surjectivite}
Soient $S$ un schéma et $f : X \fleche Y$ un morphisme de $S$-espaces algébriques. Alors
$f$ est surjectif si et seulement si pour tout morphisme de $\Spec L$ vers $S$ avec $L$ un corps
algébriquement clos, le morphisme $f_L$ obtenu par changement de base est surjectif.
\end{lem}
\begin{demo}
\'Evident.
\end{demo}

\begin{lem}
\label{lemme_alpha_surjectif}
Soient $k$ un corps algébriquement clos, $G$ un $k$-schéma en groupes localement de type fini, et $U$ un ouvert
non vide de $G$. Alors le morphisme composé
$$\xymatrix@R=0.001pc{\beta : U\times U \ar[r] & G^0\times G^0 \ar[r]^{\alpha} &G^0 \\
&(g,h) \ar@{|->}[r] &g.h^{-1} }$$
est surjectif.
\end{lem}
\begin{demo}
Soit $x$ un $k$-point de $G^0$. Comme $G^0$ est irréductible, les ouverts $U$ et $xU$ sont d'intersection
non vide. Soit $y$ un $k$-point appartenant à $U\cap xU$. On note $z$ le $k$-point $x^{-1}y$. Alors $z$ est
dans $U$, et l'on a $x=yz^{-1}$, donc $x$ est dans l'image de $\beta$. L'ouvert $Uz^{-1}$ est un voisinage
de $x$ inclus dans l'image de $\beta$. Donc $\beta(U\times U)$ contient un ouvert qui
contient tous les $k$-points de $G^0$, et comme l'ensemble des $k$-points est très dense ceci achève la démonstration.
\end{demo}

% La propreté de la composante neutre peut se lire sur les fibres, si l'on sait qu'elle est séparée.
% 
% \begin{prop} On suppose que $\piczero$ est séparé et que pour tout $s\in S$ la fibre $\Pic_{\X_s/\kappa(s)}^0$ est propre
% sur $\kappa(s)$. Alors $\piczero$ est propre sur $S$. En particulier, c'est un fermé de $P$. 
% \end{prop}
% \begin{demo}
% \end{demo}
% 
% \begin{remarque} \rm cas de la caractéristique zéro ??????
% \end{remarque}
% 
% \begin{prop}
% \end{prop}
% \begin{demo}
% \end{demo}
% 
% \begin{prop}
% \end{prop}
% \begin{demo}
% \end{demo}
% 
% \section{La composante de torsion}
% 
% Encore une fois, nous allons définir a priori un sous-foncteur du foncteur de Picard, appelé composante de torsion,
% et nous montrerons seulement ensuite que sous de bonnes hypothèses ce foncteur
% est un sous-espace ouvert et de type fini du foncteur de Picard.

\chapter{Quelques exemples}

\section{L'espace projectif}

Soit $S=\Spec \Z$ et soit $X$ l'espace projectif $\mathbb{P}_{\Z}^n$ de dimension $n$ sur  $S$. On étudie dans cette section le foncteur et le champ de Picard de $X$. Le résultat ($\Pic_{X/S}=\Z$), bien connu, très classique, est sûrement présent dans de nombreux ouvrages. Voyons rapidement comment le redémontrer. On sait déjà par exemple grâce à \cite{FGA}~V thm.~3.1 que $\Pic_{X/S}$ est un schéma. Comme $X\fleche S$ a une section, on sait aussi que ce foncteur est égal à $P_{X/S}$. On a de plus un morphisme de foncteurs
$$\varphi : \Z \flechelongue \Pic_{X/S}$$
qui à un entier $l$ associe $\O(l)$. On va montrer que $\varphi$ est un isomorphisme. L'injectivité est facile. Pour la surjectivité, nous allons en fait construire l'isomorphisme réciproque. Nous aurons besoin du résultat préliminaire suivant.

\begin{lem}
Le morphisme $\Pic_{X/S} \fleche S$ est non ramifié.
\end{lem}
\begin{demo}
On sait qu'il est de toute manière localement de présentation finie d'après (\ref{lpf}). Il faut donc montrer qu'il est formellement non ramifié, c'est-à-dire que pour tout schéma affine $S'$, tout sous-schéma fermé $S'_0$ de $S'$ défini par un idéal
$I$ de carré nul, et tout morphisme de $S'$ dans $S$, l'application canonique
$$\Hom_S(S', \Pic_{X/S}) \flechelongue \Hom_S(S_0', \Pic_{X/S})$$
est injective. On se donne donc un faisceau inversible $\Lc$ sur $X'=X\times_S S'$ dont la restriction $\Lc_0$ à $X'_0=X\times_S S'_0$ provient de la base $S'_0$. Il faut montrer que $\Lc$ provient de la base $S'$. Par hypothèse, $\Lc_0$ provient d'un faisceau $\Bc_0$ sur $S'_0$, qui lui-même provient d'un faisceau inversible $\Bc$ sur $S'$ puisque $H^2(S', I)$ est nul (on utilise~(\ref{thm_defm_fi})). Maintenant $\Lc$ et $\Bc_{|_{X'}}$ sont deux déformations de $\Lc_0$ à $X'$. D'après le théorème~(\ref{thm_defm_fi}) elles sont nécessairement isomorphes puisque $H^1(X', I)$ est nul. Donc $\Lc$ provient de la base, ce qu'il fallait démontrer.
\end{demo}

Considérons maintenant un faisceau inversible $\Lc$ sur $\mathbb{P}_{T}^n$. On peut lui associer la fonction \og degré \fg\ sur $T$ définie ainsi
$$\fonction{\deg_{\Lc}}{T}{\Z}{t}{\textrm{degré de } \Lc_t}$$
où $\Lc_t$ désigne la restriction de $\Lc$ à la fibre $\mathbb{P}_{\kappa(t)}^n$. Nous allons montrer que la fonction $\deg_{\Lc}$ est localement constante sur $T$. Soient 
$\Lc$ et $\Mc$ deux faisceaux inversibles sur $\mathbb{P}_{T}^n$ et soient $\lambda$, $\mu$, les $T$-points de $\Pic_{X/S}$ associés. On forme le produit fibré :
$$\xymatrix{U \ar[r] \ar[d] \cartesien & T \ar[d]^{(\lambda, \mu)} \\
\Pic_{X/S} \ar[r]^-{\Delta} & \Pic_{X/S}\times_S \Pic_{X/S}.}$$
Vu que $\Pic_{X/S}$ est non ramifié, la diagonale $\Delta$ est une immersion ouverte si bien que $U$ est un ouvert de $T$. Par ailleurs, il est clair qu'un point $x$ de $T$ appartient à $U$ si et seulement si $\deg_{\Lc}(x)$ et $\deg_{\Mc}(x)$ sont égaux.
On en déduit facilement que $\deg_{\Lc}$ est localement constante sur $T$.

On peut donc définir un morphisme \og degré\fg
$$\xymatrix{\deg : \Pic_{X/S} \ar[r] & \Z}$$
qui à un faisceau inversible $\Lc$ associe $\deg_{\Lc}$. Il est clair que ce morphisme est l'inverse de $\varphi$.

Enfin pour déterminer le champ de Picard de $X$ on peut appliquer~(\ref{comparaison_champ_foncteur2}) et l'on voit qu'il est isomorphe à $\Z \times \bgm$.

$$\champic(\mathbb{P}_{\Z}^n/\Spec \Z) \simeq \Z \times \bgm$$

\section{Racine $n^{\text{ième}}$ d'un faisceau inversible}

Passons maintenant à un exemple un peu plus élaboré. Soient $X$ un $S$-schéma et $\Lc$ un faisceau inversible sur $X$.
Soit $n$ un entier strictement positif. On fabrique un
champ $[\L^{\frac1n}]$ à partir de ces données de la manière suivante.
Si $U$ est un objet de $\aff$, $[\Lc^{\frac1n}]_U$ est la catégorie des triplets
$(x,\Mc,\varphi)$ où
$$\left\{ \begin{array}{l}
x : U\fleche X \text{ est un élément de }X(U) \\
\Mc \text{ est un faisceau inversible sur }$U$\\
\varphi : \Mc^{\otimes n} \fleche x^*\Lc \text{ est un isomorphisme de faisceaux inversibles.}
\end{array} \right.$$
L'ensemble des morphismes de $(x,\Mc,\varphi)$ vers $(x',\Mc',\varphi')$ est vide si $x\neq x'$, et sinon c'est l'ensemble
des isomorphismes $\psi : \Mc \fleche \Mc'$ tels que $\psi^{\otimes n}$ soit compatible avec
$\varphi$ et $\varphi'$.

Remarquons que nous avons un morphisme canonique $\pi : [\L^{\frac1n}] \fleche X$. Si l'on regarde $[\L^{\frac1n}]$ comme
un groupoïde sur $X$, alors pour tout $U\in\ob{\rm (Aff/}X)$, la catégorie fibre $[\L^{\frac1n}]_U$ est simplement la catégorie des couples
$(\Mc,\varphi)$.

Si $U\in\ob\aff$ et si $\alpha$ est un objet de $[\L^{\frac1n}]_U$, un calcul rapide montre que le foncteur $\fAut_U(\alpha)$ est représentable par $\mun$. Plus généralement, si $\alpha_1,\alpha_2$ sont deux objets de $[\L^{\frac1n}]_U$, alors le foncteur $\fIsom(\alpha_1,\alpha_2)$
est représentable par un schéma fini sur $U$ (localement ce schéma est de la forme $\Spec(A[X]/(X^n-\gamma))$ où $\Spec A$ est un ouvert de $U$ qui trivialise les objets $\alpha_1$ et $\alpha_2$ et $\gamma$ est un élément de $A^{\times}$). En d'autres termes le morphisme diagonal
$$\Delta_{\pi} : [\L^{\frac1n}] \flechelongue [\L^{\frac1n}]\times_X [\L^{\frac1n}]$$
est schématique et fini. En particulier $[\L^{\frac1n}]$ est un $S$-préchamp. Il est clair que c'est même un $S$-champ.

\begin{remarque}\rm
\label{chgt_base_racine_nieme}
Si $f : Y \fleche X$ est un morphisme de $S$-schémas, alors le $S$-champ $[\L^{\frac1n}]\times_X Y$ est canoniquement
1-isomorphe à $[(f^*\Lc)^{\frac1n}]$ (évident).
\end{remarque}

\begin{prop}
Le champ $[\L^{\frac1n}]$ est une gerbe $\emph{fppf}$ sur $X$. Si $S$ est un $\Z[\frac1n]$-schéma, alors $[\L^{\frac1n}]$
est même une gerbe étale sur $X$.
\end{prop}
\begin{demo}
Il est clair que $[\L^{\frac1n}]$ a des objets partout localement pour la topologie de Zariski, puisque pour tout $U$ le
faisceau $\O_U$ a une racine $n^{\text{ième}}$ évidente. A fortiori $\pi$ est un épimorphisme (au sens que l'on veut, étale ou \emph{fppf}).

Montrons que le morphisme diagonal est un épimorphisme \emph{fppf}, et que c'est même un épimorphisme étale
lorsque $n$ est inversible. Soient $U\in\ob\aff$ et $x\in X(U)$. Il s'agit de montrer que deux objets
quelconques de $[\L^{\frac1n}]_U$ au-dessus de $x$ sont isomorphes localement pour la topologie \emph{fppf} (resp. étale).
Quitte à localiser pour la topologie de Zariski, on peut supposer que les faisceaux sous-jacents à $\alpha_1$ et
$\alpha_2$ sont triviaux. La question revient alors à montrer que tout élément de $\gm(U)$ admet une racine
$n^{\text{ième}}$ localement pour la topologie considérée. Soit $\gamma\in\gm(U)$. On note $U=\Spec A$ et $B=A[X]/(X^n-\gamma)$.
Il est clair que $B$ est fini et plat sur $A$. D'après le going-up theorem, le morphisme $\Spec B \fleche \Spec A$
est surjectif, donc c'est une famille couvrante pour la topologie \emph{fppf} qui répond au problème posé.
Si de plus $n$ est inversible alors c'est même une famille couvrante pour la topologie étale.
\end{demo}

\begin{prop}
Si $\L$ a une racine $n^{\text{ième}}$ sur $X$, i.e. s'il existe un faisceau inversible $\Mc$ sur $X$ tel que $\Mc^{\otimes n}$
soit isomorphe à $\L$, alors $[\L^{\frac1n}]$ est canoniquement (une fois qu'on a fixé $\Mc$ et un isomorphisme entre
$\Mc^{\otimes n}$ et $\L$) 1-isomorphe au champ classifiant $\bmunfppf$ du groupe $\mun$ pour la topologie \emph{fppf}.
En particulier c'est un champ algébrique (\cite{LMB}~(10.6) et~(10.13.1)).
\end{prop}
\begin{remarque}\rm
Si $n$ est inversible, alors $\mun$ est étale et les champs $\bmunfppf$ et $\bmun$ coïncident (\cite{LMB}~(9.6)). Dans
ce cas ce sont des champs de Deligne-Mumford. Notons que si $n$ n'est pas inversible, alors $\mun$ n'est pas
lisse, et le champ $\bmun$ qui classifie les $\mun$-torseurs étales n'a aucune raison \emph{a priori} d'être algébrique.
\end{remarque}
\begin{demo}
La donnée d'un faisceau inversible $\Mc$ sur $X$ et d'un isomorphisme
$$\varphi : \flechen{\Mc^{\otimes n}}{\sim}{\L}$$
définit une section $s : X \fleche [\L^{\frac1n}]$ du morphisme structural $\pi$. Donc la gerbe \emph{fppf}
$[\L^{\frac1n}]$ sur $X$ est une gerbe neutre (au sens de \cite{LMB}~(3.20)). Le résultat découle donc de l'analogue \emph{fppf}
de \cite{LMB}~(3.21). 
\end{demo}

\begin{remarque}\rm
On a déjà vu que $[\L^{\frac1n}]$ est un $S$-champ \emph{fppf}, donc \emph{a fortiori} un $S$-champ (étale), et que sa diagonale
est représentable, séparée, et quasi-compacte (puisque $\Delta_{\pi}$ est schématique et fini et que le morphisme
$[\L^{\frac1n}]\times_X [\L^{\frac1n}] \fleche [\L^{\frac1n}]\times_S [\L^{\frac1n}]$, obtenu par changement
de base à partir de la diagonale de $X/S$, est une immersion quasi-compacte). Soit $X' \fleche X$ une famille
couvrante (pour la topologie de Zariski) telle que $\L_{|X'}$ soit trivial. Alors d'après la proposition
précédente et la remarque~(\ref{chgt_base_racine_nieme}), on a un diagramme 2-cartésien
$$\xymatrix{ \bmunfppf \ar[r] \ar[d] \cartesien& [\L^{\frac1n}] \ar[d] \\
X'\ar[r]& X}$$
ce qui prouve que $[\L^{\frac1n}]$ est un champ algébrique.
\end{remarque}

\begin{remarque}\rm
On suppose que $S$ est un $\Z[\frac1n]$-schéma et que $X$ est noethérien. Alors $\pi$ est propre, lisse, de
présentation finie, et cohomologiquement plat en dimension zéro. En particulier si $X/S$ vérifie ces
propriétés, le morphisme $[\L^{\frac1n}] \fleche S$ les vérifie aussi.
\end{remarque}

\vskip 5mm
\noindent
{\sc Calcul du groupe de Picard de $[\L^{\frac1n}]$}

\begin{lem}
Soient $X$ un schéma et $A$ un schéma en groupes commutatifs sur $X$. Soit $\Fc$ un
faisceau inversible sur une $A$-gerbe (\emph{fppf}) $\pi : \X \fleche X$.
Il existe un unique $X$-morphisme de schémas en groupes
$$\chi_{\Fc} : A \flechelongue \gm$$
tel que l'action naturelle de $A$ sur $\Fc$ soit induite par $\chi_{\Fc}$ et par la multiplication $\Fc\times \gm \fleche \Fc$
induite par la structure de $\O_{\X}$-module de $\Fc$, autrement dit tel que le diagramme suivant soit commutatif :
$$\shorthandoff{!;:?}
\xymatrix@!0 @R=1.5pc @C=3pc{A\times \Fc \ar[rd] \ar[dd] &\\ &\Fc\\ \gm\times\Fc \ar[ur]&}$$
\end{lem}
\begin{demo} Un faisceau inversible sur la gerbe $\X$ est la donnée, pour tout $x\in\ob\X_U$, d'un faisceau inversible $\Fc_x$
sur $U$, et pour tout morphisme $\varphi : x\fleche x'$ dans $\X$, d'un isomorphisme
$$L_{\Fc}(\varphi) : \Fc_x \fleche \pi(\varphi)^*\Fc_{x'}$$
ces isomorphismes vérifiant de plus une condition de compatibilité évidente.

Construisons d'abord $\chi_{\Fc}(U)$ pour un $U\in\ob{\rm (Aff/}X)$ sur lequel $\X$ a des objets.
Soit $x\in\ob\X_U$ et soit $g\in A(U)$. Via l'identification entre $A(U)$ et $\Aut(x)$, $g$ correspond à un
automorphisme $\varphi$ de $x$, et induit de ce fait un automorphisme $L_{\Fc}(\varphi)$ de $\Fc_x$. Cet
automorphisme correspond à la multiplication par un unique élément de $\gm(U)$, que l'on note pour l'instant $\chi_{\Fc}(U)(x,g)$.
Maintenant si $x$ et $x'$ sont deux objets de $\X_U$, on vérifie facilement que $\chi_{\Fc}(U)(x,g)=\chi_{\Fc}(U)(x',g)$
en utilisant la condition de compatibilité entre les $L_{\Fc}(\varphi)$, le fait que $x$ et $x'$ sont localement isomorphes
pour la topologie \emph{fppf} et le fait que $\gm$ est un faisceau pour cette même topologie. D'où le morphisme
$\chi_{\Fc}(U) : A(U) \fleche \gm(U)$. Il est clair, vu sa construction, que ce morphisme est déterminé de manière
unique par les actions naturelles de $A$ et de $\gm$ sur $\Fc$.

\'Etant donné que $\X$ a des objets partout localement pour la topologie \emph{fppf}, cette collection de morphismes
se prolonge de manière unique en un caractère $\chi_{\Fc}$ de $A$ vérifiant les propriétés annoncées.
\end{demo}

\begin{remarque}\rm
Si $\chi$ est un caractère fixé de $A$, un faisceau inversible $\Fc$ sur $\X$ est un faisceau $\chi$-tordu
de degré $d$ (au sens de Lieblich, \cite{Lieblich}~2.1.2.2) si et seulement si $\chi_{\Fc}=\chi^d$.
\end{remarque}

\begin{pptes}
\label{racine_pptes_caractere}
\begin{itemize}
\item[(1)] La construction de $\chi_{\Fc}$ est compatible au changement de base en un sens évident.
\item[(2)] Si $\Fc$ et $\Gc$ sont des faisceaux inversibles sur $\X$, alors $$\chi_{\Fc\otimes \Gc} = \chi_{\Fc} . \chi_{\Gc}\, .$$
\item[(3)] Un faisceau inversible $\Fc$ sur $\X$ provient de $X$ si et seulement si $\chi_{\Fc}$ est trivial.
\end{itemize}
\end{pptes}
\begin{demo}
La première propriété est évidente : il suffit de l'écrire. La deuxième résulte immédiatement du fait que
si l'automorphisme $L_{\Fc}(\varphi)$ (resp. $L_{\Gc}(\varphi)$) de $\Fc_x$ (resp. $\Gc_x$) est la multiplication
par $\chi_{\Fc}(\varphi)$ (resp. $\chi_{\Gc}(\varphi)$), alors l'automorphisme $L_{\Fc}(\varphi)\otimes L_{\Gc}(\varphi)$
de $\Fc_x\otimes \Gc_x$ est la multiplication par $\chi_{\Fc}(\varphi).\chi_{\Gc}(\varphi)$. Montrons maintenant le dernier
point.

Supposons tout d'abord que $\Fc$ soit isomorphe à un faisceau de la forme $\pi^*\Mc$, où $\Mc$ est un faisceau
inversible sur $X$. Il est clair que $\chi_{\Fc}$ est égal à $\chi_{\pi^*\Mc}$ donc il suffit de montrer que
$\chi_{\pi^*\Mc}$ est trivial, c'est-à-dire que pour tout objet $x$ de $\X$ et pour tout automorphisme
$\varphi$ de $x$, l'automorphisme $L_{\pi^*\Mc}(\varphi)$ de $(\pi^*\Mc)_x$ est l'identité. C'est évident par
construction de l'image inverse.

Réciproquement, supposons $\chi_{\Fc}$ trivial, i.e. supposons que pour tout objet $x$ de $\X$ et tout
automorphisme $\varphi$ de $x$, $L_{\Fc}(\varphi)$ soit l'identité de $\Fc_x$. En utilisant le lemme~(\ref{critere_relative_trivialite})
on voit que provenir de la base est une question locale sur $X$ pour la topologie \emph{fppf}.
On peut donc supposer, grâce à la propriété (1) ci-dessus, que $\X$ est une gerbe neutre, c'est-à-dire que
le morphisme structural $\pi : \X \fleche X$ a une section $s : X \fleche \X$.

On va alors montrer que $\Fc$ est isomorphe à $\pi^*s^*\Fc$. Le faisceau $\pi^*s^*\Fc$ est celui qui à tout objet $x$ de $\X$
associe $\Fc_{s(\pi(x))}$, les isomorphismes de changement de base étant simplement les isomorphismes canoniques.
Nous allons construire une collection d'isomorphismes $\rho_x$ de $\Fc_x$ dans $\Fc_{s(\pi(x))}$ (pour chaque objet $x$ de $\X$), 
compatibles avec les $L_{\Fc}(\varphi)$ et les $L_{\pi^*s^*\Fc}(\varphi)$.

Si $x$ et $s(\pi(x))$ sont isomorphes dans $\X_U$ on choisit un isomorphisme $\varphi$ de $x$ dans $s(\pi(x))$ et on pose
$\rho_x=L_{\Fc}(\varphi)$. Si $\varphi_1$ et $\varphi_2$ sont deux tels isomorphismes, alors $(\varphi_2)^{-1}\circ \varphi_1$
est un automorphisme de $x$, donc d'après l'hypothèse sur $\Fc$ on a $L_{\Fc}((\varphi_2)^{-1}\circ \varphi_1)=\Id_{\Fc_x}$
de sorte que $\rho_x$ est bien défini et ne dépend pas du choix de $\varphi$.

Dans le cas général, on sait que $x$ et $s(\pi(x))$ sont localement isomorphes pour la topologie $\emph{fppf}$ puisque
$\X$ est une gerbe. Vu l'unicité dans la construction de $\rho_x$ lorsque $x$ est isomorphe à $s(\pi(x))$, il est clair
qu'il existe un unique isomorphisme $\rho_x : \Fc_x \fleche \Fc_{s(\pi(x))}$ compatible avec ceux construits dans le cas
précédent. La collection de tous les $\rho_x$ ainsi construits répond au problème posé.
\end{demo}

\begin{exemple}[groupe de Picard de B$G$] \rm
En utilisant cette construction, on retrouve facilement le groupe de Picard du champ classifiant B$G$, où $G$ est un
$X$-schéma en groupes abéliens. En effet, le morphisme structural $\pi : \text{B}G \fleche X$ a une section, donc
$\pi^* : \Pic(X) \fleche \Pic(\text{B}G)$ a une rétraction et en particulier il est injectif. D'après la proposition
précédente, l'application $\Fc \mapsto \chi_{\Fc}$ induit un morphisme de groupes de $\Pic(\text{B}G)$ dans $\widehat{G}$
dont $\Pic(X)$ est le noyau. Ce morphisme est naturellement scindé : si $\chi : G \fleche \gm$ est un caractère
de $G$, on lui associe la classe du faisceau inversible $\Lc(\chi)$ construit de la manière suivante. Pour tout
$U\in\ob\aff$ et tout $G$-torseur $\widetilde{U}$ on définit $\Lc(\chi)_{\widetilde{U}}$ comme étant le
faisceau inversible correspondant au $\gm$-torseur sur $U$ obtenu à partir de $\widetilde{U}$ par extension
du groupe structural via le caractère $\chi$. On a donc une suite exacte courte scindée :
$$\xymatrix{1 \ar[r]& \Pic(X) \ar[r] &\Pic(\text{B}G)\ar[r] & \widehat{G}\ar[r] &1}$$
de sorte que $\Pic(\text{B}G)$ est naturellement isomorphe au produit $\Pic(X)\times \widehat{G}$.
\end{exemple}

Dans le cas du champ $[\L^{\frac1n}]$, on a un faisceau inversible \og canonique\fg, que nous noterons $\Omega$, et que
l'on peut construire de la manière suivante. Pour tout $U\in\ob\aff$ et tout objet $\alpha=(x,\Mc,\varphi)$ de
$[\L^{\frac1n}]_U$, on pose $\Omega_{\alpha}=\Mc$. Les isomorphismes de changement de base sont définis de manière
évidente. Il est clair que le caractère $\chi_{\Omega}$ associé à $\Omega$ est simplement l'injection canonique
$$\chi : \mun \flechelongue \gm.$$
En particulier, vu que $\chi$ n'est pas le caractère trivial, on peut en déduire que le faisceau $\Omega$ ne
provient pas de la base $X$ ! Enfin, pour tout objet $\alpha=(x,\Mc,\varphi)$ de
$[\L^{\frac1n}]_U$ (où $U\in\ob\aff$), le faisceau $(\Omega^{\otimes n})_{\alpha}$, qui n'est autre que $\Mc^{\otimes n}$, 
est canoniquement isomorphe à $x^*\Lc$ (via $\varphi$ !) de sorte que l'on a un isomorphisme canonique
$$\Phi : \flechen{\Omega^{\otimes n}}{\sim}{\pi^*\Lc}.$$
Nous avons maintenant presque tous les éléments en main pour décrire de manière complète le groupe de Picard
de $[\L^{\frac1n}]$. Nous terminons le travail dans la proposition ci-dessous.

\begin{prop} On note $l$ la classe du faisceau $\Lc$ dans $\Pic(X)$ et $\omega$ celle de $\Omega$ dans $\Pic([\L^{\frac1n}])$.
\begin{itemize}
\item[(1)] Le morphisme $\pi^* : \Pic(X) \fleche \Pic([\L^{\frac1n}])$ est injectif, et l'on a une suite exacte courte :
$$\xymatrix{1 \ar[r]& \Pic(X) \ar[r]& \Pic([\L^{\frac1n}]) \ar[r]& \widehat{\mun} \ar[r]& 1.}$$
\item[(2)] Le groupe $\Pic([\L^{\frac1n}])$ est isomorphe au quotient du groupe $\Pic(X) \times H^0(X,\Z)$ par le sous $H^0(X,\Z)$-module engendré par $(l^{-1},n)$ (autrement dit par la relation $\omega^n=l$).
\end{itemize}
\end{prop}
\begin{demo}
La propriété~(\ref{racine_pptes_caractere})~(2) montre que l'application $\Fc \mapsto \chi_{\Fc}$ induit un morphisme
de groupes de $\Pic([\L^{\frac1n}])$ dans $\widehat{\mun}$. Il est clair que ce morphisme est surjectif, vu que
$\widehat{\mun}$ est isomorphe au groupe $H^0(X,\Z/n\Z)$,
engendré par l'injection canonique $\chi : \mun \fleche \gm$, et que $\chi=\chi_{\Omega}$.
La propriété~(\ref{racine_pptes_caractere})~(3) montre que la suite ci-dessus est exacte en $\Pic([\L^{\frac1n}])$.
Pour en finir avec le premier point il nous reste donc juste à montrer l'injectivité de $\pi^*$.

Soit $\Nc$ un faisceau inversible sur $X$ et soit $f$ un isomorphisme de $\pi^*\Nc$ dans $\O_{[\L^{\frac1n}]}$. Il s'agit
de montrer que $\Nc$ est trivial. L'isomorphisme $f$ est donné par une collection d'isomorphismes
$$f_{\alpha} : \flechen{(\pi^*\Nc)_{\alpha}=x^*\Nc}{\sim}{(\O_{[\L^{\frac1n}]})_{\alpha}=\O_U}$$ pour
chaque objet $\alpha=(x,\Mc,\varphi)$ de $[\L^{\frac1n}]_U$, ces isomorphismes vérifiant de plus une condition de
compatibilité que nous nous dispensons d'expliciter, mais qui, en particulier, entraîne que les isomorphismes $f_{\alpha}$
et $f_{\alpha'}$ associés à deux objets $\alpha$ et $\alpha'$ au-dessus d'un même élément $x$ de $X(U)$
sont égaux dès que $\alpha$ et $\alpha'$ sont isomorphes. Si $x\in X(U)$ est un objet de $X$ au-dessus duquel $[\L^{\frac1n}]_U$
a des objets, alors vu que deux objets $\alpha$ et $\alpha'$ sont localement isomorphes, les $f_{\alpha}$ pour
$\alpha\in\ob [\L^{\frac1n}]_U$ sont tous égaux et définissent donc une section $f_x$ du faisceau $\fIsom(x^*\Nc,\O_U)$.
En fait, vu que $[\L^{\frac1n}]$ a des objets partout localement pour la topologie \emph{fppf}, la collection
des $f_{\alpha}$ définit de manière unique un élément $f_x\in \fIsom(x^*\Nc,\O_U)(U)$ pour tout $x\in X(U)$, que
$[\L^{\frac1n}]$ ait des objets au-dessus de $x$ ou non. On vérifie facilement que les isomorphismes ainsi
construits
$$f_x : \flechen{x^*\Nc}{\sim}{\O_U}$$
forment un système compatible d'isomorphismes et définissent donc un
isomorphisme de $\Nc$ dans $\O_X$, ce qui achève la démonstration.

Pour le point (2), notons $G$ le quotient du groupe $\Pic(X) \times H^0(X,\Z)$ par le sous $H^0(X,\Z)$-module engendré par la relation $\omega^n=l$. On a clairement un morphisme de $G$ dans $\Pic([\L^{\frac1n}])$ qui envoie $(0,1)$ sur
$\omega$. En utilisant les propriétés précédentes, et le fait que $\widehat{\mun}$ est isomorphe à $H^0(X,\Z/n\Z)$
et engendré par $\chi_{\Omega}$, on vérifie très facilement que ce morphisme est un isomorphisme.
\end{demo}

\begin{exemple}\rm
Prenons pour $X$ l'espace projectif $\P^k$ sur $\Spec \Z$. Alors $\Pic(X)$ est isomorphe à $\Z$. On fixe un entier relatif $l$ et on pose $\Lc=\O(l)$. La proposition précédente permet de calculer $\Pic([\L^{\frac1n}])$ pour tout $n$ appartenant à $\N^*$. Par exemple si $l=1$, on trouve $\frac1n\Z$. Si $l$ est un multiple de $n$, on est dans le cas où $\Lc$ a une racine \iem{n} et l'on trouve $\Z\times \Z/n\Z$. Dans le cas général, le groupe $\Pic([\L^{\frac1n}])$ est isomorphe (de manière non canonique) à $\frac{d}{n}\Z\times \Z/d\Z$ où $d$ est le pgcd de $n$ et $l$.
\end{exemple}

\vskip 5mm
\noindent
{\sc Foncteur de Picard relatif de $[\L^{\frac1n}]/S$}

Notons $\X=[\L^{\frac1n}]$. Pour tout schéma $U$ sur $S$, on a une suite exacte courte :
$$\xymatrix{1 \ar[r]& \Pic(X\times_S U) \ar[r]& \Pic(\X\times_S U) \ar[r]& H^0(X\times_S U, \Z/n\Z) \ar[r]& 1.}$$
Elle induit la suite exacte 
$$\xymatrix{1 \ar[r]& \frac{\Pic(X\times_S U)}{\Pic(U)} \ar[r]& \frac{\Pic(\X\times_S U)}{\Pic(U)} \ar[r]& H^0(X\times_S U, \Z/n\Z) \ar[r]& 1,}$$
d'où une suite exacte de préfaisceaux
$$\xymatrix{1 \ar[r]& P_{X/S} \ar[r]& P_{\X/S} \ar[r]& f_*\Z/n\Z \ar[r]& 1.}$$
En appliquant à cette suite exacte le foncteur \og faisceau \'etale associé\fg\ on obtient une suite exacte de faisceaux étales :
\begin{equation}
\label{sec_racine_nieme}
\xymatrix{1 \ar[r]& \Pic_{X/S} \ar[r]^{\varphi_0}& \Pic_{\X/S} \ar[r]^{\chi}& f_*\Z/n\Z \ar[r]& 1.}
\end{equation}

\begin{remarque}\rm
Si $\Lc$ a une racine ${n}^{\textrm{ième}}$ $\Rc$, la suite exacte (\ref{sec_racine_nieme}) est scindée par $i \mapsto (\omega r^{-1})^i$ où $r$ est la classe de $\Rc$ dans $\Pic(X)$, si bien que $\pic$ s'identifie au produit $\Pic_{X/S}\times_S f_*\Z/n\Z$.
\end{remarque}

\begin{remarque}\rm
Le faisceau $f_*\Z/n\Z$ n'est \emph{a priori} pas représentable. En conséquence, dans le cas général, il ne suffit pas que $\Pic_{X/S}$ soit représentable pour que $\pic$ le soit, même lorsque la gerbe $\X$ est triviale. Cependant si $f$ est ouvert, dominant et à fibres géométriquement connexes, alors $f_*\Z/n\Z= \Z/n\Z$. C'est le cas par exemple lorsque $f$ est localement de type fini, plat et cohomologiquement plat en dimension zéro.
\end{remarque}

On suppose maintenant que $f_*\Z/n\Z$ coïncide avec $\Z/n\Z$ et l'on considère le produit $\Pic_{X/S}\times_S \Z$ de $\Pic_{X/S}$ par
le groupe constant $\Z$. On va voir que $\pic$ est isomorphe au quotient de
$\Pic_{X/S}\times_S \Z$ par la relation $\omega^n=l$. On note $H$ le sous-groupe engendré par $(l^{-1},n)$. Le foncteur $\Pic_{X/S}\times_S \Z$
s'identifie à une union disjointe $\disp \coprod_{i\in \Z} (\Pic_{X/S})_i$ de
copies de $\Pic_{X/S}$ indexées par $\Z$. Pour tout $i$ appartenant à $\Z$, on note
$\mu_{\omega^i}$ le morphisme de multiplication par $\omega^i$ de $\pic$ dans 
lui-même,
%\footnote{La notation $\omega^i$ désigne la classe du faisceau
%$\Omega^{\otimes i}$, ce qui a un sens bien défini puisque la fonction $i$ est 
%localement constante sur $X$.}
et on note $\varphi_i$ le morphisme composé
$$\xymatrix{\Pic_{X/S} \ar[r]^{\varphi_0}& \pic \ar[r]^{\mu_{\omega^i}}& \pic.}$$
La collection des $\varphi_i$ définit donc un morphisme
$$\xymatrix{\varphi : \Pic_{X/S}\times_S \Z \ar[r]& \pic}$$
dont il est clair qu'il est invariant sous $H$. On vérifie facilement avec les suites exactes précédentes qu'il est universel pour les morphismes invariants sous $H$ de $\Pic_{X/S}\times_S \Z$ à valeurs dans un $S$-schéma $T$. Le foncteur $\pic$ s'identifie donc bien au quotient évoqué ci-dessus. On peut construire ce quotient \og à la main \fg\ comme suit (voir figure ci-dessous). Pour tout couple d'entiers $(i,k)$ on identifie les copies de $\Pic_{X/S}$ numéro $i$ et $i+nk$ via l'isomorphisme de translation 
$$\flechen{(\Pic_{X/S})_{i+nk}}{\mu_{l^k}}{(\Pic_{X/S})_i}.$$
La loi de groupe est induite naturellement par celle de $\Pic_{X/S}$ et par la relation $\omega^n=l$.

%%%%%%%%%%%%%%%%%%%%%%%%%%%%%%%%%%%%%%%%%%%%%%%%%%%%%%%%%%%%%%%%%%%%%%%%%%%%%%%
% This is fig4tex.tex, version 1.8, Feb 6, 2007
%
%           Authors : Yvon Lafranche, Daniel Martin
%                     IRMAR, Universite de Rennes 1 - France
%
%           E-mail  : Yvon.Lafranche@univ-rennes1.fr
%                     Daniel.Martin@univ-rennes1.fr
%
%           Web     : http://perso.univ-rennes1.fr/yvon.lafranche/
%                     http://perso.univ-rennes1.fr/daniel.martin/
%%%%%%%%%%%%%%%%%%%%%%%%%%%%%%%%%%%%%%%%%%%%%%%%%%%%%%%%%%%%%%%%%%%%%%%%%%%%%%%
\ifx\figfortexisloaded\relax \else\let\figfortexisloaded=\relax\fi
\message{version 1.8}
% CAUTION :
%   Some internal commands depend on the running environment.
%   This package has been tested with TeX and LaTeX on Mac OS X, UNIX, Windows,
%   and with TeXtures on MacOS 9.
%   In the latter case, the following \iftextures control sequence should be
%   initialized to "true".
\newif\iftextures
\catcode`\@=11
%%%%%%%%%%%%%%%%%%%%%%%%%%%%%%%%%%%%%%%%%%%%%%%%%%%%%%%%%%%%%%%%%%%%%%%%%%%%%%%
% Points with numbers >= 0 are devoted to the user.
% Points with numbers <  0 are reserved to internal use.
%%%%%%%%%%%%%%%%%%%%%%%%%%%%%%%%%%%%%%%%%%%%%%%%%%%%%%%%%%%%%%%%%%%%%%%%%%%%%%%
\newdimen\epsil@n\epsil@n=0.00005pt
\newdimen\Cepsil@n\Cepsil@n=0.005pt
\newdimen\dcq@\dcq@=254pt
\newdimen\PI@\PI@=3.141592pt
\newdimen\DemiPI@deg\DemiPI@deg=90pt
\newdimen\PI@deg\PI@deg=180pt
\newdimen\DePI@deg\DePI@deg=360pt
\chardef\t@n=10
\chardef\c@nt=100
\chardef\@lxxiv=74
\chardef\@xci=91
\mathchardef\@nMnCQn=9949
\chardef\@vi=6
\chardef\@xxx=30
\chardef\@lvi=56
\chardef\@lxxi=71
\chardef\@lxxxv=85
\mathchardef\@mmmmlxviii=4068
\mathchardef\@ccclx=360
\mathchardef\@dccxx=720
\newcount\p@rtent \newcount\f@ctech \newcount\result@tent
\newdimen\v@lmin \newdimen\v@lmax \newdimen\v@leur
\newdimen\result@t\newdimen\result@@t
\newdimen\mili@u \newdimen\c@rre \newdimen\delt@
\def\degT@rd{0.017453 }  % pi/180
\def\rdT@deg{57.295779 } % 180/pi
{\catcode`p=12 \catcode`t=12 \gdef\v@leurseule#1pt{#1}}
\def\repdecn@mb#1{\expandafter\v@leurseule\the#1\space}
\def\arct@n#1(#2,#3){{\v@lmin=#2\v@lmax=#3%
    \maxim@m{\mili@u}{-\v@lmin}{\v@lmin}\maxim@m{\c@rre}{-\v@lmax}{\v@lmax}%
    \delt@=\mili@u\m@ech\mili@u%
    \ifdim\c@rre>\@nMnCQn\mili@u\divide\v@lmax\tw@\c@lATAN\v@leur(\z@,\v@lmax)% DY > 9949 DX
    \else%
    \maxim@m{\mili@u}{-\v@lmin}{\v@lmin}\maxim@m{\c@rre}{-\v@lmax}{\v@lmax}%
    \m@ech\c@rre%
    \ifdim\mili@u>\@nMnCQn\c@rre\divide\v@lmin\tw@% DX > 9949 DY
    \maxim@m{\mili@u}{-\v@lmin}{\v@lmin}\c@lATAN\v@leur(\mili@u,\z@)%
    \else\c@lATAN\v@leur(\delt@,\v@lmax)\fi\fi%
    \ifdim\v@lmin<\z@\v@leur=-\v@leur\ifdim\v@lmax<\z@\advance\v@leur-\PI@%
    \else\advance\v@leur\PI@\fi\fi%
    \global\result@t=\v@leur}#1=\result@t}
\def\m@ech#1{\ifdim#1>1.646pt\divide\mili@u\t@n\divide\c@rre\t@n\m@ech#1\fi}
\def\c@lATAN#1(#2,#3){{\v@lmin=#2\v@lmax=#3\v@leur=\z@\delt@=\tw@ pt%
    \un@iter{0.785398}{\v@lmax<}%
    \un@iter{0.463648}{\v@lmax<}%
    \un@iter{0.244979}{\v@lmax<}%
    \un@iter{0.124355}{\v@lmax<}%
    \un@iter{0.062419}{\v@lmax<}%
    \un@iter{0.031240}{\v@lmax<}%
    \un@iter{0.015624}{\v@lmax<}%
    \un@iter{0.007812}{\v@lmax<}%
    \un@iter{0.003906}{\v@lmax<}%
    \un@iter{0.001953}{\v@lmax<}%
    \un@iter{0.000976}{\v@lmax<}%
    \un@iter{0.000488}{\v@lmax<}%
    \un@iter{0.000244}{\v@lmax<}%
    \un@iter{0.000122}{\v@lmax<}%
    \un@iter{0.000061}{\v@lmax<}%
    \un@iter{0.000030}{\v@lmax<}%
    \un@iter{0.000015}{\v@lmax<}%
    \global\result@t=\v@leur}#1=\result@t}
\def\un@iter#1#2{%
    \divide\delt@\tw@\edef\dpmn@{\repdecn@mb{\delt@}}%
    \mili@u=\v@lmin%
    \ifdim#2\z@%
      \advance\v@lmin-\dpmn@\v@lmax\advance\v@lmax\dpmn@\mili@u%
      \advance\v@leur-#1pt%
    \else%
      \advance\v@lmin\dpmn@\v@lmax\advance\v@lmax-\dpmn@\mili@u%
      \advance\v@leur#1pt%
    \fi}
\def\c@ssin#1#2#3{\expandafter\ifx\csname COS@\number#3\endcsname\relax\c@lCS{#3pt}%
    \expandafter\xdef\csname COS@\number#3\endcsname{\repdecn@mb\result@t}%
    \expandafter\xdef\csname SIN@\number#3\endcsname{\repdecn@mb\result@@t}\fi%
    \edef#1{\csname COS@\number#3\endcsname}\edef#2{\csname SIN@\number#3\endcsname}}
\def\c@lCS#1{{\mili@u=#1\p@rtent=\@ne%
    \relax\ifdim\mili@u<\z@\red@ng<-\else\red@ng>+\fi\f@ctech=\p@rtent%
    \relax\ifdim\mili@u<\z@\mili@u=-\mili@u\f@ctech=-\f@ctech\fi\c@@lCS}}
\def\c@@lCS{\v@lmin=\mili@u\c@rre=-\mili@u\advance\c@rre\DemiPI@deg\v@lmax=\c@rre%
    \mili@u\@lxxi\mili@u\divide\mili@u\@mmmmlxviii%
    \edef\v@larg{\repdecn@mb{\mili@u}}\mili@u=-\v@larg\mili@u%
    \edef\v@lmxde{\repdecn@mb{\mili@u}}%
    \c@rre\@lxxi\c@rre\divide\c@rre\@mmmmlxviii%
    \edef\v@largC{\repdecn@mb{\c@rre}}\c@rre=-\v@largC\c@rre%
    \edef\v@lmxdeC{\repdecn@mb{\c@rre}}%
    \fctc@s\mili@u\v@lmin\global\result@t\p@rtent\v@leur%
    \let\t@mp=\v@larg\let\v@larg=\v@largC\let\v@largC=\t@mp%
    \let\t@mp=\v@lmxde\let\v@lmxde=\v@lmxdeC\let\v@lmxdeC=\t@mp%
    \fctc@s\c@rre\v@lmax\global\result@@t\f@ctech\v@leur}
\def\fctc@s#1#2{\v@leur=#1\relax\ifdim#2<\@lxxxv\p@\cosser@h\else\sinser@t\fi}
\def\cosser@h{\advance\v@leur\@lvi\p@\divide\v@leur\@lvi%
    \v@leur=\v@lmxde\v@leur\advance\v@leur\@xxx\p@%
    \v@leur=\v@lmxde\v@leur\advance\v@leur\@ccclx\p@%
    \v@leur=\v@lmxde\v@leur\advance\v@leur\@dccxx\p@\divide\v@leur\@dccxx}
\def\sinser@t{\v@leur=\v@lmxdeC\p@\advance\v@leur\@vi\p@%
    \v@leur=\v@largC\v@leur\divide\v@leur\@vi}
\def\red@ng#1#2{\relax\ifdim\mili@u#1#2\DemiPI@deg\advance\mili@u#2-\PI@deg%
    \p@rtent=-\p@rtent\red@ng#1#2\fi}
\def\invers@#1#2{{\v@leur=#2\maxim@m{\v@lmax}{-\v@leur}{\v@leur}%
    \f@ctech=\@ne\m@inv@rs%
    \multiply\v@leur\f@ctech\edef\v@lv@leur{\repdecn@mb{\v@leur}}%
    \p@rtentiere{\p@rtent}{\v@leur}\v@lmin=\p@\divide\v@lmin\p@rtent%
    \inv@rs@\multiply\v@lmax\f@ctech\global\result@t=\v@lmax}#1=\result@t}
\def\m@inv@rs{\ifdim\v@lmax<\p@\multiply\v@lmax\t@n\multiply\f@ctech\t@n\m@inv@rs\fi}
\def\inv@rs@{\v@lmax=-\v@lmin\v@lmax=\v@lv@leur\v@lmax%
    \advance\v@lmax\tw@ pt\v@lmax=\repdecn@mb{\v@lmin}\v@lmax%
    \delt@=\v@lmax\advance\delt@-\v@lmin\ifdim\delt@<\z@\delt@=-\delt@\fi%
    \ifdim\delt@>\epsil@n\v@lmin=\v@lmax\inv@rs@\fi}
\def\minim@m#1#2#3{\relax\ifdim#2<#3#1=#2\else#1=#3\fi}
\def\maxim@m#1#2#3{\relax\ifdim#2>#3#1=#2\else#1=#3\fi}
\def\p@rtentiere#1#2{#1=#2\divide#1by65536 }
\def\r@undint#1#2{{\v@leur=#2\divide\v@leur\t@n\p@rtentiere{\p@rtent}{\v@leur}%
    \v@leur=\p@rtent pt\global\result@t=\t@n\v@leur}#1=\result@t}
\def\sqrt@#1#2{{\v@leur=#2%
    \minim@m{\v@lmin}{\p@}{\v@leur}\maxim@m{\v@lmax}{\p@}{\v@leur}%
    \f@ctech=\@ne\m@sqrt@\sqrt@@%
    \mili@u=\v@lmin\advance\mili@u\v@lmax\divide\mili@u\tw@\multiply\mili@u\f@ctech%
    \global\result@t=\mili@u}#1=\result@t}
\def\m@sqrt@{\ifdim\v@leur>\dcq@\divide\v@leur\c@nt\v@lmax=\v@leur%
    \multiply\f@ctech\t@n\m@sqrt@\fi}
\def\sqrt@@{\mili@u=\v@lmin\advance\mili@u\v@lmax\divide\mili@u\tw@%
    \c@rre=\repdecn@mb{\mili@u}\mili@u%
    \ifdim\c@rre<\v@leur\v@lmin=\mili@u\else\v@lmax=\mili@u\fi%
    \delt@=\v@lmax\advance\delt@-\v@lmin\ifdim\delt@>\epsil@n\sqrt@@\fi}
\def\extrairelepremi@r#1\de#2{\expandafter\lepremi@r#2@#1#2}
\def\lepremi@r#1,#2@#3#4{\def#3{#1}\def#4{#2}\ignorespaces}
\def\@cfor#1:=#2\do#3{%
  \edef\@fortemp{#2}%
  \ifx\@fortemp\empty\else\@cforloop#2,\@nil,\@nil\@@#1{#3}\fi}
\def\@cforloop#1,#2\@@#3#4{%
  \def#3{#1}%
  \ifx#3\Fig@nnil\let\@nextwhile=\Fig@fornoop\else#4\relax\let\@nextwhile=\@cforloop\fi%
  \@nextwhile#2\@@#3{#4}}

\def\@ecfor#1:=#2\do#3{%
  \def\@@cfor{\@cfor#1:=}%
  \edef\@@@cfor{#2}%
  \expandafter\@@cfor\@@@cfor\do{#3}}
\def\Fig@nnil{\@nil}
\def\Fig@fornoop#1\@@#2#3{}
\def\trtlis@rg#1#2{\def\list@@rg{#1}%
    \@ecfor\p@rv@l:=\list@@rg\do{\expandafter#2\p@rv@l|}}
\newbox\b@xvisu
\newtoks\let@xte
\newif\ifitis@K
\newcount\s@mme
\newcount\l@mbd@un \newcount\l@mbd@de
\newcount\superc@ntr@l\superc@ntr@l=\@ne        % Controle impose
\newcount\typec@ntr@l\typec@ntr@l=\superc@ntr@l % Controle souhaite
\newdimen\v@lX  \newdimen\v@lY  \newdimen\v@lZ
\newdimen\v@lXa \newdimen\v@lYa \newdimen\v@lZa
\newdimen\unit@\unit@=\p@ % Initialisation a la valeur par defaut.
\def\unit@util{pt}
\def\ptT@ptps{0.996264 }
\def\ptpsT@pt{1.00375 }
\def\ptT@unit@{1} % Initialisation correspondant a la valeur par defaut de \unit@
\def\setunit@#1{\def\unit@util{#1}\setunit@@#1:\invers@{\result@t}{\unit@}%
    \edef\ptT@unit@{\repdecn@mb\result@t}}
\def\setunit@@#1#2:{\ifcat#1a\unit@=\@ne#1#2\else\unit@=#1#2\fi}
\def\d@fm@cdim#1#2{{\v@leur=#2\v@leur=\ptT@unit@\v@leur\xdef#1{\repdecn@mb\v@leur}}}
\newif\ifBdingB@x\BdingB@xtrue
\newdimen\c@@rdXmin \newdimen\c@@rdYmin  % Dimensions de la BoundingBox
\newdimen\c@@rdXmax \newdimen\c@@rdYmax
\def\b@undb@x#1#2{\ifBdingB@x%
    \relax\ifdim#1<\c@@rdXmin\global\c@@rdXmin=#1\fi%
    \relax\ifdim#2<\c@@rdYmin\global\c@@rdYmin=#2\fi%
    \relax\ifdim#1>\c@@rdXmax\global\c@@rdXmax=#1\fi%
    \relax\ifdim#2>\c@@rdYmax\global\c@@rdYmax=#2\fi\fi}
\def\b@undb@xP#1{{\Figg@tXY{#1}\b@undb@x{\v@lX}{\v@lY}}}
\def\ellBB@x#1;#2,#3(#4,#5,#6){{\s@uvc@ntr@l\et@tellBB@x%
    \setc@ntr@l{2}\figptell-2::#1;#2,#3(#4,#6)\b@undb@xP{-2}%
    \figptell-2::#1;#2,#3(#5,#6)\b@undb@xP{-2}%
    \c@ssin{\C@}{\S@}{#6}\v@lmin=\C@ pt\v@lmax=\S@ pt%
    \mili@u=#3\v@lmin\delt@=#2\v@lmax\arct@n\v@leur(\delt@,\mili@u)%
    \mili@u=-#3\v@lmax\delt@=#2\v@lmin\arct@n\c@rre(\delt@,\mili@u)%
    \v@leur=\rdT@deg\v@leur\advance\v@leur-\DePI@deg%
    \c@rre=\rdT@deg\c@rre\advance\c@rre-\DePI@deg%
    \v@lmin=#4pt\v@lmax=#5pt%
    \loop\ifdim\v@leur<\v@lmax\ifdim\v@leur>\v@lmin%
    \edef\@ngle{\repdecn@mb\v@leur}\figptell-2::#1;#2,#3(\@ngle,#6)%
    \b@undb@xP{-2}\fi\advance\v@leur\PI@deg\repeat%
    \loop\ifdim\c@rre<\v@lmax\ifdim\c@rre>\v@lmin%
    \edef\@ngle{\repdecn@mb\c@rre}\figptell-2::#1;#2,#3(\@ngle,#6)%
    \b@undb@xP{-2}\fi\advance\c@rre\PI@deg\repeat%
    \resetc@ntr@l\et@tellBB@x}\ignorespaces}
\def\initb@undb@x{\c@@rdXmin=\maxdimen\c@@rdYmin=\maxdimen%
    \c@@rdXmax=-\maxdimen\c@@rdYmax=-\maxdimen}
\def\c@ntr@lnum#1{%
    \relax\ifnum\typec@ntr@l=\@ne%
    \ifnum#1<\z@%
    \immediate\write16{*** Forbidden point number (#1). Abort.}\end\fi\fi%
    \set@bjc@de{#1}}
\def\set@bjc@de#1{\edef\objc@de{@BJ\ifnum#1<\z@ M\romannumeral-#1\else\romannumeral#1\fi}}
\def\setc@ntr@l#1{\ifnum\superc@ntr@l>#1\typec@ntr@l=\superc@ntr@l%
    \else\typec@ntr@l=#1\fi}
\def\resetc@ntr@l#1{\global\superc@ntr@l=#1\setc@ntr@l{#1}}
\def\s@uvc@ntr@l#1{\edef#1{\the\superc@ntr@l}}
\def\c@lproscalDD#1[#2,#3]{{\Figg@tXY{#2}%
    \edef\Xu@{\repdecn@mb{\v@lX}}\edef\Yu@{\repdecn@mb{\v@lY}}\Figg@tXY{#3}%
    \global\result@t=\Xu@\v@lX\global\advance\result@t\Yu@\v@lY}#1=\result@t}
\def\c@lproscalTD#1[#2,#3]{{\Figg@tXY{#2}\edef\Xu@{\repdecn@mb{\v@lX}}%
    \edef\Yu@{\repdecn@mb{\v@lY}}\edef\Zu@{\repdecn@mb{\v@lZ}}%
    \Figg@tXY{#3}\global\result@t=\Xu@\v@lX\global\advance\result@t\Yu@\v@lY%
    \global\advance\result@t\Zu@\v@lZ}#1=\result@t}
\def\c@lprovec#1{%
    \det@rmC\v@lZa(\v@lX,\v@lY,\v@lmin,\v@lmax)%
    \det@rmC\v@lXa(\v@lY,\v@lZ,\v@lmax,\v@leur)%
    \det@rmC\v@lYa(\v@lZ,\v@lX,\v@leur,\v@lmin)%
    \Figv@ctCreg#1(\v@lXa,\v@lYa,\v@lZa)}
\def\det@rm#1[#2,#3]{{\Figg@tXY{#2}\Figg@tXYa{#3}%
    \delt@=\repdecn@mb{\v@lX}\v@lYa\advance\delt@-\repdecn@mb{\v@lY}\v@lXa%
    \global\result@t=\delt@}#1=\result@t}
\def\det@rmC#1(#2,#3,#4,#5){{\global\result@t=\repdecn@mb{#2}#5%
    \global\advance\result@t-\repdecn@mb{#3}#4}#1=\result@t}
\def\getredf@ctDD#1(#2,#3){{\maxim@m{\v@lXa}{-#2}{#2}\maxim@m{\v@lYa}{-#3}{#3}%
    \maxim@m{\v@lXa}{\v@lXa}{\v@lYa}% \v@lXa = ||X||inf
    \ifdim\v@lXa>\@xci pt\divide\v@lXa\@xci%
    \p@rtentiere{\p@rtent}{\v@lXa}\advance\p@rtent\@ne\else\p@rtent=\@ne\fi%
    \global\result@tent=\p@rtent}#1=\result@tent\ignorespaces}
\def\getredf@ctTD#1(#2,#3,#4){{\maxim@m{\v@lXa}{-#2}{#2}\maxim@m{\v@lYa}{-#3}{#3}%
    \maxim@m{\v@lZa}{-#4}{#4}\maxim@m{\v@lXa}{\v@lXa}{\v@lYa}%
    \maxim@m{\v@lXa}{\v@lXa}{\v@lZa}% \v@lXa = ||X||inf
    \ifdim\v@lXa>\@lxxiv pt\divide\v@lXa\@lxxiv%
    \p@rtentiere{\p@rtent}{\v@lXa}\advance\p@rtent\@ne\else\p@rtent=\@ne\fi%
    \global\result@tent=\p@rtent}#1=\result@tent\ignorespaces}
\def\FigptintercircB@zDD#1:#2:#3,#4[#5,#6,#7,#8]{{\s@uvc@ntr@l\et@tfigptintercircB@zDD%
    \setc@ntr@l{2}\figvectPDD-1[#5,#8]\Figg@tXY{-1}\getredf@ctDD\f@ctech(\v@lX,\v@lY)%
    \mili@u=#4\unit@\divide\mili@u\f@ctech\c@rre=\repdecn@mb{\mili@u}\mili@u%
    \figptBezierDD-5::#3[#5,#6,#7,#8]%
    \v@lmin=#3\p@\v@lmax=\v@lmin\advance\v@lmax0.1\p@%
    \loop\edef\T@{\repdecn@mb{\v@lmax}}\figptBezierDD-2::\T@[#5,#6,#7,#8]%
    \figvectPDD-1[-5,-2]\n@rmeucCDD{\delt@}{-1}\ifdim\delt@<\c@rre\v@lmin=\v@lmax%
    \advance\v@lmax0.1\p@\repeat%
    \loop\mili@u=\v@lmin\advance\mili@u\v@lmax%
    \divide\mili@u\tw@\edef\T@{\repdecn@mb{\mili@u}}\figptBezierDD-2::\T@[#5,#6,#7,#8]%
    \figvectPDD-1[-5,-2]\n@rmeucCDD{\delt@}{-1}\ifdim\delt@>\c@rre\v@lmax=\mili@u%
    \else\v@lmin=\mili@u\fi\v@leur=\v@lmax\advance\v@leur-\v@lmin%
    \ifdim\v@leur>\epsil@n\repeat\figptcopyDD#1:#2/-2/%
    \resetc@ntr@l\et@tfigptintercircB@zDD}\ignorespaces}
\def\inters@cDD#1:#2[#3,#4;#5,#6]{{\s@uvc@ntr@l\et@tinters@cDD%
    \setc@ntr@l{2}\vecunit@{-1}{#4}\vecunit@{-2}{#6}%
    \Figg@tXY{-1}\setc@ntr@l{1}\Figg@tXYa{#3}%
    \edef\A@{\repdecn@mb{\v@lX}}\edef\B@{\repdecn@mb{\v@lY}}%
    \v@lmin=\B@\v@lXa\advance\v@lmin-\A@\v@lYa%
    \Figg@tXYa{#5}\setc@ntr@l{2}\Figg@tXY{-2}%
    \edef\C@{\repdecn@mb{\v@lX}}\edef\D@{\repdecn@mb{\v@lY}}%
    \v@lmax=\D@\v@lXa\advance\v@lmax-\C@\v@lYa%
    \delt@=\A@\v@lY\advance\delt@-\B@\v@lX%
    \invers@{\v@leur}{\delt@}\edef\v@ldelta{\repdecn@mb{\v@leur}}%
    \v@lXa=\A@\v@lmax\advance\v@lXa-\C@\v@lmin%
    \v@lYa=\B@\v@lmax\advance\v@lYa-\D@\v@lmin%
    \v@lXa=\v@ldelta\v@lXa\v@lYa=\v@ldelta\v@lYa%
    \setc@ntr@l{1}\Figp@intregDD#1:{#2}(\v@lXa,\v@lYa)%
    \resetc@ntr@l\et@tinters@cDD}\ignorespaces}
\def\inters@cTD#1:#2[#3,#4;#5,#6]{{\s@uvc@ntr@l\et@tinters@cTD%
    \setc@ntr@l{2}\figvectNVTD-1[#4,#6]\figvectNVTD-2[#6,-1]\figvectPTD-1[#3,#5]%
    \r@pPSTD\v@leur[-2,-1,#4]\edef\v@lcoef{\repdecn@mb{\v@leur}}%
    \figpttraTD#1:{#2}=#3/\v@lcoef,#4/\resetc@ntr@l\et@tinters@cTD}\ignorespaces}
\def\r@pPSTD#1[#2,#3,#4]{{\Figg@tXY{#2}\edef\Xu@{\repdecn@mb{\v@lX}}%
    \edef\Yu@{\repdecn@mb{\v@lY}}\edef\Zu@{\repdecn@mb{\v@lZ}}%
    \Figg@tXY{#3}\v@lmin=\Xu@\v@lX\advance\v@lmin\Yu@\v@lY\advance\v@lmin\Zu@\v@lZ%
    \Figg@tXY{#4}\v@lmax=\Xu@\v@lX\advance\v@lmax\Yu@\v@lY\advance\v@lmax\Zu@\v@lZ%
    \invers@{\v@leur}{\v@lmax}\global\result@t=\repdecn@mb{\v@leur}\v@lmin}%
    #1=\result@t}
\def\n@rminfDD#1#2{{\Figg@tXY{#2}\maxim@m{\v@lX}{\v@lX}{-\v@lX}%
    \maxim@m{\v@lY}{\v@lY}{-\v@lY}\maxim@m{\global\result@t}{\v@lX}{\v@lY}}%
    #1=\result@t}
\def\n@rminfTD#1#2{{\Figg@tXY{#2}\maxim@m{\v@lX}{\v@lX}{-\v@lX}%
    \maxim@m{\v@lY}{\v@lY}{-\v@lY}\maxim@m{\v@lZ}{\v@lZ}{-\v@lZ}%
    \maxim@m{\v@lX}{\v@lX}{\v@lY}\maxim@m{\global\result@t}{\v@lX}{\v@lZ}}%
    #1=\result@t}
\def\n@rmeucCDD#1#2{\Figg@tXY{#2}\divide\v@lX\f@ctech\divide\v@lY\f@ctech%
    #1=\repdecn@mb{\v@lX}\v@lX\v@lX=\repdecn@mb{\v@lY}\v@lY\advance#1\v@lX}
\def\n@rmeucCTD#1#2{\Figg@tXY{#2}%
    \divide\v@lX\f@ctech\divide\v@lY\f@ctech\divide\v@lZ\f@ctech%
    #1=\repdecn@mb{\v@lX}\v@lX\v@lX=\repdecn@mb{\v@lY}\v@lY\advance#1\v@lX%
    \v@lX=\repdecn@mb{\v@lZ}\v@lZ\advance#1\v@lX}
\def\n@rmeucSVDD#1#2{{\Figg@tXY{#2}%
    \v@lXa=\repdecn@mb{\v@lX}\v@lX\v@lYa=\repdecn@mb{\v@lY}\v@lY%
    \advance\v@lXa\v@lYa\sqrt@{\global\result@t}{\v@lXa}}#1=\result@t}
\def\n@rmeucSVTD#1#2{{\Figg@tXY{#2}\v@lXa=\repdecn@mb{\v@lX}\v@lX%
    \v@lYa=\repdecn@mb{\v@lY}\v@lY\v@lZa=\repdecn@mb{\v@lZ}\v@lZ%
    \advance\v@lXa\v@lYa\advance\v@lXa\v@lZa\sqrt@{\global\result@t}{\v@lXa}}#1=\result@t}
\def\n@rmeucDD#1#2{{\Figg@tXY{#2}\getredf@ctDD\f@ctech(\v@lX,\v@lY)%
    \divide\v@lX\f@ctech\divide\v@lY\f@ctech%
    \v@lXa=\repdecn@mb{\v@lX}\v@lX\v@lYa=\repdecn@mb{\v@lY}\v@lY%
    \advance\v@lXa\v@lYa\sqrt@{\global\result@t}{\v@lXa}%
    \global\multiply\result@t\f@ctech}#1=\result@t}
\def\n@rmeucTD#1#2{{\Figg@tXY{#2}\getredf@ctTD\f@ctech(\v@lX,\v@lY,\v@lZ)%
    \divide\v@lX\f@ctech\divide\v@lY\f@ctech\divide\v@lZ\f@ctech%
    \v@lXa=\repdecn@mb{\v@lX}\v@lX%
    \v@lYa=\repdecn@mb{\v@lY}\v@lY\v@lZa=\repdecn@mb{\v@lZ}\v@lZ%
    \advance\v@lXa\v@lYa\advance\v@lXa\v@lZa\sqrt@{\global\result@t}{\v@lXa}%
    \global\multiply\result@t\f@ctech}#1=\result@t}
\def\vecunit@DD#1#2{{\Figg@tXY{#2}\getredf@ctDD\f@ctech(\v@lX,\v@lY)%
    \divide\v@lX\f@ctech\divide\v@lY\f@ctech%
    \Figv@ctCreg#1(\v@lX,\v@lY)\n@rmeucSV{\v@lYa}{#1}%
    \invers@{\v@lXa}{\v@lYa}\edef\v@lv@lXa{\repdecn@mb{\v@lXa}}%
    \v@lX=\v@lv@lXa\v@lX\v@lY=\v@lv@lXa\v@lY%
    \Figv@ctCreg#1(\v@lX,\v@lY)\multiply\v@lYa\f@ctech\global\result@t=\v@lYa}}
\def\vecunit@TD#1#2{{\Figg@tXY{#2}\getredf@ctTD\f@ctech(\v@lX,\v@lY,\v@lZ)%
    \divide\v@lX\f@ctech\divide\v@lY\f@ctech\divide\v@lZ\f@ctech%
    \Figv@ctCreg#1(\v@lX,\v@lY,\v@lZ)\n@rmeucSV{\v@lYa}{#1}%
    \invers@{\v@lXa}{\v@lYa}\edef\v@lv@lXa{\repdecn@mb{\v@lXa}}%
    \v@lX=\v@lv@lXa\v@lX\v@lY=\v@lv@lXa\v@lY\v@lZ=\v@lv@lXa\v@lZ%
    \Figv@ctCreg#1(\v@lX,\v@lY,\v@lZ)\multiply\v@lYa\f@ctech\global\result@t=\v@lYa}}
\def\vecunitC@TD[#1,#2]{\Figg@tXYa{#1}\Figg@tXY{#2}%
    \advance\v@lX-\v@lXa\advance\v@lY-\v@lYa\advance\v@lZ-\v@lZa\c@lvecunitTD}
\def\vecunitCV@TD#1{\Figg@tXY{#1}\c@lvecunitTD}
\def\c@lvecunitTD{\getredf@ctTD\f@ctech(\v@lX,\v@lY,\v@lZ)%
    \divide\v@lX\f@ctech\divide\v@lY\f@ctech\divide\v@lZ\f@ctech%
    \v@lXa=\repdecn@mb{\v@lX}\v@lX%
    \v@lYa=\repdecn@mb{\v@lY}\v@lY\v@lZa=\repdecn@mb{\v@lZ}\v@lZ%
    \advance\v@lXa\v@lYa\advance\v@lXa\v@lZa\sqrt@{\v@lYa}{\v@lXa}%
    \invers@{\v@lXa}{\v@lYa}\edef\v@lv@lXa{\repdecn@mb{\v@lXa}}%
    \v@lX=\v@lv@lXa\v@lX\v@lY=\v@lv@lXa\v@lY\v@lZ=\v@lv@lXa\v@lZ}
\def\figgetangleDD#1[#2,#3,#4]{\ifps@cri{\s@uvc@ntr@l\et@tfiggetangleDD\setc@ntr@l{2}%
    \figvectPDD-1[#2,#3]\figvectPDD-2[#2,#4]\vecunit@{-1}{-1}%
    \c@lproscalDD\delt@[-2,-1]\figvectNVDD-1[-1]\c@lproscalDD\v@leur[-2,-1]%
    \arct@n\v@lmax(\delt@,\v@leur)\v@lmax=\rdT@deg\v@lmax%
    \ifdim\v@lmax<\z@\advance\v@lmax\DePI@deg\fi\xdef#1{\repdecn@mb{\v@lmax}}%
    \resetc@ntr@l\et@tfiggetangleDD}\ignorespaces\fi}
\def\figgetangleTD#1[#2,#3,#4,#5]{\ifps@cri{\s@uvc@ntr@l\et@tfiggetangleTD\setc@ntr@l{2}%
    \figvectPTD-1[#2,#3]\figvectPTD-2[#2,#5]\figvectNVTD-3[-1,-2]%
    \figvectPTD-2[#2,#4]\figvectNVTD-4[-3,-1]%
    \vecunit@{-1}{-1}\c@lproscalTD\delt@[-2,-1]\c@lproscalTD\v@leur[-2,-4]%
    \arct@n\v@lmax(\delt@,\v@leur)\v@lmax=\rdT@deg\v@lmax%
    \ifdim\v@lmax<\z@\advance\v@lmax\DePI@deg\fi\xdef#1{\repdecn@mb{\v@lmax}}%
    \resetc@ntr@l\et@tfiggetangleTD}\ignorespaces\fi}    
\def\figgetdist#1[#2,#3]{\ifps@cri{\s@uvc@ntr@l\et@tfiggetdist\setc@ntr@l{2}%
    \figvectP-1[#2,#3]\n@rmeuc{\v@lX}{-1}\v@lX=\ptT@unit@\v@lX\xdef#1{\repdecn@mb{\v@lX}}%
    \resetc@ntr@l\et@tfiggetdist}\ignorespaces\fi}
\def\Figg@tT#1{\c@ntr@lnum{#1}%
    {\expandafter\expandafter\expandafter\extr@ctT\csname\objc@de\endcsname:%
     \ifnum\B@@ltxt=\z@\ptn@me{#1}\else\csname\objc@de T\endcsname\fi}}
\def\extr@ctT#1,#2,#3/#4:{\def\B@@ltxt{#3}}
\def\Figg@tXY#1{\c@ntr@lnum{#1}%
    \expandafter\expandafter\expandafter\extr@ctC\csname\objc@de\endcsname:}
\def\extr@ctCDD#1/#2,#3,#4:{\v@lX=#2\v@lY=#3}
\def\extr@ctCTD#1/#2,#3,#4:{\v@lX=#2\v@lY=#3\v@lZ=#4}
\def\Figg@tXYa#1{\c@ntr@lnum{#1}%
    \expandafter\expandafter\expandafter\extr@ctCa\csname\objc@de\endcsname:}
\def\extr@ctCaDD#1/#2,#3,#4:{\v@lXa=#2\v@lYa=#3}
\def\extr@ctCaTD#1/#2,#3,#4:{\v@lXa=#2\v@lYa=#3\v@lZa=#4}
\def\figinit#1{\initpr@lim\Figinit@#1,:\initpss@ttings\ignorespaces}
\def\Figinit@#1,#2:{\setunit@{#1}\def\t@xt@{#2}\ifx\t@xt@\empty\else\Figinit@@#2:\fi}
\def\Figinit@@#1#2:{\if#12 \else\Figs@tproj{#1}\initTD@\fi}
\newif\ifTr@isDim
\def\UnD@fined{UNDEFINED}
\def\ifundefined#1{\expandafter\ifx\csname#1\endcsname\relax}
\def\initpr@lim{\initb@undb@x\figsetmark{}\figsetptname{$A_{##1}$}\def\Sc@leFact{1}%
    \initDD@\figsetroundcoord{yes}\ps@critrue\expandafter\setupd@te\defaultupdate:%
    \edef\disob@unit{\UnD@fined}\edef\t@rgetpt{\UnD@fined}}
\def\initDD@{\Tr@isDimfalse%
    \ifPDFm@ke%
     \let\Ps@rcerc=\Ps@rcercBz%
     \let\Ps@rell=\Ps@rellBz%
    \fi
    \let\c@lDCUn=\c@lDCUnDD%
    \let\c@lDCDeux=\c@lDCDeuxDD%
    \let\c@ldefproj=\relax%
    \let\c@lproscal=\c@lproscalDD%
    \let\c@lprojSP=\relax%
    \let\extr@ctC=\extr@ctCDD%
    \let\extr@ctCa=\extr@ctCaDD%
    \let\extr@ctCF=\extr@ctCFDD%
    \let\Figp@intreg=\Figp@intregDD%
    \let\Figpts@xes=\Figpts@xesDD%
    \let\n@rmeucSV=\n@rmeucSVDD\let\n@rmeuc=\n@rmeucDD\let\n@rminf=\n@rminfDD%
    \let\pr@dMatV=\pr@dMatVDD%
    \let\ps@xes=\ps@xesDD%
    \let\vecunit@=\vecunit@DD%
    \let\figcoord=\figcoordDD%
    \let\figgetangle=\figgetangleDD%
    \let\figpt=\figptDD%
    \let\figptBezier=\figptBezierDD%
    \let\figptbary=\figptbaryDD%
    \let\figptcirc=\figptcircDD%
    \let\figptcircumcenter=\figptcircumcenterDD%
    \let\figptcopy=\figptcopyDD%
    \let\figptcurvcenter=\figptcurvcenterDD%
    \let\figptell=\figptellDD%
    \let\figptendnormal=\figptendnormalDD%
    \def\figptinterlineplane{\un@v@ilable{figptinterlineplane}}%
    \let\figptinterlines=\inters@cDD%
    \let\figptorthocenter=\figptorthocenterDD%
    \let\figptorthoprojline=\figptorthoprojlineDD%
    \def\figptorthoprojplane{\un@v@ilable{figptorthoprojplane}}%
    \let\figptrot=\figptrotDD%
    \let\figptscontrol=\figptscontrolDD%
    \let\figptsintercirc=\figptsintercircDD%
    \let\figptsinterlinell=\figptsinterlinellDD%
    \let\figptsorthoprojline=\figptsorthoprojlineDD%
    \let\figptsrot=\figptsrotDD%
    \let\figptssym=\figptssymDD%
    \let\figptstra=\figptstraDD%
    \let\figptsym=\figptsymDD%
    \let\figpttraC=\figpttraCDD%
    \let\figpttra=\figpttraDD%
    \def\figptvisilimSL{\un@v@ilable{figptvisilimSL}}%
    \def\figsetobdist{\un@v@ilable{figsetobdist}}%
    \def\figsettarget{\un@v@ilable{figsettarget}}%
    \def\figsetview{\un@v@ilable{figsetview}}%
    \let\figvectDBezier=\figvectDBezierDD%
    \let\figvectN=\figvectNDD%
    \let\figvectNV=\figvectNVDD%
    \let\figvectP=\figvectPDD%
    \let\figvectU=\figvectUDD%
    \let\psarccircP=\psarccircPDD%
    \let\psarccirc=\psarccircDD%
    \let\psarcell=\psarcellDD%
    \let\psarcellPA=\psarcellPADD%
    \let\psarrowBezier=\psarrowBezierDD%
    \let\psarrowcircP=\psarrowcircPDD%
    \let\psarrowcirc=\psarrowcircDD%
    \let\psarrowhead=\psarrowheadDD%
    \let\psarrow=\psarrowDD%
    \let\psBezier=\psBezierDD%
    \let\pscirc=\pscircDD%
    \let\pscurve=\pscurveDD%
    \let\psnormal=\psnormalDD%
    }
\def\initTD@{\Tr@isDimtrue\initb@undb@xTD\newt@rgetptfalse\newdis@bfalse%
    \let\c@lDCUn=\c@lDCUnTD%
    \let\c@lDCDeux=\c@lDCDeuxTD%
    \let\c@ldefproj=\c@ldefprojTD%
    \let\c@lproscal=\c@lproscalTD%
    \let\extr@ctC=\extr@ctCTD%
    \let\extr@ctCa=\extr@ctCaTD%
    \let\extr@ctCF=\extr@ctCFTD%
    \let\Figp@intreg=\Figp@intregTD%
    \let\Figpts@xes=\Figpts@xesTD%
    \let\n@rmeucSV=\n@rmeucSVTD\let\n@rmeuc=\n@rmeucTD\let\n@rminf=\n@rminfTD%
    \let\pr@dMatV=\pr@dMatVTD%
    \let\ps@xes=\ps@xesTD%
    \let\vecunit@=\vecunit@TD%
    \let\figcoord=\figcoordTD%
    \let\figgetangle=\figgetangleTD%
    \let\figpt=\figptTD%
    \let\figptBezier=\figptBezierTD%
    \let\figptbary=\figptbaryTD%
    \let\figptcirc=\figptcircTD%
    \let\figptcircumcenter=\figptcircumcenterTD%
    \let\figptcopy=\figptcopyTD%
    \let\figptcurvcenter=\figptcurvcenterTD%
    \let\figptinterlineplane=\figptinterlineplaneTD%
    \let\figptinterlines=\inters@cTD%
    \let\figptorthocenter=\figptorthocenterTD%
    \let\figptorthoprojline=\figptorthoprojlineTD%
    \let\figptorthoprojplane=\figptorthoprojplaneTD%
    \let\figptrot=\figptrotTD%
    \let\figptscontrol=\figptscontrolTD%
    \let\figptsintercirc=\figptsintercircTD%
    \let\figptsorthoprojline=\figptsorthoprojlineTD%
    \let\figptsorthoprojplane=\figptsorthoprojplaneTD%
    \let\figptsrot=\figptsrotTD%
    \let\figptssym=\figptssymTD%
    \let\figptstra=\figptstraTD%
    \let\figptsym=\figptsymTD%
    \let\figpttraC=\figpttraCTD%
    \let\figpttra=\figpttraTD%
    \let\figptvisilimSL=\figptvisilimSLTD%
    \let\figsetobdist=\figsetobdistTD%
    \let\figsettarget=\figsettargetTD%
    \let\figsetview=\figsetviewTD%
    \let\figvectDBezier=\figvectDBezierTD%
    \let\figvectN=\figvectNTD%
    \let\figvectNV=\figvectNVTD%
    \let\figvectP=\figvectPTD%
    \let\figvectU=\figvectUTD%
    \let\psarccircP=\psarccircPTD%
    \let\psarccirc=\psarccircTD%
    \let\psarcell=\psarcellTD%
    \let\psarcellPA=\psarcellPATD%
    \let\psarrowBezier=\psarrowBezierTD%
    \let\psarrowcircP=\psarrowcircPTD%
    \let\psarrowcirc=\psarrowcircTD%
    \let\psarrowhead=\psarrowheadTD%
    \let\psarrow=\psarrowTD%
    \let\psBezier=\psBezierTD%
    \let\pscirc=\pscircTD%
    \let\pscurve=\pscurveTD%
    }
\def\un@v@ilable#1{\immediate\write16{*** The macro #1 is not available in the current context.}}
\def\figinsert#1{{\def\t@xt@{#1}\relax\ifx\t@xt@\empty\Figinsert@\DefGIfilen@me,:%
    \else\expandafter\FiginsertNu@#1 :\fi}\ignorespaces}
\def\FiginsertNu@#1 #2:{\def\t@xt@{#1}\relax\ifx\t@xt@\empty\def\t@xt@{#2}%
    \ifx\t@xt@\empty\Figinsert@\DefGIfilen@me,:\else\FiginsertNu@#2:\fi%
    \else\expandafter\FiginsertNd@#1 #2:\fi}
\def\FiginsertNd@#1#2:{\ifcat#1a\Figinsert@#1#2,:\else\Figinsert@\DefGIfilen@me,#1#2,:\fi}
\def\Figinsert@#1,#2:{\def\t@xt@{#2}\ifx\t@xt@\empty\xdef\Sc@leFact{1}\else%
    \def\Xarg@##1,{\def\@rgdeux{##1}}\Xarg@#2\xdef\Sc@leFact{\@rgdeux}\fi\@psfgetbb{#1}%
    \v@lX=\@psfllx\p@\v@lX=\ptpsT@pt\v@lX\v@lX=\Sc@leFact\v@lX%
    \v@lY=\@psflly\p@\v@lY=\ptpsT@pt\v@lY\v@lY=\Sc@leFact\v@lY%
    \b@undb@x{\v@lX}{\v@lY}%
    \v@lX=\@psfurx\p@\v@lX=\ptpsT@pt\v@lX\v@lX=\Sc@leFact\v@lX%
    \v@lY=\@psfury\p@\v@lY=\ptpsT@pt\v@lY\v@lY=\Sc@leFact\v@lY%
    \b@undb@x{\v@lX}{\v@lY}%
    \ifPDFm@ke\Figinclud@PDF{#1}{\Sc@leFact}\else%
    \v@lX=\c@nt pt\v@lX=\Sc@leFact\v@lX\edef\F@ct{\repdecn@mb{\v@lX}}%
    \iftextures\special{postscriptfile #1 vscale=\F@ct\space hscale=\F@ct}%
    \else\includegraphics{#1}\fi\fi% 
    \message{[#1]}\ignorespaces}
\def\figinsertE#1{\FiginsertE@#1,:\ignorespaces}
\def\FiginsertE@#1,#2:{{\def\t@xt@{#2}\ifx\t@xt@\empty\xdef\Sc@leFact{1}\else%
    \def\Xarg@##1,{\def\@rgdeux{##1}}\Xarg@#2\xdef\Sc@leFact{\@rgdeux}\fi%
    \pdfximage{#1}\setbox\Gb@x=\hbox{\pdfrefximage\pdflastximage}%
    \v@lX=\z@\v@lY=-\Sc@leFact\dp\Gb@x\b@undb@x{\v@lX}{\v@lY}%
    \advance\v@lX\Sc@leFact\wd\Gb@x\advance\v@lY\Sc@leFact\dp\Gb@x%
    \advance\v@lY\Sc@leFact\ht\Gb@x\b@undb@x{\v@lX}{\v@lY}%
    \v@lX=\Sc@leFact\wd\Gb@x\pdfximage width \v@lX {#1}%
    \rlap{\pdfrefximage\pdflastximage}\message{[#1]}}\ignorespaces}
\def\figptDD#1:#2(#3,#4){\ifps@cri\c@ntr@lnum{#1}%
    {\v@lX=#3\unit@\v@lY=#4\unit@\Fig@dmpt{#2}{\z@}}\ignorespaces\fi}
\def\Fig@dmpt#1#2{\def\t@xt@{#1}\ifx\t@xt@\empty\def\B@@ltxt{\z@}%
    \else\expandafter\gdef\csname\objc@de T\endcsname{#1}\def\B@@ltxt{\@ne}\fi%
    \expandafter\xdef\csname\objc@de\endcsname{\ifitis@vect@r\C@dCl@svect%
    \else\C@dCl@spt\fi,\z@,\B@@ltxt/\the\v@lX,\the\v@lY,#2}}
\def\C@dCl@spt{P}
\def\C@dCl@svect{V}
\def\figptTD#1:#2(#3,#4){\ifps@cri\c@ntr@lnum{#1}%
    \def\c@@rdYZ{#4,0,0}\extrairelepremi@r\c@@rdY\de\c@@rdYZ%
    \extrairelepremi@r\c@@rdZ\de\c@@rdYZ%
    {\v@lX=#3\unit@\v@lY=\c@@rdY\unit@\v@lZ=\c@@rdZ\unit@\Fig@dmpt{#2}{\the\v@lZ}%
    \b@undb@xTD{\v@lX}{\v@lY}{\v@lZ}}\ignorespaces\fi}
\def\Figp@intregDD#1:#2(#3,#4){\c@ntr@lnum{#1}%
    {\result@t=#4\v@lX=#3\v@lY=\result@t\Fig@dmpt{#2}{\z@}}\ignorespaces}
\def\Figp@intregTD#1:#2(#3,#4){\c@ntr@lnum{#1}%
    \def\c@@rdYZ{#4,\z@,\z@}\extrairelepremi@r\c@@rdY\de\c@@rdYZ%
    \extrairelepremi@r\c@@rdZ\de\c@@rdYZ%
    {\v@lX=#3\v@lY=\c@@rdY\v@lZ=\c@@rdZ\Fig@dmpt{#2}{\the\v@lZ}%
    \b@undb@xTD{\v@lX}{\v@lY}{\v@lZ}}\ignorespaces}
\def\figptBezierDD#1:#2:#3[#4,#5,#6,#7]{\ifps@cri{\s@uvc@ntr@l\et@tfigptBezierDD%
    \FigptBezier@#3[#4,#5,#6,#7]\Figp@intregDD#1:{#2}(\v@lX,\v@lY)%
    \resetc@ntr@l\et@tfigptBezierDD}\ignorespaces\fi}
\def\figptBezierTD#1:#2:#3[#4,#5,#6,#7]{\ifps@cri{\s@uvc@ntr@l\et@tfigptBezierTD%
    \FigptBezier@#3[#4,#5,#6,#7]\Figp@intregTD#1:{#2}(\v@lX,\v@lY,\v@lZ)%
    \resetc@ntr@l\et@tfigptBezierTD}\ignorespaces\fi}
\def\FigptBezier@#1[#2,#3,#4,#5]{\setc@ntr@l{2}%
    \edef\T@{#1}\v@leur=\p@\advance\v@leur-#1pt\edef\UNmT@{\repdecn@mb{\v@leur}}%
    \figptcopy-4:/#2/\figptcopy-3:/#3/\figptcopy-2:/#4/\figptcopy-1:/#5/%
    \l@mbd@un=-4 \l@mbd@de=-\thr@@\p@rtent=\m@ne\c@lDecast%
    \l@mbd@un=-4 \l@mbd@de=-\thr@@\p@rtent=-\tw@\c@lDecast%
    \l@mbd@un=-4 \l@mbd@de=-\thr@@\p@rtent=-\thr@@\c@lDecast\Figg@tXY{-4}}
\def\c@lDCUnDD#1#2{\Figg@tXY{#1}\v@lX=\UNmT@\v@lX\v@lY=\UNmT@\v@lY%
    \Figg@tXYa{#2}\advance\v@lX\T@\v@lXa\advance\v@lY\T@\v@lYa%
    \Figp@intregDD#1:(\v@lX,\v@lY)}
\def\c@lDCUnTD#1#2{\Figg@tXY{#1}\v@lX=\UNmT@\v@lX\v@lY=\UNmT@\v@lY\v@lZ=\UNmT@\v@lZ%
    \Figg@tXYa{#2}\advance\v@lX\T@\v@lXa\advance\v@lY\T@\v@lYa\advance\v@lZ\T@\v@lZa%
    \Figp@intregTD#1:(\v@lX,\v@lY,\v@lZ)}
\def\c@lDecast{\relax\ifnum\l@mbd@un<\p@rtent\c@lDCUn{\l@mbd@un}{\l@mbd@de}%
    \advance\l@mbd@un\@ne\advance\l@mbd@de\@ne\c@lDecast\fi}
\def\figptmap#1:#2=#3/#4/#5/{\ifps@cri{\s@uvc@ntr@l\et@tfigptmap%
    \setc@ntr@l{2}\figvectP-1[#4,#3]\Figg@tXY{-1}%
    \pr@dMatV/#5/\figpttra#1:{#2}=#4/1,-1/%
    \resetc@ntr@l\et@tfigptmap}\ignorespaces\fi}
\def\pr@dMatVDD/#1,#2;#3,#4/{\v@lXa=#1\v@lX\advance\v@lXa#2\v@lY%
    \v@lYa=#3\v@lX\advance\v@lYa#4\v@lY\Figv@ctCreg-1(\v@lXa,\v@lYa)}
\def\pr@dMatVTD/#1,#2,#3;#4,#5,#6;#7,#8,#9/{%
    \v@lXa=#1\v@lX\advance\v@lXa#2\v@lY\advance\v@lXa#3\v@lZ%
    \v@lYa=#4\v@lX\advance\v@lYa#5\v@lY\advance\v@lYa#6\v@lZ%
    \v@lZa=#7\v@lX\advance\v@lZa#8\v@lY\advance\v@lZa#9\v@lZ%
    \Figv@ctCreg-1(\v@lXa,\v@lYa,\v@lZa)}
\def\figptbaryDD#1:#2[#3;#4]{\ifps@cri{\edef\list@num{#3}\extrairelepremi@r\p@int\de\list@num%
    \s@mme=\z@\@ecfor\c@ef:=#4\do{\advance\s@mme\c@ef}%
    \edef\listec@ef{#4,0}\extrairelepremi@r\c@ef\de\listec@ef%
    \Figg@tXY{\p@int}\divide\v@lX\s@mme\divide\v@lY\s@mme%
    \multiply\v@lX\c@ef\multiply\v@lY\c@ef%
    \@ecfor\p@int:=\list@num\do{\extrairelepremi@r\c@ef\de\listec@ef%
           \Figg@tXYa{\p@int}\divide\v@lXa\s@mme\divide\v@lYa\s@mme%
           \multiply\v@lXa\c@ef\multiply\v@lYa\c@ef%
           \advance\v@lX\v@lXa\advance\v@lY\v@lYa}%
    \Figp@intregDD#1:{#2}(\v@lX,\v@lY)}\ignorespaces\fi}
\def\figptbaryTD#1:#2[#3;#4]{\ifps@cri{\edef\list@num{#3}\extrairelepremi@r\p@int\de\list@num%
    \s@mme=\z@\@ecfor\c@ef:=#4\do{\advance\s@mme\c@ef}%
    \edef\listec@ef{#4,0}\extrairelepremi@r\c@ef\de\listec@ef%
    \Figg@tXY{\p@int}\divide\v@lX\s@mme\divide\v@lY\s@mme\divide\v@lZ\s@mme%
    \multiply\v@lX\c@ef\multiply\v@lY\c@ef\multiply\v@lZ\c@ef%
    \@ecfor\p@int:=\list@num\do{\extrairelepremi@r\c@ef\de\listec@ef%
           \Figg@tXYa{\p@int}\divide\v@lXa\s@mme\divide\v@lYa\s@mme\divide\v@lZa\s@mme%
           \multiply\v@lXa\c@ef\multiply\v@lYa\c@ef\multiply\v@lZa\c@ef%
           \advance\v@lX\v@lXa\advance\v@lY\v@lYa\advance\v@lZ\v@lZa}%
    \Figp@intregTD#1:{#2}(\v@lX,\v@lY,\v@lZ)}\ignorespaces\fi}
\def\figptbaryR#1:#2[#3;#4]{\ifps@cri{%
    \v@leur=\z@\@ecfor\c@ef:=#4\do{\maxim@m{\v@lmax}{\c@ef pt}{-\c@ef pt}%
    \ifdim\v@lmax>\v@leur\v@leur=\v@lmax\fi}%
    \ifdim\v@leur<\p@\f@ctech=\@M\else\ifdim\v@leur<\t@n\p@\f@ctech=\@m\else%
    \ifdim\v@leur<\c@nt\p@\f@ctech=\c@nt\else\ifdim\v@leur<\@m\p@\f@ctech=\t@n\else%
    \f@ctech=\@ne\fi\fi\fi\fi%
    \def\listec@ef{0}%
    \@ecfor\c@ef:=#4\do{\sc@lec@nvRI{\c@ef pt}\edef\listec@ef{\listec@ef,\the\s@mme}}%
    \extrairelepremi@r\c@ef\de\listec@ef\figptbary#1:#2[#3;\listec@ef]}\ignorespaces\fi}
\def\sc@lec@nvRI#1{\v@leur=#1\p@rtentiere{\s@mme}{\v@leur}\advance\v@leur-\s@mme\p@%
    \multiply\v@leur\f@ctech\p@rtentiere{\p@rtent}{\v@leur}%
    \multiply\s@mme\f@ctech\advance\s@mme\p@rtent}
\def\figptcircDD#1:#2:#3;#4(#5){\ifps@cri{\s@uvc@ntr@l\et@tfigptcircDD%
    \c@lptellDD#1:{#2}:#3;#4,#4(#5)\resetc@ntr@l\et@tfigptcircDD}\ignorespaces\fi}
\def\figptcircTD#1:#2:#3,#4,#5;#6(#7){\ifps@cri{\s@uvc@ntr@l\et@tfigptcircTD%
    \setc@ntr@l{2}\c@lExtAxes#3,#4,#5(#6)\figptellP#1:{#2}:#3,-4,-5(#7)%
    \resetc@ntr@l\et@tfigptcircTD}\ignorespaces\fi}
\def\figptcircumcenterDD#1:#2[#3,#4,#5]{\ifps@cri{\s@uvc@ntr@l\et@tfigptcircumcenterDD%
    \setc@ntr@l{2}\figvectNDD-5[#3,#4]\figptbaryDD-3:[#3,#4;1,1]%
                  \figvectNDD-6[#4,#5]\figptbaryDD-4:[#4,#5;1,1]%
    \resetc@ntr@l{2}\inters@cDD#1:{#2}[-3,-5;-4,-6]%
    \resetc@ntr@l\et@tfigptcircumcenterDD}\ignorespaces\fi}
\def\figptcircumcenterTD#1:#2[#3,#4,#5]{\ifps@cri{\s@uvc@ntr@l\et@tfigptcircumcenterTD%
    \setc@ntr@l{2}\figvectNTD-1[#3,#4,#5]%
    \figvectPTD-3[#3,#4]\figvectNVTD-5[-1,-3]\figptbaryTD-3:[#3,#4;1,1]%
    \figvectPTD-4[#4,#5]\figvectNVTD-6[-1,-4]\figptbaryTD-4:[#4,#5;1,1]%
    \resetc@ntr@l{2}\inters@cTD#1:{#2}[-3,-5;-4,-6]%
    \resetc@ntr@l\et@tfigptcircumcenterTD}\ignorespaces\fi}
\def\figptcopyDD#1:#2/#3/{\ifps@cri{\Figg@tXY{#3}%
    \Figp@intregDD#1:{#2}(\v@lX,\v@lY)}\ignorespaces\fi}
\def\figptcopyTD#1:#2/#3/{\ifps@cri{\Figg@tXY{#3}%
    \Figp@intregTD#1:{#2}(\v@lX,\v@lY,\v@lZ)}\ignorespaces\fi}
\def\figptcurvcenterDD#1:#2:#3[#4,#5,#6,#7]{\ifps@cri{\s@uvc@ntr@l\et@tfigptcurvcenterDD%
    \setc@ntr@l{2}\c@lcurvradDD#3[#4,#5,#6,#7]\edef\Sprim@{\repdecn@mb{\result@t}}%
    \figptBezierDD-1::#3[#4,#5,#6,#7]\figpttraDD#1:{#2}=-1/\Sprim@,-5/%
    \resetc@ntr@l\et@tfigptcurvcenterDD}\ignorespaces\fi}
\def\figptcurvcenterTD#1:#2:#3[#4,#5,#6,#7]{\ifps@cri{\s@uvc@ntr@l\et@tfigptcurvcenterTD%
    \setc@ntr@l{2}\figvectDBezierTD -5:1,#3[#4,#5,#6,#7]%
    \figvectDBezierTD -6:2,#3[#4,#5,#6,#7]\vecunit@TD{-5}{-5}%
    \edef\Sprim@{\repdecn@mb{\result@t}}\figvectNVTD-1[-6,-5]%
    \figvectNVTD-5[-5,-1]\c@lproscalTD\v@leur[-6,-5]%
    \invers@{\v@leur}{\v@leur}\v@leur=\Sprim@\v@leur\v@leur=\Sprim@\v@leur%
    \figptBezierTD-1::#3[#4,#5,#6,#7]\edef\Sprim@{\repdecn@mb{\v@leur}}%
    \figpttraTD#1:{#2}=-1/\Sprim@,-5/\resetc@ntr@l\et@tfigptcurvcenterTD}\ignorespaces\fi}
\def\c@lcurvradDD#1[#2,#3,#4,#5]{{\figvectDBezierDD -5:1,#1[#2,#3,#4,#5]%
    \figvectDBezierDD -6:2,#1[#2,#3,#4,#5]\vecunit@DD{-5}{-5}%
    \edef\Sprim@{\repdecn@mb{\result@t}}\figvectNVDD-5[-5]\c@lproscalDD\v@leur[-6,-5]%
    \invers@{\v@leur}{\v@leur}\v@leur=\Sprim@\v@leur\v@leur=\Sprim@\v@leur%
    \global\result@t=\v@leur}}
\def\figptellDD#1:#2:#3;#4,#5(#6,#7){\ifps@cri{\s@uvc@ntr@l\et@tfigptell%
    \c@lptellDD#1::#3;#4,#5(#6)\figptrotDD#1:{#2}=#1/#3,#7/%
    \resetc@ntr@l\et@tfigptell}\ignorespaces\fi}
\def\c@lptellDD#1:#2:#3;#4,#5(#6){\c@ssin{\C@}{\S@}{#6}\v@lmin=\C@ pt\v@lmax=\S@ pt%
    \v@lmin=#4\v@lmin\v@lmax=#5\v@lmax%
    \edef\Xc@mp{\repdecn@mb{\v@lmin}}\edef\Yc@mp{\repdecn@mb{\v@lmax}}%
    \setc@ntr@l{2}\figvectC-1(\Xc@mp,\Yc@mp)\figpttraDD#1:{#2}=#3/1,-1/}
\def\figptellP#1:#2:#3,#4,#5(#6){\ifps@cri{\s@uvc@ntr@l\et@tfigptellP%
    \setc@ntr@l{2}\figvectP-1[#3,#4]\figvectP-2[#3,#5]%
    \v@leur=#6pt\c@lptellP{#3}{-1}{-2}\figptcopy#1:{#2}/-3/%
    \resetc@ntr@l\et@tfigptellP}\ignorespaces\fi}
\def\c@lptellP#1#2#3{\edef\@ngle{\repdecn@mb\v@leur}\c@ssin{\C@}{\S@}{\@ngle}%
    \figpttra-3:=#1/\C@,#2/\figpttra-3:=-3/\S@,#3/}
\def\figptendnormalDD#1:#2:#3,#4[#5,#6]{\ifps@cri{\s@uvc@ntr@l\et@tfigptendnormal%
    \Figg@tXYa{#5}\Figg@tXY{#6}%
    \advance\v@lX-\v@lXa\advance\v@lY-\v@lYa%
    \setc@ntr@l{2}\Figv@ctCreg-1(\v@lX,\v@lY)\vecunit@{-1}{-1}\Figg@tXY{-1}%
    \delt@=#3\unit@\maxim@m{\delt@}{\delt@}{-\delt@}\edef\l@ngueur{\repdecn@mb{\delt@}}%
    \v@lX=\l@ngueur\v@lX\v@lY=\l@ngueur\v@lY%
    \delt@=\p@\advance\delt@-#4pt\edef\l@ngueur{\repdecn@mb{\delt@}}%
    \figptbaryR-1:[#5,#6;#4,\l@ngueur]\Figg@tXYa{-1}%
    \advance\v@lXa\v@lY\advance\v@lYa-\v@lX%
    \setc@ntr@l{1}\Figp@intregDD#1:{#2}(\v@lXa,\v@lYa)\resetc@ntr@l\et@tfigptendnormal}%
    \ignorespaces\fi}
\def\figptexcenter#1:#2[#3,#4,#5]{\ifps@cri{\let@xte={-}%
    \Figptexinsc@nter#1:#2[#3,#4,#5]}\ignorespaces\fi}
\def\figptincenter#1:#2[#3,#4,#5]{\ifps@cri{\let@xte={}%
    \Figptexinsc@nter#1:#2[#3,#4,#5]}\ignorespaces\fi}
\let\figptinscribedcenter=\figptincenter% pour compatibilite avec anciennes versions
\def\Figptexinsc@nter#1:#2[#3,#4,#5]{%
    \figgetdist\LA@[#4,#5]\figgetdist\LB@[#3,#5]\figgetdist\LC@[#3,#4]%
    \figptbaryR#1:{#2}[#3,#4,#5;\the\let@xte\LA@,\LB@,\LC@]}
\def\figptinterlineplaneTD#1:#2[#3,#4;#5,#6]{\ifps@cri{\s@uvc@ntr@l\et@tfigptinterlineplane%
    \setc@ntr@l{2}\figvectPTD-1[#3,#5]\vecunit@TD{-2}{#6}%
    \r@pPSTD\v@leur[-2,-1,#4]\edef\v@lcoef{\repdecn@mb{\v@leur}}%
    \figpttraTD#1:{#2}=#3/\v@lcoef,#4/\resetc@ntr@l\et@tfigptinterlineplane}\ignorespaces\fi}
\def\figptorthocenterDD#1:#2[#3,#4,#5]{\ifps@cri{\s@uvc@ntr@l\et@tfigptorthocenterDD%
    \setc@ntr@l{2}\figvectNDD-3[#3,#4]\figvectNDD-4[#4,#5]%
    \resetc@ntr@l{2}\inters@cDD#1:{#2}[#5,-3;#3,-4]%
    \resetc@ntr@l\et@tfigptorthocenterDD}\ignorespaces\fi}
\def\figptorthocenterTD#1:#2[#3,#4,#5]{\ifps@cri{\s@uvc@ntr@l\et@tfigptorthocenterTD%
    \setc@ntr@l{2}\figvectNTD-1[#3,#4,#5]%
    \figvectPTD-2[#3,#4]\figvectNVTD-3[-1,-2]%
    \figvectPTD-2[#4,#5]\figvectNVTD-4[-1,-2]%
    \resetc@ntr@l{2}\inters@cTD#1:{#2}[#5,-3;#3,-4]%
    \resetc@ntr@l\et@tfigptorthocenterTD}\ignorespaces\fi}
\def\figptorthoprojlineDD#1:#2=#3/#4,#5/{\ifps@cri{\s@uvc@ntr@l\et@tfigptorthoprojlineDD%
    \setc@ntr@l{2}\figvectPDD-3[#4,#5]\figvectNVDD-4[-3]\resetc@ntr@l{2}%
    \inters@cDD#1:{#2}[#3,-4;#4,-3]\resetc@ntr@l\et@tfigptorthoprojlineDD}\ignorespaces\fi}
\def\figptorthoprojlineTD#1:#2=#3/#4,#5/{\ifps@cri{\s@uvc@ntr@l\et@tfigptorthoprojlineTD%
    \setc@ntr@l{2}\figvectPTD-1[#4,#3]\figvectPTD-2[#4,#5]\vecunit@TD{-2}{-2}%
    \c@lproscalTD\v@leur[-1,-2]\edef\v@lcoef{\repdecn@mb{\v@leur}}%
    \figpttraTD#1:{#2}=#4/\v@lcoef,-2/\resetc@ntr@l\et@tfigptorthoprojlineTD}\ignorespaces\fi}
\def\figptorthoprojplaneTD#1:#2=#3/#4,#5/{\ifps@cri{\s@uvc@ntr@l\et@tfigptorthoprojplane%
    \setc@ntr@l{2}\figvectPTD-1[#3,#4]\vecunit@TD{-2}{#5}%
    \c@lproscalTD\v@leur[-1,-2]\edef\v@lcoef{\repdecn@mb{\v@leur}}%
    \figpttraTD#1:{#2}=#3/\v@lcoef,-2/\resetc@ntr@l\et@tfigptorthoprojplane}\ignorespaces\fi}
\def\figpthom#1:#2=#3/#4,#5/{\ifps@cri{\s@uvc@ntr@l\et@tfigpthom%
    \setc@ntr@l{2}\figvectP-1[#4,#3]\figpttra#1:{#2}=#4/#5,-1/%
    \resetc@ntr@l\et@tfigpthom}\ignorespaces\fi}
\def\figptrotDD#1:#2=#3/#4,#5/{\ifps@cri{\s@uvc@ntr@l\et@tfigptrotDD%
    \c@ssin{\C@}{\S@}{#5}\setc@ntr@l{2}\figvectPDD-1[#4,#3]\Figg@tXY{-1}%
    \v@lXa=\C@\v@lX\advance\v@lXa-\S@\v@lY%
    \v@lYa=\S@\v@lX\advance\v@lYa\C@\v@lY%
    \Figv@ctCreg-1(\v@lXa,\v@lYa)\figpttraDD#1:{#2}=#4/1,-1/%
    \resetc@ntr@l\et@tfigptrotDD}\ignorespaces\fi}
\def\figptrotTD#1:#2=#3/#4,#5,#6/{\ifps@cri{\s@uvc@ntr@l\et@tfigptrotTD%
    \c@ssin{\C@}{\S@}{#5}%
    \setc@ntr@l{2}\figptorthoprojplaneTD-3:=#4/#3,#6/\figvectPTD-2[-3,#3]%
    \n@rmeucTD\v@leur{-2}\ifdim\v@leur<\Cepsil@n\Figg@tXYa{#3}\else%
    \edef\v@lcoef{\repdecn@mb{\v@leur}}\figvectNVTD-1[#6,-2]%
    \Figg@tXYa{-1}\v@lXa=\v@lcoef\v@lXa\v@lYa=\v@lcoef\v@lYa\v@lZa=\v@lcoef\v@lZa%
    \v@lXa=\S@\v@lXa\v@lYa=\S@\v@lYa\v@lZa=\S@\v@lZa\Figg@tXY{-2}%
    \advance\v@lXa\C@\v@lX\advance\v@lYa\C@\v@lY\advance\v@lZa\C@\v@lZ%
    \Figg@tXY{-3}\advance\v@lXa\v@lX\advance\v@lYa\v@lY\advance\v@lZa\v@lZ\fi%
    \Figp@intregTD#1:{#2}(\v@lXa,\v@lYa,\v@lZa)\resetc@ntr@l\et@tfigptrotTD}\ignorespaces\fi}
\def\figptsymDD#1:#2=#3/#4,#5/{\ifps@cri{\s@uvc@ntr@l\et@tfigptsymDD%
    \resetc@ntr@l{2}\figptorthoprojlineDD-5:=#3/#4,#5/\figvectPDD-2[#3,-5]%
    \figpttraDD#1:{#2}=#3/2,-2/\resetc@ntr@l\et@tfigptsymDD}\ignorespaces\fi}
\def\figptsymTD#1:#2=#3/#4,#5/{\ifps@cri{\s@uvc@ntr@l\et@tfigptsymTD%
    \resetc@ntr@l{2}\figptorthoprojplaneTD-3:=#3/#4,#5/\figvectPTD-2[#3,-3]%
    \figpttraTD#1:{#2}=#3/2,-2/\resetc@ntr@l\et@tfigptsymTD}\ignorespaces\fi}
\def\figpttraDD#1:#2=#3/#4,#5/{\ifps@cri{\Figg@tXYa{#5}\v@lXa=#4\v@lXa\v@lYa=#4\v@lYa%
    \Figg@tXY{#3}\advance\v@lX\v@lXa\advance\v@lY\v@lYa%
    \Figp@intregDD#1:{#2}(\v@lX,\v@lY)}\ignorespaces\fi}
\def\figpttraTD#1:#2=#3/#4,#5/{\ifps@cri{\Figg@tXYa{#5}\v@lXa=#4\v@lXa\v@lYa=#4\v@lYa%
    \v@lZa=#4\v@lZa\Figg@tXY{#3}\advance\v@lX\v@lXa\advance\v@lY\v@lYa%
    \advance\v@lZ\v@lZa\Figp@intregTD#1:{#2}(\v@lX,\v@lY,\v@lZ)}\ignorespaces\fi}
\def\figpttraCDD#1:#2=#3/#4,#5/{\ifps@cri{\v@lXa=#4\unit@\v@lYa=#5\unit@%
    \Figg@tXY{#3}\advance\v@lX\v@lXa\advance\v@lY\v@lYa%
    \Figp@intregDD#1:{#2}(\v@lX,\v@lY)}\ignorespaces\fi}
\def\figpttraCTD#1:#2=#3/#4,#5,#6/{\ifps@cri{\v@lXa=#4\unit@\v@lYa=#5\unit@\v@lZa=#6\unit@%
    \Figg@tXY{#3}\advance\v@lX\v@lXa\advance\v@lY\v@lYa\advance\v@lZ\v@lZa%
    \Figp@intregTD#1:{#2}(\v@lX,\v@lY,\v@lZ)}\ignorespaces\fi}
\def\figptsaxes#1:#2(#3){\ifps@cri{\an@lys@xes#3,:\ifx\t@xt@\empty%
    \ifTr@isDim\Figpts@xes#1:#2(0,#3,0,#3,0,#3)\else\Figpts@xes#1:#2(0,#3,0,#3)\fi%
    \else\Figpts@xes#1:#2(#3)\fi}\ignorespaces\fi}
\def\Figpts@xesDD#1:#2(#3,#4,#5,#6){%
    \s@mme=#1\figpttraC\the\s@mme:$x$=#2/#4,0/%
    \advance\s@mme\@ne\figpttraC\the\s@mme:$y$=#2/0,#6/}
\def\Figpts@xesTD#1:#2(#3,#4,#5,#6,#7,#8){%
    \s@mme=#1\figpttraC\the\s@mme:$x$=#2/#4,0,0/%
    \advance\s@mme\@ne\figpttraC\the\s@mme:$y$=#2/0,#6,0/%
    \advance\s@mme\@ne\figpttraC\the\s@mme:$z$=#2/0,0,#8/}
\def\figptsmap#1=#2/#3/#4/{\ifps@cri{\s@uvc@ntr@l\et@tfigptsmap%
    \setc@ntr@l{2}\def\list@num{#2}\s@mme=#1%
    \@ecfor\p@int:=\list@num\do{\figvectP-1[#3,\p@int]\Figg@tXY{-1}%
    \pr@dMatV/#4/\figpttra\the\s@mme:=#3/1,-1/\advance\s@mme\@ne}%
    \resetc@ntr@l\et@tfigptsmap}\ignorespaces\fi}
\def\figptscontrolDD#1[#2,#3,#4,#5]{\ifps@cri{\s@uvc@ntr@l\et@tfigptscontrolDD\setc@ntr@l{2}%
    \v@lX=\z@\v@lY=\z@\Figtr@nptDD{-5}{#2}\Figtr@nptDD{2}{#5}%
    \divide\v@lX\@vi\divide\v@lY\@vi%
    \Figtr@nptDD{3}{#3}\Figtr@nptDD{-1.5}{#4}\Figp@intregDD-1:(\v@lX,\v@lY)%
    \v@lX=\z@\v@lY=\z@\Figtr@nptDD{2}{#2}\Figtr@nptDD{-5}{#5}%
    \divide\v@lX\@vi\divide\v@lY\@vi\Figtr@nptDD{-1.5}{#3}\Figtr@nptDD{3}{#4}%
    \s@mme=#1\advance\s@mme\@ne\Figp@intregDD\the\s@mme:(\v@lX,\v@lY)%
    \figptcopyDD#1:/-1/\resetc@ntr@l\et@tfigptscontrolDD}\ignorespaces\fi}
\def\figptscontrolTD#1[#2,#3,#4,#5]{\ifps@cri{\s@uvc@ntr@l\et@tfigptscontrolTD\setc@ntr@l{2}%
    \v@lX=\z@\v@lY=\z@\v@lZ=\z@\Figtr@nptTD{-5}{#2}\Figtr@nptTD{2}{#5}%
    \divide\v@lX\@vi\divide\v@lY\@vi\divide\v@lZ\@vi%
    \Figtr@nptTD{3}{#3}\Figtr@nptTD{-1.5}{#4}\Figp@intregTD-1:(\v@lX,\v@lY,\v@lZ)%
    \v@lX=\z@\v@lY=\z@\v@lZ=\z@\Figtr@nptTD{2}{#2}\Figtr@nptTD{-5}{#5}%
    \divide\v@lX\@vi\divide\v@lY\@vi\divide\v@lZ\@vi\Figtr@nptTD{-1.5}{#3}\Figtr@nptTD{3}{#4}%
    \s@mme=#1\advance\s@mme\@ne\Figp@intregTD\the\s@mme:(\v@lX,\v@lY,\v@lZ)%
    \figptcopyTD#1:/-1/\resetc@ntr@l\et@tfigptscontrolTD}\ignorespaces\fi}
\def\Figtr@nptDD#1#2{\Figg@tXYa{#2}\v@lXa=#1\v@lXa\v@lYa=#1\v@lYa%
    \advance\v@lX\v@lXa\advance\v@lY\v@lYa}
\def\Figtr@nptTD#1#2{\Figg@tXYa{#2}\v@lXa=#1\v@lXa\v@lYa=#1\v@lYa\v@lZa=#1\v@lZa%
    \advance\v@lX\v@lXa\advance\v@lY\v@lYa\advance\v@lZ\v@lZa}
\def\figptscontrolcurve#1,#2[#3]{\ifps@cri{\s@uvc@ntr@l\et@tfigptscontrolcurve%
    \def\list@num{#3}\extrairelepremi@r\Ak@\de\list@num%
    \extrairelepremi@r\Ai@\de\list@num\extrairelepremi@r\Aj@\de\list@num%
    \s@mme=#1\figptcopy\the\s@mme:/\Ai@/%
    \setc@ntr@l{2}\figvectP -1[\Ak@,\Aj@]%
    \@ecfor\Ak@:=\list@num\do{\advance\s@mme\@ne\figpttra\the\s@mme:=\Ai@/\curv@roundness,-1/%
       \figvectP -1[\Ai@,\Ak@]\advance\s@mme\@ne\figpttra\the\s@mme:=\Aj@/-\curv@roundness,-1/%
       \advance\s@mme\@ne\figptcopy\the\s@mme:/\Aj@/%
       \edef\Ai@{\Aj@}\edef\Aj@{\Ak@}}\advance\s@mme-#1\divide\s@mme\thr@@%
       \xdef#2{\the\s@mme}%
    \resetc@ntr@l\et@tfigptscontrolcurve}\ignorespaces\fi}
\def\figptsintercircDD#1[#2,#3;#4,#5]{\ifps@cri{\s@uvc@ntr@l\et@tfigptsintercircDD%
    \setc@ntr@l{2}\let\c@lNVintc=\c@lNVintcDD\Figptsintercirc@#1[#2,#3;#4,#5]%    
    \resetc@ntr@l\et@tfigptsintercircDD}\ignorespaces\fi}
\def\figptsintercircTD#1[#2,#3;#4,#5;#6]{\ifps@cri{\s@uvc@ntr@l\et@tfigptsintercircTD%
    \setc@ntr@l{2}\let\c@lNVintc=\c@lNVintcTD\vecunitC@TD[#2,#6]%
    \Figv@ctCreg-3(\v@lX,\v@lY,\v@lZ)\Figptsintercirc@#1[#2,#3;#4,#5]%
    \resetc@ntr@l\et@tfigptsintercircTD}\ignorespaces\fi}
\def\Figptsintercirc@#1[#2,#3;#4,#5]{\figvectP-1[#2,#4]%
    \vecunit@{-1}{-1}\delt@=\result@t\f@ctech=\result@tent%
    \s@mme=#1\advance\s@mme\@ne\figptcopy#1:/#2/\figptcopy\the\s@mme:/#4/%
    \ifdim\delt@=\z@\else%
    \v@lmin=#3\unit@\v@lmax=#5\unit@\v@leur=\v@lmin\advance\v@leur\v@lmax%
    \ifdim\v@leur>\delt@%
    \v@leur=\v@lmin\advance\v@leur-\v@lmax\maxim@m{\v@leur}{\v@leur}{-\v@leur}%
    \ifdim\v@leur<\delt@%
    \divide\v@lmin\f@ctech\divide\v@lmax\f@ctech\divide\delt@\f@ctech%
    \v@lmin=\repdecn@mb{\v@lmin}\v@lmin\v@lmax=\repdecn@mb{\v@lmax}\v@lmax%
    \invers@{\v@leur}{\delt@}\advance\v@lmax-\v@lmin%
    \v@lmax=-\repdecn@mb{\v@leur}\v@lmax\advance\delt@\v@lmax\delt@=.5\delt@%
    \v@lmax=\delt@\multiply\v@lmax\f@ctech%
    \edef\t@ille{\repdecn@mb{\v@lmax}}\figpttra-2:=#2/\t@ille,-1/%
    \delt@=\repdecn@mb{\delt@}\delt@\advance\v@lmin-\delt@%
    \sqrt@{\v@leur}{\v@lmin}\multiply\v@leur\f@ctech\edef\t@ille{\repdecn@mb{\v@leur}}%
    \c@lNVintc\figpttra#1:=-2/-\t@ille,-1/\figpttra\the\s@mme:=-2/\t@ille,-1/\fi\fi\fi}
\def\c@lNVintcDD{\Figg@tXY{-1}\Figv@ctCreg-1(-\v@lY,\v@lX)} % <=> \figvectNVDD-1[-1]
\def\c@lNVintcTD{{\Figg@tXY{-3}\v@lmin=\v@lX\v@lmax=\v@lY\v@leur=\v@lZ%
    \Figg@tXY{-1}\c@lprovec{-3}\vecunit@{-3}{-3}% <=> \figvectNVTD-3[-1,-3]\vecunit@{-3}{-3}
    \Figg@tXY{-1}\v@lmin=\v@lX\v@lmax=\v@lY%
    \v@leur=\v@lZ\Figg@tXY{-3}\c@lprovec{-1}}} % <=> \figvectNVTD-1[-3,-1]
\def\figptsinterlinellDD#1[#2,#3,#4,#5;#6,#7]{\ifps@cri{\s@uvc@ntr@l\et@tfigptsinterlinellDD%
    \figptcopy#1:/#6/\s@mme=#1\advance\s@mme\@ne\figptcopy\the\s@mme:/#7/%
    \v@lmin=#3\unit@\v@lmax=#4\unit@% a, b
    \setc@ntr@l{2}\figptbaryDD-4:[#6,#7;1,1]\figptsrotDD-3=-4,#7/#2,-#5/% D et rotation
    \Figg@tXY{-3}\Figg@tXYa{#2}\advance\v@lX-\v@lXa\advance\v@lY-\v@lYa% alpha, beta
    \figvectP-1[-3,-2]\Figg@tXYa{-1}\figvectP-3[-4,#7]\Figptsint@rLE{#1}% u1, u2
    \resetc@ntr@l\et@tfigptsinterlinellDD}\ignorespaces\fi}
\def\figptsinterlinellP#1[#2,#3,#4;#5,#6]{\ifps@cri{\s@uvc@ntr@l\et@tfigptsinterlinellP%
    \figptcopy#1:/#5/\s@mme=#1\advance\s@mme\@ne\figptcopy\the\s@mme:/#6/\setc@ntr@l{2}%
    \figvectP-1[#2,#3]\vecunit@{-1}{-1}\v@lmin=\result@t% a
    \figvectP-2[#2,#4]\vecunit@{-2}{-2}\v@lmax=\result@t% b
    \figptbary-4:[#5,#6;1,1]% D
    \figvectP-3[#2,-4]\c@lproscal\v@lX[-3,-1]\c@lproscal\v@lY[-3,-2]% alpha, beta
    \figvectP-3[-4,#6]\c@lproscal\v@lXa[-3,-1]\c@lproscal\v@lYa[-3,-2]% u1, u2
    \Figptsint@rLE{#1}\resetc@ntr@l\et@tfigptsinterlinellP}\ignorespaces\fi}
\def\Figptsint@rLE#1{%
    \getredf@ctDD\f@ctech(\v@lmin,\v@lmax)%
    \getredf@ctDD\p@rtent(\v@lX,\v@lY)\ifnum\p@rtent>\f@ctech\f@ctech=\p@rtent\fi%
    \getredf@ctDD\p@rtent(\v@lXa,\v@lYa)\ifnum\p@rtent>\f@ctech\f@ctech=\p@rtent\fi%
    \divide\v@lmin\f@ctech\divide\v@lmax\f@ctech\divide\v@lX\f@ctech\divide\v@lY\f@ctech%
    \divide\v@lXa\f@ctech\divide\v@lYa\f@ctech%
    \c@rre=\repdecn@mb\v@lXa\v@lmax\mili@u=\repdecn@mb\v@lYa\v@lmin%
    \getredf@ctDD\f@ctech(\c@rre,\mili@u)%
    \c@rre=\repdecn@mb\v@lX\v@lmax\mili@u=\repdecn@mb\v@lY\v@lmin%
    \getredf@ctDD\p@rtent(\c@rre,\mili@u)\ifnum\p@rtent>\f@ctech\f@ctech=\p@rtent\fi%
    \divide\v@lmin\f@ctech\divide\v@lmax\f@ctech\divide\v@lX\f@ctech\divide\v@lY\f@ctech%
    \divide\v@lXa\f@ctech\divide\v@lYa\f@ctech%
    \v@lmin=\repdecn@mb{\v@lmin}\v@lmin\v@lmax=\repdecn@mb{\v@lmax}\v@lmax%
    \edef\G@xde{\repdecn@mb\v@lmin}\edef\P@xde{\repdecn@mb\v@lmax}%
    \c@rre=-\v@lmax\v@leur=\repdecn@mb\v@lY\v@lY\advance\c@rre\v@leur\c@rre=\G@xde\c@rre%
    \v@leur=\repdecn@mb\v@lX\v@lX\v@leur=\P@xde\v@leur\advance\c@rre\v@leur% C
    \v@lmin=\repdecn@mb\v@lYa\v@lmin\v@lmax=\repdecn@mb\v@lXa\v@lmax%
    \mili@u=\repdecn@mb\v@lX\v@lmax\advance\mili@u\repdecn@mb\v@lY\v@lmin% B
    \v@lmax=\repdecn@mb\v@lXa\v@lmax\advance\v@lmax\repdecn@mb\v@lYa\v@lmin% A
    \ifdim\v@lmax>\epsil@n%
    \maxim@m{\v@leur}{\c@rre}{-\c@rre}\maxim@m{\v@lmin}{\mili@u}{-\mili@u}%
    \maxim@m{\v@leur}{\v@leur}{\v@lmin}\maxim@m{\v@lmin}{\v@lmax}{-\v@lmax}%
    \maxim@m{\v@leur}{\v@leur}{\v@lmin}\p@rtentiere{\p@rtent}{\v@leur}\advance\p@rtent\@ne%
    \divide\c@rre\p@rtent\divide\mili@u\p@rtent\divide\v@lmax\p@rtent%
    \delt@=\repdecn@mb{\mili@u}\mili@u\v@leur=\repdecn@mb{\v@lmax}\c@rre%
    \advance\delt@-\v@leur\ifdim\delt@<\z@\else\sqrt@\delt@\delt@%
    \invers@\v@lmax\v@lmax\edef\Uns@rAp{\repdecn@mb\v@lmax}%
    \v@leur=-\mili@u\advance\v@leur-\delt@\v@leur=\Uns@rAp\v@leur%
    \edef\t@ille{\repdecn@mb\v@leur}\figpttra#1:=-4/\t@ille,-3/\s@mme=#1\advance\s@mme\@ne%
    \v@leur=-\mili@u\advance\v@leur\delt@\v@leur=\Uns@rAp\v@leur%
    \edef\t@ille{\repdecn@mb\v@leur}\figpttra\the\s@mme:=-4/\t@ille,-3/\fi\fi}
\def\figptsorthoprojlineDD#1=#2/#3,#4/{\ifps@cri{\s@uvc@ntr@l\et@tfigptsorthoprojlineDD%
    \setc@ntr@l{2}\figvectPDD-3[#3,#4]\figvectNVDD-4[-3]\resetc@ntr@l{2}%
    \def\list@num{#2}\s@mme=#1\@ecfor\p@int:=\list@num\do{%
    \inters@cDD\the\s@mme:[\p@int,-4;#3,-3]\advance\s@mme\@ne}%
    \resetc@ntr@l\et@tfigptsorthoprojlineDD}\ignorespaces\fi}
\def\figptsorthoprojlineTD#1=#2/#3,#4/{\ifps@cri{\s@uvc@ntr@l\et@tfigptsorthoprojlineTD%
    \setc@ntr@l{2}\figvectPTD-2[#3,#4]\vecunit@TD{-2}{-2}%
    \def\list@num{#2}\s@mme=#1\@ecfor\p@int:=\list@num\do{%
    \figvectPTD-1[#3,\p@int]\c@lproscalTD\v@leur[-1,-2]%
    \edef\v@lcoef{\repdecn@mb{\v@leur}}\figpttraTD\the\s@mme:=#3/\v@lcoef,-2/%
    \advance\s@mme\@ne}\resetc@ntr@l\et@tfigptsorthoprojlineTD}\ignorespaces\fi}
\def\figptsorthoprojplaneTD#1=#2/#3,#4/{\ifps@cri{\s@uvc@ntr@l\et@tfigptsorthoprojplane%
    \setc@ntr@l{2}\vecunit@TD{-2}{#4}%
    \def\list@num{#2}\s@mme=#1\@ecfor\p@int:=\list@num\do{\figvectPTD-1[\p@int,#3]%
    \c@lproscalTD\v@leur[-1,-2]\edef\v@lcoef{\repdecn@mb{\v@leur}}%
    \figpttraTD\the\s@mme:=\p@int/\v@lcoef,-2/\advance\s@mme\@ne}%
    \resetc@ntr@l\et@tfigptsorthoprojplane}\ignorespaces\fi}
\def\figptshom#1=#2/#3,#4/{\ifps@cri{\s@uvc@ntr@l\et@tfigptshom%
    \setc@ntr@l{2}\def\list@num{#2}\s@mme=#1%
    \@ecfor\p@int:=\list@num\do{\figvectP-1[#3,\p@int]%
    \figpttra\the\s@mme:=#3/#4,-1/\advance\s@mme\@ne}%
    \resetc@ntr@l\et@tfigptshom}\ignorespaces\fi}
\def\figptsrotDD#1=#2/#3,#4/{\ifps@cri{\s@uvc@ntr@l\et@tfigptsrotDD%
    \c@ssin{\C@}{\S@}{#4}\setc@ntr@l{2}\def\list@num{#2}\s@mme=#1%
    \@ecfor\p@int:=\list@num\do{\figvectPDD-1[#3,\p@int]\Figg@tXY{-1}%
    \v@lXa=\C@\v@lX\advance\v@lXa-\S@\v@lY%
    \v@lYa=\S@\v@lX\advance\v@lYa\C@\v@lY%
    \Figv@ctCreg-1(\v@lXa,\v@lYa)\figpttraDD\the\s@mme:=#3/1,-1/\advance\s@mme\@ne}%
    \resetc@ntr@l\et@tfigptsrotDD}\ignorespaces\fi}
\def\figptsrotTD#1=#2/#3,#4,#5/{\ifps@cri{\s@uvc@ntr@l\et@tfigptsrotTD%
    \c@ssin{\C@}{\S@}{#4}%
    \setc@ntr@l{2}\def\list@num{#2}\s@mme=#1%
    \@ecfor\p@int:=\list@num\do{\figptorthoprojplaneTD-3:=#3/\p@int,#5/%
    \figvectPTD-2[-3,\p@int]%
    \figvectNVTD-1[#5,-2]\n@rmeucTD\v@leur{-2}\edef\v@lcoef{\repdecn@mb{\v@leur}}%
    \Figg@tXYa{-1}\v@lXa=\v@lcoef\v@lXa\v@lYa=\v@lcoef\v@lYa\v@lZa=\v@lcoef\v@lZa%
    \v@lXa=\S@\v@lXa\v@lYa=\S@\v@lYa\v@lZa=\S@\v@lZa\Figg@tXY{-2}%
    \advance\v@lXa\C@\v@lX\advance\v@lYa\C@\v@lY\advance\v@lZa\C@\v@lZ%
    \Figg@tXY{-3}\advance\v@lXa\v@lX\advance\v@lYa\v@lY\advance\v@lZa\v@lZ%
    \Figp@intregTD\the\s@mme:(\v@lXa,\v@lYa,\v@lZa)\advance\s@mme\@ne}%
    \resetc@ntr@l\et@tfigptsrotTD}\ignorespaces\fi}
\def\figptssymDD#1=#2/#3,#4/{\ifps@cri{\s@uvc@ntr@l\et@tfigptssymDD%
    \setc@ntr@l{2}\figvectPDD-3[#3,#4]\Figg@tXY{-3}\Figv@ctCreg-4(-\v@lY,\v@lX)%
    \resetc@ntr@l{2}\def\list@num{#2}\s@mme=#1%
    \@ecfor\p@int:=\list@num\do{\inters@cDD-5:[#3,-3;\p@int,-4]\figvectPDD-2[\p@int,-5]%
    \figpttraDD\the\s@mme:=\p@int/2,-2/\advance\s@mme\@ne}%
    \resetc@ntr@l\et@tfigptssymDD}\ignorespaces\fi}
\def\figptssymTD#1=#2/#3,#4/{\ifps@cri{\s@uvc@ntr@l\et@tfigptssymTD%
    \setc@ntr@l{2}\vecunit@TD{-2}{#4}\def\list@num{#2}\s@mme=#1%
    \@ecfor\p@int:=\list@num\do{\figvectPTD-1[\p@int,#3]%
    \c@lproscalTD\v@leur[-1,-2]\v@leur=2\v@leur\edef\v@lcoef{\repdecn@mb{\v@leur}}%
    \figpttraTD\the\s@mme:=\p@int/\v@lcoef,-2/\advance\s@mme\@ne}%
    \resetc@ntr@l\et@tfigptssymTD}\ignorespaces\fi}
\def\figptstraDD#1=#2/#3,#4/{\ifps@cri{\Figg@tXYa{#4}\v@lXa=#3\v@lXa\v@lYa=#3\v@lYa%
    \def\list@num{#2}\s@mme=#1\@ecfor\p@int:=\list@num\do{\Figg@tXY{\p@int}%
    \advance\v@lX\v@lXa\advance\v@lY\v@lYa%
    \Figp@intregDD\the\s@mme:(\v@lX,\v@lY)\advance\s@mme\@ne}}\ignorespaces\fi}
\def\figptstraTD#1=#2/#3,#4/{\ifps@cri{\Figg@tXYa{#4}\v@lXa=#3\v@lXa\v@lYa=#3\v@lYa%
    \v@lZa=#3\v@lZa\def\list@num{#2}\s@mme=#1\@ecfor\p@int:=\list@num\do{\Figg@tXY{\p@int}%
    \advance\v@lX\v@lXa\advance\v@lY\v@lYa\advance\v@lZ\v@lZa%
    \Figp@intregTD\the\s@mme:(\v@lX,\v@lY,\v@lZ)\advance\s@mme\@ne}}\ignorespaces\fi}
\def\figptvisilimSLTD#1:#2[#3,#4;#5,#6]{\ifps@cri{\s@uvc@ntr@l\et@tfigptvisilimSLTD%
    \setc@ntr@l{2}\figvectP-1[#3,#4]\n@rminf{\delt@}{-1}%
    \ifcase\curr@ntproj\v@lX=\cxa@\p@\v@lY=-\p@\v@lZ=\cxb@\p@% Proj cav
    \Figv@ctCreg-2(\v@lX,\v@lY,\v@lZ)\figvectP-3[#5,#6]\figvectNV-1[-2,-3]%
    \or\figvectP-1[#5,#6]\vecunitCV@TD{-1}\v@lmin=\v@lX\v@lmax=\v@lY% Proj ortho
    \v@leur=\v@lZ\v@lX=\cza@\p@\v@lY=\czb@\p@\v@lZ=\czc@\p@\c@lprovec{-1}%
    \or\c@ley@pt{-2}\figvectN-1[#5,#6,-2]\fi% Proj rea
    \edef\Ai@{#3}\edef\Aj@{#4}\figvectP-2[#5,\Ai@]\c@lproscal\v@leur[-1,-2]%
    \ifdim\v@leur>\z@\p@rtent=\@ne\else\p@rtent=\m@ne\fi%
    \figvectP-2[#5,\Aj@]\c@lproscal\v@leur[-1,-2]%
    \ifdim\p@rtent\v@leur>\z@\figptcopy#1:#2/#3/%
    \message{*** \BS@ figptvisilimSL: points are on the same side.}\else%
    \figptcopy-3:/#3/\figptcopy-4:/#4/%
    \loop\figptbary-5:[-3,-4;1,1]\figvectP-2[#5,-5]\c@lproscal\v@leur[-1,-2]%
    \ifdim\p@rtent\v@leur>\z@\figptcopy-3:/-5/\else\figptcopy-4:/-5/\fi%
    \divide\delt@\tw@\ifdim\delt@>\epsil@n\repeat%
    \figptbary#1:#2[-3,-4;1,1]\fi\resetc@ntr@l\et@tfigptvisilimSLTD}\ignorespaces\fi}
\def\c@ley@pt#1{\t@stp@r\ifitis@K\v@lX=\cza@\p@\v@lY=\czb@\p@\v@lZ=\czc@\p@%
    \Figv@ctCreg-1(\v@lX,\v@lY,\v@lZ)\Figp@intreg-2:(\wd\Bt@rget,\ht\Bt@rget,\dp\Bt@rget)%
    \figpttra#1:=-2/-\disob@intern,-1/\else\end\fi}
\def\t@stp@r{\itis@Ktrue\ifnewt@rgetpt\else\itis@Kfalse%
    \message{*** \BS@ figptvisilimXX: target point undefined.}\fi\ifnewdis@b\else%
    \itis@Kfalse\message{*** \BS@ figptvisilimXX: observation distance undefined.}\fi%
    \ifitis@K\else\message{*** This macro must be called after \BS@ psbeginfig or after
    having set the missing parameter(s) with \BS@ figset proj()}\fi}
\def\figscan#1(#2,#3){{\s@uvc@ntr@l\et@tfigscan\@psfgetbb{#1}\if@psfbbfound\else%
    \def\@psfllx{0}\def\@psflly{20}\def\@psfurx{540}\def\@psfury{640}\fi%
    \unit@=\@ne bp\setc@ntr@l{2}\figsetmark{}%
    \def\minst@p{20pt}%
    \v@lX=\@psfllx\p@\v@lX=\Sc@leFact\v@lX\r@undint\v@lX\v@lX%
    \v@lY=\@psflly\p@\v@lY=\Sc@leFact\v@lY\ifdim\v@lY>\z@\r@undint\v@lY\v@lY\fi%
    \delt@=\@psfury\p@\delt@=\Sc@leFact\delt@%
    \advance\delt@-\v@lY\v@lXa=\@psfurx\p@\v@lXa=\Sc@leFact\v@lXa\v@leur=\minst@p%
    \edef\valv@lY{\repdecn@mb{\v@lY}}\edef\LgTr@it{\the\delt@}%
    \loop\ifdim\v@lX<\v@lXa\edef\valv@lX{\repdecn@mb{\v@lX}}%
    \figptDD -1:(\valv@lX,\valv@lY)\figwriten -1:\hbox{\vrule height\LgTr@it}(0)%
    \ifdim\v@leur<\minst@p\else\figsetmark{\raise-8bp\hbox{$\scriptscriptstyle\triangle$}}%
    \figwrites -1:\@ffichnb{0}{\valv@lX}(6)\v@leur=\z@\figsetmark{}\fi%
    \advance\v@leur#2pt\advance\v@lX#2pt\repeat%
    \def\minst@p{10pt}%
    \v@lX=\@psfllx\p@\v@lX=\Sc@leFact\v@lX\ifdim\v@lX>\z@\r@undint\v@lX\v@lX\fi%
    \v@lY=\@psflly\p@\v@lY=\Sc@leFact\v@lY\r@undint\v@lY\v@lY%
    \delt@=\@psfurx\p@\delt@=\Sc@leFact\delt@%
    \advance\delt@-\v@lX\v@lYa=\@psfury\p@\v@lYa=\Sc@leFact\v@lYa\v@leur=\minst@p%
    \edef\valv@lX{\repdecn@mb{\v@lX}}\edef\LgTr@it{\the\delt@}%
    \loop\ifdim\v@lY<\v@lYa\edef\valv@lY{\repdecn@mb{\v@lY}}%
    \figptDD -1:(\valv@lX,\valv@lY)\figwritee -1:\vbox{\hrule width\LgTr@it}(0)%
    \ifdim\v@leur<\minst@p\else\figsetmark{$\triangleright$\kern4bp}%
    \figwritew -1:\@ffichnb{0}{\valv@lY}(6)\v@leur=\z@\figsetmark{}\fi%
    \advance\v@leur#3pt\advance\v@lY#3pt\repeat%
    \resetc@ntr@l\et@tfigscan}\ignorespaces}
\def\figshowpts[#1,#2]{{\figsetmark{$\bullet$}\figsetptname{\bf ##1}%
    \p@rtent=#2\relax\ifnum\p@rtent<\z@\p@rtent=\z@\fi%
    \s@mme=#1\relax\ifnum\s@mme<\z@\s@mme=\z@\fi%
    \loop\ifnum\s@mme<\p@rtent\pt@rvect{\s@mme}%
    \ifitis@K\figwriten{\the\s@mme}:(4pt)\fi\advance\s@mme\@ne\repeat%
    \pt@rvect{\s@mme}\ifitis@K\figwriten{\the\s@mme}:(4pt)\fi}\ignorespaces}
\def\pt@rvect#1{\set@bjc@de{#1}%
    \expandafter\expandafter\expandafter\inqpt@rvec\csname\objc@de\endcsname:}
\def\inqpt@rvec#1#2:{\if#1\C@dCl@spt\itis@Ktrue\else\itis@Kfalse\fi}
\def\figshowsettings{{%
    \immediate\write16{====================================================================}%
    \immediate\write16{ Current settings about:}%
    \immediate\write16{ --- GENERAL ---}%
    \immediate\write16{Scale factor and Unit = \unit@util\space (\the\unit@)
     \space -> \BS@ figinit{ScaleFactorUnit}}%
    \immediate\write16{Update mode = \ifpstestm@de yes\else no\fi
     \space-> \BS@ psset(update=yes/no) or \BS@ pssetdefault(update=yes/no)}%
    \immediate\write16{ --- PRINTING ---}%
    \immediate\write16{Implicit point name = \ptn@me{i} \space-> \BS@ figsetptname{Name}}%
    \immediate\write16{Point marker = \the\c@nsymb \space -> \BS@ figsetmark{Mark}}%
    \immediate\write16{Print rounded coordinates = \ifr@undcoord yes\else no\fi
     \space-> \BS@ figsetroundcoord{yes/no}}%
    \immediate\write16{ --- GRAPHICAL (general) ---}%
    \immediate\write16{First-level (or primary) settings:}%
    \immediate\write16{ Color = \curr@ntcolor \space-> \BS@ psset(color=ColorDefinition)}%
    \immediate\write16{ Filling mode = \iffillm@de yes\else no\fi
     \space-> \BS@ psset(fillmode=yes/no)}%
    \immediate\write16{ Line join = \curr@ntjoin \space-> \BS@ psset(join=miter/round/bevel)}%
    \immediate\write16{ Line style = \curr@ntdash \space-> \BS@ psset(dash=Index/Pattern)}%
    \immediate\write16{ Line width = \curr@ntwidth
     \space-> \BS@ psset(width=real in PostScript units)}%
    \immediate\write16{Second-level (or secondary) settings:}%
    \immediate\write16{ Color = \sec@ndcolor \space-> \BS@ psset second(color=ColorDefinition)}%
    \immediate\write16{ Line style = \curr@ntseconddash
     \space-> \BS@ psset second(dash=Index/Pattern)}%
    \immediate\write16{ Line width = \curr@ntsecondwidth
     \space-> \BS@ psset second(width=real in PostScript units)}%
    \immediate\write16{Third-level (or ternary) settings:}%
    \immediate\write16{ Color = \th@rdcolor \space-> \BS@ psset third(color=ColorDefinition)}%
    \immediate\write16{ --- GRAPHICAL (specific) ---}%
    \immediate\write16{Arrow-head:}%
    \immediate\write16{ (half-)Angle = \@rrowheadangle
     \space-> \BS@ psset arrowhead(angle=real in degrees)}%
    \immediate\write16{ Filling mode = \if@rrowhfill yes\else no\fi
     \space-> \BS@ psset arrowhead(fillmode=yes/no)}%
    \immediate\write16{ "Outside" = \if@rrowhout yes\else no\fi
     \space-> \BS@ psset arrowhead(out=yes/no)}%
    \immediate\write16{ Length = \@rrowheadlength
     \if@rrowratio\space(not active)\else\space(active)\fi
     \space-> \BS@ psset arrowhead(length=real in user coord.)}%
    \immediate\write16{ Ratio = \@rrowheadratio
     \if@rrowratio\space(active)\else\space(not active)\fi
     \space-> \BS@ psset arrowhead(ratio=real in [0,1])}%
    \immediate\write16{Curve: Roundness = \curv@roundness
     \space-> \BS@ psset curve(roundness=real in [0,0.5])}%
    \immediate\write16{Mesh: Diagonal = \the\c@ntrolmesh
     \space-> \BS@ psset mesh(diag=integer in {-1,0,1})}%
    \immediate\write16{Flow chart:}%
    \immediate\write16{ Arrow position = \@rrowp@s
     \space-> \BS@ psset flowchart(arrowposition=real in [0,1])}%
    \immediate\write16{ Arrow reference point = \ifcase\@rrowr@fpt start\else end\fi
     \space-> \BS@ psset flowchart(arrowrefpt = start/end)}%     
    \immediate\write16{ Line type = \ifcase\fclin@typ@ curve\else polygon\fi
     \space-> \BS@ psset flowchart(line=polygon/curve)}%
    \immediate\write16{ Padding = (\Xp@dd, \Yp@dd)
     \space-> \BS@ psset flowchart(padding = real in user coord.)}%
    \immediate\write16{\space\space\space\space(or
     \BS@ psset flowchart(xpadding=real, ypadding=real) )}%
    \immediate\write16{ Radius = \fclin@r@d
     \space-> \BS@ psset flowchart(radius=positive real in user coord.)}%
    \immediate\write16{ Shape = \fcsh@pe
     \space-> \BS@ psset flowchart(shape = rectangle, ellipse or lozenge)}%
    \immediate\write16{ Thickness = \thickn@ss
     \space-> \BS@ psset flowchart(thickness = real in user coord.)}%
    \ifTr@isDim%
    \immediate\write16{ --- 3D to 2D PROJECTION ---}%
    \immediate\write16{Projection : \typ@proj \space-> \BS@ figinit{ScaleFactorUnit, ProjType}}%
    \immediate\write16{Longitude (psi) = \v@lPsi \space-> \BS@ figset proj(psi=real in degrees)}%
    \ifcase\curr@ntproj\immediate\write16{Depth coeff. (Lambda)
     \space = \v@lTheta \space-> \BS@ figset proj(lambda=real in [0,1])}%
    \else\immediate\write16{Latitude (theta)
     \space = \v@lTheta \space-> \BS@ figset proj(theta=real in degrees)}%
    \fi%
    \ifnum\curr@ntproj=\tw@%
    \immediate\write16{Observation distance = \disob@unit
     \space-> \BS@ figset proj(dist=real in user coord.)}%
    \immediate\write16{Target point = \t@rgetpt \space-> \BS@ figset proj(targetpt=pt number)}%
     \v@lX=\ptT@unit@\wd\Bt@rget\v@lY=\ptT@unit@\ht\Bt@rget\v@lZ=\ptT@unit@\dp\Bt@rget%
    \immediate\write16{ Its coordinates are
     (\repdecn@mb{\v@lX}, \repdecn@mb{\v@lY}, \repdecn@mb{\v@lZ})}%
    \fi%
    \fi%
    \immediate\write16{====================================================================}%
    \ignorespaces}}
{\catcode`\/=0 \catcode`/\=12 /gdef/BS@{\}}
\newif\ifitis@vect@r
\def\figvectC#1(#2,#3){{\itis@vect@rtrue\figpt#1:(#2,#3)}\ignorespaces}
\def\Figv@ctCreg#1(#2,#3){{\itis@vect@rtrue\Figp@intreg#1:(#2,#3)}\ignorespaces}
\def\figvectDBezierDD#1:#2,#3[#4,#5,#6,#7]{\ifps@cri{\s@uvc@ntr@l\et@tfigvectDBezierDD%
    \FigvectDBezier@#2,#3[#4,#5,#6,#7]\v@lX=\c@ef\v@lX\v@lY=\c@ef\v@lY%
    \Figv@ctCreg#1(\v@lX,\v@lY)\resetc@ntr@l\et@tfigvectDBezierDD}\ignorespaces\fi}
\def\figvectDBezierTD#1:#2,#3[#4,#5,#6,#7]{\ifps@cri{\s@uvc@ntr@l\et@tfigvectDBezierTD%
    \FigvectDBezier@#2,#3[#4,#5,#6,#7]\v@lX=\c@ef\v@lX\v@lY=\c@ef\v@lY\v@lZ=\c@ef\v@lZ%
    \Figv@ctCreg#1(\v@lX,\v@lY,\v@lZ)\resetc@ntr@l\et@tfigvectDBezierTD}\ignorespaces\fi}
\def\FigvectDBezier@#1,#2[#3,#4,#5,#6]{\setc@ntr@l{2}%
    \edef\T@{#2}\v@leur=\p@\advance\v@leur-#2pt\edef\UNmT@{\repdecn@mb{\v@leur}}%
    \ifnum#1=\tw@\def\c@ef{6}\else\def\c@ef{3}\fi%
    \figptcopy-4:/#3/\figptcopy-3:/#4/\figptcopy-2:/#5/\figptcopy-1:/#6/%
    \l@mbd@un=-4 \l@mbd@de=-\thr@@\p@rtent=\m@ne\c@lDecast%
    \ifnum#1=\tw@\c@lDCDeux{-4}{-3}\c@lDCDeux{-3}{-2}\c@lDCDeux{-4}{-3}\else%
    \l@mbd@un=-4 \l@mbd@de=-\thr@@\p@rtent=-\tw@\c@lDecast%
    \c@lDCDeux{-4}{-3}\fi\Figg@tXY{-4}}
\def\c@lDCDeuxDD#1#2{\Figg@tXY{#2}\Figg@tXYa{#1}%
    \advance\v@lX-\v@lXa\advance\v@lY-\v@lYa\Figp@intregDD#1:(\v@lX,\v@lY)}
\def\c@lDCDeuxTD#1#2{\Figg@tXY{#2}\Figg@tXYa{#1}\advance\v@lX-\v@lXa%
    \advance\v@lY-\v@lYa\advance\v@lZ-\v@lZa\Figp@intregTD#1:(\v@lX,\v@lY,\v@lZ)}
\def\figvectNDD#1[#2,#3]{\ifps@cri{\Figg@tXYa{#2}\Figg@tXY{#3}%
    \advance\v@lX-\v@lXa\advance\v@lY-\v@lYa%
    \Figv@ctCreg#1(-\v@lY,\v@lX)}\ignorespaces\fi}
\def\figvectNTD#1[#2,#3,#4]{\ifps@cri{\vecunitC@TD[#2,#4]\v@lmin=\v@lX\v@lmax=\v@lY%
    \v@leur=\v@lZ\vecunitC@TD[#2,#3]\c@lprovec{#1}}\ignorespaces\fi}
\def\figvectNVDD#1[#2]{\ifps@cri{\Figg@tXY{#2}\Figv@ctCreg#1(-\v@lY,\v@lX)}\ignorespaces\fi}
\def\figvectNVTD#1[#2,#3]{\ifps@cri{\vecunitCV@TD{#3}\v@lmin=\v@lX\v@lmax=\v@lY%
    \v@leur=\v@lZ\vecunitCV@TD{#2}\c@lprovec{#1}}\ignorespaces\fi}
\def\figvectPDD#1[#2,#3]{\ifps@cri{\Figg@tXYa{#2}\Figg@tXY{#3}%
    \advance\v@lX-\v@lXa\advance\v@lY-\v@lYa%
    \Figv@ctCreg#1(\v@lX,\v@lY)}\ignorespaces\fi}
\def\figvectPTD#1[#2,#3]{\ifps@cri{\Figg@tXYa{#2}\Figg@tXY{#3}%
    \advance\v@lX-\v@lXa\advance\v@lY-\v@lYa\advance\v@lZ-\v@lZa%
    \Figv@ctCreg#1(\v@lX,\v@lY,\v@lZ)}\ignorespaces\fi}
\def\figvectUDD#1[#2]{\ifps@cri{\n@rmeuc\v@leur{#2}\invers@\v@leur\v@leur%
    \delt@=\repdecn@mb{\v@leur}\unit@\edef\v@ldelt@{\repdecn@mb{\delt@}}%
    \Figg@tXY{#2}\v@lX=\v@ldelt@\v@lX\v@lY=\v@ldelt@\v@lY%
    \Figv@ctCreg#1(\v@lX,\v@lY)}\ignorespaces\fi}
\def\figvectUTD#1[#2]{\ifps@cri{\n@rmeuc\v@leur{#2}\invers@\v@leur\v@leur%
    \delt@=\repdecn@mb{\v@leur}\unit@\edef\v@ldelt@{\repdecn@mb{\delt@}}%
    \Figg@tXY{#2}\v@lX=\v@ldelt@\v@lX\v@lY=\v@ldelt@\v@lY\v@lZ=\v@ldelt@\v@lZ%
    \Figv@ctCreg#1(\v@lX,\v@lY,\v@lZ)}\ignorespaces\fi}
\def\figvisu#1#2#3{\c@ldefproj\initb@undb@x\setbox\b@xvisu=\hbox{\ignorespaces#3}%
    \v@lXa=-\c@@rdYmin\v@lYa=\c@@rdYmax\advance\v@lYa-\c@@rdYmin%
    \v@lX=\c@@rdXmax\advance\v@lX-\c@@rdXmin%
    \setbox#1=\hbox{#2}\v@lY=-\v@lX\maxim@m{\v@lX}{\v@lX}{\wd#1}%
    \advance\v@lY\v@lX\divide\v@lY\tw@\advance\v@lY-\c@@rdXmin%
    \setbox#1=\vbox{\parindent0mm\hsize=\v@lX\vskip\v@lYa%
    \rlap{\hskip\v@lY\smash{\raise\v@lXa\box\b@xvisu}}%
    \def\t@xt@{#2}\ifx\t@xt@\empty\else\medskip\centerline{#2}\fi}\wd#1=\v@lX}
\newbox\Bt@rget\setbox\Bt@rget=\null
\newbox\BminTD@\setbox\BminTD@=\null
\newbox\BmaxTD@\setbox\BmaxTD@=\null
\newif\ifnewt@rgetpt\newif\ifnewdis@b
\def\b@undb@xTD#1#2#3{%
    \relax\ifdim#1<\wd\BminTD@\global\wd\BminTD@=#1\fi%
    \relax\ifdim#2<\ht\BminTD@\global\ht\BminTD@=#2\fi%
    \relax\ifdim#3<\dp\BminTD@\global\dp\BminTD@=#3\fi%
    \relax\ifdim#1>\wd\BmaxTD@\global\wd\BmaxTD@=#1\fi%
    \relax\ifdim#2>\ht\BmaxTD@\global\ht\BmaxTD@=#2\fi%
    \relax\ifdim#3>\dp\BmaxTD@\global\dp\BmaxTD@=#3\fi}
\def\c@ldefdisob{{\ifdim\wd\BminTD@<\maxdimen\v@leur=\wd\BmaxTD@\advance\v@leur-\wd\BminTD@%
    \delt@=\ht\BmaxTD@\advance\delt@-\ht\BminTD@\maxim@m{\v@leur}{\v@leur}{\delt@}%
    \delt@=\dp\BmaxTD@\advance\delt@-\dp\BminTD@\maxim@m{\v@leur}{\v@leur}{\delt@}%
    \v@leur=5\v@leur\else\v@leur=800pt\fi\c@ldefdisob@{\v@leur}}}
\def\c@ldefdisob@#1{{\v@leur=#1\ifdim\v@leur<\p@\v@leur=800pt\fi%
    \xdef\disob@intern{\repdecn@mb{\v@leur}}%
    \delt@=\ptT@unit@\v@leur\xdef\disob@unit{\repdecn@mb{\delt@}}%
    \f@ctech=\@ne\loop\ifdim\v@leur>\t@n pt\divide\v@leur\t@n\multiply\f@ctech\t@n\repeat%
    \xdef\disob@{\repdecn@mb{\v@leur}}\xdef\divf@ctproj{\the\f@ctech}}%
    \global\newdis@btrue}
\def\c@ldeft@rgetpt{\newt@rgetpttrue\def\t@rgetpt{CenterBoundBox}{%
    \delt@=\wd\BmaxTD@\advance\delt@-\wd\BminTD@\divide\delt@\tw@%
    \v@leur=\wd\BminTD@\advance\v@leur\delt@\global\wd\Bt@rget=\v@leur%
    \delt@=\ht\BmaxTD@\advance\delt@-\ht\BminTD@\divide\delt@\tw@%
    \v@leur=\ht\BminTD@\advance\v@leur\delt@\global\ht\Bt@rget=\v@leur%
    \delt@=\dp\BmaxTD@\advance\delt@-\dp\BminTD@\divide\delt@\tw@%
    \v@leur=\dp\BminTD@\advance\v@leur\delt@\global\dp\Bt@rget=\v@leur}}
\def\c@ldefprojTD{\ifnewt@rgetpt\else\c@ldeft@rgetpt\fi\ifnewdis@b\else\c@ldefdisob\fi}
\def\c@lprojcav{% Projection cavaliere : X = x + y L cos t, Y = z + y L sin t
    \v@lZa=\cxa@\v@lY\advance\v@lX\v@lZa%
    \v@lZa=\cxb@\v@lY\v@lY=\v@lZ\advance\v@lY\v@lZa\ignorespaces}
\def\c@lprojrea{% Projection realiste
    \advance\v@lX-\wd\Bt@rget\advance\v@lY-\ht\Bt@rget\advance\v@lZ-\dp\Bt@rget%
    \v@lZa=\cza@\v@lX\advance\v@lZa\czb@\v@lY\advance\v@lZa\czc@\v@lZ%
    \divide\v@lZa\divf@ctproj\advance\v@lZa\disob@ pt\invers@{\v@lZa}{\v@lZa}%
    \v@lZa=\disob@\v@lZa\edef\v@lcoef{\repdecn@mb{\v@lZa}}%
    \v@lXa=\cxa@\v@lX\advance\v@lXa\cxb@\v@lY\v@lXa=\v@lcoef\v@lXa%
    \v@lY=\cyb@\v@lY\advance\v@lY\cya@\v@lX\advance\v@lY\cyc@\v@lZ%
    \v@lY=\v@lcoef\v@lY\v@lX=\v@lXa\ignorespaces}
\def\c@lprojort{% Projection orthogonale
    \v@lXa=\cxa@\v@lX\advance\v@lXa\cxb@\v@lY%
    \v@lY=\cyb@\v@lY\advance\v@lY\cya@\v@lX\advance\v@lY\cyc@\v@lZ%
    \v@lX=\v@lXa\ignorespaces}
\def\Figptpr@j#1:#2/#3/{{\Figg@tXY{#3}\superc@lprojSP%
    \Figp@intregDD#1:{#2}(\v@lX,\v@lY)}\ignorespaces}
\def\figsetobdistTD(#1){{\ifcurr@ntPS%
    \immediate\write16{*** \BS@ figsetobdist is ignored inside a
     \BS@ psbeginfig-\BS@ psendfig block.}%
    \else\v@leur=#1\unit@\c@ldefdisob@{\v@leur}\fi}\ignorespaces}
\def\Figs@tproj#1{%
    \if#13 \d@faultproj\else\if#1c\d@faultproj%
    \else\if#1o\xdef\curr@ntproj{1}\xdef\typ@proj{orthogonal}%
         \figsetviewTD(\def@ultpsi,\def@ulttheta)%
         \global\let\c@lprojSP=\c@lprojort\global\let\superc@lprojSP=\c@lprojort%
    \else\if#1r\xdef\curr@ntproj{2}\xdef\typ@proj{realistic}%
         \figsetviewTD(\def@ultpsi,\def@ulttheta)%
         \global\let\c@lprojSP=\c@lprojrea\global\let\superc@lprojSP=\c@lprojrea%
    \else\d@faultproj\message{*** Unknown projection. Cavalier projection assumed.}%
    \fi\fi\fi\fi}
\def\d@faultproj{\xdef\curr@ntproj{0}\xdef\typ@proj{cavalier}\figsetviewTD(\def@ultpsi,0.5)%
         \global\let\c@lprojSP=\c@lprojcav\global\let\superc@lprojSP=\c@lprojcav}
\def\figsettargetTD[#1]{{\ifcurr@ntPS%
    \immediate\write16{*** \BS@ figsettarget is ignored inside a
     \BS@ psbeginfig-\BS@ psendfig block.}%
    \else\global\newt@rgetpttrue\xdef\t@rgetpt{#1}\Figg@tXY{#1}\global\wd\Bt@rget=\v@lX%
    \global\ht\Bt@rget=\v@lY\global\dp\Bt@rget=\v@lZ\fi}\ignorespaces}
\def\figsetviewTD(#1){\ifcurr@ntPS%
     \immediate\write16{*** \BS@ figsetview is ignored inside a
     \BS@ psbeginfig-\BS@ psendfig block.}\else\Figsetview@#1,:\fi\ignorespaces}
\def\Figsetview@#1,#2:{{\xdef\v@lPsi{#1}\def\t@xt@{#2}%
    \ifx\t@xt@\empty\def\@rgdeux{\v@lTheta}\else%
    \def\Xarg@##1,{\edef\@rgdeux{##1}}\Xarg@#2\fi%
    \c@ssin{\costhet@}{\sinthet@}{#1}\v@lmin=\costhet@ pt\v@lmax=\sinthet@ pt%
    \ifcase\curr@ntproj%
    \v@leur=\@rgdeux\v@lmin\xdef\cxa@{\repdecn@mb{\v@leur}}%
    \v@leur=\@rgdeux\v@lmax\xdef\cxb@{\repdecn@mb{\v@leur}}\v@leur=\@rgdeux pt%
    \relax\ifdim\v@leur>\p@\message{*** Lambda too large ! See \BS@ figset proj() !}\fi%
    \else%
    \v@lmax=-\v@lmax\xdef\cxa@{\repdecn@mb{\v@lmax}}\xdef\cxb@{\costhet@}%
    \ifx\t@xt@\empty\edef\@rgdeux{\def@ulttheta}\fi\c@ssin{\C@}{\S@}{\@rgdeux}%
    \v@lmax=-\S@ pt%
    \v@leur=\v@lmax\v@leur=\costhet@\v@leur\xdef\cya@{\repdecn@mb{\v@leur}}%
    \v@leur=\v@lmax\v@leur=\sinthet@\v@leur\xdef\cyb@{\repdecn@mb{\v@leur}}%
    \xdef\cyc@{\C@}\v@lmin=-\C@ pt%
    \v@leur=\v@lmin\v@leur=\costhet@\v@leur\xdef\cza@{\repdecn@mb{\v@leur}}%
    \v@leur=\v@lmin\v@leur=\sinthet@\v@leur\xdef\czb@{\repdecn@mb{\v@leur}}%
    \xdef\czc@{\repdecn@mb{\v@lmax}}\fi%
    \xdef\v@lTheta{\@rgdeux}}}
\def\def@ultpsi{40}
\def\def@ulttheta{25}
\def\figset#1(#2){\def\t@xt@{#1}\ifx\t@xt@\empty\trtlis@rg{#2}{\Figsetwr@te}% write
    \else\keln@mde#1|%
    \def\n@mref{pr}\ifx\l@debut\n@mref\ifcurr@ntPS% projection
     \immediate\write16{*** \BS@ figset proj(...) is ignored inside a
     \BS@ psbeginfig-\BS@ psendfig block.}\else\trtlis@rg{#2}{\Figsetpr@j}\fi\else%
    \def\n@mref{wr}\ifx\l@debut\n@mref\trtlis@rg{#2}{\Figsetwr@te}\else% write
    \immediate\write16{*** Unknown keyword: \BS@ figset #1(...)}%
    \fi\fi\fi\ignorespaces}
\def\Figsetpr@j#1=#2|{\keln@mtr#1|%
    \def\n@mref{dep}\ifx\l@debut\n@mref\Figsetd@p{#2}\else% depth (lambda)
    \def\n@mref{dis}\ifx\l@debut\n@mref%
     \ifnum\curr@ntproj=\tw@\figsetobdist(#2)\else\Figset@rr\fi\else% dist
    \def\n@mref{lam}\ifx\l@debut\n@mref\Figsetd@p{#2}\else% depth (lambda)
    \def\n@mref{lat}\ifx\l@debut\n@mref\Figsetth@{#2}\else% latitude (theta)
    \def\n@mref{lon}\ifx\l@debut\n@mref\figsetview(#2)\else% longitude (psi)
    \def\n@mref{psi}\ifx\l@debut\n@mref\figsetview(#2)\else% longitude (psi)
    \def\n@mref{tar}\ifx\l@debut\n@mref%
     \ifnum\curr@ntproj=\tw@\figsettarget[#2]\else\Figset@rr\fi\else% target point
    \def\n@mref{the}\ifx\l@debut\n@mref\Figsetth@{#2}\else% latitude (theta)
    \immediate\write16{*** Unknown attribute: \BS@ figset proj(..., #1=...).}%
    \fi\fi\fi\fi\fi\fi\fi\fi}
\def\Figsetd@p#1{\ifnum\curr@ntproj=\z@\figsetview(\v@lPsi,#1)\else\Figset@rr\fi}
\def\Figsetth@#1{\ifnum\curr@ntproj=\z@\Figset@rr\else\figsetview(\v@lPsi,#1)\fi}
\def\Figset@rr{\message{*** \BS@ figset proj(): Attribute "\n@mref" ignored, incompatible
    with current projection}}
\def\initb@undb@xTD{\wd\BminTD@=\maxdimen\ht\BminTD@=\maxdimen\dp\BminTD@=\maxdimen%
    \wd\BmaxTD@=-\maxdimen\ht\BmaxTD@=-\maxdimen\dp\BmaxTD@=-\maxdimen}
\newbox\Gb@x      % boite a tout faire
\newbox\Gb@xSC    % boite qui contient le point marker
\newtoks\c@nsymb  % the point marker
\newif\ifr@undcoord\newif\ifunitpr@sent
\def\unssqrttw@{0.707106 }
\def\figAst{\raise-1.15ex\hbox{$\ast$}}
\def\figBullet{\raise-1.15ex\hbox{$\bullet$}}
\def\figCirc{\raise-1.15ex\hbox{$\circ$}}
\def\figDiamond{\raise-1.15ex\hbox{$\diamond$}}%
\def\boxit#1#2{\leavevmode\hbox{\vrule\vbox{\hrule\vglue#1%
    \vtop{\hbox{\kern#1{#2}\kern#1}\vglue#1\hrule}}\vrule}}
\def\centertext#1#2{\vbox{\hsize#1\parindent0cm%
    \leftskip=0pt plus 1fil\rightskip=0pt plus 1fil\parfillskip=0pt{#2}}}
\def\lefttext#1#2{\vbox{\hsize#1\parindent0cm\rightskip=0pt plus 1fil#2}}
\def\c@nterpt{\ignorespaces%
    \kern-.5\wd\Gb@xSC%
    \raise-.5\ht\Gb@xSC\rlap{\hbox{\raise.5\dp\Gb@xSC\hbox{\copy\Gb@xSC}}}%
    \kern .5\wd\Gb@xSC\ignorespaces}
\def\b@undb@xSC#1#2{{\v@lXa=#1\v@lYa=#2%
    \v@leur=\ht\Gb@xSC\advance\v@leur\dp\Gb@xSC%
    \advance\v@lXa-.5\wd\Gb@xSC\advance\v@lYa-.5\v@leur\b@undb@x{\v@lXa}{\v@lYa}%
    \advance\v@lXa\wd\Gb@xSC\advance\v@lYa\v@leur\b@undb@x{\v@lXa}{\v@lYa}}}
\def\@keldist#1#2{\edef\Dist@n{#2}\y@tiunit{\Dist@n}%
    \ifunitpr@sent#1=\Dist@n\else#1=\Dist@n\unit@\fi}
\def\y@tiunit#1{\unitpr@sentfalse\expandafter\y@tiunit@#1:}
\def\y@tiunit@#1#2:{\ifcat#1a\unitpr@senttrue\else\def\l@suite{#2}%
    \ifx\l@suite\empty\else\y@tiunit@#2:\fi\fi}
\def\figcoordDD#1{{\v@lX=\ptT@unit@\v@lX\v@lY=\ptT@unit@\v@lY%
    \ifr@undcoord\ifcase#1\v@leur=0.5pt\or\v@leur=0.05pt\or\v@leur=0.005pt%
    \or\v@leur=0.0005pt\else\v@leur=\z@\fi%
    \ifdim\v@lX<\z@\advance\v@lX-\v@leur\else\advance\v@lX\v@leur\fi%
    \ifdim\v@lY<\z@\advance\v@lY-\v@leur\else\advance\v@lY\v@leur\fi\fi%
    (\@ffichnb{#1}{\repdecn@mb{\v@lX}},\ifmmode\else\thinspace\fi%
    \@ffichnb{#1}{\repdecn@mb{\v@lY}})}}
\def\@ffichnb#1#2{\def\@@ffich{\@ffich#1(}\edef\n@mbre{#2}%
    \expandafter\@@ffich\n@mbre)}
\def\@ffich#1(#2.#3){{#2\ifnum#1>\z@.\fi\def\dig@ts{#3}\s@mme=\z@%
    \loop\ifnum\s@mme<#1\expandafter\@ffichdec\dig@ts:\advance\s@mme\@ne\repeat}}
\def\@ffichdec#1#2:{\relax#1\def\dig@ts{#20}}
\def\figcoordTD#1{{\v@lX=\ptT@unit@\v@lX\v@lY=\ptT@unit@\v@lY\v@lZ=\ptT@unit@\v@lZ%
    \ifr@undcoord\ifcase#1\v@leur=0.5pt\or\v@leur=0.05pt\or\v@leur=0.005pt%
    \or\v@leur=0.0005pt\else\v@leur=\z@\fi%
    \ifdim\v@lX<\z@\advance\v@lX-\v@leur\else\advance\v@lX\v@leur\fi%
    \ifdim\v@lY<\z@\advance\v@lY-\v@leur\else\advance\v@lY\v@leur\fi%
    \ifdim\v@lZ<\z@\advance\v@lZ-\v@leur\else\advance\v@lZ\v@leur\fi\fi%
    (\@ffichnb{#1}{\repdecn@mb{\v@lX}},\ifmmode\else\thinspace\fi%
     \@ffichnb{#1}{\repdecn@mb{\v@lY}},\ifmmode\else\thinspace\fi%
     \@ffichnb{#1}{\repdecn@mb{\v@lZ}})}}
\def\figsetroundcoord#1{\expandafter\Figsetr@undcoord#1:\ignorespaces}
\def\Figsetr@undcoord#1#2:{\if#1n\r@undcoordfalse\else\r@undcoordtrue\fi}
\def\Figsetwr@te#1=#2|{\keln@mun#1|%
    \def\n@mref{m}\ifx\l@debut\n@mref\figsetmark{#2}\else% mark
    \immediate\write16{*** Unknown attribute: \BS@ figset (..., #1=...)}%
    \fi}
\def\figsetmark#1{\c@nsymb={#1}\setbox\Gb@xSC=\hbox{\the\c@nsymb}\ignorespaces}
\def\figsetptname#1{\def\ptn@me##1{#1}\ignorespaces}
\def\FigWrit@L#1:#2(#3,#4){\ignorespaces\@keldist\v@leur{#3}\@keldist\delt@{#4}%
    \C@rp@r@m\def\list@num{#1}\@ecfor\p@int:=\list@num\do{\FigWrit@pt{\p@int}{#2}}}
\def\FigWrit@pt#1#2{\FigWp@r@m{#1}{#2}\Vc@rrect\figWp@si%
    \ifdim\wd\Gb@xSC>\z@\b@undb@xSC{\v@lX}{\v@lY}\fi\figWBB@x}
\def\FigWp@r@m#1#2{\Figg@tXY{#1}%
    \setbox\Gb@x=\hbox{\def\t@xt@{#2}\ifx\t@xt@\empty\Figg@tT{#1}\else#2\fi}\c@lprojSP}
\let\Vc@rrect=\relax
\let\C@rp@r@m=\relax
\def\figwrite[#1]#2{{\ignorespaces\def\list@num{#1}\@ecfor\p@int:=\list@num\do{%
    \setbox\Gb@x=\hbox{\def\t@xt@{#2}\ifx\t@xt@\empty\Figg@tT{\p@int}\else#2\fi}%
    \Figwrit@{\p@int}}}\ignorespaces}
\def\Figwrit@#1{\Figg@tXY{#1}\c@lprojSP%
    \rlap{\kern\v@lX\raise\v@lY\hbox{\unhcopy\Gb@x}}\v@leur=\v@lY%
    \advance\v@lY\ht\Gb@x\b@undb@x{\v@lX}{\v@lY}\advance\v@lX\wd\Gb@x%
    \v@lY=\v@leur\advance\v@lY-\dp\Gb@x\b@undb@x{\v@lX}{\v@lY}}
\def\figwritec[#1]#2{{\ignorespaces\def\list@num{#1}%
    \@ecfor\p@int:=\list@num\do{\Figwrit@c{\p@int}{#2}}}\ignorespaces}
\def\Figwrit@c#1#2{\FigWp@r@m{#1}{#2}%
    \rlap{\kern\v@lX\raise\v@lY\hbox{\rlap{\kern-.5\wd\Gb@x%
    \raise-.5\ht\Gb@x\hbox{\raise.5\dp\Gb@x\hbox{\unhcopy\Gb@x}}}}}%
    \v@leur=\ht\Gb@x\advance\v@leur\dp\Gb@x%
    \advance\v@lX-.5\wd\Gb@x\advance\v@lY-.5\v@leur\b@undb@x{\v@lX}{\v@lY}%
    \advance\v@lX\wd\Gb@x\advance\v@lY\v@leur\b@undb@x{\v@lX}{\v@lY}}
\def\figwritep[#1]{{\ignorespaces\def\list@num{#1}\setbox\Gb@x=\hbox{\c@nterpt}%
    \@ecfor\p@int:=\list@num\do{\Figwrit@{\p@int}}}\ignorespaces}
\def\figwritew#1:#2(#3){\figwritegcw#1:{#2}(#3,0pt)}
\def\figwritee#1:#2(#3){\figwritegce#1:{#2}(#3,0pt)}
\def\figwriten#1:#2(#3){{\def\Vc@rrect{\v@lZ=\v@leur\advance\v@lZ\dp\Gb@x}%
    \Figwrit@NS#1:{#2}(#3)}\ignorespaces}
\def\figwrites#1:#2(#3){{\def\Vc@rrect{\v@lZ=-\v@leur\advance\v@lZ-\ht\Gb@x}%
    \Figwrit@NS#1:{#2}(#3)}\ignorespaces}
\def\Figwrit@NS#1:#2(#3){\let\figWp@si=\FigWp@siNS\let\figWBB@x=\FigWBB@xNS%
    \FigWrit@L#1:{#2}(#3,0pt)}
\def\FigWp@siNS{\rlap{\kern\v@lX\raise\v@lY\hbox{\rlap{\kern-.5\wd\Gb@x%
    \raise\v@lZ\hbox{\unhcopy\Gb@x}}\c@nterpt}}}
\def\FigWBB@xNS{\advance\v@lY\v@lZ%
    \advance\v@lY-\dp\Gb@x\advance\v@lX-.5\wd\Gb@x\b@undb@x{\v@lX}{\v@lY}%
    \advance\v@lY\ht\Gb@x\advance\v@lY\dp\Gb@x%
    \advance\v@lX\wd\Gb@x\b@undb@x{\v@lX}{\v@lY}}
\def\figwritenw#1:#2(#3){{\let\figWp@si=\FigWp@sigW\let\figWBB@x=\FigWBB@xgWE%
    \def\C@rp@r@m{\v@leur=\unssqrttw@\v@leur\delt@=\v@leur%
    \ifdim\delt@=\z@\delt@=\epsil@n\fi}\let@xte={-}\FigWrit@L#1:{#2}(#3,0pt)}\ignorespaces}
\def\figwritesw#1:#2(#3){{\let\figWp@si=\FigWp@sigW\let\figWBB@x=\FigWBB@xgWE%
    \def\C@rp@r@m{\v@leur=\unssqrttw@\v@leur\delt@=-\v@leur%
    \ifdim\delt@=\z@\delt@=-\epsil@n\fi}\let@xte={-}\FigWrit@L#1:{#2}(#3,0pt)}\ignorespaces}
\def\figwritene#1:#2(#3){{\let\figWp@si=\FigWp@sigE\let\figWBB@x=\FigWBB@xgWE%
    \def\C@rp@r@m{\v@leur=\unssqrttw@\v@leur\delt@=\v@leur%
    \ifdim\delt@=\z@\delt@=\epsil@n\fi}\let@xte={}\FigWrit@L#1:{#2}(#3,0pt)}\ignorespaces}
\def\figwritese#1:#2(#3){{\let\figWp@si=\FigWp@sigE\let\figWBB@x=\FigWBB@xgWE%
    \def\C@rp@r@m{\v@leur=\unssqrttw@\v@leur\delt@=-\v@leur%
    \ifdim\delt@=\z@\delt@=-\epsil@n\fi}\let@xte={}\FigWrit@L#1:{#2}(#3,0pt)}\ignorespaces}
\def\figwritegw#1:#2(#3,#4){{\let\figWp@si=\FigWp@sigW\let\figWBB@x=\FigWBB@xgWE%
    \let@xte={-}\FigWrit@L#1:{#2}(#3,#4)}\ignorespaces}
\def\figwritege#1:#2(#3,#4){{\let\figWp@si=\FigWp@sigE\let\figWBB@x=\FigWBB@xgWE%
    \let@xte={}\FigWrit@L#1:{#2}(#3,#4)}\ignorespaces}
\def\FigWp@sigW{\v@lXa=\z@\v@lYa=\ht\Gb@x\advance\v@lYa\dp\Gb@x%
    \ifdim\delt@>\z@\relax%
    \rlap{\kern\v@lX\raise\v@lY\hbox{\rlap{\kern-\wd\Gb@x\kern-\v@leur%
          \raise\delt@\hbox{\raise\dp\Gb@x\hbox{\unhcopy\Gb@x}}}\c@nterpt}}%
    \else\ifdim\delt@<\z@\relax\v@lYa=-\v@lYa%
    \rlap{\kern\v@lX\raise\v@lY\hbox{\rlap{\kern-\wd\Gb@x\kern-\v@leur%
          \raise\delt@\hbox{\raise-\ht\Gb@x\hbox{\unhcopy\Gb@x}}}\c@nterpt}}%
    \else\v@lXa=-.5\v@lYa%
    \rlap{\kern\v@lX\raise\v@lY\hbox{\rlap{\kern-\wd\Gb@x\kern-\v@leur%
          \raise-.5\ht\Gb@x\hbox{\raise.5\dp\Gb@x\hbox{\unhcopy\Gb@x}}}\c@nterpt}}%
    \fi\fi}
\def\FigWp@sigE{\v@lXa=\z@\v@lYa=\ht\Gb@x\advance\v@lYa\dp\Gb@x%
    \ifdim\delt@>\z@\relax%
    \rlap{\kern\v@lX\raise\v@lY\hbox{\c@nterpt\kern\v@leur%
          \raise\delt@\hbox{\raise\dp\Gb@x\hbox{\unhcopy\Gb@x}}}}%
    \else\ifdim\delt@<\z@\relax\v@lYa=-\v@lYa%
    \rlap{\kern\v@lX\raise\v@lY\hbox{\c@nterpt\kern\v@leur%
          \raise\delt@\hbox{\raise-\ht\Gb@x\hbox{\unhcopy\Gb@x}}}}%
    \else\v@lXa=-.5\v@lYa%
    \rlap{\kern\v@lX\raise\v@lY\hbox{\c@nterpt\kern\v@leur%
          \raise-.5\ht\Gb@x\hbox{\raise.5\dp\Gb@x\hbox{\unhcopy\Gb@x}}}}%
    \fi\fi}
\def\FigWBB@xgWE{\advance\v@lY\delt@%
    \advance\v@lX\the\let@xte\v@leur\advance\v@lY\v@lXa\b@undb@x{\v@lX}{\v@lY}%
    \advance\v@lX\the\let@xte\wd\Gb@x\advance\v@lY\v@lYa\b@undb@x{\v@lX}{\v@lY}}
\def\figwritegcw#1:#2(#3,#4){{\let\figWp@si=\FigWp@sigcW\let\figWBB@x=\FigWBB@xgcWE%
    \let@xte={-}\FigWrit@L#1:{#2}(#3,#4)}\ignorespaces}
\def\figwritegce#1:#2(#3,#4){{\let\figWp@si=\FigWp@sigcE\let\figWBB@x=\FigWBB@xgcWE%
    \let@xte={}\FigWrit@L#1:{#2}(#3,#4)}\ignorespaces}
\def\FigWp@sigcW{\rlap{\kern\v@lX\raise\v@lY\hbox{\rlap{\kern-\wd\Gb@x\kern-\v@leur%
     \raise-.5\ht\Gb@x\hbox{\raise\delt@\hbox{\raise.5\dp\Gb@x\hbox{\unhcopy\Gb@x}}}}%
     \c@nterpt}}}
\def\FigWp@sigcE{\rlap{\kern\v@lX\raise\v@lY\hbox{\c@nterpt\kern\v@leur%
    \raise-.5\ht\Gb@x\hbox{\raise\delt@\hbox{\raise.5\dp\Gb@x\hbox{\unhcopy\Gb@x}}}}}}
\def\FigWBB@xgcWE{\v@lZ=\ht\Gb@x\advance\v@lZ\dp\Gb@x%
    \advance\v@lX\the\let@xte\v@leur\advance\v@lY\delt@\advance\v@lY.5\v@lZ%
    \b@undb@x{\v@lX}{\v@lY}%
    \advance\v@lX\the\let@xte\wd\Gb@x\advance\v@lY-\v@lZ\b@undb@x{\v@lX}{\v@lY}}
\def\figwritebn#1:#2(#3){{\def\Vc@rrect{\v@lZ=\v@leur}\Figwrit@NS#1:{#2}(#3)}\ignorespaces}
\def\figwritebs#1:#2(#3){{\def\Vc@rrect{\v@lZ=-\v@leur}\Figwrit@NS#1:{#2}(#3)}\ignorespaces}
\def\figwritebw#1:#2(#3){{\let\figWp@si=\FigWp@sibW\let\figWBB@x=\FigWBB@xbWE%
    \let@xte={-}\FigWrit@L#1:{#2}(#3,0pt)}\ignorespaces}
\def\figwritebe#1:#2(#3){{\let\figWp@si=\FigWp@sibE\let\figWBB@x=\FigWBB@xbWE%
    \let@xte={}\FigWrit@L#1:{#2}(#3,0pt)}\ignorespaces}
\def\FigWp@sibW{\rlap{\kern\v@lX\raise\v@lY\hbox{\rlap{\kern-\wd\Gb@x\kern-\v@leur%
          \hbox{\unhcopy\Gb@x}}\c@nterpt}}}
\def\FigWp@sibE{\rlap{\kern\v@lX\raise\v@lY\hbox{\c@nterpt\kern\v@leur%
          \hbox{\unhcopy\Gb@x}}}}
\def\FigWBB@xbWE{\v@lZ=\ht\Gb@x\advance\v@lZ\dp\Gb@x%
    \advance\v@lX\the\let@xte\v@leur\advance\v@lY\ht\Gb@x\b@undb@x{\v@lX}{\v@lY}%
    \advance\v@lX\the\let@xte\wd\Gb@x\advance\v@lY-\v@lZ\b@undb@x{\v@lX}{\v@lY}}
\newread\frf@g  \newwrite\fwf@g
\newif\ifcurr@ntPS
\newif\ifps@cri
\newif\ifUse@llipse
\newif\ifpsdebugmode \psdebugmodefalse 
\newif\ifPDFm@ke
\ifx\pdfliteral\undefined\else\ifnum\pdfoutput>\z@\PDFm@ketrue\fi\fi
\ifPDFm@ke
 \def\c@mcurveto{c}
 \def\c@mfill{f}
 \def\c@mgsave{q}
 \def\c@mgrestore{Q}
 \def\c@mlineto{l}
 \def\c@mmoveto{m}
 \def\c@msetgray{g}     \def\c@msetgrayStroke{G}
 \def\c@msetcmykcolor{k}\def\c@msetcmykcolorStroke{K}
 \def\c@msetrgbcolor{rg}\def\c@msetrgbcolorStroke{RG}
 \def\d@fprimarC@lor{\curr@ntcolor\space\curr@ntcolorc@md%
               \space\curr@ntcolor\space\curr@ntcolorc@mdStroke}
 \def\d@fsecondC@lor{\sec@ndcolor\space\sec@ndcolorc@md%
               \space\sec@ndcolor\space\sec@ndcolorc@mdStroke}
 \def\d@fthirdC@lor{\th@rdcolor\space\th@rdcolorc@md%
              \space\th@rdcolor\space\th@rdcolorc@mdStroke}
 \def\c@msetdash{d}
 \def\c@msetlinejoin{j}
 \def\c@msetlinewidth{w}
 \def\f@gclosestroke{\immediate\write\fwf@g{s}}
 \def\f@gfill{\immediate\write\fwf@g{\fillc@md}}% Voir la def de \fillc@md ****
 \def\f@gnewpath{}
 \def\f@gstroke{\immediate\write\fwf@g{S}}
\else
 \let\figinsertE=\figinsert
 \def\c@mcurveto{curveto}
 \def\c@mfill{fill}
 \def\c@mgsave{gsave}
 \def\c@mgrestore{grestore}
 \def\c@mlineto{lineto}
 \def\c@mmoveto{moveto}
 \def\c@msetgray{setgray}          \def\c@msetgrayStroke{}
 \def\c@msetcmykcolor{setcmykcolor}\def\c@msetcmykcolorStroke{}
 \def\c@msetrgbcolor{setrgbcolor}  \def\c@msetrgbcolorStroke{}
 \def\d@fprimarC@lor{\curr@ntcolor\space\curr@ntcolorc@md}
 \def\d@fsecondC@lor{\sec@ndcolor\space\sec@ndcolorc@md}
 \def\d@fthirdC@lor{\th@rdcolor\space\th@rdcolorc@md}
 \def\c@msetdash{setdash}
 \def\c@msetlinejoin{setlinejoin}
 \def\c@msetlinewidth{setlinewidth}
 \def\f@gclosestroke{\immediate\write\fwf@g{closepath\space stroke}}
 \def\f@gfill{\immediate\write\fwf@g{\fillc@md}}
 \def\f@gnewpath{\immediate\write\fwf@g{newpath}}
 \def\f@gstroke{\immediate\write\fwf@g{stroke}}
\fi
\def\c@pypsfile#1#2{\c@pyfil@{\immediate\write#1}{#2}}
\def\Figinclud@PDF#1#2{\openin\frf@g=#1\pdfliteral{q #2 0 0 #2 0 0 cm}%
    \c@pyfil@{\pdfliteral}{\frf@g}\pdfliteral{Q}\closein\frf@g}
\newif\ifmored@ta
\def\c@pyfil@#1#2{\def\blankline{\par}{\catcode`\%=12
    \loop\ifeof#2\mored@tafalse\else\mored@tatrue\immediate\read#2 to\tr@c
    \ifx\tr@c\blankline\else#1{\tr@c}\fi\fi\ifmored@ta\repeat}}
\def\keln@mun#1#2|{\def\l@debut{#1}\def\l@suite{#2}}
\def\keln@mde#1#2#3|{\def\l@debut{#1#2}\def\l@suite{#3}}
\def\keln@mtr#1#2#3#4|{\def\l@debut{#1#2#3}\def\l@suite{#4}}
\def\keln@mqu#1#2#3#4#5|{\def\l@debut{#1#2#3#4}\def\l@suite{#5}}
\let\@psffilein=\frf@g % file to \read
\newif\if@psffileok    % continue looking for the bounding box?
\newif\if@psfbbfound   % success?
\newif\if@psfverbose   % report what you're making?
\@psfverbosetrue
\def\@psfgetbb#1{\global\@psfbbfoundfalse%
\global\def\@psfllx{0}\global\def\@psflly{0}%
\global\def\@psfurx{30}\global\def\@psfury{30}%
\openin\@psffilein=#1
\ifeof\@psffilein\errmessage{I couldn't open #1, will ignore it}\else
   \edef\setcolonc@tcode{\catcode`\noexpand\:\the\catcode`\:\relax}%
   {\@psffileoktrue \chardef\other=12
    \def\do##1{\catcode`##1=\other}\dospecials \catcode`\ =10 \setcolonc@tcode
    \loop
       \read\@psffilein to \@psffileline
       \ifeof\@psffilein\@psffileokfalse\else
          \expandafter\@psfaux\@psffileline:. \\%
       \fi
   \if@psffileok\repeat
   \if@psfbbfound\else
    \if@psfverbose\message{No bounding box comment in #1; using defaults}\fi\fi
   }\closein\@psffilein\fi}%
{\catcode`\%=12 \global\let\@psfpercent=%\global\def\@psfbblit{%BoundingBox}}%
\long\def\@psfaux#1#2:#3\\{\ifx#1\@psfpercent
   \def\testit{#2}\ifx\testit\@psfbblit
      \@psfgrab #3 . . . \\%
      \@psffileokfalse
      \global\@psfbbfoundtrue
   \fi\else\ifx#1\par\else\@psffileokfalse\fi\fi}%
\def\@psfempty{}%
\def\@psfgrab #1 #2 #3 #4 #5\\{%
\global\def\@psfllx{#1}\ifx\@psfllx\@psfempty
      \@psfgrab #2 #3 #4 #5 .\\\else
   \global\def\@psflly{#2}%
   \global\def\@psfurx{#3}\global\def\@psfury{#4}\fi}%
\def\PSwrit@cmd#1#2#3{{\Figg@tXY{#1}\c@lprojSP\b@undb@x{\v@lX}{\v@lY}%
    \v@lX=\ptT@ptps\v@lX\v@lY=\ptT@ptps\v@lY%
    \immediate\write#3{\repdecn@mb{\v@lX}\space\repdecn@mb{\v@lY}\space#2}}}
\def\PSwrit@cmdS#1#2#3#4#5{{\Figg@tXY{#1}\c@lprojSP\b@undb@x{\v@lX}{\v@lY}%
    \global\result@t=\v@lX\global\result@@t=\v@lY%
    \v@lX=\ptT@ptps\v@lX\v@lY=\ptT@ptps\v@lY%
    \immediate\write#3{\repdecn@mb{\v@lX}\space\repdecn@mb{\v@lY}\space#2}}%
    \edef#4{\the\result@t}\edef#5{\the\result@@t}}
\def\psaltitude#1[#2,#3,#4]{{\ifcurr@ntPS\ifps@cri%
    \PSc@mment{psaltitude Square Dim=#1, Triangle=[#2 / #3,#4]}%
    \s@uvc@ntr@l\et@tpsaltitude\resetc@ntr@l{2}\figptorthoprojline-5:=#2/#3,#4/%
    \figvectP -1[#3,#4]\n@rminf{\v@leur}{-1}\vecunit@{-3}{-1}%
    \figvectP -1[-5,#3]\n@rminf{\v@lmin}{-1}\figvectP -2[-5,#4]\n@rminf{\v@lmax}{-2}%
    \ifdim\v@lmin<\v@lmax\s@mme=#3\else\v@lmax=\v@lmin\s@mme=#4\fi%
    \figvectP -4[-5,#2]\vecunit@{-4}{-4}\delt@=#1\unit@%
    \edef\t@ille{\repdecn@mb{\delt@}}\figpttra-1:=-5/\t@ille,-3/%
    \figptstra-3=-5,-1/\t@ille,-4/\psline[#2,-5]\psline[-1,-2,-3]%
    \ifdim\v@leur<\v@lmax\Pss@tsecondSt\psline[-5,\the\s@mme]\Psrest@reSt\fi%
    \PSc@mment{End psaltitude}\resetc@ntr@l\et@tpsaltitude\fi\fi}}
\def\Ps@rcerc#1;#2(#3,#4){\ellBB@x#1;#2,#2(#3,#4,0)%
    \f@gnewpath{\delt@=#2\unit@\delt@=\ptT@ptps\delt@%
    \BdingB@xfalse%
    \PSwrit@cmd{#1}{\repdecn@mb{\delt@}\space #3\space #4\space arc}{\fwf@g}}}
\def\psarccircDD#1;#2(#3,#4){\ifcurr@ntPS\ifps@cri%
    \PSc@mment{psarccircDD Center=#1 ; Radius=#2 (Ang1=#3, Ang2=#4)}%
    \iffillm@de\Ps@rcerc#1;#2(#3,#4)%
    \f@gfill%
    \else\Ps@rcerc#1;#2(#3,#4)\f@gstroke\fi%
    \PSc@mment{End psarccircDD}\fi\fi}
\def\psarccircTD#1,#2,#3;#4(#5,#6){{\ifcurr@ntPS\ifps@cri\s@uvc@ntr@l\et@tpsarccircTD%
    \PSc@mment{psarccircTD Center=#1,P1=#2,P2=#3 ; Radius=#4 (Ang1=#5, Ang2=#6)}%
    \setc@ntr@l{2}\c@lExtAxes#1,#2,#3(#4)\psarcellPATD#1,-4,-5(#5,#6)%
    \PSc@mment{End psarccircTD}\resetc@ntr@l\et@tpsarccircTD\fi\fi}}
\def\c@lExtAxes#1,#2,#3(#4){%
    \figvectPTD-5[#1,#2]\vecunit@{-5}{-5}\figvectNTD-4[#1,#2,#3]\vecunit@{-4}{-4}%
    \figvectNVTD-3[-4,-5]\delt@=#4\unit@\edef\r@yon{\repdecn@mb{\delt@}}%
    \figpttra-4:=#1/\r@yon,-5/\figpttra-5:=#1/\r@yon,-3/}
\def\psarccircPDD#1;#2[#3,#4]{{\ifcurr@ntPS\ifps@cri\s@uvc@ntr@l\et@tpsarccircPDD%
    \PSc@mment{psarccircPDD Center=#1; Radius=#2, [P1=#3, P2=#4]}%
    \Ps@ngleparam#1;#2[#3,#4]\ifdim\v@lmin>\v@lmax\advance\v@lmax\DePI@deg\fi%
    \edef\@ngdeb{\repdecn@mb{\v@lmin}}\edef\@ngfin{\repdecn@mb{\v@lmax}}%
    \psarccirc#1;\r@dius(\@ngdeb,\@ngfin)%
    \PSc@mment{End psarccircPDD}\resetc@ntr@l\et@tpsarccircPDD\fi\fi}}
\def\psarccircPTD#1;#2[#3,#4,#5]{{\ifcurr@ntPS\ifps@cri\s@uvc@ntr@l\et@tpsarccircPTD%
    \PSc@mment{psarccircPTD Center=#1; Radius=#2, [P1=#3, P2=#4, P3=#5]}%
    \setc@ntr@l{2}\c@lExtAxes#1,#3,#5(#2)\psarcellPP#1,-4,-5[#3,#4]%
    \PSc@mment{End psarccircPTD}\resetc@ntr@l\et@tpsarccircPTD\fi\fi}}
\def\Ps@ngleparam#1;#2[#3,#4]{\setc@ntr@l{2}%
    \figvectPDD-1[#1,#3]\vecunit@{-1}{-1}\Figg@tXY{-1}\arct@n\v@lmin(\v@lX,\v@lY)%
    \figvectPDD-2[#1,#4]\vecunit@{-2}{-2}\Figg@tXY{-2}\arct@n\v@lmax(\v@lX,\v@lY)%
    \v@lmin=\rdT@deg\v@lmin\v@lmax=\rdT@deg\v@lmax%
    \v@leur=#2pt\maxim@m{\mili@u}{-\v@leur}{\v@leur}%
    \edef\r@dius{\repdecn@mb{\mili@u}}}
\def\Ps@rcercBz#1;#2(#3,#4){\Ps@rellBz#1;#2,#2(#3,#4,0)}
\def\Ps@rellBz#1;#2,#3(#4,#5,#6){%
    \ellBB@x#1;#2,#3(#4,#5,#6)\BdingB@xfalse%
    \c@lNbarcs{#4}{#5}\v@leur=#4pt\setc@ntr@l{2}\figptell-13::#1;#2,#3(#4,#6)%
    \f@gnewpath\PSwrit@cmd{-13}{\c@mmoveto}{\fwf@g}%
    \s@mme=\z@\bcl@rellBz#1;#2,#3(#6)\BdingB@xtrue}
\def\bcl@rellBz#1;#2,#3(#4){\relax%
    \ifnum\s@mme<\p@rtent\advance\s@mme\@ne%
    \advance\v@leur\delt@\edef\@ngle{\repdecn@mb\v@leur}\figptell-14::#1;#2,#3(\@ngle,#4)%
    \advance\v@leur\delt@\edef\@ngle{\repdecn@mb\v@leur}\figptell-15::#1;#2,#3(\@ngle,#4)%
    \advance\v@leur\delt@\edef\@ngle{\repdecn@mb\v@leur}\figptell-16::#1;#2,#3(\@ngle,#4)%
    \figptscontrolDD-18[-13,-14,-15,-16]%
    \PSwrit@cmd{-18}{}{\fwf@g}\PSwrit@cmd{-17}{}{\fwf@g}%
    \PSwrit@cmd{-16}{\c@mcurveto}{\fwf@g}%
    \figptcopyDD-13:/-16/\bcl@rellBz#1;#2,#3(#4)\fi}
\def\Ps@rell#1;#2,#3(#4,#5,#6){\ellBB@x#1;#2,#3(#4,#5,#6)%
    \f@gnewpath{\v@lmin=#2\unit@\v@lmin=\ptT@ptps\v@lmin%
    \v@lmax=#3\unit@\v@lmax=\ptT@ptps\v@lmax\BdingB@xfalse%
    \PSwrit@cmd{#1}%
    {#6\space\repdecn@mb{\v@lmin}\space\repdecn@mb{\v@lmax}\space #4\space #5\space ellipse}{\fwf@g}}%
    \global\Use@llipsetrue}
\def\psarcellDD#1;#2,#3(#4,#5,#6){{\ifcurr@ntPS\ifps@cri%
    \PSc@mment{psarcellDD Center=#1 ; XRad=#2, YRad=#3 (Ang1=#4, Ang2=#5, Inclination=#6)}%
    \iffillm@de\Ps@rell#1;#2,#3(#4,#5,#6)%
    \f@gfill%
    \else\Ps@rell#1;#2,#3(#4,#5,#6)\f@gstroke\fi%
    \PSc@mment{End psarcellDD}\fi\fi}}
\def\psarcellTD#1;#2,#3(#4,#5,#6){{\ifcurr@ntPS\ifps@cri\s@uvc@ntr@l\et@tpsarcellTD%
    \PSc@mment{psarcellTD Center=#1 ; XRad=#2, YRad=#3 (Ang1=#4, Ang2=#5, Inclination=#6)}%
    \setc@ntr@l{2}\figpttraC -8:=#1/#2,0,0/\figpttraC -7:=#1/0,#3,0/%
    \figvectC -4(0,0,1)\figptsrot -8=-8,-7/#1,#6,-4/\psarcellPATD#1,-8,-7(#4,#5)%
    \PSc@mment{End psarcellTD}\resetc@ntr@l\et@tpsarcellTD\fi\fi}}
\def\psarcellPADD#1,#2,#3(#4,#5){{\ifcurr@ntPS\ifps@cri\s@uvc@ntr@l\et@tpsarcellPADD%
    \PSc@mment{psarcellPADD Center=#1,PtAxis1=#2,PtAxis2=#3 (Ang1=#4, Ang2=#5)}%
    \setc@ntr@l{2}\figvectPDD-1[#1,#2]\vecunit@DD{-1}{-1}\v@lX=\ptT@unit@\result@t%
    \edef\XR@d{\repdecn@mb{\v@lX}}\Figg@tXY{-1}\arct@n\v@lmin(\v@lX,\v@lY)%
    \v@lmin=\rdT@deg\v@lmin\edef\Inclin@{\repdecn@mb{\v@lmin}}%
    \figgetdist\YR@d[#1,#3]\psarcellDD#1;\XR@d,\YR@d(#4,#5,\Inclin@)%
    \PSc@mment{End psarcellPADD}\resetc@ntr@l\et@tpsarcellPADD\fi\fi}}
\def\psarcellPATD#1,#2,#3(#4,#5){{\ifcurr@ntPS\ifps@cri\s@uvc@ntr@l\et@tpsarcellPATD%
    \PSc@mment{psarcellPATD Center=#1,PtAxis1=#2,PtAxis2=#3 (Ang1=#4, Ang2=#5)}%
    \iffillm@de\Ps@rellPATD#1,#2,#3(#4,#5)%
    \f@gfill%
    \else\Ps@rellPATD#1,#2,#3(#4,#5)\f@gstroke\fi%
    \PSc@mment{End psarcellPATD}\resetc@ntr@l\et@tpsarcellPATD\fi\fi}}
\def\Ps@rellPATD#1,#2,#3(#4,#5){\let\c@lprojSP=\relax%
    \setc@ntr@l{2}\figvectPTD-1[#1,#2]\figvectPTD-2[#1,#3]\c@lNbarcs{#4}{#5}%
    \v@leur=#4pt\c@lptellP{#1}{-1}{-2}\Figptpr@j-5:/-3/%
    \f@gnewpath\PSwrit@cmdS{-5}{\c@mmoveto}{\fwf@g}{\X@un}{\Y@un}%
    \edef\C@nt@r{#1}\s@mme=\z@\bcl@rellPATD}
\def\bcl@rellPATD{\relax%
    \ifnum\s@mme<\p@rtent\advance\s@mme\@ne%
    \advance\v@leur\delt@\c@lptellP{\C@nt@r}{-1}{-2}\Figptpr@j-4:/-3/%
    \advance\v@leur\delt@\c@lptellP{\C@nt@r}{-1}{-2}\Figptpr@j-6:/-3/%
    \advance\v@leur\delt@\c@lptellP{\C@nt@r}{-1}{-2}\Figptpr@j-3:/-3/%
    \v@lX=\z@\v@lY=\z@\Figtr@nptDD{-5}{-5}\Figtr@nptDD{2}{-3}%
    \divide\v@lX\@vi\divide\v@lY\@vi%
    \Figtr@nptDD{3}{-4}\Figtr@nptDD{-1.5}{-6}\v@lmin=\v@lX\v@lmax=\v@lY%
    \v@lX=\z@\v@lY=\z@\Figtr@nptDD{2}{-5}\Figtr@nptDD{-5}{-3}%
    \divide\v@lX\@vi\divide\v@lY\@vi\Figtr@nptDD{-1.5}{-4}\Figtr@nptDD{3}{-6}%
    \BdingB@xfalse%
    \Figp@intregDD-4:(\v@lmin,\v@lmax)\PSwrit@cmdS{-4}{}{\fwf@g}{\X@de}{\Y@de}%
    \Figp@intregDD-4:(\v@lX,\v@lY)\PSwrit@cmdS{-4}{}{\fwf@g}{\X@tr}{\Y@tr}%
    \BdingB@xtrue\PSwrit@cmdS{-3}{\c@mcurveto}{\fwf@g}{\X@qu}{\Y@qu}%
    \B@zierBB@x{1}{\Y@un}(\X@un,\X@de,\X@tr,\X@qu)%
    \B@zierBB@x{2}{\X@un}(\Y@un,\Y@de,\Y@tr,\Y@qu)%
    \edef\X@un{\X@qu}\edef\Y@un{\Y@qu}\figptcopyDD-5:/-3/\bcl@rellPATD\fi}
\def\c@lNbarcs#1#2{%
    \delt@=#2pt\advance\delt@-#1pt\maxim@m{\v@lmax}{\delt@}{-\delt@}%
    \v@leur=\v@lmax\divide\v@leur45 \p@rtentiere{\p@rtent}{\v@leur}\advance\p@rtent\@ne%
    \s@mme=\p@rtent\multiply\s@mme\thr@@\divide\delt@\s@mme}
\def\psarcellPP#1,#2,#3[#4,#5]{{\ifcurr@ntPS\ifps@cri\s@uvc@ntr@l\et@tpsarcellPP%
    \PSc@mment{psarcellPP Center=#1,PtAxis1=#2,PtAxis2=#3 [Point1=#4, Point2=#5]}%
    \setc@ntr@l{2}\figvectP-2[#1,#3]\vecunit@{-2}{-2}\v@lmin=\result@t%
    \invers@{\v@lmax}{\v@lmin}%
    \figvectP-1[#1,#2]\vecunit@{-1}{-1}\v@leur=\result@t%
    \v@leur=\repdecn@mb{\v@lmax}\v@leur\edef\AsB@{\repdecn@mb{\v@leur}}% a/b
    \c@lAngle{#1}{#4}{\v@lmin}\edef\@ngdeb{\repdecn@mb{\v@lmin}}%
    \c@lAngle{#1}{#5}{\v@lmax}\ifdim\v@lmin>\v@lmax\advance\v@lmax\DePI@deg\fi%
    \edef\@ngfin{\repdecn@mb{\v@lmax}}\psarcellPA#1,#2,#3(\@ngdeb,\@ngfin)%
    \PSc@mment{End psarcellPP}\resetc@ntr@l\et@tpsarcellPP\fi\fi}}
\def\c@lAngle#1#2#3{\figvectP-3[#1,#2]%
    \c@lproscal\delt@[-3,-1]\c@lproscal\v@leur[-3,-2]%
    \v@leur=\AsB@\v@leur\arct@n#3(\delt@,\v@leur)#3=\rdT@deg#3}
\newif\if@rrowratio\@rrowratiotrue
\newif\if@rrowhfill
\newif\if@rrowhout
\def\Psset@rrowhe@d#1=#2|{\keln@mun#1|%
    \def\n@mref{a}\ifx\l@debut\n@mref\pssetarrowheadangle{#2}\else% angle
    \def\n@mref{f}\ifx\l@debut\n@mref\pssetarrowheadfill{#2}\else% fillmode
    \def\n@mref{l}\ifx\l@debut\n@mref\pssetarrowheadlength{#2}\else% length
    \def\n@mref{o}\ifx\l@debut\n@mref\pssetarrowheadout{#2}\else% out
    \def\n@mref{r}\ifx\l@debut\n@mref\pssetarrowheadratio{#2}\else% ratio
    \immediate\write16{*** Unknown attribute: \BS@ psset arrowhead(..., #1=...)}%
    \fi\fi\fi\fi\fi}
\def\pssetarrowheadangle#1{\edef\@rrowheadangle{#1}{\c@ssin{\C@}{\S@}{#1}%
    \xdef\C@AHANG{\C@}\xdef\S@AHANG{\S@}\v@lmax=\S@ pt%
    \invers@{\v@leur}{\v@lmax}\maxim@m{\v@leur}{\v@leur}{-\v@leur}%
    \xdef\UNSS@N{\the\v@leur}}}
\def\pssetarrowheadfill#1{\expandafter\set@rrowhfill#1:}
\def\set@rrowhfill#1#2:{\if#1n\@rrowhfillfalse\else\@rrowhfilltrue\fi}
\def\pssetarrowheadout#1{\expandafter\set@rrowhout#1:}
\def\set@rrowhout#1#2:{\if#1n\@rrowhoutfalse\else\@rrowhouttrue\fi}
\def\pssetarrowheadlength#1{\edef\@rrowheadlength{#1}\@rrowratiofalse}
\def\pssetarrowheadratio#1{\edef\@rrowheadratio{#1}\@rrowratiotrue}
\def\psresetarrowhead{%
    \pssetarrowheadangle{\defaultarrowheadangle}%
    \pssetarrowheadfill{\defaultarrowheadfill}%
    \pssetarrowheadout{\defaultarrowheadout}%
    \pssetarrowheadratio{\defaultarrowheadratio}%
    \d@fm@cdim\defaultarrowheadlength{\defaulth@rdahlength}% Valeur par defaut...
    \pssetarrowheadlength{\defaultarrowheadlength}}
\def\defaultarrowheadratio{0.1}
\def\defaultarrowheadangle{20}
\def\defaultarrowheadfill{no}
\def\defaultarrowheadout{no}
\def\defaulth@rdahlength{8pt}
\def\psarrowDD[#1,#2]{{\ifcurr@ntPS\ifps@cri\s@uvc@ntr@l\et@tpsarrow%
    \PSc@mment{psarrowDD [Pt1,Pt2]=[#1,#2]}\pssetfillmode{no}%
    \psarrowheadDD[#1,#2]\setc@ntr@l{2}\psline[#1,-3]%
    \PSc@mment{End psarrowDD}\resetc@ntr@l\et@tpsarrow\fi\fi}}
\def\psarrowTD[#1,#2]{{\ifcurr@ntPS\ifps@cri\s@uvc@ntr@l\et@tpsarrowTD%
    \PSc@mment{psarrowTD [Pt1,Pt2]=[#1,#2]}\resetc@ntr@l{2}%
    \Figptpr@j-5:/#1/\Figptpr@j-6:/#2/\let\c@lprojSP=\relax\psarrowDD[-5,-6]%
    \PSc@mment{End psarrowTD}\resetc@ntr@l\et@tpsarrowTD\fi\fi}}
\def\psarrowheadDD[#1,#2]{{\ifcurr@ntPS\ifps@cri\s@uvc@ntr@l\et@tpsarrowheadDD%
    \if@rrowhfill\def\@hangle{-\@rrowheadangle}\else\def\@hangle{\@rrowheadangle}\fi%
    \if@rrowratio%
    \if@rrowhout\def\@hratio{-\@rrowheadratio}\else\def\@hratio{\@rrowheadratio}\fi%
    \PSc@mment{psarrowheadDD Ratio=\@hratio, Angle=\@hangle, [Pt1,Pt2]=[#1,#2]}%
    \Ps@rrowhead\@hratio,\@hangle[#1,#2]%
    \else%
    \if@rrowhout\def\@hlength{-\@rrowheadlength}\else\def\@hlength{\@rrowheadlength}\fi%
    \PSc@mment{psarrowheadDD Length=\@hlength, Angle=\@hangle, [Pt1,Pt2]=[#1,#2]}%
    \Ps@rrowheadfd\@hlength,\@hangle[#1,#2]%
    \fi%
    \PSc@mment{End psarrowheadDD}\resetc@ntr@l\et@tpsarrowheadDD\fi\fi}}
\def\psarrowheadTD[#1,#2]{{\ifcurr@ntPS\ifps@cri\s@uvc@ntr@l\et@tpsarrowheadTD%
    \PSc@mment{psarrowheadTD [Pt1,Pt2]=[#1,#2]}\resetc@ntr@l{2}%
    \Figptpr@j-5:/#1/\Figptpr@j-6:/#2/\let\c@lprojSP=\relax\psarrowheadDD[-5,-6]%
    \PSc@mment{End psarrowheadTD}\resetc@ntr@l\et@tpsarrowheadTD\fi\fi}}
\def\Ps@rrowhead#1,#2[#3,#4]{\v@leur=#1\p@\maxim@m{\v@leur}{\v@leur}{-\v@leur}%
    \ifdim\v@leur>\Cepsil@n{% Arrow is not degenerated
    \PSc@mment{ps@rrowhead Ratio=#1, Angle=#2, [Pt1,Pt2]=[#3,#4]}\v@leur=\UNSS@N%
    \v@leur=\curr@ntwidth\v@leur\v@leur=\ptpsT@pt\v@leur\delt@=.5\v@leur% = width / (2 sin(Angle))
    \setc@ntr@l{2}\figvectPDD-3[#4,#3]%
    \Figg@tXY{-3}\v@lX=#1\v@lX\v@lY=#1\v@lY\Figv@ctCreg-3(\v@lX,\v@lY)%
    \vecunit@{-4}{-3}\mili@u=\result@t%
    \ifdim#2pt>\z@\v@lXa=-\C@AHANG\delt@%
     \edef\c@ef{\repdecn@mb{\v@lXa}}\figpttraDD-3:=-3/\c@ef,-4/\fi%
    \edef\c@ef{\repdecn@mb{\delt@}}%
    \v@lXa=\mili@u\v@lXa=\C@AHANG\v@lXa%
    \v@lYa=\ptpsT@pt\p@\v@lYa=\curr@ntwidth\v@lYa\v@lYa=\sDcc@ngle\v@lYa%
    \advance\v@lXa-\v@lYa\gdef\sDcc@ngle{0}%
    \ifdim\v@lXa>\v@leur\edef\c@efendpt{\repdecn@mb{\v@leur}}%
    \else\edef\c@efendpt{\repdecn@mb{\v@lXa}}\fi%
    \Figg@tXY{-3}\v@lmin=\v@lX\v@lmax=\v@lY%
    \v@lXa=\C@AHANG\v@lmin\v@lYa=\S@AHANG\v@lmax\advance\v@lXa\v@lYa%
    \v@lYa=-\S@AHANG\v@lmin\v@lX=\C@AHANG\v@lmax\advance\v@lYa\v@lX%
    \setc@ntr@l{1}\Figg@tXY{#4}\advance\v@lX\v@lXa\advance\v@lY\v@lYa%
    \setc@ntr@l{2}\Figp@intregDD-2:(\v@lX,\v@lY)%
    \v@lXa=\C@AHANG\v@lmin\v@lYa=-\S@AHANG\v@lmax\advance\v@lXa\v@lYa%
    \v@lYa=\S@AHANG\v@lmin\v@lX=\C@AHANG\v@lmax\advance\v@lYa\v@lX%
    \setc@ntr@l{1}\Figg@tXY{#4}\advance\v@lX\v@lXa\advance\v@lY\v@lYa%
    \setc@ntr@l{2}\Figp@intregDD-1:(\v@lX,\v@lY)%
    \ifdim#2pt<\z@\fillm@detrue\psline[-2,#4,-1]% fill
    \else\figptstraDD-3=#4,-2,-1/\c@ef,-4/\psline[-2,-3,-1]\fi% no fill
    \ifdim#1pt>\z@\figpttraDD-3:=#4/\c@efendpt,-4/\else\figptcopyDD-3:/#4/\fi%
    \PSc@mment{End ps@rrowhead}}\fi}
\def\sDcc@ngle{0}% Initialisation
\def\Ps@rrowheadfd#1,#2[#3,#4]{{%
    \PSc@mment{ps@rrowheadfd Length=#1, Angle=#2, [Pt1,Pt2]=[#3,#4]}%
    \setc@ntr@l{2}\figvectPDD-1[#3,#4]\n@rmeucDD{\v@leur}{-1}\v@leur=\ptT@unit@\v@leur%
    \invers@{\v@leur}{\v@leur}\v@leur=#1\v@leur\edef\R@tio{\repdecn@mb{\v@leur}}%
    \Ps@rrowhead\R@tio,#2[#3,#4]\PSc@mment{End ps@rrowheadfd}}}
\def\psarrowBezierDD[#1,#2,#3,#4]{{\ifcurr@ntPS\ifps@cri\s@uvc@ntr@l\et@tpsarrowBezierDD%
    \PSc@mment{psarrowBezierDD Control points=#1,#2,#3,#4}\setc@ntr@l{2}%
    \if@rrowratio\c@larclengthDD\v@leur,10[#1,#2,#3,#4]\else\v@leur=\z@\fi%
    \Ps@rrowB@zDD\v@leur[#1,#2,#3,#4]%
    \PSc@mment{End psarrowBezierDD}\resetc@ntr@l\et@tpsarrowBezierDD\fi\fi}}
\def\psarrowBezierTD[#1,#2,#3,#4]{{\ifcurr@ntPS\ifps@cri\s@uvc@ntr@l\et@tpsarrowBezierTD%
    \PSc@mment{psarrowBezierTD Control points=#1,#2,#3,#4}\resetc@ntr@l{2}%
    \Figptpr@j-7:/#1/\Figptpr@j-8:/#2/\Figptpr@j-9:/#3/\Figptpr@j-10:/#4/%
    \let\c@lprojSP=\relax\ifnum\curr@ntproj<\tw@\psarrowBezierDD[-7,-8,-9,-10]%
    \else\f@gnewpath\PSwrit@cmd{-7}{\c@mmoveto}{\fwf@g}%
    \if@rrowratio\c@larclengthDD\mili@u,10[-7,-8,-9,-10]\else\mili@u=\z@\fi%
    \p@rtent=\NBz@rcs\advance\p@rtent\m@ne\subB@zierTD\p@rtent[#1,#2,#3,#4]%
    \f@gstroke%
    \advance\v@lmin\p@rtent\delt@% Initialized in \subB@zierTD
    \v@leur=\v@lmin\advance\v@leur0.33333 \delt@\edef\unti@rs{\repdecn@mb{\v@leur}}%
    \v@leur=\v@lmin\advance\v@leur0.66666 \delt@\edef\deti@rs{\repdecn@mb{\v@leur}}%
    \figptcopyDD-8:/-10/\c@lsubBzarc\unti@rs,\deti@rs[#1,#2,#3,#4]%
    \figptcopyDD-8:/-4/\figptcopyDD-9:/-3/\Ps@rrowB@zDD\mili@u[-7,-8,-9,-10]\fi%
    \PSc@mment{End psarrowBezierTD}\resetc@ntr@l\et@tpsarrowBezierTD\fi\fi}}
\def\c@larclengthDD#1,#2[#3,#4,#5,#6]{{\p@rtent=#2\figptcopyDD-5:/#3/%
    \delt@=\p@\divide\delt@\p@rtent\c@rre=\z@\v@leur=\z@\s@mme=\z@%
    \loop\ifnum\s@mme<\p@rtent\advance\s@mme\@ne\advance\v@leur\delt@%
    \edef\T@{\repdecn@mb{\v@leur}}\figptBezierDD-6::\T@[#3,#4,#5,#6]%
    \figvectPDD-1[-5,-6]\n@rmeucDD{\mili@u}{-1}\advance\c@rre\mili@u%
    \figptcopyDD-5:/-6/\repeat\global\result@t=\ptT@unit@\c@rre}#1=\result@t}
\def\Ps@rrowB@zDD#1[#2,#3,#4,#5]{{\pssetfillmode{no}%
    \if@rrowratio\delt@=\@rrowheadratio#1\else\delt@=\@rrowheadlength pt\fi%
    \v@leur=\C@AHANG\delt@\edef\R@dius{\repdecn@mb{\v@leur}}%
    \FigptintercircB@zDD-5::0,\R@dius[#5,#4,#3,#2]%
    \pssetarrowheadlength{\repdecn@mb{\delt@}}\psarrowheadDD[-5,#5]%
    \let\n@rmeuc=\n@rmeucDD\figgetdist\R@dius[#5,-3]%
    \FigptintercircB@zDD-6::0,\R@dius[#5,#4,#3,#2]%
    \figptBezierDD-5::0.33333[#5,#4,#3,#2]\figptBezierDD-3::0.66666[#5,#4,#3,#2]%
    \figptscontrolDD-5[-6,-5,-3,#2]\psBezierDD1[-6,-5,-4,#2]}}
\def\psarrowcircDD#1;#2(#3,#4){{\ifcurr@ntPS\ifps@cri\s@uvc@ntr@l\et@tpsarrowcircDD%
    \PSc@mment{psarrowcircDD Center=#1 ; Radius=#2 (Ang1=#3,Ang2=#4)}%
    \pssetfillmode{no}\Pscirc@rrowhead#1;#2(#3,#4)%
    \setc@ntr@l{2}\figvectPDD -4[#1,-3]\vecunit@{-4}{-4}%
    \Figg@tXY{-4}\arct@n\v@lmin(\v@lX,\v@lY)%
    \v@lmin=\rdT@deg\v@lmin\v@leur=#4pt\advance\v@leur-\v@lmin%
    \maxim@m{\v@leur}{\v@leur}{-\v@leur}%
    \ifdim\v@leur>\DemiPI@deg\relax\ifdim\v@lmin<#4pt\advance\v@lmin\DePI@deg%
    \else\advance\v@lmin-\DePI@deg\fi\fi\edef\ar@ngle{\repdecn@mb{\v@lmin}}%
    \ifdim#3pt<#4pt\psarccirc#1;#2(#3,\ar@ngle)\else\psarccirc#1;#2(\ar@ngle,#3)\fi%
    \PSc@mment{End psarrowcircDD}\resetc@ntr@l\et@tpsarrowcircDD\fi\fi}}
\def\psarrowcircTD#1,#2,#3;#4(#5,#6){{\ifcurr@ntPS\ifps@cri\s@uvc@ntr@l\et@tpsarrowcircTD%
    \PSc@mment{psarrowcircTD Center=#1,P1=#2,P2=#3 ; Radius=#4 (Ang1=#5, Ang2=#6)}%
    \resetc@ntr@l{2}\c@lExtAxes#1,#2,#3(#4)\let\c@lprojSP=\relax%
    \figvectPTD-11[#1,-4]\figvectPTD-12[#1,-5]\c@lNbarcs{#5}{#6}%
    \if@rrowratio\v@lmax=\degT@rd\v@lmax\edef\D@lpha{\repdecn@mb{\v@lmax}}\fi%
    \advance\p@rtent\m@ne\mili@u=\z@%
    \v@leur=#5pt\c@lptellP{#1}{-11}{-12}\Figptpr@j-9:/-3/%
    \f@gnewpath\PSwrit@cmdS{-9}{\c@mmoveto}{\fwf@g}{\X@un}{\Y@un}%
    \edef\C@nt@r{#1}\s@mme=\z@\bcl@rcircTD\f@gstroke%
    \advance\v@leur\delt@\c@lptellP{#1}{-11}{-12}\Figptpr@j-5:/-3/%
    \advance\v@leur\delt@\c@lptellP{#1}{-11}{-12}\Figptpr@j-6:/-3/%
    \advance\v@leur\delt@\c@lptellP{#1}{-11}{-12}\Figptpr@j-10:/-3/%
    \figptscontrolDD-8[-9,-5,-6,-10]%
    \if@rrowratio\c@lcurvradDD0.5[-9,-8,-7,-10]\advance\mili@u\result@t%
    \maxim@m{\mili@u}{\mili@u}{-\mili@u}\mili@u=\ptT@unit@\mili@u%
    \mili@u=\D@lpha\mili@u\advance\p@rtent\@ne\divide\mili@u\p@rtent\fi%
    \Ps@rrowB@zDD\mili@u[-9,-8,-7,-10]%
    \PSc@mment{End psarrowcircTD}\resetc@ntr@l\et@tpsarrowcircTD\fi\fi}}
\def\bcl@rcircTD{\relax%
    \ifnum\s@mme<\p@rtent\advance\s@mme\@ne%
    \advance\v@leur\delt@\c@lptellP{\C@nt@r}{-11}{-12}\Figptpr@j-5:/-3/%
    \advance\v@leur\delt@\c@lptellP{\C@nt@r}{-11}{-12}\Figptpr@j-6:/-3/%
    \advance\v@leur\delt@\c@lptellP{\C@nt@r}{-11}{-12}\Figptpr@j-10:/-3/%
    \figptscontrolDD-8[-9,-5,-6,-10]\BdingB@xfalse%
    \PSwrit@cmdS{-8}{}{\fwf@g}{\X@de}{\Y@de}\PSwrit@cmdS{-7}{}{\fwf@g}{\X@tr}{\Y@tr}%
    \BdingB@xtrue\PSwrit@cmdS{-10}{\c@mcurveto}{\fwf@g}{\X@qu}{\Y@qu}%
    \if@rrowratio\c@lcurvradDD0.5[-9,-8,-7,-10]\advance\mili@u\result@t\fi%
    \B@zierBB@x{1}{\Y@un}(\X@un,\X@de,\X@tr,\X@qu)%
    \B@zierBB@x{2}{\X@un}(\Y@un,\Y@de,\Y@tr,\Y@qu)%
    \edef\X@un{\X@qu}\edef\Y@un{\Y@qu}\figptcopyDD-9:/-10/\bcl@rcircTD\fi}
\def\Pscirc@rrowhead#1;#2(#3,#4){{%
    \PSc@mment{pscirc@rrowhead Center=#1 ; Radius=#2 (Ang1=#3,Ang2=#4)}%
    \v@leur=#2\unit@\edef\s@glen{\repdecn@mb{\v@leur}}\v@lY=\z@\v@lX=\v@leur%
    \resetc@ntr@l{2}\Figv@ctCreg-3(\v@lX,\v@lY)\figpttraDD-5:=#1/1,-3/%
    \figptrotDD-5:=-5/#1,#4/%
    \figvectPDD-3[#1,-5]\Figg@tXY{-3}\v@leur=\v@lX%
    \ifdim#3pt<#4pt\v@lX=\v@lY\v@lY=-\v@leur\else\v@lX=-\v@lY\v@lY=\v@leur\fi%
    \Figv@ctCreg-3(\v@lX,\v@lY)\vecunit@{-3}{-3}%
    \if@rrowratio\v@leur=#4pt\advance\v@leur-#3pt\maxim@m{\mili@u}{-\v@leur}{\v@leur}%
    \mili@u=\degT@rd\mili@u\v@leur=\s@glen\mili@u\edef\s@glen{\repdecn@mb{\v@leur}}%
    \mili@u=#2\mili@u\mili@u=\@rrowheadratio\mili@u\else\mili@u=\@rrowheadlength pt\fi%
    \figpttraDD-6:=-5/\s@glen,-3/\v@leur=#2pt\v@leur=2\v@leur%
    \invers@{\v@leur}{\v@leur}\c@rre=\repdecn@mb{\v@leur}\mili@u% = sin = L/(2R)
    \mili@u=\c@rre\mili@u=\repdecn@mb{\c@rre}\mili@u%
    \v@leur=\p@\advance\v@leur-\mili@u% \v@leur = cos*cos
    \invers@{\mili@u}{2\v@leur}\delt@=\c@rre\delt@=\repdecn@mb{\mili@u}\delt@%
    \xdef\sDcc@ngle{\repdecn@mb{\delt@}}% sin/(2*cos*cos) used in \Ps@rrowhead
    \sqrt@{\mili@u}{\v@leur}\arct@n\v@leur(\mili@u,\c@rre)%
    \v@leur=\rdT@deg\v@leur% \cor@ngle = atan(L/sqrt(4R*R-L*L))
    \ifdim#3pt<#4pt\v@leur=-\v@leur\fi%
    \if@rrowhout\v@leur=-\v@leur\fi\edef\cor@ngle{\repdecn@mb{\v@leur}}%
    \figptrotDD-6:=-6/-5,\cor@ngle/\psarrowheadDD[-6,-5]%
    \PSc@mment{End pscirc@rrowhead}}}
\def\psarrowcircPDD#1;#2[#3,#4]{{\ifcurr@ntPS\ifps@cri%
    \PSc@mment{psarrowcircPDD Center=#1; Radius=#2, [P1=#3,P2=#4]}%
    \s@uvc@ntr@l\et@tpsarrowcircPDD\Ps@ngleparam#1;#2[#3,#4]%
    \ifdim\v@leur>\z@\ifdim\v@lmin>\v@lmax\advance\v@lmax\DePI@deg\fi%
    \else\ifdim\v@lmin<\v@lmax\advance\v@lmin\DePI@deg\fi\fi%
    \edef\@ngdeb{\repdecn@mb{\v@lmin}}\edef\@ngfin{\repdecn@mb{\v@lmax}}%
    \psarrowcirc#1;\r@dius(\@ngdeb,\@ngfin)%
    \PSc@mment{End psarrowcircPDD}\resetc@ntr@l\et@tpsarrowcircPDD\fi\fi}}
\def\psarrowcircPTD#1;#2[#3,#4,#5]{{\ifcurr@ntPS\ifps@cri\s@uvc@ntr@l\et@tpsarrowcircPTD%
    \PSc@mment{psarrowcircPTD Center=#1; Radius=#2, [P1=#3,P2=#4,P3=#5]}%
    \figgetangleTD\@ngfin[#1,#3,#4,#5]\v@leur=#2pt%
    \maxim@m{\mili@u}{-\v@leur}{\v@leur}\edef\r@dius{\repdecn@mb{\mili@u}}%
    \ifdim\v@leur<\z@\v@lmax=\@ngfin pt\advance\v@lmax-\DePI@deg%
    \edef\@ngfin{\repdecn@mb{\v@lmax}}\fi\psarrowcircTD#1,#3,#5;\r@dius(0,\@ngfin)%
    \PSc@mment{End psarrowcircPTD}\resetc@ntr@l\et@tpsarrowcircPTD\fi\fi}}
\def\psaxes#1(#2){{\ifcurr@ntPS\ifps@cri\s@uvc@ntr@l\et@tpsaxes%
    \PSc@mment{psaxes Origin=#1 Range=(#2)}\an@lys@xes#2,:\resetc@ntr@l{2}%
    \ifx\t@xt@\empty\ifTr@isDim\ps@xes#1(0,#2,0,#2,0,#2)\else\ps@xes#1(0,#2,0,#2)\fi%
    \else\ps@xes#1(#2)\fi\PSc@mment{End psaxes}\resetc@ntr@l\et@tpsaxes\fi\fi}}
\def\an@lys@xes#1,#2:{\def\t@xt@{#2}}
\def\ps@xesDD#1(#2,#3,#4,#5){%
    \figpttraC-5:=#1/#2,0/\figpttraC-6:=#1/#3,0/\psarrowDD[-5,-6]%
    \figpttraC-5:=#1/0,#4/\figpttraC-6:=#1/0,#5/\psarrowDD[-5,-6]}
\def\ps@xesTD#1(#2,#3,#4,#5,#6,#7){%
    \figpttraC-7:=#1/#2,0,0/\figpttraC-8:=#1/#3,0,0/\psarrowTD[-7,-8]%
    \figpttraC-7:=#1/0,#4,0/\figpttraC-8:=#1/0,#5,0/\psarrowTD[-7,-8]%
    \figpttraC-7:=#1/0,0,#6/\figpttraC-8:=#1/0,0,#7/\psarrowTD[-7,-8]}
\edef\DefGIfilen@me{\jobname GI.anx}
\def\psbeginfig#1{\def\t@xt@{#1}\relax\ifx\t@xt@\empty\Psb@ginfig\DefGIfilen@me%
    \else\expandafter\Psb@ginfigNu@#1 :\fi}
\def\Psb@ginfigNu@#1 #2:{\def\t@xt@{#1}\relax\ifx\t@xt@\empty\def\t@xt@{#2}%
    \ifx\t@xt@\empty\Psb@ginfig\DefGIfilen@me\else\Psb@ginfigNu@#2:\fi%
    \else\Psb@ginfig{#1}\fi}
\def\Psb@ginfig#1{\ifcurr@ntPS\else%
    \edef\PSfilen@me{#1}\edef\auxfilen@me{\jobname.anx}%
    \ifpstestm@de\ps@critrue\else\openin\frf@g=\PSfilen@me\relax%
    \ifeof\frf@g\ps@critrue\else\ps@crifalse\fi\closein\frf@g\fi%
    \curr@ntPStrue\c@ldefproj\expandafter\setupd@te\defaultupdate:%
    \ifps@cri\initb@undb@x%
    \immediate\openout\fwf@g=\auxfilen@me\initpss@ttings\fi%
    \fi}
\def\initpss@ttings{\psreset{arrowhead,curve,first,flowchart,mesh,second,third}%
    \Use@llipsefalse}
\def\B@zierBB@x#1#2(#3,#4,#5,#6){{\c@rre=\t@n\epsil@n% Do not reduce this value
    \v@lmax=#4\advance\v@lmax-#5\v@lmax=\thr@@\v@lmax\advance\v@lmax#6\advance\v@lmax-#3%
    \mili@u=#4\mili@u=-\tw@\mili@u\advance\mili@u#3\advance\mili@u#5%
    \v@lmin=#4\advance\v@lmin-#3\maxim@m{\v@leur}{-\v@lmax}{\v@lmax}%
    \maxim@m{\delt@}{-\mili@u}{\mili@u}\maxim@m{\v@leur}{\v@leur}{\delt@}%
    \maxim@m{\delt@}{-\v@lmin}{\v@lmin}\maxim@m{\v@leur}{\v@leur}{\delt@}%
    \ifdim\v@leur>\c@rre\invers@{\v@leur}{\v@leur}\edef\Uns@rM@x{\repdecn@mb{\v@leur}}%
    \v@lmax=\Uns@rM@x\v@lmax\mili@u=\Uns@rM@x\mili@u\v@lmin=\Uns@rM@x\v@lmin%
    \maxim@m{\v@leur}{-\v@lmax}{\v@lmax}\ifdim\v@leur<\c@rre%
    \maxim@m{\v@leur}{-\mili@u}{\mili@u}\ifdim\v@leur<\c@rre\else%
    \invers@{\mili@u}{\mili@u}\v@leur=-0.5\v@lmin%
    \v@leur=\repdecn@mb{\mili@u}\v@leur\m@jBBB@x{\v@leur}{#1}{#2}(#3,#4,#5,#6)\fi%
    \else\delt@=\repdecn@mb{\mili@u}\mili@u\v@leur=\repdecn@mb{\v@lmax}\v@lmin%
    \advance\delt@-\v@leur\ifdim\delt@<\z@\else\invers@{\v@lmax}{\v@lmax}%
    \edef\Uns@rAp{\repdecn@mb{\v@lmax}}\sqrt@{\delt@}{\delt@}%
    \v@leur=-\mili@u\advance\v@leur\delt@\v@leur=\Uns@rAp\v@leur%
    \m@jBBB@x{\v@leur}{#1}{#2}(#3,#4,#5,#6)%
    \v@leur=-\mili@u\advance\v@leur-\delt@\v@leur=\Uns@rAp\v@leur%
    \m@jBBB@x{\v@leur}{#1}{#2}(#3,#4,#5,#6)\fi\fi\fi}}
\def\m@jBBB@x#1#2#3(#4,#5,#6,#7){{\relax\ifdim#1>\z@\ifdim#1<\p@%
    \edef\T@{\repdecn@mb{#1}}\v@lX=\p@\advance\v@lX-#1\edef\UNmT@{\repdecn@mb{\v@lX}}%
    \v@lX=#4\v@lY=#5\v@lZ=#6\v@lXa=#7\v@lX=\UNmT@\v@lX\advance\v@lX\T@\v@lY%
    \v@lY=\UNmT@\v@lY\advance\v@lY\T@\v@lZ\v@lZ=\UNmT@\v@lZ\advance\v@lZ\T@\v@lXa%
    \v@lX=\UNmT@\v@lX\advance\v@lX\T@\v@lY\v@lY=\UNmT@\v@lY\advance\v@lY\T@\v@lZ%
    \v@lX=\UNmT@\v@lX\advance\v@lX\T@\v@lY%
    \ifcase#2\or\v@lY=#3\or\v@lY=\v@lX\v@lX=#3\fi\b@undb@x{\v@lX}{\v@lY}\fi\fi}}
\def\PsB@zier#1[#2]{{\f@gnewpath%
    \s@mme=\z@\def\list@num{#2,0}\extrairelepremi@r\p@int\de\list@num%
    \PSwrit@cmdS{\p@int}{\c@mmoveto}{\fwf@g}{\X@un}{\Y@un}\p@rtent=#1\bclB@zier}}
\def\bclB@zier{\relax%
    \ifnum\s@mme<\p@rtent\advance\s@mme\@ne\BdingB@xfalse%
    \extrairelepremi@r\p@int\de\list@num\PSwrit@cmdS{\p@int}{}{\fwf@g}{\X@de}{\Y@de}%
    \extrairelepremi@r\p@int\de\list@num\PSwrit@cmdS{\p@int}{}{\fwf@g}{\X@tr}{\Y@tr}%
    \BdingB@xtrue%
    \extrairelepremi@r\p@int\de\list@num\PSwrit@cmdS{\p@int}{\c@mcurveto}{\fwf@g}{\X@qu}{\Y@qu}%
    \B@zierBB@x{1}{\Y@un}(\X@un,\X@de,\X@tr,\X@qu)%
    \B@zierBB@x{2}{\X@un}(\Y@un,\Y@de,\Y@tr,\Y@qu)%
    \edef\X@un{\X@qu}\edef\Y@un{\Y@qu}\bclB@zier\fi}
\def\psBezierDD#1[#2]{\ifcurr@ntPS\ifps@cri%
    \PSc@mment{psBezierDD N arcs=#1, Control points=#2}%
    \iffillm@de\PsB@zier#1[#2]%
    \f@gfill%
    \else\PsB@zier#1[#2]\f@gstroke\fi%
    \PSc@mment{End psBezierDD}\fi\fi}
\def\psBezierTD#1[#2]{\ifcurr@ntPS\ifps@cri\s@uvc@ntr@l\et@tpsBezierTD%
    \PSc@mment{psBezierTD N arcs=#1, Control points=#2}%
    \iffillm@de\PsB@zierTD#1[#2]%
    \f@gfill%
    \else\PsB@zierTD#1[#2]\f@gstroke\fi%
    \PSc@mment{End psBezierTD}\resetc@ntr@l\et@tpsBezierTD\fi\fi}
\def\PsB@zierTD#1[#2]{\ifnum\curr@ntproj<\tw@\PsB@zier#1[#2]\else\PsB@zier@TD#1[#2]\fi}
\def\PsB@zier@TD#1[#2]{{\f@gnewpath%
    \s@mme=\z@\def\list@num{#2,0}\extrairelepremi@r\p@int\de\list@num%
    \let\c@lprojSP=\relax\setc@ntr@l{2}\Figptpr@j-7:/\p@int/%
    \PSwrit@cmd{-7}{\c@mmoveto}{\fwf@g}%
    \loop\ifnum\s@mme<#1\advance\s@mme\@ne\extrairelepremi@r\p@intun\de\list@num%
    \extrairelepremi@r\p@intde\de\list@num\extrairelepremi@r\p@inttr\de\list@num%
    \subB@zierTD\NBz@rcs[\p@int,\p@intun,\p@intde,\p@inttr]\edef\p@int{\p@inttr}\repeat}}
\def\subB@zierTD#1[#2,#3,#4,#5]{\delt@=\p@\divide\delt@\NBz@rcs\v@lmin=\z@%
    {\Figg@tXY{-7}\edef\X@un{\the\v@lX}\edef\Y@un{\the\v@lY}%
    \s@mme=\z@\loop\ifnum\s@mme<#1\advance\s@mme\@ne%
    \v@leur=\v@lmin\advance\v@leur0.33333 \delt@\edef\unti@rs{\repdecn@mb{\v@leur}}%
    \v@leur=\v@lmin\advance\v@leur0.66666 \delt@\edef\deti@rs{\repdecn@mb{\v@leur}}%
    \advance\v@lmin\delt@\edef\trti@rs{\repdecn@mb{\v@lmin}}%
    \figptBezierTD-8::\trti@rs[#2,#3,#4,#5]\Figptpr@j-8:/-8/%
    \c@lsubBzarc\unti@rs,\deti@rs[#2,#3,#4,#5]\BdingB@xfalse%
    \PSwrit@cmdS{-4}{}{\fwf@g}{\X@de}{\Y@de}\PSwrit@cmdS{-3}{}{\fwf@g}{\X@tr}{\Y@tr}%
    \BdingB@xtrue\PSwrit@cmdS{-8}{\c@mcurveto}{\fwf@g}{\X@qu}{\Y@qu}%
    \B@zierBB@x{1}{\Y@un}(\X@un,\X@de,\X@tr,\X@qu)%
    \B@zierBB@x{2}{\X@un}(\Y@un,\Y@de,\Y@tr,\Y@qu)%
    \edef\X@un{\X@qu}\edef\Y@un{\Y@qu}\figptcopyDD-7:/-8/\repeat}}
\def\NBz@rcs{2}
\def\c@lsubBzarc#1,#2[#3,#4,#5,#6]{\figptBezierTD-5::#1[#3,#4,#5,#6]%
    \figptBezierTD-6::#2[#3,#4,#5,#6]\Figptpr@j-4:/-5/\Figptpr@j-5:/-6/%
    \figptscontrolDD-4[-7,-4,-5,-8]}
\def\pscircDD#1(#2){\ifcurr@ntPS\ifps@cri\PSc@mment{pscircDD Center=#1 (Radius=#2)}%
    \psarccircDD#1;#2(0,360)\PSc@mment{End pscircDD}\fi\fi}
\def\pscircTD#1,#2,#3(#4){\ifcurr@ntPS\ifps@cri%
    \PSc@mment{pscircTD Center=#1,P1=#2,P2=#3 (Radius=#4)}%
    \psarccircTD#1,#2,#3;#4(0,360)\PSc@mment{End pscircTD}\fi\fi}
{\catcode`\%=12\gdef\p@urcent{%}}
\def\PSc@mment#1{\ifpsdebugmode\immediate\write\fwf@g{\p@urcent\space#1}\fi}
{\catcode`\[=1\catcode`\{=12\gdef\acc@louv[{}}
{\catcode`\]=2\catcode`\}=12\gdef\acc@lfer{}]]
\def\PSdict@{\ifUse@llipse%
    \immediate\write\fwf@g{/ellipsedict 9 dict def ellipsedict /mtrx matrix put}%
    \immediate\write\fwf@g{/ellipse \acc@louv ellipsedict begin}%
    \immediate\write\fwf@g{ /endangle exch def /startangle exch def}%
    \immediate\write\fwf@g{ /yrad exch def /xrad exch def}%
    \immediate\write\fwf@g{ /rotangle exch def /y exch def /x exch def}%
    \immediate\write\fwf@g{ /savematrix mtrx currentmatrix def}%
    \immediate\write\fwf@g{ x y translate rotangle rotate xrad yrad scale}%
    \immediate\write\fwf@g{ 0 0 1 startangle endangle arc}%
    \immediate\write\fwf@g{ savematrix setmatrix end\acc@lfer def}%
    \fi\PShe@der{EndProlog}}
\def\Pssetc@rve#1=#2|{\keln@mun#1|%
    \def\n@mref{r}\ifx\l@debut\n@mref\pssetroundness{#2}\else% roundness
    \immediate\write16{*** Unknown attribute: \BS@ psset curve(..., #1=...)}%
    \fi}
\def\pssetroundness#1{\edef\curv@roundness{#1}}
\def\defaultroundness{0.2} % Valeur par defaut
\def\pscurveDD[#1]{{\ifcurr@ntPS\ifps@cri\PSc@mment{pscurveDD Points=#1}%
    \s@uvc@ntr@l\et@tpscurveDD%
    \iffillm@de\Psc@rveDD\curv@roundness[#1]%
    \f@gfill%
    \else\Psc@rveDD\curv@roundness[#1]\f@gstroke\fi%
    \PSc@mment{End pscurveDD}\resetc@ntr@l\et@tpscurveDD\fi\fi}}
\def\pscurveTD[#1]{{\ifcurr@ntPS\ifps@cri%
    \PSc@mment{pscurveTD Points=#1}\s@uvc@ntr@l\et@tpscurveTD\let\c@lprojSP=\relax%
    \iffillm@de\Psc@rveTD\curv@roundness[#1]%
    \f@gfill%
    \else\Psc@rveTD\curv@roundness[#1]\f@gstroke\fi%
    \PSc@mment{End pscurveTD}\resetc@ntr@l\et@tpscurveTD\fi\fi}}
\def\Psc@rveDD#1[#2]{%
    \def\list@num{#2}\extrairelepremi@r\Ak@\de\list@num%
    \extrairelepremi@r\Ai@\de\list@num\extrairelepremi@r\Aj@\de\list@num%
    \f@gnewpath\PSwrit@cmdS{\Ai@}{\c@mmoveto}{\fwf@g}{\X@un}{\Y@un}%
    \setc@ntr@l{2}\figvectPDD -1[\Ak@,\Aj@]%
    \@ecfor\Ak@:=\list@num\do{\figpttraDD-2:=\Ai@/#1,-1/\BdingB@xfalse%
       \PSwrit@cmdS{-2}{}{\fwf@g}{\X@de}{\Y@de}%
       \figvectPDD -1[\Ai@,\Ak@]\figpttraDD-2:=\Aj@/-#1,-1/%
       \PSwrit@cmdS{-2}{}{\fwf@g}{\X@tr}{\Y@tr}\BdingB@xtrue%
       \PSwrit@cmdS{\Aj@}{\c@mcurveto}{\fwf@g}{\X@qu}{\Y@qu}%
       \B@zierBB@x{1}{\Y@un}(\X@un,\X@de,\X@tr,\X@qu)%
       \B@zierBB@x{2}{\X@un}(\Y@un,\Y@de,\Y@tr,\Y@qu)%
       \edef\X@un{\X@qu}\edef\Y@un{\Y@qu}\edef\Ai@{\Aj@}\edef\Aj@{\Ak@}}}
\def\Psc@rveTD#1[#2]{\ifnum\curr@ntproj<\tw@\Psc@rvePPTD#1[#2]\else\Psc@rveCPTD#1[#2]\fi}
\def\Psc@rvePPTD#1[#2]{\setc@ntr@l{2}%
    \def\list@num{#2}\extrairelepremi@r\Ak@\de\list@num\Figptpr@j-5:/\Ak@/%
    \extrairelepremi@r\Ai@\de\list@num\Figptpr@j-3:/\Ai@/%
    \extrairelepremi@r\Aj@\de\list@num\Figptpr@j-4:/\Aj@/%
    \f@gnewpath\PSwrit@cmdS{-3}{\c@mmoveto}{\fwf@g}{\X@un}{\Y@un}%
    \figvectPDD -1[-5,-4]%
    \@ecfor\Ak@:=\list@num\do{\Figptpr@j-5:/\Ak@/\figpttraDD-2:=-3/#1,-1/%
       \BdingB@xfalse\PSwrit@cmdS{-2}{}{\fwf@g}{\X@de}{\Y@de}%
       \figvectPDD -1[-3,-5]\figpttraDD-2:=-4/-#1,-1/%
       \PSwrit@cmdS{-2}{}{\fwf@g}{\X@tr}{\Y@tr}\BdingB@xtrue%
       \PSwrit@cmdS{-4}{\c@mcurveto}{\fwf@g}{\X@qu}{\Y@qu}%
       \B@zierBB@x{1}{\Y@un}(\X@un,\X@de,\X@tr,\X@qu)%
       \B@zierBB@x{2}{\X@un}(\Y@un,\Y@de,\Y@tr,\Y@qu)%
       \edef\X@un{\X@qu}\edef\Y@un{\Y@qu}\figptcopyDD-3:/-4/\figptcopyDD-4:/-5/}}
\def\Psc@rveCPTD#1[#2]{\setc@ntr@l{2}%
    \def\list@num{#2}\extrairelepremi@r\Ak@\de\list@num%
    \extrairelepremi@r\Ai@\de\list@num\extrairelepremi@r\Aj@\de\list@num%
    \Figptpr@j-7:/\Ai@/%
    \f@gnewpath\PSwrit@cmd{-7}{\c@mmoveto}{\fwf@g}%
    \figvectPTD -9[\Ak@,\Aj@]%
    \@ecfor\Ak@:=\list@num\do{\figpttraTD-10:=\Ai@/#1,-9/%
       \figvectPTD -9[\Ai@,\Ak@]\figpttraTD-11:=\Aj@/-#1,-9/%
       \subB@zierTD\NBz@rcs[\Ai@,-10,-11,\Aj@]\edef\Ai@{\Aj@}\edef\Aj@{\Ak@}}}
\def\psendfig{\ifcurr@ntPS\ifps@cri\immediate\closeout\fwf@g%
    \immediate\openout\fwf@g=\PSfilen@me\relax%
    \ifPDFm@ke\PSBdingB@x\else%
    \immediate\write\fwf@g{\p@urcent\string!PS-Adobe-2.0 EPSF-2.0}%
    \PShe@der{Creator\string: TeX (fig4tex.tex)}%
    \PShe@der{Title\string: \PSfilen@me}%
    \PShe@der{CreationDate\string: \the\day/\the\month/\the\year}%
    \PSBdingB@x%
    \PShe@der{EndComments}\PSdict@\fi%
    \immediate\write\fwf@g{\c@mgsave}%
    \openin\frf@g=\auxfilen@me\c@pypsfile\fwf@g\frf@g\closein\frf@g%
    \immediate\write\fwf@g{\c@mgrestore}%
    \PSc@mment{End of file.}\immediate\closeout\fwf@g%
    \immediate\openout\fwf@g=\auxfilen@me\immediate\closeout\fwf@g%
    \immediate\write16{File \PSfilen@me\space created.}\fi\fi\curr@ntPSfalse\ps@critrue}
\def\PShe@der#1{\immediate\write\fwf@g{\p@urcent\p@urcent#1}}
\def\PSBdingB@x{{\v@lX=\ptT@ptps\c@@rdXmin\v@lY=\ptT@ptps\c@@rdYmin%
     \v@lXa=\ptT@ptps\c@@rdXmax\v@lYa=\ptT@ptps\c@@rdYmax%
     \PShe@der{BoundingBox\string: \repdecn@mb{\v@lX}\space\repdecn@mb{\v@lY}%
     \space\repdecn@mb{\v@lXa}\space\repdecn@mb{\v@lYa}}}}
\def\psfcconnect[#1]{{\ifcurr@ntPS\ifps@cri\PSc@mment{psfcconnect Points=#1}%
    \pssetfillmode{no}\s@uvc@ntr@l\et@tpsfcconnect\resetc@ntr@l{2}%
    \fcc@nnect@[#1]\resetc@ntr@l\et@tpsfcconnect\PSc@mment{End psfcconnect}\fi\fi}}
\def\fcc@nnect@[#1]{\let\N@rm=\n@rmeucDD\def\list@num{#1}%
    \extrairelepremi@r\Ai@\de\list@num\edef\pr@m{\Ai@}\v@leur=\z@\p@rtent=\@ne\c@llgtot%
    \ifcase\fclin@typ@\edef\list@num{[\pr@m,#1,\Ai@}\expandafter\pscurve\list@num]%
    \else\ifdim\fclin@r@d\p@>\z@\Pslin@conge[#1]\else\psline[#1]\fi\fi%
    \v@leur=\@rrowp@s\v@leur\edef\list@num{#1,\Ai@,0}%
    \extrairelepremi@r\Ai@\de\list@num\mili@u=\epsil@n\c@llgpart%
    \advance\mili@u-\epsil@n\advance\mili@u-\delt@\advance\v@leur-\mili@u%
    \ifcase\fclin@typ@\invers@\mili@u\delt@%
    \ifnum\@rrowr@fpt>\z@\advance\delt@-\v@leur\v@leur=\delt@\fi%
    \v@leur=\repdecn@mb\v@leur\mili@u\edef\v@lt{\repdecn@mb\v@leur}%
    \extrairelepremi@r\Ak@\de\list@num%
    \figvectPDD-1[\pr@m,\Aj@]\figpttraDD-6:=\Ai@/\curv@roundness,-1/%
    \figvectPDD-1[\Ak@,\Ai@]\figpttraDD-7:=\Aj@/\curv@roundness,-1/%
    \delt@=\@rrowheadlength\p@\delt@=\C@AHANG\delt@\edef\R@dius{\repdecn@mb{\delt@}}%
    \ifcase\@rrowr@fpt%
    \FigptintercircB@zDD-8::\v@lt,\R@dius[\Ai@,-6,-7,\Aj@]\psarrowheadDD[-5,-8]\else%
    \FigptintercircB@zDD-8::\v@lt,\R@dius[\Aj@,-7,-6,\Ai@]\psarrowheadDD[-8,-5]\fi%
    \else\advance\delt@-\v@leur%
    \p@rtentiere{\p@rtent}{\delt@}\edef\C@efun{\the\p@rtent}%
    \p@rtentiere{\p@rtent}{\v@leur}\edef\C@efde{\the\p@rtent}%
    \figptbaryDD-5:[\Ai@,\Aj@;\C@efun,\C@efde]\ifcase\@rrowr@fpt%
    \delt@=\@rrowheadlength\unit@\delt@=\C@AHANG\delt@\edef\t@ille{\repdecn@mb{\delt@}}%
    \figvectPDD-2[\Ai@,\Aj@]\vecunit@{-2}{-2}\figpttraDD-5:=-5/\t@ille,-2/\fi%
    \psarrowheadDD[\Ai@,-5]\fi}
\def\c@llgtot{\@ecfor\Aj@:=\list@num\do{\figvectP-1[\Ai@,\Aj@]\N@rm\delt@{-1}%
    \advance\v@leur\delt@\advance\p@rtent\@ne\edef\Ai@{\Aj@}}}
\def\c@llgpart{\extrairelepremi@r\Aj@\de\list@num\figvectP-1[\Ai@,\Aj@]\N@rm\delt@{-1}%
    \advance\mili@u\delt@\ifdim\mili@u<\v@leur\edef\pr@m{\Ai@}\edef\Ai@{\Aj@}\c@llgpart\fi}
\def\Pslin@conge[#1]{\ifnum\p@rtent>\tw@{\def\list@num{#1}%
    \extrairelepremi@r\Ai@\de\list@num\extrairelepremi@r\Aj@\de\list@num%
    \figptcopy-6:/\Ai@/\figvectP-3[\Ai@,\Aj@]\vecunit@{-3}{-3}\v@lmax=\result@t%
    \@ecfor\Ak@:=\list@num\do{\figvectP-4[\Aj@,\Ak@]\vecunit@{-4}{-4}%
    \minim@m\v@lmin\v@lmax\result@t\v@lmax=\result@t%
    \det@rm\delt@[-3,-4]\maxim@m\mili@u{\delt@}{-\delt@}\ifdim\mili@u>\Cepsil@n%
    \ifdim\delt@>\z@\figgetangleDD\Angl@[\Aj@,\Ak@,\Ai@]\else%
    \figgetangleDD\Angl@[\Aj@,\Ai@,\Ak@]\fi%
    \v@leur=\PI@deg\advance\v@leur-\Angl@\p@\divide\v@leur\tw@%
    \edef\Angl@{\repdecn@mb\v@leur}\c@ssin{\C@}{\S@}{\Angl@}\v@leur=\fclin@r@d\unit@%
    \v@leur=\S@\v@leur\mili@u=\C@\p@\invers@\mili@u\mili@u%
    \v@leur=\repdecn@mb{\mili@u}\v@leur%
    \minim@m\v@leur\v@leur\v@lmin\edef\t@ille{\repdecn@mb{\v@leur}}%
    \figpttra-5:=\Aj@/-\t@ille,-3/\psline[-6,-5]\figpttra-6:=\Aj@/\t@ille,-4/%
    \figvectNVDD-3[-3]\figvectNVDD-8[-4]\inters@cDD-7:[-5,-3;-6,-8]%
    \ifdim\delt@>\z@\psarccircP-7;\fclin@r@d[-5,-6]\else\psarccircP-7;\fclin@r@d[-6,-5]\fi%
    \else\psline[-6,\Aj@]\figptcopy-6:/\Aj@/\fi% Points alignes
    \edef\Ai@{\Aj@}\edef\Aj@{\Ak@}\figptcopy-3:/-4/}\psline[-6,\Aj@]}\else\psline[#1]\fi}
\def\psfcnode[#1]#2{{\ifcurr@ntPS\ifps@cri\PSc@mment{psfcnode Points=#1}%
    \s@uvc@ntr@l\et@tpsfcnode\resetc@ntr@l{2}%
    \def\t@xt@{#2}\ifx\t@xt@\empty\def\g@tt@xt{\setbox\Gb@x=\hbox{\Figg@tT{\p@int}}}%
    \else\def\g@tt@xt{\setbox\Gb@x=\hbox{#2}}\fi%
    \v@lmin=\h@rdfcXp@dd\advance\v@lmin\Xp@dd\unit@\multiply\v@lmin\tw@%
    \v@lmax=\h@rdfcYp@dd\advance\v@lmax\Yp@dd\unit@\multiply\v@lmax\tw@%
    \Figv@ctCreg-8(\unit@,-\unit@)\def\list@num{#1}%
    \delt@=\curr@ntwidth bp\divide\delt@\tw@%
    \fcn@de\PSc@mment{End psfcnode}\resetc@ntr@l\et@tpsfcnode\fi\fi}}
\def\d@butn@de{\g@tt@xt\v@lX=\wd\Gb@x%
    \v@lY=\ht\Gb@x\advance\v@lY\dp\Gb@x\advance\v@lX\v@lmin\advance\v@lY\v@lmax}
\def\fcn@deE{%
    \@ecfor\p@int:=\list@num\do{\d@butn@de\v@lX=\unssqrttw@\v@lX\v@lY=\unssqrttw@\v@lY%
    \ifdim\thickn@ss\p@>\z@% Shadow
    \v@lXa=\v@lX\advance\v@lXa\delt@\v@lXa=\ptT@unit@\v@lXa\edef\XR@d{\repdecn@mb\v@lXa}%
    \v@lYa=\v@lY\advance\v@lYa\delt@\v@lYa=\ptT@unit@\v@lYa\edef\YR@d{\repdecn@mb\v@lYa}%
    \arct@n\v@leur(\v@lXa,\v@lYa)\v@leur=\rdT@deg\v@leur\edef\@nglde{\repdecn@mb\v@leur}%
    {\c@lptellDD-2::\p@int;\XR@d,\YR@d(\@nglde)}% \v@lmin & \v@lmax modified in \c@lptellDD
    \advance\v@leur-\PI@deg\edef\@nglun{\repdecn@mb\v@leur}%
    {\c@lptellDD-3::\p@int;\XR@d,\YR@d(\@nglun)}%
    \figptstra-6=-3,-2,\p@int/\thickn@ss,-8/\pssetfillmode{yes}\us@secondC@lor%
    \psline[-2,-3,-6,-5]\psarcell-4;\XR@d,\YR@d(\@nglun,\@nglde,0)\fi% End shadow
    \v@lX=\ptT@unit@\v@lX\v@lY=\ptT@unit@\v@lY%
    \edef\XR@d{\repdecn@mb\v@lX}\edef\YR@d{\repdecn@mb\v@lY}%
    \pssetfillmode{yes}\us@thirdC@lor\psarcell\p@int;\XR@d,\YR@d(0,360,0)%
    \pssetfillmode{no}\us@primarC@lor\psarcell\p@int;\XR@d,\YR@d(0,360,0)}}
\def\fcn@deL{\delt@=\ptT@unit@\delt@\edef\t@ille{\repdecn@mb\delt@}%
    \@ecfor\p@int:=\list@num\do{\Figg@tXYa{\p@int}\d@butn@de%
    \ifdim\v@lX>\v@lY\itis@Ktrue\else\itis@Kfalse\fi%
    \advance\v@lXa-\v@lX\Figp@intreg-1:(\v@lXa,\v@lYa)%
    \advance\v@lXa\v@lX\advance\v@lYa-\v@lY\Figp@intreg-2:(\v@lXa,\v@lYa)%
    \advance\v@lXa\v@lX\advance\v@lYa\v@lY\Figp@intreg-3:(\v@lXa,\v@lYa)%
    \advance\v@lXa-\v@lX\advance\v@lYa\v@lY\Figp@intreg-4:(\v@lXa,\v@lYa)%
    \ifdim\thickn@ss\p@>\z@\Figg@tXYa{\p@int}\pssetfillmode{yes}\us@secondC@lor% Shadow
    \c@lpt@xt{-1}{-4}\c@lpt@xt@\v@lXa\v@lYa\v@lX\v@lY\c@rre\delt@%
    \Figp@intregDD-9:(\v@lZ,\v@lYa)\Figp@intregDD-11:(\v@lZa,\v@lYa)%
    \c@lpt@xt{-4}{-3}\c@lpt@xt@\v@lYa\v@lXa\v@lY\v@lX\delt@\c@rre%
    \Figp@intregDD-12:(\v@lXa,\v@lZ)\Figp@intregDD-10:(\v@lXa,\v@lZa)%
    \ifitis@K\figptstra-7=-9,-10,-11/\thickn@ss,-8/\psline[-9,-11,-5,-6,-7]\else%
    \figptstra-7=-10,-11,-12/\thickn@ss,-8/\psline[-10,-12,-5,-6,-7]\fi\fi% End shadow
    \pssetfillmode{yes}\us@thirdC@lor\psline[-1,-2,-3,-4]%
    \pssetfillmode{no}\us@primarC@lor\psline[-1,-2,-3,-4,-1]}}
\def\c@lpt@xt#1#2{\figvectN-7[#1,#2]\vecunit@{-7}{-7}\figpttra-5:=#1/\t@ille,-7/%
    \figvectP-7[#1,#2]\Figg@tXY{-7}\c@rre=\v@lX\delt@=\v@lY\Figg@tXY{-5}}
\def\c@lpt@xt@#1#2#3#4#5#6{\v@lZ=#6\invers@{\v@lZ}{\v@lZ}\v@leur=\repdecn@mb{#5}\v@lZ%
    \v@lZ=#2\advance\v@lZ-#4\mili@u=\repdecn@mb{\v@leur}\v@lZ%
    \v@lZ=#3\advance\v@lZ\mili@u\v@lZa=-\v@lZ\advance\v@lZa\tw@#1}
\def\fcn@deR{\@ecfor\p@int:=\list@num\do{\Figg@tXYa{\p@int}\d@butn@de%
    \advance\v@lXa-0.5\v@lX\advance\v@lYa-0.5\v@lY\Figp@intreg-1:(\v@lXa,\v@lYa)%
    \advance\v@lXa\v@lX\Figp@intreg-2:(\v@lXa,\v@lYa)%
    \advance\v@lYa\v@lY\Figp@intreg-3:(\v@lXa,\v@lYa)%
    \advance\v@lXa-\v@lX\Figp@intreg-4:(\v@lXa,\v@lYa)%
    \ifdim\thickn@ss\p@>\z@\pssetfillmode{yes}\us@secondC@lor% Shadow
    \Figv@ctCreg-5(-\delt@,-\delt@)\figpttra-9:=-1/1,-5/%
    \Figv@ctCreg-5(\delt@,-\delt@)\figpttra-10:=-2/1,-5/%
    \Figv@ctCreg-5(\delt@,\delt@)\figpttra-11:=-3/1,-5/%
    \figptstra-7=-9,-10,-11/\thickn@ss,-8/\psline[-9,-11,-5,-6,-7]\fi% End shadow
    \pssetfillmode{yes}\us@thirdC@lor\psline[-1,-2,-3,-4]%
    \pssetfillmode{no}\us@primarC@lor\psline[-1,-2,-3,-4,-1]}}
\def\Pssetfl@wchart#1=#2|{\keln@mtr#1|%
    \def\n@mref{arr}\ifx\l@debut\n@mref\expandafter\keln@mtr\l@suite|%
     \def\n@mref{owp}\ifx\l@debut\n@mref\edef\@rrowp@s{#2}\else% arrowposition
     \def\n@mref{owr}\ifx\l@debut\n@mref\setfcr@fpt#2|\else% arrowrefpt
     \immediate\write16{*** Unknown attribute: \BS@ psset flowchart(..., #1=...)}%
     \fi\fi\else%
    \def\n@mref{lin}\ifx\l@debut\n@mref\setfccurv@#2|\else% line
    \def\n@mref{pad}\ifx\l@debut\n@mref\edef\Xp@dd{#2}\edef\Yp@dd{#2}\else% padding
    \def\n@mref{rad}\ifx\l@debut\n@mref\edef\fclin@r@d{#2}\else% connection radius
    \def\n@mref{sha}\ifx\l@debut\n@mref\setfcshap@#2|\else% shape
    \def\n@mref{thi}\ifx\l@debut\n@mref\edef\thickn@ss{#2}\else% thickness
    \def\n@mref{xpa}\ifx\l@debut\n@mref\edef\Xp@dd{#2}\else% xpadding
    \def\n@mref{ypa}\ifx\l@debut\n@mref\edef\Yp@dd{#2}\else% ypadding
    \immediate\write16{*** Unknown attribute: \BS@ psset flowchart(..., #1=...)}%
    \fi\fi\fi\fi\fi\fi\fi\fi}
\def\setfcr@fpt#1#2|{\if#1e\def\@rrowr@fpt{1}\else\def\@rrowr@fpt{0}\fi}
\def\setfccurv@#1#2|{\if#1c\def\fclin@typ@{0}\else\def\fclin@typ@{1}\fi}
\def\setfcshap@#1#2|{%
    \if#1e\let\fcn@de=\fcn@deE\def\h@rdfcXp@dd{4pt}\def\h@rdfcYp@dd{4pt}%
     \edef\fcsh@pe{ellipse}\else%
    \if#1l\let\fcn@de=\fcn@deL\def\h@rdfcXp@dd{4pt}\def\h@rdfcYp@dd{4pt}%
     \edef\fcsh@pe{lozenge}\else%
          \let\fcn@de=\fcn@deR\def\h@rdfcXp@dd{6pt}\def\h@rdfcYp@dd{6pt}%
     \edef\fcsh@pe{rectangle}\fi\fi}
\def\psline[#1]{{\ifcurr@ntPS\ifps@cri\PSc@mment{psline Points=#1}%
    \let\pslign@=\Pslign@P\Pslin@{#1}\PSc@mment{End psline}\fi\fi}}
\def\pslineF#1{{\ifcurr@ntPS\ifps@cri\PSc@mment{pslineF Filename=#1}%
    \let\pslign@=\Pslign@F\Pslin@{#1}\PSc@mment{End pslineF}\fi\fi}}
\def\pslineC(#1){{\ifcurr@ntPS\ifps@cri\PSc@mment{pslineC}%
    \let\pslign@=\Pslign@C\Pslin@{#1}\PSc@mment{End pslineC}\fi\fi}}
\def\Pslin@#1{\iffillm@de\pslign@{#1}%
    \f@gfill%
    \else\pslign@{#1}\ifx\derp@int\premp@int%
    \f@gclosestroke%
    \else\f@gstroke\fi\fi}
\def\Pslign@P#1{\def\list@num{#1}\extrairelepremi@r\p@int\de\list@num%
    \edef\premp@int{\p@int}\f@gnewpath%
    \PSwrit@cmd{\p@int}{\c@mmoveto}{\fwf@g}%
    \@ecfor\p@int:=\list@num\do{\PSwrit@cmd{\p@int}{\c@mlineto}{\fwf@g}%
    \edef\derp@int{\p@int}}}
\def\Pslign@F#1{\s@uvc@ntr@l\et@tPslign@F\setc@ntr@l{2}\openin\frf@g=#1\relax%
    \ifeof\frf@g\message{*** File #1 not found !}\end\else%
    \read\frf@g to\tr@c\edef\premp@int{\tr@c}\expandafter\extr@ctCF\tr@c:%
    \f@gnewpath\PSwrit@cmd{-1}{\c@mmoveto}{\fwf@g}%
    \loop\read\frf@g to\tr@c\ifeof\frf@g\mored@tafalse\else\mored@tatrue\fi%
    \ifmored@ta\expandafter\extr@ctCF\tr@c:\PSwrit@cmd{-1}{\c@mlineto}{\fwf@g}%
    \edef\derp@int{\tr@c}\repeat\fi\closein\frf@g\resetc@ntr@l\et@tPslign@F}
\def\extr@ctCFDD#1 #2:{\v@lX=#1\unit@\v@lY=#2\unit@\Figp@intregDD-1:(\v@lX,\v@lY)}
\def\extr@ctCFTD#1 #2 #3:{\v@lX=#1\unit@\v@lY=#2\unit@\v@lZ=#3\unit@%
    \Figp@intregTD-1:(\v@lX,\v@lY,\v@lZ)}
\def\Pslign@C#1{\s@uvc@ntr@l\et@tPslign@C\setc@ntr@l{2}%
    \def\list@num{#1}\extrairelepremi@r\p@int\de\list@num%
    \edef\premp@int{\p@int}\f@gnewpath%
    \expandafter\Pslign@C@\p@int:\PSwrit@cmd{-1}{\c@mmoveto}{\fwf@g}%
    \@ecfor\p@int:=\list@num\do{\expandafter\Pslign@C@\p@int:%
    \PSwrit@cmd{-1}{\c@mlineto}{\fwf@g}\edef\derp@int{\p@int}}%
    \resetc@ntr@l\et@tPslign@C}
\def\Pslign@C@#1 #2:{{\def\t@xt@{#1}\ifx\t@xt@\empty\Pslign@C@#2:% Discard leading spaces
    \else\extr@ctCF#1 #2:\fi}}
\newcount\c@ntrolmesh
\def\Pssetm@sh#1=#2|{\keln@mun#1|%
    \def\n@mref{d}\ifx\l@debut\n@mref\pssetmeshdiag{#2}\else% diag
    \immediate\write16{*** Unknown attribute: \BS@ psset mesh(..., #1=...)}%
    \fi}
\def\pssetmeshdiag#1{\c@ntrolmesh=#1}
\def\defaultmeshdiag{0}    % Valeur par defaut
\def\psmesh#1,#2[#3,#4,#5,#6]{{\ifcurr@ntPS\ifps@cri%
    \PSc@mment{psmesh N1=#1, N2=#2, Quadrangle=[#3,#4,#5,#6]}%
    \s@uvc@ntr@l\et@tpsmesh\Pss@tsecondSt\setc@ntr@l{2}%
    \ifnum#1>\@ne\Psmeshp@rt#1[#3,#4,#5,#6]\fi%
    \ifnum#2>\@ne\Psmeshp@rt#2[#4,#5,#6,#3]\fi%
    \ifnum\c@ntrolmesh>\z@\Psmeshdi@g#1,#2[#3,#4,#5,#6]\fi%
    \ifnum\c@ntrolmesh<\z@\Psmeshdi@g#2,#1[#4,#5,#6,#3]\fi\Psrest@reSt%
    \psline[#3,#4,#5,#6,#3]\PSc@mment{End psmesh}\resetc@ntr@l\et@tpsmesh\fi\fi}}
\def\Psmeshp@rt#1[#2,#3,#4,#5]{{\l@mbd@un=\@ne\l@mbd@de=#1\loop%
    \ifnum\l@mbd@un<#1\advance\l@mbd@de\m@ne\figptbary-1:[#2,#3;\l@mbd@de,\l@mbd@un]%
    \figptbary-2:[#5,#4;\l@mbd@de,\l@mbd@un]\psline[-1,-2]\advance\l@mbd@un\@ne\repeat}}
\def\Psmeshdi@g#1,#2[#3,#4,#5,#6]{\figptcopy-2:/#3/\figptcopy-3:/#6/%
    \l@mbd@un=\z@\l@mbd@de=#1\loop\ifnum\l@mbd@un<#1%
    \advance\l@mbd@un\@ne\advance\l@mbd@de\m@ne\figptcopy-1:/-2/\figptcopy-4:/-3/%
    \figptbary-2:[#3,#4;\l@mbd@de,\l@mbd@un]%
    \figptbary-3:[#6,#5;\l@mbd@de,\l@mbd@un]\Psmeshdi@gp@rt#2[-1,-2,-3,-4]\repeat}
\def\Psmeshdi@gp@rt#1[#2,#3,#4,#5]{{\l@mbd@un=\z@\l@mbd@de=#1\loop%
    \ifnum\l@mbd@un<#1\figptbary-5:[#2,#5;\l@mbd@de,\l@mbd@un]%
    \advance\l@mbd@de\m@ne\advance\l@mbd@un\@ne%
    \figptbary-6:[#3,#4;\l@mbd@de,\l@mbd@un]\psline[-5,-6]\repeat}}
\def\psnormalDD#1,#2[#3,#4]{{\ifcurr@ntPS\ifps@cri%
    \PSc@mment{psnormal Length=#1, Lambda=#2 [Pt1,Pt2]=[#3,#4]}%
    \s@uvc@ntr@l\et@tpsnormal\resetc@ntr@l{2}\figptendnormal-6::#1,#2[#3,#4]%
    \figptcopyDD-5:/-1/\psarrow[-5,-6]%
    \PSc@mment{End psnormal}\resetc@ntr@l\et@tpsnormal\fi\fi}}
\def\psreset#1{\trtlis@rg{#1}{\Psreset@}}
\def\Psreset@#1|{\keln@mde#1|%
    \def\n@mref{ar}\ifx\l@debut\n@mref\psresetarrowhead\else% arrowhead
    \def\n@mref{cu}\ifx\l@debut\n@mref\psset curve(roundness=\defaultroundness)\else% curve
    \def\n@mref{fi}\ifx\l@debut\n@mref\psset (color=\defaultcolor,dash=\defaultdash,%
         fill=\defaultfill,join=\defaultjoin,width=\defaultwidth)\else% primary settings
    \def\n@mref{fl}\ifx\l@debut\n@mref\psset flowchart(arrowp=\defaultfcarrowposition,%
	arrowr=\defaultfcarrowrefpt,line=\defaultfcline,xpadd=\defaultfcxpadding,%
	ypadd=\defaultfcypadding,radius=\defaultfcradius,shape=\defaultfcshape,%
	thick=\defaultfcthickness)\else% flow chart
    \def\n@mref{me}\ifx\l@debut\n@mref\psset mesh(diag=\defaultmeshdiag)\else% mesh
    \def\n@mref{se}\ifx\l@debut\n@mref\psresetsecondsettings\else% secondary
    \def\n@mref{th}\ifx\l@debut\n@mref\psset third(color=\defaultthirdcolor)\else% ternary
    \immediate\write16{*** Unknown keyword #1 (\BS@ psreset).}%
    \fi\fi\fi\fi\fi\fi\fi}
\def\psset#1(#2){\def\t@xt@{#1}\ifx\t@xt@\empty\trtlis@rg{#2}{\Pssetf@rst}% primary settings
    \else\keln@mde#1|%
    \def\n@mref{ar}\ifx\l@debut\n@mref\trtlis@rg{#2}{\Psset@rrowhe@d}\else% arrow-head
    \def\n@mref{cu}\ifx\l@debut\n@mref\trtlis@rg{#2}{\Pssetc@rve}\else% curve
    \def\n@mref{fi}\ifx\l@debut\n@mref\trtlis@rg{#2}{\Pssetf@rst}\else% primary settings
    \def\n@mref{fl}\ifx\l@debut\n@mref\trtlis@rg{#2}{\Pssetfl@wchart}\else% flow chart
    \def\n@mref{me}\ifx\l@debut\n@mref\trtlis@rg{#2}{\Pssetm@sh}\else% mesh
    \def\n@mref{se}\ifx\l@debut\n@mref\trtlis@rg{#2}{\Pssets@cond}\else% secondary settings
    \def\n@mref{th}\ifx\l@debut\n@mref\trtlis@rg{#2}{\Pssetth@rd}\else% ternary settings
    \immediate\write16{*** Unknown keyword: \BS@ psset #1(...)}%
    \fi\fi\fi\fi\fi\fi\fi\fi}
\def\pssetdefault#1(#2){\ifcurr@ntPS\immediate\write16{*** \BS@ pssetdefault is ignored
    inside a \BS@ psbeginfig-\BS@ psendfig block.}%
    \immediate\write16{*** It must be called before \BS@ psbeginfig.}\else%
    \def\t@xt@{#1}\ifx\t@xt@\empty\trtlis@rg{#2}{\Pssd@f@rst}\else\keln@mde#1|%
    \def\n@mref{ar}\ifx\l@debut\n@mref\trtlis@rg{#2}{\Pssd@@rrowhe@d}\else% arrow-head
    \def\n@mref{cu}\ifx\l@debut\n@mref\trtlis@rg{#2}{\Pssd@c@rve}\else% curve
    \def\n@mref{fi}\ifx\l@debut\n@mref\trtlis@rg{#2}{\Pssd@f@rst}\else% primary settings
    \def\n@mref{fl}\ifx\l@debut\n@mref\trtlis@rg{#2}{\Pssd@fl@wchart}\else% flow chart
    \def\n@mref{me}\ifx\l@debut\n@mref\trtlis@rg{#2}{\Pssd@m@sh}\else% mesh
    \def\n@mref{se}\ifx\l@debut\n@mref\trtlis@rg{#2}{\Pssd@s@cond}\else% secondary settings
    \def\n@mref{th}\ifx\l@debut\n@mref\trtlis@rg{#2}{\Pssd@th@rd}\else% ternary settings
    \immediate\write16{*** Unknown keyword: \BS@ pssetdefault #1(...)}%
    \fi\fi\fi\fi\fi\fi\fi\fi\initpss@ttings\fi}
\def\Pssd@f@rst#1=#2|{\keln@mun#1|%
    \def\n@mref{c}\ifx\l@debut\n@mref\edef\defaultcolor{#2}\else% color
    \def\n@mref{d}\ifx\l@debut\n@mref\edef\defaultdash{#2}\else% dash
    \def\n@mref{f}\ifx\l@debut\n@mref\edef\defaultfill{#2}\else% fillmode
    \def\n@mref{j}\ifx\l@debut\n@mref\edef\defaultjoin{#2}\else% line join
    \def\n@mref{u}\ifx\l@debut\n@mref\edef\defaultupdate{#2}\pssetupdate{#2}\else% update
    \def\n@mref{w}\ifx\l@debut\n@mref\edef\defaultwidth{#2}\else% line width
    \immediate\write16{*** Unknown attribute: \BS@ pssetdefault (..., #1=...)}%
    \fi\fi\fi\fi\fi\fi}
\def\Pssd@@rrowhe@d#1=#2|{\keln@mun#1|%
    \def\n@mref{a}\ifx\l@debut\n@mref\edef\defaultarrowheadangle{#2}\else% angle
    \def\n@mref{f}\ifx\l@debut\n@mref\edef\defaultarrowheadangle{#2}\else% fillmode
    \def\n@mref{l}\ifx\l@debut\n@mref\y@tiunit{#2}\ifunitpr@sent%
     \edef\defaulth@rdahlength{#2}\else\edef\defaulth@rdahlength{#2pt}%
     \message{*** \BS@ pssetdefault (..., #1=#2, ...) : unit is missing, pt is assumed.}%
     \fi\else% length
    \def\n@mref{o}\ifx\l@debut\n@mref\edef\defaultarrowheadout{#2}\else% out
    \def\n@mref{r}\ifx\l@debut\n@mref\edef\defaultarrowheadratio{#2}\else% ratio
    \immediate\write16{*** Unknown attribute: \BS@ pssetdefault arrowhead(..., #1=...)}%
    \fi\fi\fi\fi\fi}
\def\Pssd@c@rve#1=#2|{\keln@mun#1|%
    \def\n@mref{r}\ifx\l@debut\n@mref\edef\defaultroundness{#2}\else%
    \immediate\write16{*** Unknown attribute: \BS@ pssetdefault curve(..., #1=...)}%
    \fi}
\def\Pssd@fl@wchart#1=#2|{\keln@mtr#1|%
    \def\n@mref{arr}\ifx\l@debut\n@mref\expandafter\keln@mtr\l@suite|%
     \def\n@mref{owp}\ifx\l@debut\n@mref\edef\defaultfcarrowposition{#2}\else% arrowposition
     \def\n@mref{owr}\ifx\l@debut\n@mref\edef\defaultfcarrowrefpt{#2}\else% arrowrefpt
     \immediate\write16{*** Unknown attribute: \BS@ pssetdefault flowchart(..., #1=...)}%
     \fi\fi\else%
    \def\n@mref{lin}\ifx\l@debut\n@mref\edef\defaultfcline{#2}\else% line
    \def\n@mref{pad}\ifx\l@debut\n@mref\edef\defaultfcxpadding{#2}%
                    \edef\defaultfcypadding{#2}\else% padding
    \def\n@mref{rad}\ifx\l@debut\n@mref\edef\defaultfcradius{#2}\else% connection radius
    \def\n@mref{sha}\ifx\l@debut\n@mref\edef\defaultfcshape{#2}\else% shape
    \def\n@mref{thi}\ifx\l@debut\n@mref\edef\defaultfcthickness{#2}\else% thickness
    \def\n@mref{xpa}\ifx\l@debut\n@mref\edef\defaultfcxpadding{#2}\else% xpadding
    \def\n@mref{ypa}\ifx\l@debut\n@mref\edef\defaultfcypadding{#2}\else% ypadding
    \immediate\write16{*** Unknown attribute: \BS@ pssetdefault flowchart(..., #1=...)}%
    \fi\fi\fi\fi\fi\fi\fi\fi}
\def\defaultfcarrowposition{0.5}\let\defaultfcarrowpos=\defaultfcarrowposition
\def\defaultfcarrowrefpt{start}
\def\defaultfcline{polygon}
\def\defaultfcradius{0}
\def\defaultfcshape{rectangle}
\def\defaultfcthickness{0}\let\defaultfcthick=\defaultfcthickness
\def\defaultfcxpadding{0}\let\defaultfcxpad=\defaultfcxpadding
\def\defaultfcypadding{0}\let\defaultfcypad=\defaultfcypadding
\def\Pssd@m@sh#1=#2|{\keln@mun#1|%
    \def\n@mref{d}\ifx\l@debut\n@mref\edef\defaultmeshdiag{#2}\else%
    \immediate\write16{*** Unknown attribute: \BS@ pssetdefault mesh(..., #1=...)}%
    \fi}
\def\Pssd@s@cond#1=#2|{\keln@mun#1|%
    \def\n@mref{c}\ifx\l@debut\n@mref\edef\defaultsecondcolor{#2}\else%
    \def\n@mref{d}\ifx\l@debut\n@mref\edef\defaultseconddash{#2}\else%
    \def\n@mref{w}\ifx\l@debut\n@mref\edef\defaultsecondwidth{#2}\else%
    \immediate\write16{*** Unknown attribute: \BS@ pssetdefault second(..., #1=...)}%
    \fi\fi\fi}
\def\Pssd@th@rd#1=#2|{\keln@mun#1|%
    \def\n@mref{c}\ifx\l@debut\n@mref\edef\defaultthirdcolor{#2}\else%
    \immediate\write16{*** Unknown attribute: \BS@ pssetdefault third(..., #1=...)}%
    \fi}
\newif\iffillm@de
\def\pssetfillmode#1{\expandafter\setfillm@de#1:}
\def\setfillm@de#1#2:{\if#1n\fillm@defalse\else\fillm@detrue\fi}
\def\defaultfill{no}     % Valeur par defaut
\newif\ifpstestm@de
\def\pssetupdate#1{\ifcurr@ntPS\immediate\write16{*** \BS@ pssetupdate is ignored inside a
     \BS@ psbeginfig-\BS@ psendfig block.}%
    \immediate\write16{*** It must be called before \BS@ psbeginfig.}%
    \else\expandafter\setupd@te#1:\fi}
\def\setupd@te#1#2:{\if#1n\pstestm@defalse\else\pstestm@detrue\fi}
\def\defaultupdate{no}     % Valeur par defaut
\def\Pssetc@lor#1{\ifps@cri\result@tent=\@ne\expandafter\c@lnbV@l#1 :%
    \def\curr@ntcolor{}\def\curr@ntcolorc@md{}%
    \ifcase\result@tent\or\pssetgray{#1}\or\or\pssetrgb{#1}\or\pssetcmyk{#1}\fi\fi}
\def\pssetcmyk#1{\ifps@cri\def\curr@ntcolor{#1}\def\curr@ntcolorc@md{\c@msetcmykcolor}%
    \def\curr@ntcolorc@mdStroke{\c@msetcmykcolorStroke}%
    \ifcurr@ntPS\PSc@mment{pssetcmyk Color=#1}\us@primarC@lor\fi\fi}
\def\pssetrgb#1{\ifps@cri\def\curr@ntcolor{#1}\def\curr@ntcolorc@md{\c@msetrgbcolor}%
    \def\curr@ntcolorc@mdStroke{\c@msetrgbcolorStroke}%
    \ifcurr@ntPS\PSc@mment{pssetrgb Color=#1}\us@primarC@lor\fi\fi}
\def\pssetgray#1{\ifps@cri\def\curr@ntcolor{#1}\def\curr@ntcolorc@md{\c@msetgray}%
    \def\curr@ntcolorc@mdStroke{\c@msetgrayStroke}%
    \ifcurr@ntPS\PSc@mment{pssetgray Gray level=#1}\us@primarC@lor\fi\fi}
\def\us@primarC@lor{\immediate\write\fwf@g{\d@fprimarC@lor}%
    \let\fillc@md=\prfillc@md}
\def\prfillc@md{\d@fprimarC@lor\space\c@mfill}
\def\defaultcolor{0}       % Valeur par defaut
\def\c@lnbV@l#1 #2:{\def\t@xt@{#1}\relax\ifx\t@xt@\empty\c@lnbV@l#2:% Discard leading spaces
    \else\c@lnbV@l@#1 #2:\fi}
\def\c@lnbV@l@#1 #2:{\def\t@xt@{#2}\ifx\t@xt@\empty%
    \def\t@xt@{#1}\ifx\t@xt@\empty\advance\result@tent\m@ne\fi% Discard trailing spaces
    \else\advance\result@tent\@ne\c@lnbV@l@#2:\fi}
\def\Blackcmyk{0 0 0 1}
\def\Whitecmyk{0 0 0 0}
\def\Cyancmyk{1 0 0 0}
\def\Magentacmyk{0 1 0 0}
\def\Yellowcmyk{0 0 1 0}
\def\Redcmyk{0 1 1 0}
\def\Greencmyk{1 0 1 0}
\def\Bluecmyk{1 1 0 0}
\def\Graycmyk{0 0 0 0.50}
\def\BrickRedcmyk{0 0.89 0.94 0.28} % PANTONE 1805
\def\Browncmyk{0 0.81 1 0.60} % PANTONE 1615
\def\ForestGreencmyk{0.91 0 0.88 0.12} % PANTONE 349
\def\Goldenrodcmyk{ 0 0.10 0.84 0} % PANTONE 109
\def\Marooncmyk{0 0.87 0.68 0.32} % PANTONE 201
\def\Orangecmyk{0 0.61 0.87 0} % PANTONE ORANGE-021
\def\Purplecmyk{0.45 0.86 0 0} % PANTONE PURPLE
\def\RoyalBluecmyk{1. 0.50 0 0} % No PANTONE match
\def\Violetcmyk{0.79 0.88 0 0} % PANTONE VIOLET
\def\Blackrgb{0 0 0}
\def\Whitergb{1 1 1}
\def\Redrgb{1 0 0}
\def\Greenrgb{0 1 0}
\def\Bluergb{0 0 1}
\def\Cyanrgb{0 1 1}
\def\Magentargb{1 0 1}
\def\Yellowrgb{1 1 0}
\def\Grayrgb{0.5 0.5 0.5}
\def\Chocolatergb{0.824 0.412 0.118}
\def\DarkGoldenrodrgb{0.722 0.525 0.043}
\def\DarkOrangergb{1 0.549 0}
\def\Firebrickrgb{0.698 0.133 0.133}
\def\ForestGreenrgb{0.133 0.545 0.133}
\def\Goldrgb{1 0.843 0}
\def\HotPinkrgb{1 0.412 0.706}
\def\Maroonrgb{0.690 0.188 0.376}
\def\Pinkrgb{1 0.753 0.796}
\def\RoyalBluergb{0.255 0.412 0.882}
\def\Pssetf@rst#1=#2|{\keln@mun#1|%
    \def\n@mref{c}\ifx\l@debut\n@mref\Pssetc@lor{#2}\else% color
    \def\n@mref{d}\ifx\l@debut\n@mref\pssetdash{#2}\else% dash
    \def\n@mref{f}\ifx\l@debut\n@mref\pssetfillmode{#2}\else% fillmode
    \def\n@mref{j}\ifx\l@debut\n@mref\pssetjoin{#2}\else% line join
    \def\n@mref{u}\ifx\l@debut\n@mref\pssetupdate{#2}\else% update
    \def\n@mref{w}\ifx\l@debut\n@mref\pssetwidth{#2}\else% line width
    \immediate\write16{*** Unknown attribute: \BS@ psset (..., #1=...)}%
    \fi\fi\fi\fi\fi\fi}
\def\s@uvdash#1{\edef#1{\curr@ntdash}}
\def\defaultdash{1}        % Valeur par defaut
\def\pssetdash#1{\ifps@cri\edef\curr@ntdash{#1}\ifcurr@ntPS\expandafter\Pssetd@sh#1 :\fi\fi}
\def\Pssetd@shI#1{\PSc@mment{pssetdash Index=#1}\ifcase#1%
    \or\immediate\write\fwf@g{[] 0 \c@msetdash}%         Index=1
    \or\immediate\write\fwf@g{[6 2] 0 \c@msetdash}%      Index=2
    \or\immediate\write\fwf@g{[4 2] 0 \c@msetdash}%      Index=3
    \or\immediate\write\fwf@g{[2 2] 0 \c@msetdash}%      Index=4
    \or\immediate\write\fwf@g{[1 2] 0 \c@msetdash}%      Index=5
    \or\immediate\write\fwf@g{[2 4] 0 \c@msetdash}%      Index=6
    \or\immediate\write\fwf@g{[3 5] 0 \c@msetdash}%      Index=7
    \or\immediate\write\fwf@g{[3 3] 0 \c@msetdash}%      Index=8
    \or\immediate\write\fwf@g{[3 5 1 5] 0 \c@msetdash}%  Index=9
    \or\immediate\write\fwf@g{[6 4 2 4] 0 \c@msetdash}%  Index=10
    \fi}
\def\Pssetd@sh#1 #2:{{\def\t@xt@{#1}\ifx\t@xt@\empty\Pssetd@sh#2:% Discard leading spaces
    \else\def\t@xt@{#2}\ifx\t@xt@\empty\Pssetd@shI{#1}\else\s@mme=\@ne\def\debutp@t{#1}%
    \an@lysd@sh#2:\ifodd\s@mme\edef\debutp@t{\debutp@t\space\finp@t}\def\finp@t{0}\fi%
    \PSc@mment{pssetdash Pattern=#1 #2}%
    \immediate\write\fwf@g{[\debutp@t] \finp@t\space\c@msetdash}\fi\fi}}
\def\an@lysd@sh#1 #2:{\def\t@xt@{#2}\ifx\t@xt@\empty\def\finp@t{#1}\else%
    \edef\debutp@t{\debutp@t\space#1}\advance\s@mme\@ne\an@lysd@sh#2:\fi}
\def\s@uvwidth#1{\edef#1{\curr@ntwidth}}
\def\defaultwidth{0.4}     % Valeur par defaut
\def\pssetwidth#1{\ifps@cri\edef\curr@ntwidth{#1}\ifcurr@ntPS%
    \PSc@mment{pssetwidth Width=#1}\immediate\write\fwf@g{#1 \c@msetlinewidth}\fi\fi}
\def\pssetjoin#1{\ifps@cri\edef\curr@ntjoin{#1}\ifcurr@ntPS\expandafter\Pssetj@in#1:\fi\fi}
\def\Pssetj@in#1#2:{\PSc@mment{pssetjoin join=#1}%
    \if#1r\def\t@xt@{1}\else\if#1b\def\t@xt@{2}\else\def\t@xt@{0}\fi\fi%
    \immediate\write\fwf@g{\t@xt@\space\c@msetlinejoin}}
\def\defaultjoin{miter}   % Valeur par defaut
\def\Pssets@cond#1=#2|{\keln@mun#1|%
    \def\n@mref{c}\ifx\l@debut\n@mref\Pssets@condcolor{#2}\else%
    \def\n@mref{d}\ifx\l@debut\n@mref\pssetseconddash{#2}\else%
    \def\n@mref{w}\ifx\l@debut\n@mref\pssetsecondwidth{#2}\else%
    \immediate\write16{*** Unknown attribute: \BS@ psset second(..., #1=...)}%
    \fi\fi\fi}
\def\pssetseconddash#1{\edef\curr@ntseconddash{#1}}
\def\defaultseconddash{4}  % Valeur par defaut
\def\pssetsecondwidth#1{\edef\curr@ntsecondwidth{#1}}
\edef\defaultsecondwidth{\defaultwidth} % Valeur par defaut
\def\psresetsecondsettings{%
    \pssetseconddash{\defaultseconddash}\pssetsecondwidth{\defaultsecondwidth}%
    \Pssets@condcolor{\defaultsecondcolor}}
\def\Pssets@condcolor#1{\ifps@cri\result@tent=\@ne\expandafter\c@lnbV@l#1 :%
    \def\sec@ndcolor{}\def\sec@ndcolorc@md{}%
    \ifcase\result@tent\or\pssetsecondgray{#1}\or\or\pssetsecondrgb{#1}%
    \or\pssetsecondcmyk{#1}\fi\fi}
\def\pssetsecondcmyk#1{\def\sec@ndcolor{#1}\def\sec@ndcolorc@md{\c@msetcmykcolor}%
    \def\sec@ndcolorc@mdStroke{\c@msetcmykcolorStroke}}
\def\pssetsecondrgb#1{\def\sec@ndcolor{#1}\def\sec@ndcolorc@md{\c@msetrgbcolor}%
    \def\sec@ndcolorc@mdStroke{\c@msetrgbcolorStroke}}
\def\pssetsecondgray#1{\def\sec@ndcolor{#1}\def\sec@ndcolorc@md{\c@msetgray}%
    \def\sec@ndcolorc@mdStroke{\c@msetgrayStroke}}
\def\us@secondC@lor{\immediate\write\fwf@g{\d@fsecondC@lor}%
    \let\fillc@md=\sdfillc@md}
\def\sdfillc@md{\d@fsecondC@lor\space\c@mfill}
\edef\defaultsecondcolor{\defaultcolor} % Valeur par defaut
\def\Pss@tsecondSt{%
    \s@uvdash{\typ@dash}\pssetdash{\curr@ntseconddash}%
    \s@uvwidth{\typ@width}\pssetwidth{\curr@ntsecondwidth}\us@secondC@lor}
\def\Psrest@reSt{\pssetwidth{\typ@width}\pssetdash{\typ@dash}\us@primarC@lor}
\def\Pssetth@rd#1=#2|{\keln@mun#1|%
    \def\n@mref{c}\ifx\l@debut\n@mref\Pssetth@rdcolor{#2}\else%
    \immediate\write16{*** Unknown attribute: \BS@ psset third(..., #1=...)}%
    \fi}
\def\Pssetth@rdcolor#1{\ifps@cri\result@tent=\@ne\expandafter\c@lnbV@l#1 :%
    \def\th@rdcolor{}\def\th@rdcolorc@md{}%
    \ifcase\result@tent\or\Pssetth@rdgray{#1}\or\or\Pssetth@rdrgb{#1}%
    \or\Pssetth@rdcmyk{#1}\fi\fi}
\def\Pssetth@rdcmyk#1{\def\th@rdcolor{#1}\def\th@rdcolorc@md{\c@msetcmykcolor}%
    \def\th@rdcolorc@mdStroke{\c@msetcmykcolorStroke}}
\def\Pssetth@rdrgb#1{\def\th@rdcolor{#1}\def\th@rdcolorc@md{\c@msetrgbcolor}%
    \def\th@rdcolorc@mdStroke{\c@msetrgbcolorStroke}}
\def\Pssetth@rdgray#1{\def\th@rdcolor{#1}\def\th@rdcolorc@md{\c@msetgray}%
    \def\th@rdcolorc@mdStroke{\c@msetgrayStroke}}
\def\us@thirdC@lor{\immediate\write\fwf@g{\d@fthirdC@lor}%
    \let\fillc@md=\thfillc@md}
\def\thfillc@md{\d@fthirdC@lor\space\c@mfill}
\def\defaultthirdcolor{1}  % Valeur par defaut
\def\pstrimesh#1[#2,#3,#4]{{\ifcurr@ntPS\ifps@cri%
    \PSc@mment{pstrimesh Type=#1, Triangle=[#2,#3,#4]}%
    \s@uvc@ntr@l\et@tpstrimesh\ifnum#1>\@ne\Pss@tsecondSt\setc@ntr@l{2}%
    \Pstrimeshp@rt#1[#2,#3,#4]\Pstrimeshp@rt#1[#3,#4,#2]%
    \Pstrimeshp@rt#1[#4,#2,#3]\Psrest@reSt\fi\psline[#2,#3,#4,#2]%
    \PSc@mment{End pstrimesh}\resetc@ntr@l\et@tpstrimesh\fi\fi}}
\def\Pstrimeshp@rt#1[#2,#3,#4]{{\l@mbd@un=\@ne\l@mbd@de=#1\loop\ifnum\l@mbd@de>\@ne%
    \advance\l@mbd@de\m@ne\figptbary-1:[#2,#3;\l@mbd@de,\l@mbd@un]%
    \figptbary-2:[#2,#4;\l@mbd@de,\l@mbd@un]\psline[-1,-2]%
    \advance\l@mbd@un\@ne\repeat}}
\initpr@lim\initpss@ttings% Initialisation preliminaire
\catcode`\@=12

\pssetdefault (update=yes)

\newbox\dessin

\def\MyPSfile{FigDom.ps}
% 1. Definition of characteristic points (Interpolation points)
\figinit{pt}
\figpt 1:(-120, 60)\figpt 2:(-100, 65)\figpt 3:(-30, 55)\figpt 4:(50, 60)\figpt 5:(70, 57)
\figpt 6:(-120,30) \figpt 7:(-100, 37)\figpt 8:(-30, 23)\figpt 9:(50, 30)\figpt 10:(70, 26)
\figpt 11:(-120,-30)\figpt 12:(-100, -23)\figpt 13:(-30, -34)\figpt 14:(50, -30)\figpt 15:(70, -35)
\figpt 16:(-120,-60) \figpt 17:(-100, -57) \figpt 18:(-30, -67)\figpt 19:(50,-60)
\figpt 20:(70, -65)
\figptbary 21:c.g.[8,13 ; 1,1]
% 2. Creation of the postscript file
\psbeginfig{\MyPSfile}
\pscurve[1,1,2,3,4,5,5]
\pscurve[6,6,7,8,9,10,10]
\pscurve[11,11,12,13,14,15,15]
\pscurve[16,16,17,18,19,20,20]
\psset arrowhead(fillmode=yes,ratio=0.07)\psarrow [2,13]
%\psset arrowhead(out=yes)
\psarrow [7,18]
\psarrow [8,19]
\psarrow [3,14]
\psendfig
% 3. Writing text on the figure
\figvisu{\dessin}{Le foncteur de Picard de $[\L^{\frac1n}]$}{
\figinsert{\MyPSfile}
\figsetmark{$\figBullet$}
\figwritep[2,3,4,7,8,9,13,14,18,19]
%\figsetmark{}
\figwriten 3:$0$(4)
\figwriten 8:$0$(4)
\figwritene 13:$0$(4)
\figwritene 18:$0$(4)
\figwriten 4:$l$(4)
\figwriten 9:$l$(4)
\figwritene 14:$l$(4)
\figwritene 19:$l$(4)
\figwriten 2:$-l$(4)
\figwriten 7:$-l$(4)
\figsetmark{}
\figwritee 5:${(\Pic_{X/S})_n}$(4)
\figwritee 10:${(\Pic_{X/S})_{n-1}}$(4)
\figwritee 15:${(\Pic_{X/S})_0}$(4)
\figwritee 20:${(\Pic_{X/S})_{-1}}$(4)
\figwriten 21:$\vdots$(2)
\figwriten 3:$\vdots$(20)
\figwrites 18:$\vdots$(10)
}

\centerline{\box\dessin}

Ceci montre en particulier que si $\Pic_{X/S}$ est représentable, alors $\pic$ l'est aussi\footnote{Mais bien sûr, on le savait déjà dans le cas où $f$ est propre, plat et cohomologiquement plat en dimension zéro.}.

\vskip 5mm
\noindent
{\sc Description du champ de Picard de $[\L^{\frac1n}]$}
\vskip 5mm

On a une \og suite exacte \fg\ de champs de Picard :
\begin{equation}
\label{sec_racine_nieme_champ}
\xymatrix{1 \ar[r]& \champic(X/S) \ar[r]^{\pi^*}& \champic(\X/S) \ar[r]^{\chi}& f_*\Z/n\Z \ar[r]& 1.}
\end{equation}

Autrement dit, $\pi^*$ est pleinement fidèle, $\chi$ est un épimorphisme, et si $\Fc$ est un objet de $\champic(\X/S)$, il provient de $\champic(X/S)$ si et seulement si son caractère $\chi_{\Fc}$ est nul. Tout ceci a déjà été prouvé. De même que précédemment, si l'on suppose que $f_*\Z/n\Z=\Z/n\Z$, alors le champ $\champic(\X/S)$ s'identifie au champ obtenu à partir de $\champic(X/S)\times_S \Z$ en recollant les copies numéro $i$ et $i+nk$ le long de l'isomorphisme
$$\flechen{(\champic(X/S))_{i+nk}}{\mu_{l^k}}{(\champic(X/S))_i}$$
pour tous $i,k$ appartenant à $\Z$. En particulier il suffit dans ce cas que $\champic(X/S)$ soit algébrique pour que $\champic(\X/S)$ le soit aussi.

Dans le cas où $\Lc $ a une racine \iem{n} $\Rc$ sur $X$, la suite exacte~(\ref{sec_racine_nieme_champ}) est scindée et $\champic(\X/S)$ s'identifie au produit $\champic(X/S) \times_S f_*\Z/n\Z$.

\section{Courbes tordues d'Abramovich et Vistoli}
\label{courbes_tordues}

Abramovich et Vistoli ont mis au jour dans \cite{Abramovich_Vistoli_note}, \cite{Abramovich_Vistoli_CMFFS} et \cite{Abramovich_Vistoli_CSSM} une classe de courbes \og tordues\fg\ qui apparaissent naturellement lorsque l'on cherche à compactifier certains espaces de modules. Ces courbes sont des courbes nodales munies d'une \og structure champêtre\fg\ supplémentaire aux points singuliers ou en certains points marqués. Nous nous proposons de décrire le foncteur de Picard des courbes tordues \emph{lisses}\footnote{Nous espérons traiter le cas général dans un avenir proche.}. Nous allons voir que la structure supplémentaire modifie le foncteur de Picard de la courbe d'une manière très analogue à ce que nous avons pu observer dans la section précédente. Nous commençons par quelques rappels sur les courbes tordues.

\subsection{Rappel des définitions et propriétés élémentaires}

\begin{sousdefi}[\cite{Abramovich_Vistoli_CSSM} 4.1.2 ou \cite{Olsson_log_twisted_curves} 1.2]
Soit $S$ un schéma. Une courbe tordue sur $S$ est un champ de Deligne-Mumford $f : \Cc \fleche S$ modéré, propre, plat et de présentation finie sur $S$ dont les fibres sont purement de dimension~1, géométriquement connexes et ont au plus des singularités nodales, vérifiant de plus les propriétés suivantes :
\begin{itemize}
\item[1)] Si $\pi : \Cc \fleche C$ est l'espace de modules grossier de $\Cc$ et si $C_{\rm{lis}}$ est le lieu lisse de $C$ sur $S$, alors le sous-champ ouvert $\Cc\times_C C_{\rm{lis}}$ est le lieu lisse de $\Cc$ sur $S$.
\item[2)] Pour tout point géométrique $\overline{s} \fleche S$ le morphisme induit $\Cc_{\overline{s}} \fleche C_{\overline{s}}$ est un isomorphisme au-dessus d'un ouvert dense de $C_{\overline{s}}$.
\end{itemize}

Une courbe tordue $n$-pointée est une courbe tordue munie d'une collection $\{\Sigma_i\}_{i=1}^n$ de sous-champs fermés de $\Cc$ deux à deux disjoints tels que :
\begin{itemize}
\item[(i)] Pour tout $i$, le sous-champ fermé $\Sigma_i$ est dans le lieu lisse de $\Cc$.
\item[(ii)] Pour tout $i$, le morphisme $\Sigma_i\fleche \Cc\fleche S$ est une gerbe étale sur $S$.
\item[(iii)] Si $\Cc_{\emph{gén}}$ est l'ouvert complémentaire des $\Sigma_i$ dans
$\Cc_{\rm{lis}}$, alors $\Cc_{\emph{gén}}$ est un schéma.
\end{itemize}
\end{sousdefi}

\begin{sousremarque} \rm
On rappelle qu'un champ de Deligne-Mumford est dit modéré si pour tout corps algébriquement clos $k$ et tout morphisme $x : \Spec k \fleche \Cc$ le groupe $\Aut(x)$ est d'ordre inversible dans $k$.
\end{sousremarque}

\begin{sousremarque} \rm 
D'après la proposition~4.1.1 de \cite{Abramovich_Vistoli_CSSM}, l'espace de modules grossier $C$ est une courbe nodale propre et plate sur $S$, de présentation finie et à fibres géométriquement connexes. Si de plus la courbe $\Cc$ est $n$-pointée, alors l'espace de modules grossier $D_i$ du sous-champ fermé $\Sigma_i$ est naturellement un sous-schéma fermé de $C$, et le morphisme composé $D_i \fleche C \fleche S$ est un isomorphisme. Les $D_i$ définissent donc des sections de $C \fleche S$ qui en font une courbe nodale $n$-pointée au sens usuel.
\end{sousremarque}

\begin{sousthm}[\cite{Abramovich_Vistoli_CSSM} 3.2.3 ou \cite{Olsson_log_twisted_curves} 2.2]
Au voisinage (étale) d'un point marqué, la courbe $\Cc \fleche S$ est de la forme $[U/\mu_r]\fleche \Spec A$ où $U=\Spec A[x]$ et où un générateur de $\mu_r$ agit sur $U$ par $x\mapsto \xi.x$ avec $\xi$ une racine primitive \iem{r} de l'unité.

Au voisinage d'un n\oe ud, la courbe $\Cc \fleche S$ est de la forme $[U/\mu_r]\fleche \Spec A$ où $U=\Spec (A[x,y]/(xy-t))$ pour un certain $t\in A$ et où un générateur de $\mu_r$ agit sur $U$ par $(x,y)\mapsto (\xi.x,\xi'.y)$ avec $\xi$ et $\xi'$ des racines primitives ${r}^{\textrm{ièmes}}$ de l'unité.
\end{sousthm}

\begin{sousremarque}\rm
Si $p\in C$ est un n\oe ud fixé, et si $\xi$ et $\xi'$ sont les racines primitives ${r}^{\textrm{ièmes}}$ de l'unité qui apparaissent ci-dessus, on dit que le n\oe ud $p$ est \og balancé\fg\ si le produit $\xi.\xi'$ est égal à 1. On dit que la courbe $\Cc$ est balancée si tous ses n\oe uds sont balancés. Signalons que si $p$ n'est pas balancé, 
on a nécessairement $t=0$ dans la description locale ci-dessus. Autrement dit on ne peut pas faire disparaître le n\oe ud en déformant la courbe (\cite{Olsson_log_twisted_curves}~2.2).
\end{sousremarque}

Signalons qu'Olsson montre dans \cite{Olsson_log_twisted_curves} que se donner une courbe tordue revient à se donner une courbe nodale classique munie d'une certaine \og structure logarithmique\fg. Pour les courbes tordues lisses, Cadman donne une autre description, plus élémentaire.

\begin{sousthm}[\cite{Cadman_USTITCOC} 2.2.4 et 4.1]
\label{Cadman_equiv_AV}
Se donner une courbe tordue $(\Cc, \{\Sigma_i\}_{i=1}^n)$ $n$-pointée lisse sur un schéma $S$ noethérien et connexe est équivalent à se donner une courbe $n$-pointée $(C,\{\sigma_i\}_{i=1}^n)$ lisse sur $S$ et un $n$-uplet $\overrightarrow{r}=(r_1,\dots, r_n)$ d'entiers strictement positifs inversibles sur $S$. La courbe tordue $\Cc$ est alors isomorphe au champ $C_{\D,\overrightarrow{r}}$ défini de la manière suivante.

Chaque section $\sigma_i$ définit un diviseur de Cartier effectif $D_i$ de $C$. On note $s_{D_i}$ la section canonique de $\Oc(D_i)$ qui s'annule sur $D_i$. La collection des $\Oc(D_i)$ et des $s_{D_i}$ correspond à un morphisme $C \fleche [\A^n/\gm^n]$. On note aussi $\theta_{\overrightarrow{r}}$ le morphisme de $[\A^n/\gm^n]$ dans lui-même qui envoie un $n$-uplet de faisceaux inversibles $(L_1, \dots, L_n)$ muni de sections $(t_1, \dots, t_n)$ sur le $n$-uplet $(L_1^{r_1}, \dots, L_n^{r_n})$ muni de $(t_1^{r_1}, \dots, t_n^{r_n})$. On définit alors $C_{\D,\overrightarrow{r}}$ comme étant le produit fibré
$$C\times_{[\A^n/\gm^n],\theta_{\overrightarrow{r}}} [\A^n/\gm^n].$$
\end{sousthm}

\subsection{Description du foncteur de Picard des courbes tordues lisses}

Cadman décrit dans \cite{Cadman_USTITCOC} les faisceaux inversibles sur une courbe tordue lisse sur une base connexe et noethérienne (corollaire~3.2.1) : un faisceau inversible sur $C_{\D,\overrightarrow{r}}$ s'écrit de manière unique sous la forme $\pi^*L \otimes \prod_{i=1}^n \Tc_i^{\otimes k_i}$ où $\pi : \Cc \fleche C$ est la projection de $\Cc$ sur son espace de modules grossier, $L$ est un faisceau inversible sur $C$, $\Tc_i$ est le faisceau inversible $\Oc_{\Cc}(\Sigma_i)$ et $k_i$ est un entier compris entre $0$ et $r_i-1$. Il est clair que les hypothèses noethériennes ne sont pas essentielles pour ce résultat. Par ailleurs on peut aussi supprimer l'hypothèse de connexité sur la base : il faut alors remplacer les entiers $k_i$ par des fonctions localement constantes à valeurs dans $\Z$. On obtient ainsi le théorème suivant.

\begin{sousthm}
Soient $S$ un schéma, $C$ une courbe lisse $n$-pointée sur $S$, $\overrightarrow{r}$ un $n$-uplet d'entiers strictement positifs et $C_{\D,\overrightarrow{r}}$ la courbe tordue associée par la construction de Cadman. Soit $\Lc$ un faisceau inversible sur $C_{\D,\overrightarrow{r}}$. Alors il existe un faisceau inversible $L$ sur $C$ et des fonctions localement constantes $k_i$ appartenant à $H^0(S, \Z)$ prenant leurs valeurs dans $\{0, \dots, r_i-1\}$ tels que
$$\Lc \simeq \pi^* L \otimes \prod_{i=1}^n \Tc_i^{k_i}.$$
De plus les fonctions $k_i$ sont uniques, $L$ est unique à isomorphisme près, et $\Tc_i^{r_i}$ est isomorphe à $\pi^*\Oc(D_i)$. $\square$
\end{sousthm}

En particulier on a une suite exacte courte :
$$\xymatrix{0 \ar[r]& \Pic(C) \ar[r]& \Pic(\Cc) \ar[r]&
\disp \prod_{i=1}^n H^0(S, \Z/r_i\Z) \ar[r]& 0.}$$
Signalons que si $S$ est le spectre d'un corps algébriquement clos, Chiodo (\cite{Chiodo_twisted_curves_spin_structures}) obtient cette suite exacte d'une autre manière pour une courbe tordue quelconque.

Soit $(\Cc,\{\Sigma_i\}_{i=1}^{n})$ une courbe tordue $n$-pointée lisse sur une base $S$ noethérienne et connexe. D'après le théorème~\ref{Cadman_equiv_AV}, $\Cc$ est isomorphe au champ $C_{\D,\overrightarrow{r}}$ où $C$ est l'espace de modules grossier de $\Cc$, $\overrightarrow{r}$ est un $n$-uplet d'entiers positifs et $\D=(D_1,\dots, D_n)$ est le $n$-uplet de diviseurs effectifs de Cartier de $C$ correspondant aux $\Sigma_i$. Pour tout $i$ le morphisme de $D_i$ vers $S$ est un isomorphisme. On note $\Tc_i$ le faisceau $\Oc_{\Cc}(\Sigma_i)$. Alors toutes ces données sont compatibles au changement de base. Plus précisément, si $T\fleche S$ est un morphisme de changement de base, le produit fibré $\Cc\times_S T$ est isomorphe au champ $C'_{\D',\overrightarrow{r}}$ où $C'$ (resp. $D'_i$) est le produit fibré $C\times_S T$ (resp. $D_i\times_S T$) et $\D'=(D'_1,\dots, D'_n)$. De plus le faisceau inversible canonique $\Tc'_i$ n'est autre que $\Phi^*\Tc_i$ où $\Phi$ est la projection de $\Cc'$ sur $\Cc$. En appliquant le théorème précédent à la courbe $\Cc'$, on obtient pour tout $T$ une suite exacte courte
$$\xymatrix{0 \ar[r]& \Pic(C\times_S T) \ar[r]& \Pic(\Cc\times_S T) \ar[r]&
\disp \prod_{i=1}^n H^0(T, \Z/r_i\Z) \ar[r]& 0,}$$
autrement dit on a une suite exacte courte de préfaisceaux
$$\xymatrix{0 \ar[r]& \Pic_{C} \ar[r]& \Pic_{\Cc} \ar[r]&
\disp \prod_{i=1}^n \Z/r_i\Z \ar[r]& 0.}$$
En appliquant le foncteur \og faisceau associé pour la topologie étale\fg\ on en déduit une suite exacte courte de faisceaux étales :
$$\xymatrix{0 \ar[r]& \Pic_{C/S} \ar[r]& \Pic_{\Cc/S} \ar[r]&
\disp \prod_{i=1}^n \Z/r_i\Z \ar[r]& 0.}$$

Notons $l_i$ la classe de $\Oc_C(D_i)$ dans $\Pic_{C/S}(S)$ et $t_i$ la classe de $\Tc_i$ dans $\Pic_{\Cc/S}(S)$. On a dans $\Pic_{\Cc/S}$ la relation $t_i^{r_i}=l_i$. En procédant comme pour le cas du champ $[\L^{\frac1n}]$, on voit que le foncteur $\Pic_{\Cc/S}$ s'identifie au foncteur quotient de $\Pic_{C/S}\times_S \Z^n$ par les relations $t_i^{r_i}=l_i$ (où par abus les $t_1, \dots, t_n$ désignent aussi les générateurs canoniques de $\Z^n$). On peut construire ce quotient à la main comme suit. Le produit $\Pic_{C/S}\times_S \Z^n$ est une union disjointe de copies de $\Pic_{C/S}$ indexées par les $n$-uplets $\underline{\alpha}=(\alpha_1, \dots, \alpha_n)$ appartenant à $\Z^n$. Alors $\Pic_{\Cc/S}$ est obtenu en identifiant pour tout $\underline{\alpha}$, pour tout entier $k$ appartenant à $\Z$ et pour tout entier $i$ compris entre 1 et $n$, les copies $(\Pic_{C/S})_{\underline{\alpha}}$ et $(\Pic_{C/S})_{(\alpha_1, \dots, \alpha_i+kr_i, \dots,\alpha_n)}$ via l'isomorphisme de multiplication par $l_i$ :
$$\xymatrix{(\Pic_{C/S})_{(\alpha_1, \dots, \alpha_i+kr_i, \dots,\alpha_n)} \ar[r]^-{\mu_{l_i}}& (\Pic_{C/S})_{\underline{\alpha}}}$$
La loi de groupe est évidente.

\begin{sousremarque} \rm
Si l'on ne tient pas compte de la structure de groupe, on voit que $\Pic_{\Cc/S}$ s'identifie à une union disjointe de $r_1\dots r_n$ copies de $\Pic_{C/S}$.
\end{sousremarque}

\begin{sousremarque} \rm
La composante neutre de $\Pic_{\Cc/S}$ est la même que celle de $\Pic_{C/S}$. Autrement dit, le morphisme de $\Pic_{C/S}$ vers $\Pic_{\Cc/S}$ induit un isomorphisme naturel :
$$\xymatrix{\Pic^0_{C/S} \ar[r]^{\sim} & \Pic^0_{\Cc/S}.}$$
\end{sousremarque}
\appendix
\chaptermark{Annexes}
\chapter*{Annexes}
\addcontentsline{toc}{chapter}{Annexes}
\setcounter{chapter}{1}

Ainsi qu'il a été dit en introduction, la présente annexe rassemble les résultats relatifs à la cohomologie des faisceaux sur les champs algébriques nécessaires au texte principal. Voici en résumé la liste des sujets qui y sont abordés.
\begin{itemize}
\item[$\bullet$] Les deux premières sections rappellent les définitions du topos lisse-étale et du couple de foncteurs adjoints $(f^{-1},f_*)$ associé à un morphisme $f$ de champs algébriques. La seule nouveauté est la vérification du fait que sur un champ de Deligne-Mumford, les groupes de cohomologie lisse-étale coïncident avec les groupes de cohomologie étale.
\item[$\bullet$ Section \ref{site_llc} :] Nous introduisons un nouveau site, le site lisse-lisse champêtre, dont les objets sont les morphismes représentables et lisses $\Uc \fleche \X$ de \emph{champs algébriques}. Il définit le même topos que le site lisse-étale mais présente au moins deux avantages : il se comporte mieux vis-à-vis des images directes et il a un objet final.
\item[$\bullet$ Section~\ref{annexe_desc_coh} :] C'est au départ la nécessité de disposer de techniques de descente cohomologique à la Deligne-Saint-Donat qui a motivé ce travail. Nous avions en particulier besoin d'un analogue pour les champs algébriques de la suite spectrale de descente relative à un morphisme lisse et surjectif de schémas. Il s'est finalement avéré que l'introduction du site lisse-lisse champêtre rendait ce résultat presque trivial (voir proposition~\ref{suite_spectrale_de_descente}).
\item[$\bullet$ Sections~\ref{Faisceaux_acycliques} et~\ref{paragraphe_ss_Leray} :] Nous décrivons ici une classe de faisceaux acycliques adaptée aux particularités du site lisse-étale. Ces faisceaux \og $\gll$-acycliques\fg\ nous sont surtout utiles pour obtenir la suite spectrale de Leray relative à un morphisme de champs algébriques (thm.~\ref{ss_Leray}). Nous montrons aussi que les images directes supérieures d'un faisceau lisse-étale abélien peuvent être calculées comme l'on imagine (cf. prop.~\ref{prop_image_directe_sup}).
\item[$\bullet$ Section~\ref{coh_et_chgt_de_base} :] La formation des groupes de cohomologie et des images directes supérieures commute au changement de base plat.
\item[$\bullet$ Section \ref{par_coh_et_ext_inf} :] Il est d'usage, lorsque $i : X \fleche \widetilde{X}$ est une extension infinitésimale, d'identifier les catégories de faisceaux Zariski sur $X$ et sur $\widetilde{X}$. Ceci est tout à fait légitime puisque $X$ et $\widetilde{X}$ ont le même espace topologique sous-jacent. Mieux : le foncteur qui à un ouvert étale $\widetilde{U}$ de $\widetilde{X}$ associe l'ouvert étale $\widetilde{U}\times_{\widetilde{X}} X$ de $X$ définit une équivalence entre les sites étales de $\widetilde{X}$ et de $X$, ce qui permet d'identifier aussi les faisceaux étales. Il faut faire nettement plus attention avec la topologie lisse-étale. On peut en effet vérifier facilement que le foncteur ci-dessus n'est même pas fidèle. Heureusement, on peut tout de même identifier les groupes de cohomologie des faisceaux abéliens sur $X$ et sur $\widetilde{X}$ via le foncteur $i_*$ (cf.~\ref{coh_et_ext_inf}). Cette section contient également un résultat analogue pour les images directes supérieures.
\item[$\bullet$ Section \ref{Un_résultat_de_descente} :] On y trouvera un résultat technique utilisé dans la section suivante. Plus précisément, soit $\Fc$ un faisceau sur le site lisse-lisse champêtre d'un champ algébrique $\X$. On suppose qu'il existe une présentation $X \fleche \X$ de $\X$ telle que la restriction de $\Fc$ au site lisse de $X$ soit représentable par un espace algébrique lisse sur $X$. Alors $\Fc$ est lui-même représentable par un unique champ algébrique lisse sur $\X$.
\item[$\bullet$ Section \ref{Cohomologie_et_torseurs} :] On y décrit les torseurs du topos lisse-étale et l'on constate dans le cas particulier d'un groupe lisse sur la base $S$ que le $H^1$ au sens des foncteurs dérivés coïncide avec le groupe des classes de torseurs.
\item[$\bullet$ Section \ref{Cohomologie_plate} :] Après quelques généralités sur la cohomologie plate sur les champs algébriques, nous donnons principalement deux résultats. D'une part la suite spectrale qui relie la cohomologie plate à la cohomologie lisse-étale, et d'autre part la généralisation aux champs algébriques du théorème de Grothendieck suivant lequel dans le cas d'un groupe lisse, la cohomologie plate coïncide avec la cohomologie étale (cf. \cite{Dix}, exposé~VI,
paragraphe~11). 
\end{itemize}

\section{\hskip-2.3pt Cohomologie lisse-\'etale sur les champs alg\'ebriques}

\subsection{Cohomologie des faisceaux}
Rappelons bri\`evement, pour la
commodit\'e du lecteur, la d\'efinition du site lisse-\'etale d'un champ
alg\'ebrique donn\'ee au chapitre 12 de \cite{LMB}.

\begin{sousdefi}
\label{def_lisse_etale}
Soit $\X$ un $S$-champ alg\'ebrique. On appelle site \emph{lisse-\'etale} de
$\X$ et on note \emph{Lis-\'et}$(\X)$ le site d\'efini comme suit.

Les ouverts lisses-\'etales de $\X$ sont les couples $(U,u)$ o\`u $U$ est un
$S$-espace alg\'ebrique et $u:U\fleche \X$ un morphisme repr\'esentable et
lisse. Une fl\`eche entre deux tels ouverts $(U,u)$ et $(V,v)$ est un couple
$(\varphi,\alpha)$ faisant 2-commuter le diagramme suivant :
$$\shorthandoff{!;:?}
\xymatrix@R=0.9pc@C=0.9pc{U \ar[rr]^{\varphi} \ar[rd]_u  & \raisebox{-3ex}{$^{\alpha} \FlecheNE$} & V \ar[ld]^v\\
& \X&}$$
Une famille couvrante de $(U,u)$ est une collection de morphismes
$$((\varphi_i,\alpha_i):(U_i,u_i)\flechelongue (U,u))_{i\in I}$$ telle que le
1-morphisme d'espaces alg\'ebriques
$$\coprod_{i\in I} \varphi_i : \coprod_{i\in I} U_i \flechelongue U$$
soit \'etale et surjectif.

Le site \'etale de $\X$, not\'e $\et(\X)$, est la sous-cat\'egorie pleine de
$\liset(\X)$ dont les
objets sont les couples $(U,u)$ o\`u $u$ est un morphisme \'etale, munie de la
topologie induite par celle de $\liset(\X)$.
\end{sousdefi}

\begin{sousremarque}\rm
\label{remarque_faisceau_etale_associe}
Notons $(.)_{\textrm{\'et}}$ l'op\'eration de restriction des pr\'efaisceaux ou
faisceaux lisses-\'etales au site \'etale de $\X$. Notons encore $\underline{a}$
le foncteur \og faisceau associ\'e \fg\ (le contexte devant lever toute
ambigu\"{i}t\'e sur le site concern\'e, \'etale ou lisse-\'etale).
Maintenant si $\Fc$ est un pr\'efaisceau sur le site lisse-\'etale de $\X$,
on a un isomorphisme canonique :
$$\flechen{\underline{a}(\Fc_{\textrm{\'et}})}{\sim}
{(\underline{a}\Fc)_{\textrm{\'et}}}.$$
Ceci est essentiellement d\^u au fait que si $(U,u)$ est un ouvert \'etale de
$\X$ et si $((\varphi_i,\alpha_i) : (U_i,u_i) \fleche (U,u))_i$ est une famille
couvrante dans $\liset(\X)$, alors les $U_i$ sont \'etales sur $\X$ et la
famille consid\'er\'ee est couvrante dans $\et(\X)$.
\end{sousremarque}

On notera encore $\O_{\X}$ le faisceau structural de $\X$ d\'efini en (12.7.1)
dans~\cite{LMB}. Si $\Ac$ est un faisceau d'anneaux sur $\liset(\X)$, on notera
$\Mod_{\Ac}(\X)$ la cat\'egorie des faisceaux de $\Ac$-modules sur le site
lisse-\'etale de $\X$, ou plus simplement $\Mod(\X)$ lorsque $\Ac=\O_{\X}$.
La cat\'egorie des faisceaux ab\'eliens sera not\'ee $\Ab(\X)$.
\begin{sousprop}
La cat\'egorie $\Mod_{\Ac}(\X)$ est une cat\'egorie ab\'elienne avec suffisamment
d'objets injectifs. En particulier il en est ainsi des cat\'egories $\Mod(\X)$
et $\Ab(\X)$.
\end{sousprop}
\begin{demo}
\cite{SGA4_1} II (6.7) et \cite{Tohoku} th\'eor\`eme (1.10.1).
\end{demo}

On rappelle que le foncteur \og sections globales \fg\  est d\'efini de la
mani\`ere suivante. Si $\Fc$ est un faisceau lisse-\'etale sur $\X$, l'ensemble
$\Gamma(\X, \Fc)$ est l'ensemble des familles $(s_{(U,u)})$ de sections de $\Fc$
sur les $(U,u)\in \ob\liset(\X)$ qui sont compatibles aux fl\`eches de
restriction en un sens \'evident (cf. \cite{LMB} (12.5.3)). On v\'erifie
imm\'ediatement que le foncteur $\Gamma(\X, .) : \Ab(\X)\fleche \Ab$ est exact
\`a gauche. On d\'efinit alors $H^i(\X, .)$ comme \'etant le
$i^{\textrm{\`eme}}$ foncteur d\'eriv\'e \`a droite de $\Gamma(\X, .)$.

\begin{sousprop}
\label{cohomologie_modules_egale_cohomologie_abelienne}
Sur la cat\'egorie $\Mod(\X)$, le foncteur $H^i(\X,.)$ co\"{i}ncide avec le
$i^{\textrm{\`eme}}$ foncteur d\'eriv\'e \`a droite de $\Gamma(\X, .) : \Mod(\X)
\fleche \Ab$.
\end{sousprop}
\begin{demo}
\cite{SGA4_2} V (3.5).
\end{demo}

Nous esp\'erons ici que le lecteur a pr\'esents \`a l'esprit les r\'esultats
(12.2.1) et (12.3.3) de \cite{LMB}, ainsi que les notations qui y sont
introduites. Si tel n'\'etait pas le cas, nous l'invitons \`a se les
rem\'emorer.
\begin{souslem} Soit $\X$ un $S$-champ de Deligne-Mumford et
soient $\Fc$, $\Gc$ des faisceaux lisses-étales d'ensembles (resp. de groupes
abéliens) sur $\X$. On suppose que $\Fc$ est cartésien. Alors l'application
$$\Hom(\Fc, \Gc) \flechelongue \Hom(\Fc_{\text{\rm ét}},\Gc_{\text{\rm ét}})$$
induite par le foncteur d'inclusion $\xymatrix@C=1pc{\et(\X) \ar@{^{(}->}[r] &
\liset(\X)}$  est bijective.
\label{lemme_de_prolongement_des_morphismes_etales}
\end{souslem}
\begin{demo} Soit $g : \Fc_{\text{ét}} \fleche \Gc_{\text{ét}}$ un morphisme de
faisceaux (resp. de faisceaux abéliens). Il est donné par la
collection des $g_{U,u} : \Fc_{U,u} \fleche \Gc_{U,u}$ pour
$(U,u)\in\ob\et(\X)$, ces morphismes $g_{U,u}$ satisfaisant de plus à une
condition de compatibilité évidente avec les morphismes de changement de base
$\theta_{\varphi,\alpha}$ pour tout morphisme $(\varphi,\alpha)$ dans $\et(\X)$. Il
faut montrer que $g$ se prolonge de manière unique en un morphisme de faisceaux
lisses-étales. Soit $(U,u)\in\ob\liset(\X)$. Soit $P: X \fleche \X$ une
présentation étale de $\X$. On note
\begin{equation*}
\xymatrix{V \ar[r]^{\psi} \ar[d]_{\varphi}\cartesien & X \ar[d]^P \\
U\ar[r]^u & \X,}
\end{equation*}
et on pose $v=P\circ \psi$. Un morphisme $h$ qui prolonge $g$ fait
nécessairement commuter le diagramme suivant\footnote{Afin de ne pas alourdir
l'exposé, nous avons cru bon de ne pas préciser les 2-isomorphismes dans les
morphismes entre ouverts lisses-étales. Nous notons donc $\theta_{\varphi}$ au
lieu de $\theta_{\varphi,\alpha}$. Une notation comme $\theta_{\psi,\Fc}$ désigne alors le $\theta_{\varphi}$ relatif à $\Fc$.} :
$$\xymatrix{\psi^{-1}\Fc_{X,P} \ar[r]^{\psi^{-1}(g_{X,P})}
\ar[d]_{\theta_{\psi,\Fc}^{\natural}} &
\psi^{-1}\Gc_{X,P} \ar[d]^{\theta_{\psi,\Gc}^{\natural}}\\
\Fc_{V,v} \ar[r]^{h_{V,v}} & \Gc_{V,v}}$$
Comme $\Fc$ est cartésien, $\theta_{\psi,\Fc}^{\natural}$ est un isomorphisme,
donc $h_{V,v}$ est déterminé de manière unique par
$h_{V,v}=\theta_{\psi,\Gc}^{\natural}\circ (\psi^{-1}g_{X,P}) \circ
(\theta_{\psi,\Fc}^{\natural})^{-1}$. Le diagramme analogue obtenu à partir du
morphisme $\varphi$ nous donne $\varphi^{-1}(h_{U,u})=
(\theta_{\varphi,\Gc}^{\natural})^{-1}\circ
h_{V,v} \circ \theta_{\varphi,\Fc}^{\natural}$
(o\`u $\theta_{\varphi,\Gc}^{\natural}$ est un isomorphisme parce que $\varphi$ est
étale). Comme $\varphi$ est étale surjectif, le foncteur $\varphi^{-1}$ est
pleinement fidèle de sorte que la relation précédente détermine de manière
unique le morphisme $h_{U,u}$. On vérifie facilement que si $(U,u)$ et $(U',u')$
sont deux ouverts lisses-étales et si $\varphi : U'
\fleche U$ est un morphisme entre ces ouverts, alors le diagramme suivant
commute :
$$\xymatrix{\varphi^{-1}\Fc_{U,u} \ar[r]^{\varphi^{-1}h_{U,u}}
\ar[d]_{\theta_{\varphi,\Fc}^{\natural}} &
\varphi^{-1}\Gc_{U,u} \ar[d]^{\theta_{\varphi,\Gc}^{\natural}}\\
\Fc_{U',u'} \ar[r]^{h_{U',u'}} & \Gc_{U',u'}.}$$
En d'autres termes, la collection des
$h_{U,u}$ ainsi construits définit bien un morphisme de faisceaux lisses-étales
de $\Fc$ dans $\Gc$ qui prolonge $g$. 
\end{demo}

\begin{sousremarque}\rm
Le résultat est faux si on ne suppose pas $\Fc$ cartésien. Il suffit de considérer un faisceau lisse-étale
abélien $\Fc$ non nul tel que $\Fc_{\text{ét}}$ soit nul. Alors, dans $\Hom(\Fc, \Fc)$ on a au moins deux
éléments distincts, à savoir l'identité et le morphisme nul, tandis que $\Hom(\Fc_{\text{ét}}, \Fc_{\text{ét}})$
est réduit à zéro. Pour exhiber un tel faisceau, on peut par exemple prendre $\X=\Spec k$, le spectre d'un corps,
et poser $\Fc(U,u)=\Omega_{U/k}$.
\end{sousremarque}

\begin{souscor}
Soient $\X$ un $S$-champ algébrique et $\Fc$ un faisceau lisse-étale sur $\X$. Alors le morphisme canonique
$$\Gamma_{\text{\rm lis-ét}}(\X,\Fc) \flechelongue \Gamma_{\text{\rm ét}}(\X,\Fc)$$
est un isomorphisme.
\end{souscor}
\begin{demo}
Il suffit juste d'appliquer (\ref{lemme_de_prolongement_des_morphismes_etales}) en prenant pour $\Fc$ le
faisceau final, i.e. le faisceau qui à tout ouvert lisse-étale associe un singleton. Il est clair que
ce faisceau est bien cartésien.
\end{demo}

\begin{souscor}
Soit $\X$ un $S$-champ de Deligne-Mumford, et soit $\Fc$ un objet injectif de la
catégorie des faisceaux lisses-étales abéliens sur $\X$. Alors la restriction
$\Fc_{\text{\rm ét}}$ de $\Fc$ au site étale de $\X$ est un objet injectif de la
catégorie des faisceaux étales abéliens.
\end{souscor}
\begin{demo}
Soit $0 \fleche \Mc \fleche \Mc'$ une suite exacte de faisceaux étales, et soit
$f : \Mc \fleche \Fc_{\text{ét}}$ un morphisme donné. Il s'agit de montrer que
$f$ se prolonge en un morphisme $\Mc' \fleche \Fc_{\text{ét}}$. D'après
\cite{LMB} (12.3.3), le foncteur d'inclusion $\et(\X) \fleche \liset(\X)$ induit
une équivalence de catégories $\X_{\textrm{lis-ét,cart}} \fleche
\X_{\textrm{ét}}$. De plus cette équivalence préserve les faisceaux abéliens,
de sorte qu'elle induit une équivalence de la sous-catégorie pleine de $\Ab(\X)$
dont les objets sont les faisceaux abéliens cartésiens, sur la catégorie
$\Ab_{\textrm{ét}}(\X)$ des faisceaux étales abéliens. Il existe donc un
morphisme de faisceaux lisses-\'etales cartésiens
$\widetilde{\Mc} \fleche \widetilde{\Mc'}$ qui par restriction au site
étale de $\X$ induit $\Mc \fleche \Mc'$. De plus si $(U,u)$ est un ouvert
lisse-étale de $\X$, le morphisme $\widetilde{\Mc}_{U,u} \fleche
\widetilde{\Mc'}_{U,u}$ est par construction égal au morphisme $u^{-1}\Mc
\fleche u^{-1}\Mc'$, où $u$ est considéré comme un morphisme de topos de
$U_{\textrm{ét}}$ dans $\X_{\textrm{ét}}$. Le foncteur $u^{-1}$ étant exact, on
en déduit que $\widetilde{\Mc} \fleche \widetilde{\Mc'}$ est un monomorphisme.
Par ailleurs, il existe en vertu du lemme
(\ref{lemme_de_prolongement_des_morphismes_etales}) un morphisme
$\tilde{f} :  \widetilde{\Mc} \fleche \Fc$ qui induit $f$ par restriction au
site étale. Comme $\Fc$ est injectif, $\widetilde{f}$ se prolonge à
$\widetilde{\Mc'}$ et la restriction au site étale de $\X$ de ce prolongement
fournit un prolongement de $f$ à  ${\Mc'}$, ce qui achève la démonstration.
\end{demo}

\begin{sousprop}
\label{coh_et_egale_coh_liset}
Soient $\X$ un $S$-champ de Deligne-Mumford et
$\Fc\in\Ab(\X)$ un faisceau lisse-\'etale ab\'elien sur $\X$.
On note $\Fc_{\text{\rm \'et}}$ la restriction de $\Fc$
au site \'etale de $\X$. Alors pour tout $q\geq 0$ on a un isomorphisme
canonique :
$$\flechen{H^q_{\text{\rm lis-\'et}}(\X,\Fc)}{\sim}
{H^q_{\text{\rm \'et}}(\X,\Fc_{\text{\rm \'et}}).}$$
\end{sousprop}
\begin{demo}
C'est évident compte tenu des deux résultats précédents.
\begin{comment}
Compte tenu du résultat précédent, il suffit clairement de montrer que le
morphisme $$\Gamma_{\textrm{lis-ét}}(\X,\Fc) \fleche
\Gamma_{\textrm{ét}}(\X,\Fc)$$ qui associe à une famille compatible $(s_{U,u})$
de sections de $\Fc$ sur les ouverts lisses-\'etales $(U,u)$ de $\X$, la famille
form\'ee des $s_{U,u}$ pour $u$ \emph{\'etale}, est un isomorphisme.
Donnons-nous donc une famille $(s_{U,u})$ d\'efinie sur le site \'etale, et
montrons qu'elle se prolonge de mani\`ere unique au site lisse-\'etale. Soit
$(U,u)$ un ouvert lisse-\'etale de $\X$. Reprenant les notations du diagramme
(\ref{lemme_diagramme_notations}), on a n\'ecessairement
$s_{V,v}=(s_{X,P})|_{(V,v)}$. On note aussi $w$ le morphisme $u\circ
\varphi\circ \pr_1$ (qui est \'egal \`a $u\circ \varphi\circ \pr_2$) de
$V\times_U V$ dans $\X$. Vu que le morphisme $V\fleche U$ est
\'etale surjectif, le diagramme
$$\xymatrix{\Fc(U,u)\ar[r]& \Fc(V,v) \ar@<0.5ex>[r] \ar@<-0.5ex>[r] &
\Fc(V\times_U V,w)}$$
est exact. Comme la section $s_{V,v}$ provient de $s_{X,P}$, ses deux images
r\'eciproques dans $\Fc(V\times_U V,w)$ co\"{i}ncident, donc il existe un
unique $s_{U,u}\in \Fc(U,u)$ qui induit $s_{V,v}$. On v\'erifie facilement que
la famille ainsi construite est compatible au changement de base.
\end{comment}
\end{demo}

\subsection{Fonctorialité du topos lisse-étale}
\label{fonctorialite_lisse_etale}
Soit $f: \X \fleche \Y$ un 1-morphisme de $S$-champs alg\'ebriques. On rappelle
(cf. \cite{LMB}, (12.5)) que l'on peut lui associer un couple de foncteurs
adjoints $(f^{-1},f_*)$ de la mani\`ere suivante. Pour tout faisceau
lisse-\'etale $\Fc$ sur $\X$ et tout ouvert lisse-\'etale $(V,v)$ de $\Y$, on
pose
$$(f_*\Fc)(V,v)=\lpro \Fc(U,u)$$
o\`u la limite projective est prise sur les carr\'es 2-commutatifs de $S$-champs
alg\'ebriques
$$\xymatrix{U\ar[r]^g \ar[d]_u & V\ar[d]^v\\ \X\ar[r]^f&\Y}$$
avec $(U,u)\in\ob\liset(\X)$. Les flèches de restriction sont définies de
manière évidente. Le foncteur image directe est décrit de manière légèrement
différente dans \cite{LMB}, mais il est très facile de vérifier qu'il s'agit
bien là du même foncteur. Ainsi défini, il est évident qu'il commute aux limites
projectives arbitraires, de sorte que l'on sait \emph{a priori} qu'il admet un adjoint
à gauche, déterminé de manière unique à unique isomorphisme près. Nous noterons
$f^{-1}$ cet adjoint à gauche. On peut en donner une description simple.
Pour tout faisceau lisse-\'etale $\Gc$ sur $\Y$,
on d\'efinit d'abord un pr\'efaisceau $\widehat{f^{-1}}\Gc$ sur $\liset(\X)$ en
posant pour tout ouvert lisse-\'etale $(U,u)$ de $\X$,
$$(\widehat{f^{-1}}\Gc)(U,u)=\lind \Gc(V,v)$$
o\`u la limite inductive est prise sur les carr\'es 2-commutatifs de $S$-champs
alg\'ebriques
$$\xymatrix{U\ar[r]^g \ar[d]_u & V\ar[d]^v\\ \X\ar[r]^f&\Y}$$
avec $(V,v)\in\ob\liset(\Y)$. Les flèches de restriction du préfaisceau
$\widehat{f^{-1}}\Gc$ sont définies de
manière évidente. On d\'efinit alors $f^{-1}\Gc$ comme \'etant le
faisceau lisse-\'etale associ\'e au pr\'efaisceau $\widehat{f^{-1}}\Gc$ .
Il est clair que ces foncteurs induisent un couple de foncteurs adjoints entre
les cat\'egories des faisceaux ab\'eliens sur $\X$ et sur $\Y$.
Rappelons aussi que ces foncteurs s'expriment plus simplement dans les cas
particuliers suivants (\cite{LMB} (12.5.1) et (12.5.2)). Si $f$ est lisse, alors
pour tout faisceau lisse-\'etale $\Gc$ sur $\Y$ le faisceau lisse-\'etale
$f^{-1}\Gc$ n'est autre que la restriction de $\Gc$ au site lisse-\'etale de
$\X$ par le foncteur $(U,u)\mapsto (U,f\circ u)$.
Si $f$ est repr\'esentable,
le foncteur image directe est induit par une application continue $\liset(\X)
\fleche \liset(\Y)$. Pour tout faisceau lisse-\'etale $\Fc$ sur $\X$ le faisceau
lisse-\'etale $f_*\Fc$ est donn\'e par $f_*\Fc(V,v)=\Fc(\X\times_{\Y} V,
\pr_{\X})$.

\begin{sousremarque}\rm
Comme on pourra le lire bient\^ot dans la prochaine \'edition de \cite{LMB}, ou
d\`es aujourd'hui dans \cite{Olsson_Sheaves_on_Artin_stacks}, il serait erron\'e
de penser que le couple $(f^{-1},f_*)$ est toujours un
morphisme de topos. Il peut en effet arriver que le foncteur $f^{-1}$ ne soit
pas exact, même lorsque $\X$ et $\Y$ sont des schémas. Pour s'en convaincre on
consultera les r\'ef\'erences cit\'ees. C'est toutefois le cas d\`es que $f$
est lisse.
\end{sousremarque}

\begin{comment}
On d\'efinit enfin le foncteur $f^* : \Mod(\Y) \fleche \Mod(\X)$ de la mani\`ere
suivante. (...)
Il est clair que ce foncteur est un adjoint \`a gauche du foncteur $f_*:
\Mod(\X) \fleche \Mod(\Y)$.
\end{comment}

Le foncteur $f_* : \Ab(\X)\fleche \Ab(\Y)$, qui a un adjoint
\`a gauche, commute aux limites projectives et en particulier est exact \`a
gauche. On peut donc d\'efinir les foncteurs d\'eriv\'es $R^qf_* : \Ab(\X)
\fleche \Ab(\Y)$. Nous les étudierons plus en détail ultérieurement : nous ne disposons
pas à l'heure actuelle de tous les outils techniques nécessaires.

\subsection{Le site lisse-lisse champêtre d'un champ algébrique}
\label{site_llc}
Pour un certain nombre de considérations techniques, le site lisse-étale
défini ci-dessus ne contient pas suffisamment d'ouverts pour être vraiment
commode. En effet, lorsque $f$ n'est pas repr\'esentable, le foncteur $f_*$ n'est pas
induit par une application continue $\liset(\X) \fleche \liset(\Y)$, ce qui pose problème
par exemple lorsque l'on essaye de calculer les foncteurs images directes supérieures $R^qf_*$
(voir le paragraphe (\ref{images_directes_supérieures})). C'est la
raison pour laquelle nous introduisons un site un peu plus gros, qui ne
pr\'esentera plus
les m\^emes inconv\'enients. Nous démontrons ensuite (\ref{topos_llc_equiv_topos_liset}) que le topos qu'il définit est équivalent
au topos lisse-étale. Le choix de la topologie lisse plut\^ot qu'\'etale
pour ce site est essentiellement d\^u au fait que pour la topologie
\'etale, les \og ouverts lisses champêtres\fg\ 
ne sont pas toujours recouverts par un espace alg\'ebrique.

\medskip

Soit $\X$ un $S$-champ alg\'ebrique. On d\'efinit la 2-cat\'egorie des
ouverts lisses champêtres de la mani\`ere suivante. Les objets sont les couples $(\Uc,u)$
o\`u $\Uc$ est un $S$-champ alg\'ebrique et $u : \Uc \fleche \X$ est
un morphisme repr\'esentable et lisse. Un 1-morphisme entre deux tels ouverts
$(\Uc,u)$ et $(\Vc,v)$ est un couple $(\varphi,\alpha)$ o\`u $\varphi : \Uc
\fleche \Vc$ est un 1-morphisme (automatiquement représentable !) de $S$-champs alg\'ebriques et $\alpha : u
\Rightarrow v\circ \varphi$ est un 2-isomorphisme.
$$\shorthandoff{!;:?}
\xymatrix@R=0.9pc@C=0.9pc{\Uc \ar[rr]^{\varphi} \ar[rd]_u  & \raisebox{-3ex}{$^{\alpha} \FlecheNE$} & \Vc \ar[ld]^v\\
& \X&}$$
Si $(\varphi,\alpha)$ et $(\psi,\beta)$ sont deux 1-morphismes de $(\U,u)$ dans
$(\Vc,v)$, un 2-morphisme entre $(\varphi,\alpha)$ et $(\psi,\beta)$ est un
2-isomorphisme $\gamma : \varphi \Rightarrow \psi$ tel que
$\beta=(v_*\gamma)\circ \alpha$.

\begin{souslem}
\label{lemme_site_champetre}
Soient $(\U,u)$ et $(\Vc,v)$ deux ouverts lisses champêtres. Alors la cat\'egorie
des morphismes de $(\U,u)$ dans $(\Vc,v)$ est \'equivalente \`a une cat\'egorie discr\`ete.
\end{souslem}
\begin{demo}
Il suffit de montrer que si $(\varphi,\alpha)$ est un objet de
$\Hom((\U,u),(\Vc,v))$, alors le groupe des automorphismes de $(\varphi,\alpha)$
est r\'eduit \`a l'identit\'e. Soit $\gamma\in\Aut(\varphi,\alpha)$. Alors
$\gamma$ est un 2-automorphisme de $\varphi$ qui induit l'identit\'e de $v\circ
\varphi$. Or d'apr\`es \cite{LMB} (8.1.2), (2.2) et (2.3), un 1-morphisme de
$S$-champs alg\'ebriques $v : \Vc \fleche \X$ est repr\'esentable si
et seulement si pour tout $U\in\ob\aff$, le foncteur $v_U : \Vc_U
\fleche \X_U$ est fid\`ele. En particulier le morphisme $\Aut(\varphi) \fleche
\Aut(v\circ \varphi)$ est injectif, donc $\gamma$ est l'identit\'e de $\varphi$.
\end{demo}

\`A l'avenir on identifiera $\Hom((\U,u),(\Vc,v))$ \`a une cat\'egorie
discr\`ete \'equivalente et on parlera de l'\emph{ensemble} des morphismes de
$(\Uc, u)$ dans $(\Vc,v)$, et de la \emph{cat\'egorie} des ouverts lisses champêtres.
Il est clair que cette cat\'egorie admet des produits fibr\'es.

\begin{sousdefi}
On appelle site lisse-lisse champêtre, et on note $\gll(\X)$, la cat\'egorie des
ouverts lisses champêtres de $\X$ munie de la topologie engendr\'ee par la pr\'etopologie
pour laquelle les familles couvrantes sont les familles de morphismes
$$((\Uc_i,u_i)\flechelongue (\Uc,u))_{i\in I}$$ telles que le morphisme (automatiquement représentable)
$$\xymatrix{\disp \coprod_{i\in I} u_i : \coprod_{i\in I} \Uc_i \ar[r] &\Uc}$$
soit lisse et surjectif.
\end{sousdefi}

\begin{sousprop}
\label{topos_llc_equiv_topos_liset}
\begin{itemize}
\item[1)] La topologie de $\liset(\X)$ est \'egale \`a la topologie engendr\'ee par la
pr\'etopologie dite lisse, pour laquelle les familles couvrantes sont les
familles de morphismes $((U_i,u_i)\fleche (U,u))_{i\in I}$ telles que le
morphisme $\coprod_{i\in I} u_i : \coprod_{i\in I} U_i \fleche U$
soit lisse et surjectif.
\item[2)] Le foncteur d'inclusion $\xymatrix@C=1pc{\liset(\X) \ar@{^(->}[r]& \gll(\X)}$
induit une \'equivalence
de topos de la cat\'egorie des faisceaux sur le site lisse-lisse champêtre vers la
cat\'egorie des faisceaux sur $\liset(\X)$.
\end{itemize}
\end{sousprop}
\begin{demo}
1) Si $V\fleche U$ est un morphisme lisse et surjectif de $S$-espaces
alg\'ebriques,
il existe un $S$-espace alg\'ebrique $U'$ et un morphisme $U'\fleche U$ \'etale
surjectif qui se factorise par $V$. On en d\'eduit facilement que toute famille
couvrante de la pr\'etopologie lisse admet un raffinement qui est une famille
couvrante de la pr\'etopologie \'etale. En particulier, si $R$ est le crible
engendr\'e par une famille couvrante de la pr\'etopologie lisse, il contient un
crible engendr\'e par une famille couvrante de la pr\'etopologie \'etale.
D'apr\`es la proposition~II~(1.4) de \cite{SGA4_1}, si $T$ est la topologie
engendr\'ee par une pr\'etopologie $E$, et si $(U,u)$ est un ouvert, pour qu'un
crible $R$ de $(U,u)$ soit un crible couvrant de $(U,u)$ pour $T$, il faut et il
suffit qu'il contienne un crible $R'$ engendr\'e par une famille couvrante de la
pr\'etopologie $E$. Notre assertion en r\'esulte imm\'ediatement.

2) Compte tenu de 1) il est clair que la topologie de $\liset(\X)$ est bien la
topologie induite par celle de $\gll(\X)$ via le foncteur d'inclusion
$\xymatrix@C=1pc{\liset(\X) \ar@{^(->}[r]& \gll(\X)}$. Il suffit alors d'appliquer le
lemme de comparaison~\cite{SGA4_1}~III~(4.1).
\end{demo}

\begin{sousremarque}
\label{rem_foncteurs_compatibles}\rm
Notons $\X_{\gll}$ le topos des faisceaux sur $\gll(\X)$. Si $f : \X \fleche \Y$
est un 1-morphisme de $S$-champs alg\'ebriques, il induit un couple de foncteurs
adjoints $(f^{-1}_{\text{llc}},f_{*}^{\text{llc}})$ d\'efinis de mani\`ere
\'evidente. En particulier le foncteur $f_*^{\text{llc}}$ est simplement donn\'e
par $(f_{*}^{\text{llc}}\Fc)(\Uc,u)=\Fc(\X\times_{\Y} \Uc, \pr_{\X})$ pour tout
ouvert lisse champêtre $u : \Uc \fleche \Y$ de $\Y$. On peut aussi
d\'efinir un foncteur \og sections globales \fg\ en posant simplement 
$\Gamma_{\text{llc}}(\X,\Fc)=\Fc(\X, \Id_{\X})$. Il est alors clair que ces
foncteurs co\"{i}ncident avec leurs analogues sur les faisceaux lisses-\'etales
via l'\'equivalence de cat\'egories ci-dessus. Nous supprimerons \`a l'avenir la
mention \og llc \fg\ dans les notations qui d\'esignent ces foncteurs.
\end{sousremarque}

\begin{sousremarque}\rm
\label{restriction_faisceaux_injectifs}
Si $u : \Uc \fleche \X$ est un morphisme repr\'esentable et lisse de $S$-champs
alg\'ebriques, et si $\Ac$ est un anneau du topos $\X_{\gll}$, le foncteur
$$\fonction{u^*}{\Mod_{\Ac}(\X)}{\Mod_{\Ac_{(\Uc,u)}}(\Uc)}
{\Fc}{\Fc_{(\Uc,u)}}$$
o\`u $\Fc_{(\Uc,u)}$ (resp. $\Ac_{(\Uc,u)}$) d\'esigne la restriction de $\Fc$
(resp. $\Ac$) au site lisse-lisse champêtre de $\Uc$, commute
aux limites projectives et inductives arbitraires. Il a donc un adjoint \`a
gauche $u_!$, et de plus cet adjoint est exact (tout ceci r\'esulte de
\cite{SGA4_1}~IV~(11.3.1)). Par suite $u^*$ transforme les objets injectifs en
objets injectifs (voir aussi \cite{SGA4_2} V (2.2)).
\end{sousremarque}

\subsection{La suite spectrale relative à un recouvrement}
\label{annexe_desc_coh}
Soit $P : \Uc^0 \fleche \X$ un morphisme représentable, lisse et surjectif de
$S$-champs algébriques.
On note
\begin{eqnarray*}
\Uc^1 &=& \Uc^0\times_{\X}\Uc^0\\
\Uc^2 &=& \Uc^0\times_{\X}\Uc^0\times_{\X}\Uc^0\\
&\vdots& \\
\Uc^n &=& \Uc^0\times_{\X} \dots \times_{\X}\Uc^0 \quad \text{($n+1$ fois)}\\
&\vdots&
\end{eqnarray*}
Les $\Uc^n$ forment avec les diagonales partielles et les projections un
$S$-champ
alg\'ebrique simplicial muni d'une augmentation vers $\X$ :
$$\xymatrix{\dots & \Uc^n & \dots & \Uc^2 \ar[r] \ar@<2ex>[r] \ar@<-2ex>[r] &
  \Uc^1 \ar@<1ex>[l] \ar@<-1ex>[l] \ar@<1ex>[r] \ar@<-1ex>[r] &
  \Uc^0 \ar[l] \ar[r] & \X}$$

Soit $\Fc$ un faisceau ab\'elien sur le site lisse-\'etale de $\X$. On
note $\Fc^i$ la restriction de $\Fc$ au site lisse-étale de $\Uc^i$.
On cherche alors \`a calculer
la cohomologie de $\Fc$ en fonction de la cohomologie des $\Fc^i$ sur les
$\Uc^i$.
On peut associer au champ alg\'ebrique simplicial ci-dessus un complexe de
\v{C}ech de la mani\`ere suivante. Pour $n\geq 2$, on note $p_{1\dots\hat{l}\dots n}$, o\`u la
notation $\hat{l}$ signifie que l'indice $l$ est omis, la projection
$\Uc^{n-1}\fleche \Uc^{n-2}$ qui \og oublie \fg\ le facteur d'indice $l$ de
$\Uc^{n-1}$. Par exemple $p_1$ et $p_2$ d\'esignent respectivement les premi\`ere
et seconde projections de $\Uc^0\times_{\X} \Uc^0$ sur $\Uc^0$. On d\'efinit alors le
complexe de \v{C}ech $S(H^q)$ comme \'etant le complexe :
$$\xymatrix{H^q(\Uc^0,\Fc^0) \ar[r]^{d^0}& H^q(\Uc^1,\Fc^1) \ar[r]^-{d^1}&
  \dots \ar[r] & H^q(\Uc^p,\Fc^p) \ar[r]^-{d^p} & \dots}$$
avec
\begin{eqnarray*}
d^0 &=& p_2^*-p_1^*\\
d^1 &=& p_{23}^*-p_{13}^*+p_{12}^*\\
&\vdots& \\
d^p &=& \sum_{l=1}^{p+2} (-1)^{l+1} p_{1\dots \hat{l}\dots (p+2)}^*.
\end{eqnarray*}

On d\'esigne par $\check{H}^p(H^q(\Uc^{\bullet},\Fc^{\bullet}))$ l'homologie
en degr\'e $p$ de ce complexe :
$$\check{H}^p(H^q(\Uc^{\bullet},\Fc^{\bullet}))=\frac{\Ker d^p}{\Im d^{p-1}}.$$

Le r\'esultat suivant, qui donne la suite
spectrale reliant la \og cohomologie de \v{C}ech \fg\ relativement \`a la
famille couvrante $P : \Uc^0 \fleche \X$, \`a la cohomologie lisse-\'etale de $\Fc$
sur $\X$, est essentiellement trivial. Il ne fait pas appel aux techniques de descente cohomologique de Deligne
présentées par Saint-Donat dans l'exposé V~bis de \cite{SGA4_2} : il s'agit juste de la \og suite spectrale
de Cartan-Leray relative à un recouvrement\fg. Il nous rendra cependant de précieux services.
Le cas qui nous intéressera le plus en pratique sera celui où
$\Uc$ est un espace algébrique, de sorte que $P$ est une présentation de $\X$.

\begin{sousprop}[\cite{SGA4_2} V (3.3)]
\label{suite_spectrale_de_descente}
Reprenons les hypoth\`eses et notations pr\'ec\'edentes. Il existe une suite
spectrale (fonctorielle) :
$$E_2^{p,q}=\check{H}^p(H^q(\Uc^{\bullet},\Fc^{\bullet})) \Rightarrow H^{p+q}(\X,\Fc).$$
\end{sousprop}
\begin{demo}
Pour pouvoir appliquer le corollaire V~(3.3) de \cite{SGA4_2}, on regarde $\Fc$ comme un faisceau
sur le site lisse-lisse champêtre de $\X$. 
\end{demo}

\subsection{Faisceaux acycliques}
\label{Faisceaux_acycliques}

Soit $f : \X \fleche \Y$ un morphisme de $S$-champs algébriques. Dans la mesure où le couple de foncteurs
adjoints $(f^{-1},f_*)$ n'est pas un morphisme de topos, il n'est pas évident \emph{a priori} que le foncteur
$f_*$ transforme les faisceaux injectifs en faisceaux injectifs. C'est essentiellement à ce défaut que le présent
paragraphe doit son existence. En effet, pour obtenir la suite spectrale de Leray relative à un morphisme
de champs algébriques (cf. paragraphe (\ref{paragraphe_ss_Leray})), il nous faudra montrer que le foncteur $f_*$
transforme les faisceaux injectifs en faisceaux acycliques pour le foncteur \og sections globales\fg.
L'usage de faisceaux \og flasques\fg\ en un certain sens va nous permettre de résoudre ce problème, mais il faut bien choisir la
classe de faisceaux que l'on considère. En effet, vu que le topos $\X_{\text{lis-ét}}$ n'est pas engendré par ses ouverts
(la catégorie des ouverts ne contient pas assez de monomorphismes) une définition naïve des faisceaux flasques,
comme celle adoptée par Godement (\cite{Godement_faisceaux}) dans le cas de faisceaux sur un espace topologique, ne convient pas
puisque de tels faisceaux ne sont pas nécessairement acycliques (cf. \cite{SGA4_2}~V, exercice~4.16). Par
ailleurs, les faisceaux flasques proposés dans SGA~4 ne répondent pas non plus à nos besoins : pour
montrer qu'ils sont préservés par le foncteur $f_*$, on utilise là-aussi réellement l'exactitude du foncteur
$f^{-1}$ (\cite{SGA4_2}~V~4.9). Rappelons la définition ci-dessous.

\begin{sousdefi}[\cite{SGA4_2} V 4.2]
Soit $\X$ un $S$-champ algébrique et soit $\Fc$ un faisceau lisse-étale abélien. On dit que $\Fc$ est
$\gll$-acyclique si pour tout morphisme représentable et lisse $u : \Uc \fleche \X$ et pour tout $q>0$, le
groupe $H^q(\Uc, \Fc_{(\Uc,u)})$ est nul.
\end{sousdefi}

\begin{sousremarque} \rm S'il est évident que les faisceaux flasques au sens de SGA~4 sont $\gll$-acycliques, il n'y a aucune
raison \emph{a priori} pour que la réciproque soit vraie. Pour quelques commentaires sur ce genre de questions, on pourra
consulter \cite{SGA4_2}~V~4.6 et~4.13. Par ailleurs il est évident que les faisceaux injectifs sont flasques,
donc aussi $\gll$-acycliques.
\end{sousremarque}

La proposition suivante, quoique fortement inspirée de la proposition~V~4.3 de \cite{SGA4_2}, en diffère légèrement.

\begin{sousprop}
\label{prop_acyclicite}
Soient $\X$ un $S$-champ algébrique, et $\Fc$ un faisceau sur $\gll(\X)$. Les propositions suivantes sont équivalentes :
\begin{enumerate}
\item[(1)] $\Fc$ est $\gll$-acyclique ;
\item[(2)] pour toute famille couvrante $((\Uc',u') \fleche (\Uc,u))$ dans $\gll(\X)$, et pour tout $q>0$,
le groupe $H^q(\Uc'/\Uc,\Fc)$ est nul (où $H^q(\Uc'/\Uc,\Fc)$ désigne le $q^{\text{ième}}$ groupe de cohomologie de
\v{C}ech de $\Fc$, aussi noté $\check{H}^q(H^0(\Uc'^{\bullet},\Fc_{(\Uc,u)}^{\bullet}))$ dans la section précédente).
\end{enumerate}
\end{sousprop}
\begin{demo}
L'implication $(1)\Rightarrow (2)$ résulte de l'implication $(i)\Rightarrow (ii)$ de \cite{SGA4_2}~V~4.3. Notons
qu'il est très facile, si l'on préfère, de la redémontrer directement à partir de la suite
spectrale~(\ref{suite_spectrale_de_descente}). Pour la réciproque, on suppose que la condition (2) est vérifiée,
et on va montrer par récurrence sur $n$ que pour tout $n>0$ et pour tout $(\Uc,u)\in\ob\gll(\X)$,
le groupe $H^q(\Uc, \Fc_{(\Uc,u)})$ est nul. Soit $x\in H^1(\Uc, \Fc_{(\Uc,u)})$. On sait que pour
toute famille couvrante à un élément $\Uc'\fleche \U$, on a une suite spectrale
$$E_2^{p,q}=\check{H}^p(H^q(\Uc'^{\bullet},\Fc_{(\Uc,u)}^{\bullet})) \Rightarrow H^{p+q}(\Uc, \Fc_{(\Uc,u)}).$$
Comme le terme $E_2^{1,0}$ est nul par hypothèse, la suite exacte des termes de bas degré nous montre
que le morphisme $H^1(\Uc, \Fc_{(\Uc,u)}) \fleche E_2^{0,1}$ est injectif. \emph{A fortiori} le morphisme
$H^1(\Uc, \Fc_{(\Uc,u)})\fleche H^1(\Uc', \Fc_{(\Uc',u')})$ est lui aussi injectif. Or il existe une
famille couvrante $\Uc'\fleche \U$ telle que $x_{|_{\Uc'}}$ soit nul, d'où la nullité de $x$.
On suppose maintenant la condition vérifiée jusqu'au rang $n-1$. Si $x$ est un élément de
$H^n(\Uc, \Fc_{(\Uc,u)})$, on choisit une famille couvrante $\Uc'\fleche \U$ telle que $x_{|_{\Uc'}}$ soit nul.
Un examen attentif des zéros de la suite spectrale relative à cette famille couvrante
nous montre que le morphisme $H^n(\Uc, \Fc_{(\Uc,u)})\fleche H^n(\Uc', \Fc_{(\Uc',u')})$ est
injectif, ce qui prouve que $x$ est nul.
\end{demo}

\subsection{Images directes supérieures et suite spectrale de Leray relative à un morphisme de champs algébriques}
\label{paragraphe_ss_Leray}
\label{images_directes_supérieures}

Nous allons \`a pr\'esent essayer de d\'ecrire un peu mieux le foncteur
$R^qf_*$. Pour d\'emontrer la proposition (\ref{prop_image_directe_sup})
ci-dessous, nous
aimerions appliquer directement la proposition~V~(5.1) de \cite{SGA4_2}.
Telle quelle, elle ne permet malheureusement pas d'obtenir une description
satisfaisante des images directes sup\'erieures,
pour deux raisons. La premi\`ere est qu'elle suppose que le couple
$(f^{-1},f_*)$ est un morphisme de topos, ce qui, comme on l'a vu ci-dessus,
n'est pas vrai en g\'en\'eral. Cet obstacle est assez inoffensif, puisqu'un
examen attentif de la
d\'emonstration r\'ev\`ele que cette hypoth\`ese n'est pas essentielle, et que
le r\'esultat de \emph{loc. cit.} est valable dans notre cas (voir la d\'emonstration
de (\ref{prop_image_directe_sup})~1)). Le second obstacle est le manque d'ouverts
dans la topologie lisse-étale, résolu par l'usage du site lisse-lisse champêtre.

\begin{sousprop}
\label{prop_image_directe_sup}
\begin{enumerate}
\item[1)] Soit $\Fc$ un faisceau ab\'elien sur $\liset(\X)$. Alors le faisceau
$R^qf_*\Fc$ est le faisceau associ\'e au pr\'efaisceau qui \`a tout ouvert
lisse-\'etale $(U,u)$ de
$\Y$ associe $H^q(\X\times_{\Y} U, \Fc_{(\X\times_{\Y} U, \pr_{\X})})$.
$$\xymatrix{\X\times_{\Y}U \ar[r] \ar[d]_{\pr_{\X}} \cartesien& U\ar[d]^u \\
\X \ar[r]^f & \Y}$$
\item[2)] La restriction
de $(R^qf_*)\Fc$ au site \'etale de $\Y$ est le faisceau \'etale associ\'e au
pr\'efaisceau qui \`a tout ouvert \'etale $(U,u)$ de $\Y$ associe
$H^q(\X\times_{\Y} U, \Fc_{(\X\times_{\Y} U, \pr_{\X})})$.
\item[3)] Le faisceau \'evoqu\'e en 2) est encore la valeur en $\Fc$ du
$q^{\textrm{i\`eme}}$ foncteur d\'eriv\'e \`a droite du foncteur qui \`a un
faisceau lisse-\'etale ab\'elien 
$\Fc$ associe la restriction de $f_*\Fc$ au site \'etale de $\Y$.
\item[4)] Si $\Y=\Spec A$ est un schéma affine et si $\X$ est quasi-compact, alors pour
tout faisceau quasi-cohérent $\Fc$ sur $\X$ :
$$H^0(\Spec A, R^qf_*\Fc)\simeq H^q(\X,\Fc).$$
\end{enumerate}
\end{sousprop}
\begin{demo}
1) Vu la remarque (\ref{rem_foncteurs_compatibles}), on peut remplacer le site
lisse-\'etale par le site lisse-lisse champêtre.
Nous reprenons ici les arguments de \cite{SGA4_2}~V~(5.1).
Notons $\hat{f_*}$ le foncteur image directe de la cat\'egorie des
pr\'efaisceaux sur $\gll(\X)$ vers celle des pr\'efaisceaux sur $\gll(\Y)$,
$\underline{a}$ le foncteur \og faisceau associ\'e \fg, et $\Hc^0$ le foncteur
\og pr\'efaisceau sous-jacent \fg. Il est clair que l'on a un isomorphisme de
foncteurs :
$$f_* \simeq \underline{a}\hat{f_*} \Hc^0.$$
Comme les foncteurs $\underline{a}$ et $\hat{f_*}$ sont exacts, il en r\'esulte
un isomorphisme $R^qf_* \simeq \underline{a}\hat{f_*}R^q\Hc^0$. Il nous reste
juste \`a calculer $R^q\Hc^0$. Soit $\Fc$ un faisceau ab\'elien sur $\gll(\X)$.
Le pr\'efaisceau $R^q\Hc^0\Fc$ est l'homologie en degr\'e $q$ (calcul\'ee dans
la cat\'egorie des pr\'efaisceaux !) du complexe $0 \fleche \Ic^{\bullet}$ o\`u
$0\fleche \Fc \fleche \Ic^{\bullet}$ est une r\'esolution injective de $\Fc$.
Donc c'est le pr\'efaisceau qui \`a un gros ouvert lisse $(\Uc,u)$ de $\X$
associe l'homologie en degr\'e $q$ du complexe :
$$0\flechelongue \Gamma(\Uc,\Ic^0_{(\Uc,u)}) \flechelongue \Gamma(\Uc,\Ic^1_{(\Uc,u)})
\flechelongue \dots$$
Or les faisceaux $\Ic^n_{(\Uc,u)}$ sont injectifs (remarque
(\ref{restriction_faisceaux_injectifs})), donc ils forment une r\'esolution
injective de $\Fc_{(\Uc,u)}$, d'o\`u il r\'esulte que le $q^{\textrm{i\`eme}}$
groupe d'homologie du complexe ci-dessus n'est autre que
$H^q(\Uc,\Fc_{(\Uc,u)})$, ce qui prouve notre assertion.

2) r\'esulte imm\'ediatement de 1) et de la remarque
(\ref{remarque_faisceau_etale_associe}). Pour 3), il suffit de constater que
le foncteur $(.)_{\textrm{\'et}}$ de
restriction des faisceaux au site \'etale est exact, de sorte que
$R^q((.)_{\textrm{\'et}}\circ f_*) \simeq
(.)_{\textrm{\'et}}\circ (R^qf_{*})$.

4) On considère la suite spectrale de Leray relative au morphisme $f$. (On n'utilisera pas 4) avant (\ref{ss_Leray}).)
$$H^p(\Spec A, R^qf_*\Fc) \Rightarrow H^{p+q}(\X,\Fc)$$
Or $R^qf_*\Fc$ est quasi-cohérent (\cite{LMB}~(13.2.6)), donc pour tout $p>0$ et
pour tout entier~$q$, $H^p(\Spec A, R^qf_*\Fc)=0$. On en déduit les
isomorphismes escomptés (\cite{Cartan_Eilenberg}, chap. XV, théorème~5.12
p.~328).
\end{demo}

\begin{sousprop}
\label{images_directes_sup_compatibles}
La restriction du foncteur $R^qf_*$ \`a la cat\'egorie
$\Mod(\X)$ co\"{i}ncide avec le $q^{\textrm{i\`eme}}$ foncteur d\'eriv\'e \`a
droite de $f_* : \Mod(\X) \fleche \Mod(\Y)$.
\end{sousprop}
\begin{demo}
Il suffit de montrer que les objets injectifs de $\Mod(\X)$ sont acycliques pour
le foncteur $f_* : \Ab(\X) \fleche \Ab(\Y)$. Soit $\Ic$ un faisceau de
$\O_{\X}$-modules injectif et soit $q$ un entier strictement positif. Pour
montrer que $R^qf_* \Ic$ est nul, il suffit d'apr\`es la proposition
pr\'ec\'edente de montrer que pour tout ouvert lisse-\'etale $(U,u)$ de $\Y$, le
groupe $H^q(\X\times_{\Y} U, \Ic_{(\X\times_{\Y} U,\pr_{\X})})$ est nul. Or
les groupes de cohomologie du faisceau $\Ic_{(\X\times_{\Y} U,\pr_{\X})}$ vu
comme faisceau de modules ou comme faisceau ab\'elien co\"{i}ncident
(\ref{cohomologie_modules_egale_cohomologie_abelienne}), et ce faisceau est un
objet injectif de la cat\'egorie $\Mod(\X\times_{\Y} U)$ (cf. remarque
(\ref{restriction_faisceaux_injectifs})), ce qui prouve que sa cohomologie est
nulle.
\end{demo}

\begin{souslem}
\label{faisceaux_gll_acycliques_preserves}
Le foncteur $f_* : \Ab(\X)\fleche \Ab(\Y)$ préserve les faisceaux
$\gll$-acy\-cli\-ques. En particulier il transforme les faisceaux injectifs en faisceaux acycliques pour
le foncteur $\Gamma(\Y,.)$. 
\end{souslem}
\begin{demo}
Soit $\Fc$ un faisceau $\gll$-acyclique sur $\X$. Soit $((\Uc^0,u^0)\fleche (\Uc,u))$ une famille
couvrante à un élément de $\gll(\Y)$. On veut montrer que pour tout entier $p$ strictement positif le groupe $H^p(\Uc^0/\Uc,f_*\Fc)$ est nul.
On note $\Uc^{\bullet}$ le champ algébrique simplicial obtenu en prenant le cosquelette
du morphisme $\Uc^0 \fleche \Uc$. On note encore $\Vc$ (resp. $\Vc^0$) le produit fibré
$\X\times_{\Y} \Uc$ (resp. $\X\times_{\Y} \Uc^0$) et $\Vc^{\bullet}$ le champ algébrique simplicial obtenu en prenant le cosquelette
du morphisme $\Vc^0 \fleche \Vc$.  Il est canoniquement isomorphe à $\X\times_{\Y} \Uc^{\bullet}$.
$$\xymatrix{\dots & \Vc^1 \ar@<1ex>[r] \ar@<-1ex>[r]\ar[d]_{f_{\Uc}^1}\cartesien
   & \Vc^0 \ar[l] \ar[r]\ar[d]_{f_{\Uc}^0} \cartesien& \Vc \ar[r]^v\ar[d]_{f_{\Uc}}\cartesien& \X\ar[d]^{f}\\
\dots & \Uc^1 \ar@<1ex>[r] \ar@<-1ex>[r] & \Uc^0 \ar[l] \ar[r] & \Uc \ar[r]^u& \Y}$$
Par définition, le groupe $H^p(\Uc^0/\Uc,f_*\Fc)$ est l'homologie en degré $p$ du complexe
$$0\flechelongue H^0(\Uc^0,(f_*\Fc)^0) \flechelongue H^0(\Uc^1,(f_*\Fc)^1) \flechelongue \dots$$
où $(f_*\Fc)^i$ désigne la restriction à $\Uc^i$ du faisceau $f_*\Fc$.
Or pour tout $i$, le faisceau $(f_*\Fc)^i$ est exactement le faisceau $f_{\Uc*}^i(\Fc^i)$,
où $\Fc^i$ est la restriction à $\Vc^i$ de $\Fc$. On en déduit que le complexe ci-dessus n'est autre
que
$$0\flechelongue H^0(\Vc^0,\Fc^0) \flechelongue H^0(\Vc^1,\Fc^1) \flechelongue \dots$$
d'où un isomorphisme
$$H^p(\Uc^0/\Uc,f_*\Fc) \simeq H^p(\Vc^0/\Vc,\Fc).$$
L'assertion est maintenant claire puisque $\Vc^0\fleche \Vc$ est une famille couvrante de $\gll(\X)$.
\end{demo}

\begin{sousthm}[suite spectrale de Leray relative à $f : \X \fleche \Y$]
\label{ss_Leray}
Soit $$f : \X\flechelongue \Y$$ un morphisme de $S$-champs algébriques, et soit $\Fc$ un faisceau lisse-étale abélien sur
$\X$. Il existe une suite spectrale :
$$H^p(\Y,R^qf_*\Fc) \Rightarrow H^{p+q}(\X,\Fc).$$
\end{sousthm}
\begin{demo}
Compte tenu du lemme précédent, c'est juste la suite spectrale d'un foncteur composé
(cf. \cite{Lang_Algebra}~XX, théorème~(9.6)). Le fait que le foncteur $\Gamma(\X,.)$ soit bien le foncteur
composé $\Gamma(\Y,.)\circ f_*$ est évident lorsqu'on regarde $\Fc$ comme un faisceau sur le site lisse-lisse
champêtre de $\X$.
\end{demo}

\begin{sousprop}
Soient $f : \X \fleche \Y$ et $g : \Y \fleche \Zc$ des morphismes de $S$-champs algébriques et soit
$\Fc$ un faisceau lisse-étale abélien sur $\X$. Alors on a une suite spectrale :
$$R^pg_*R^qf_* \Fc \Rightarrow R^{p+q}(g\circ f)_* \Fc.$$
\end{sousprop}
\begin{demo}
C'est encore une fois la suite spectrale d'un foncteur composé. Il faut juste montrer que si $\Ic$ est
un faisceau injectif sur $\X$, alors $f_*\Ic$ est acyclique pour le foncteur $g_*$, autrement dit que
le faisceau $R^qg_*(f_*\Ic)$ est nul pour tout $q>0$. C'est évident en utilisant la proposition~(\ref{prop_image_directe_sup}),
le lemme~(\ref{faisceaux_gll_acycliques_preserves}) et la remarque~(\ref{restriction_faisceaux_injectifs}).
\end{demo}

\subsection{Cohomologie et changement de base}
\label{coh_et_chgt_de_base}

Nous aurons besoin par la suite d'hyperrecouvrements. La définition que, pour la commodité du lecteur,
nous rappelons ci-dessous, est plus simple en apparence que celle présentée dans \cite{SGA4_2}. Elle lui est
équivalente dès que la catégorie sous-jacente à $C$ admet des sommes directes arbitraires
(en effet les objets semi-représentables sont alors représentables), ce que nous supposons
dès à présent.

\begin{sousdefi}[cf. \cite{SGA4_2}~V (7.3.1.2)]
Soient $C$ un site (admettant des sommes directes arbitraires), $X$ un objet de $C$ et $U^{\bullet}\fleche X$
un objet simplicial de $C$ muni d'une augmentation vers $X$, ou de manière équivalente un objet simplicial
de $C/X$. On dit que $U^{\bullet}$ est un hyperrecouvrement de $X$ s'il possède les propriétés suivantes :
\begin{enumerate}
\item[(1)] Pour tout entier $n\geq 0$ le morphisme canonique
$$U^{n+1} \flechelongue (\cosq_n\,(\sq_n\,U^{\bullet}))_{n+1}$$
est un morphisme couvrant.
\item[(2)] Le morphisme $U^0 \fleche X$ est couvrant.
\end{enumerate}
\end{sousdefi}

\begin{sousremarque}\rm
Dans la définition précédente, $\sq_n$ désigne le foncteur qui à un objet simplicial associe son tronqué à l'ordre
$n$. Le foncteur $\cosq_n$ est l'adjoint à droite de $\sq_n$. Enfin $(\cosq_n\,(\sq_n\, U^{\bullet}))_{n+1}$ désigne
le terme d'indice $n+1$ de l'objet simplicial $\cosq_n\,(\sq_n\,U^{\bullet})$.
\end{sousremarque}

Si $U^{\bullet}\fleche X$ est un hyperrecouvrement de $X$ et si $\Fc$ est un faisceau abélien sur $X$, on lui
associe un complexe
$$0 \flechelongue H^q(U^0,\Fc^0) \flechelongue H^q(U^1,\Fc^1) \flechelongue  \dots$$
de manière tout à fait analogue à ce qui a été fait dans la section (\ref{annexe_desc_coh}) (où $\Fc^i$ désigne
la restriction de $\Fc$ à $U^i$). On note encore $\check{H}^p(H^q(U^{\bullet},\Fc^{\bullet}))$ l'homologie
en degr\'e $p$ de ce complexe. On a alors une suite spectrale (\cite{SGA4_2}~V~(7.4.0.3)) :
$$\check{H}^p(H^q(U^{\bullet},\Fc^{\bullet}))\Rightarrow H^{p+q}(X,\Fc).$$

\begin{sousprop}Soit $\X$ un $S$-champ alg\'ebrique quasi-compact, avec $S=\Spec A$ affine, et
soit $\Fc$ un faisceau quasi-coh\'erent sur $\X$.
Soient $A'$ une $A$-alg\`ebre plate, $S'=\Spec A'$, et $\X'=\X\times_S S'$.
On note $\Fc'$ l'image inverse de $\Fc$ sur $\X'$.
Alors pour tout $q\geq 0$ le morphisme naturel
$$\flechen{H^q(\X,\Fc)\otimes_A A'}{}{ H^q(\X', \Fc')}$$
est un isomorphisme.
\label{lemme_coh_et_chgt_base}
\end{sousprop}
\begin{demo}
C'est évident lorsque $\X$ est un schéma affine et $q=0$. Dans le cas général, soit $U^0\fleche \X$ une présentation
de $\X$, telle que $U^0$ soit un schéma affine. On pose $V^1=U^0\times_{\X}U^0$. Soit $W^1\fleche V^1$ une
présentation de l'espace algébrique $V^1$ dont la source $W^1$ est un schéma affine. On obtient alors un 
hyperrecouvrement tronqué à l'ordre 1 en posant $U^1=U^0 \coprod W^1$.
$$\xymatrix{U^1  \ar@<1ex>[r] \ar@<-1ex>[r] & U^0 \ar[l] \ar[r] & \X}$$
On note $U^{\bullet}$ le 1-cosquelette de ce diagramme. C'est clairement un
hyperrecouvrement\footnote{C'est même un hyperrecouvrement de type 1 au sens de \cite{SGA4_2}~V~(7.3.1.1).} de $\X$. De plus,
on peut voir facilement dans la construction du cosquelette (cf. \cite{Duskin} (0.8)) que pour tout $n\geq 0$ le champ algébrique
$U^{n+2}$ s'exprime en termes de produits fibrés obtenus à partir du diagramme :
$$\xymatrix{U^{n+1} \ar@<1.5ex>[r] \ar@<-1.5ex>[r]^{\vdots}  &U^n.}$$
On en déduit que pour tout $n\geq 0$, $U^n$ est un schéma affine. D'après \cite{SGA4_2}~V~(7.4.0.3)
on a une suite spectrale :
$$E_2^{p,q}=\check{H}^p(H^q(U^{\bullet},\Fc^{\bullet}))\Rightarrow H^{p+q}(X,\Fc).$$
Comme $\Fc$ est quasi-cohérent, on a $H^q(U^i,\Fc^i)=0$ pour tout $q>0$ et pour tout $i$, donc $E_2^{p,q}=0$ pour
tout $q>0$. On en déduit pour tout $p$ un isomorphisme :
$$\xymatrix{\check{H}^p(H^0(U^{\bullet},\Fc^{\bullet})) \ar[r]^-{\sim} & H^p(\X,\Fc).}$$
Par ailleurs l'objet simplicial $U^{\bullet}\times_S S'$ obtenu par changement de base est un hyperrecouvrement
de $\X'$, et les objets qui le composent sont des schémas affines. Donc on a aussi un isomorphisme :
$$\xymatrix{\check{H}^p(H^0(U^{\bullet}\times_S S',\Fc'^{\bullet})) \ar[r]^-{\sim} & H^p(\X',\Fc').}$$
Vu que $A'$ est plat sur $A$, l'opération qui consiste à prendre l'homologie en degré $p$ d'un complexe de
$A$-modules donné commute avec l'extension des scalaires à $A'$. On en déduit que dans le diagramme commutatif
$$\xymatrix{\check{H}^p(H^0(U^{\bullet},\Fc^{\bullet}))\otimes_A A' \ar[d] \ar[r]^-{\sim} & H^p(\X,\Fc)\otimes_A A'\ar[d]\\
    \check{H}^p(H^0(U^{\bullet}\times_S S',\Fc'^{\bullet})) \ar[r]^-{\sim} & H^p(\X',\Fc')}$$
la flèche verticale de gauche est un isomorphisme, donc celle de droite aussi.
\end{demo}

\begin{sousprop}
\label{images_directes_commutent_chgt_base_plat}
Soit $f : \X \fleche \Y$ un morphisme quasi-compact de $S$-champs alg\'ebriques, et
soit $\Fc$ un faisceau quasi-coh\'erent sur $\X$. Soit $u : \Y' \fleche \Y$ un
morphisme plat de changement de base.
$$\xymatrix{\X' \ar[r]^v \ar[d]_g \cartesien& \X\ar[d]^f \\
\Y' \ar[r]^u & \Y}$$
Alors pour tout $q\geq 0$ le morphisme naturel
$$u^*R^qf_*\Fc \flechelongue (R^qg_*)(v^*\Fc)$$
est un isomorphisme.
\end{sousprop}
\begin{demo}
Le cas où $u$ est lisse découle formellement de la proposition~(\ref{prop_image_directe_sup}).
On se ramène alors facilement au cas où $\Y$ et $\Y'$ sont des schémas affines. On les note
respectivement $S=\Spec A$ et $S'=\Spec A'$. Par \cite{LMB}~(13.2.6), le
faisceau $R^qf_*\Fc$ est quasi-cohérent. De plus on a des équivalences de
catégories :
$$\xymatrix{(A\text{-}\Mod) \simeq \Mod_{\textrm{qcoh}}(\O_{S_{\textrm{Zar}}})
\ar[r]^-{\sim}& \Mod_{\textrm{qcoh}}(\O_{S_{\textrm{lis-ét}}})}$$
$$\xymatrix{(A'\text{-}\Mod) \simeq \Mod_{\textrm{qcoh}}(\O_{S'_{\textrm{Zar}}})
\ar[r]^-{\sim}& \Mod_{\textrm{qcoh}}(\O_{S'_{\textrm{lis-ét}}})}$$
induites par les équivalences de \cite{LMB}~(13.2.3) et (13.1.2). Via ces
équivalences de catégories, le morphisme naturel
$$u^*R^qf_*\Fc \flechelongue (R^qg_*)(v^*\Fc)$$
correspond au morphisme
$$H^0(S,R^qf_*\Fc)\otimes_A A' \flechelongue H^0(S',R^qg_*(v^*\Fc)).$$
Il suffit alors d'appliquer~(\ref{prop_image_directe_sup})~4)
et~(\ref{lemme_coh_et_chgt_base}).
\end{demo}

\subsection{Cohomologie et extensions infinit\'esimales}
\label{par_coh_et_ext_inf}
\begin{souslem}
%\label{defm_sites_equivalents}
\label{coh_et_ext_inf}
Soit $i : \X \fleche \Xt$ une immersion ferm\'ee de champs alg\'ebriques
d\'efinie par un id\'eal quasi-coh\'erent de $\Xt$ de carr\'e nul, et soit $\Fc$
un faisceau ab\'elien sur $\X$. Alors pour tout $q>0$ le faisceau $R^qi_*\Fc$ est nul. En particulier
pour tout $n$ le morphisme naturel :
$$ H^n(\Xt, i_*\Fc)\flechelongue H^n(\X,\Fc) $$
est un isomorphisme.
\end{souslem}
\begin{demo}
La seconde assertion résulte clairement de la première grâce à la suite spectrale de Leray.
Pour obtenir la première, en vertu du lemme (\ref{lemme_nullite_faisceau_associe}) ci-dessous,
et vu la description du faisceau $R^qi_*\Fc$ donnée en (\ref{prop_image_directe_sup}), il suffit
de montrer que pour tout ouvert lisse-étale $(\widetilde{U}, \tilde{u})$ de $\Xt$ et pour tout $\xi\in H^q_{\text{lis-ét}}(U, \Fc_{(U,u)})$,
où $(U,u)$ désigne l'ouvert lisse-étale de $\X$ obtenu par changement de base, il existe un morphisme étale et surjectif
$\widetilde{V}\fleche \widetilde{U}$ tel que $\xi_{|_{\widetilde{V}}}$ soit nul. De plus on peut supposer que
$\widetilde{U}$ est un schéma.
Notons $j : U \fleche \widetilde{U}$ l'immersion fermée obtenue par changement de base à partir de $i$, et $j_*$ le
foncteur image directe de la catégorie des faisceaux étales abéliens sur $U$ vers la catégorie des faisceaux étales abéliens
sur $\widetilde{U}$. Comme $j$ est définie par un idéal quasi-cohérent de carré nul, elle induit une équivalence entre les sites
étales de $U$ et de $\widetilde{U}$, donc les faisceaux $R^qj_*\Fc_{U,u}$ sont nuls pour $q>0$ (où l'on rappelle que $\Fc_{U,u}$ désigne
la restriction de $\Fc_{(U,u)}$ au site \emph{étale} de $U$). Or $R^qj_*\Fc_{U,u}$ est le faisceau étale associé au préfaisceau
sur $\widetilde{U}_{\text{ét}}$ qui à $(\widetilde{V},\tilde{v})$ associe
\begin{eqnarray*}
H^q_{\text{ét}}(V,(\Fc_{U,u})_{V,v}) &=& H^q_{\text{ét}}(V,\Fc_{V,u\circ v})\\
&=&H^q_{\text{lis-ét}}(V,\Fc_{(V,u\circ v)})
\end{eqnarray*}
où les notations sont celles du diagramme suivant
$$\xymatrix{V\ar[r]^v \ar[d]\cartesien& U \cartesien\ar[r]^u \ar[d]_j & \X\ar[d]^i\\
\widetilde{V}\ar[r]_{\tilde{v}} & \widetilde{U}\ar[r]_{\tilde{u}} &\Xt.}$$
En particulier pour tout $\xi\in H^q_{\text{lis-ét}}(U, \Fc_{(U,u)})$, il existe un morphisme étale et surjectif
$\widetilde{V}\fleche \widetilde{U}$ tel que $\xi_{|_{\widetilde{V}}}$ soit nul, ce qui achève la preuve.
\end{demo}

\begin{souslem}
\label{lemme_nullite_faisceau_associe}
Soient $\X$ un $S$-champ algébrique, $F$ un préfaisceau sur $\liset(\X)$, et $\underline{a}F$ le faisceau associé. Alors
$\underline{a}F$ est nul si et seulement si pour tout $(U,u)\in\ob\liset(\X)$ (resp. pour tout $(U,u)\in\ob\liset(\X)$ avec $U$ un schéma)
et pour tout $\xi\in F(U,u)$ il existe un schéma $U'$ et un morphisme $U'\fleche U$ étale et surjectif tel que $\xi_{|_{U'}}$ soit nul.
\end{souslem}
\begin{demo}
Cela résulte trivialement de la construction du faisceau associé à un préfaisceau. La restriction au cas où $U$ est un schéma provient 
de l'équivalence entre les catégories de faisceaux sur $\liset(\X)$ et sur $\lisets(\X)$ (cf. \cite{LMB}~(12.1.2)).
\end{demo}

\begin{souscor}
Soit
$$\xymatrix{\X \ar[r]^i \ar[d]_f \cartesien &\Xt\ar[d]^{\widetilde{f}}\\ \Y \ar[r]^j &\Yt}$$
un carré 2-cartésien de $S$-champs algébriques, dans lequel on suppose que $i$ et $j$ sont des immersions fermées
définies par des idéaux quasi-cohérents de carrés nuls. Soit $\Fc$ un faisceau abélien sur $\X$. Alors on 
a pour tout $p$ un isomorphisme canonique :
$$\xymatrix{R^p\widetilde{f}_*(i_*\Fc) \ar[r]^{\sim} &j_*R^pf_*\Fc.}$$
\end{souscor}
\begin{demo}
On a une suite spectrale :
$$R^p\widetilde{f}_*(R^qi_*\Fc) \Rightarrow R^{p+q}(\widetilde{f}\circ i)_*\Fc.$$
Or pour tout $q>0$ le faisceau $(R^qi_*\Fc)$ est nul, donc on obtient pour tout $p$ un isomorphisme :
$$\xymatrix{R^p\widetilde{f}_*(i_*\Fc) \ar[r]^{\sim} &R^p(\widetilde{f}\circ i)_*\Fc.}$$
Maintenant le foncteur $(\widetilde{f}\circ i)_*$ est canoniquement isomorphe à $j_*\circ f_*$, et $j_*$
est exact grâce au lemme précédent, donc $R^p(\widetilde{f}\circ i)_* \simeq j_*R^pf_*$, ce qui achève la preuve.
\end{demo}

\subsection{Un résultat de descente}
\label{Un_résultat_de_descente}

L'objet de ce paragraphe est de démontrer le théorème suivant, qui n'est autre qu'un résultat de descente pour les champs algébriques.
\begin{sousthm}
\label{descente_pour_les_champs}
Soit $\Fc$ un faisceau sur $\gll(\X)$ qui, localement pour la topologie lisse sur $\X$, est représentable.
Alors $\Fc$ est représentable par un unique objet $(\Pc, p)$ de $\gll(\X)$. Autrement dit il existe un morphisme représentable et lisse
de $S$-champs algébriques $p : \Pc \fleche \X$ qui représente $\Fc$.
\end{sousthm}
\begin{demo}
L'hypothèse sur $\Fc$ signifie qu'il existe une présentation $x_0 : X_0\fleche \X$ de $\X$ 
telle que la restriction de $\Fc$ au site lisse de $X_0$ soit représentable par
un $X_0$-espace algébrique lisse $F_0$.  On note $X_1$ le produit fibré $X_0\times_{\X} X_0$ et $p_1, p_2$ les projections
canoniques de $X_1$ sur $X_0$. L'unicité de $\Pc$ est claire en vertu du lemme de Yoneda. Pour l'existence, nous allons construire en premier lieu une collection de faisceaux qui seront les fibres
de $\Pc$ au-dessus de chaque objet de $\X$.

\begin{souslem}
\label{lemme_descente_pour_les_champs}
Il existe une collection de $U$-espaces $F_{U,u}$ indexée par les couples
$(U,u)$ où $U\in\ob\aff$ et $u\in\ob\X_U$ et une collection d'isomorphismes
de transition
$$\theta_{\varphi,\alpha} : F_{V,v} \flechelongue F_{U,u}\times_U V$$
indexée par les diagrammes 2-commutatifs
$$\shorthandoff{!;:?}
\xymatrix@R=0.9pc@C=0.9pc{V \ar[rr]^{\varphi} \ar[rd]_v  & \raisebox{-3ex}{$^{\alpha} \FlecheNE$} & U \ar[ld]^u\\
& \X&}$$
avec les propriétés suivantes :
\begin{itemize}
\item[(i)] les $\theta_{\varphi,\alpha}$ vérifient une condition de cocycle
évidente ;
\item[(ii)] si $u$ est lisse (i.e. si $(U,u)$ est un objet du
site lisse-étale de $\X$) alors $F_{U,u}(U)$ est canoniquement isomorphe
à $\Fc(U,u)$.
\end{itemize}
\end{souslem}
\begin{demo}

\begin{etape}{Cas de $X_0$ et $X_1$.}
On pose $F_{X_0,x_0}=F_0$, $x_1=x_0\circ p_1$ et
$F_{X_1, x_1}=F_0\times_{X_0,p_1} X_1$.
$$\xymatrix{F_{X_1, x_1} \ar[r] \ar[d] \cartesien & F_{X_0,x_0}\ar[d] \\
X_1 \ar[r]_{p_1} &X_0 }$$
On a un diagramme 2-commutatif :
$$\shorthandoff{!;:?}
\xymatrix@R=0.9pc@C=0.9pc{X_1 \ar[rr]^{p_2} \ar[rd]_{x_1} & \raisebox{-3ex}{$^{\gamma} \FlecheNE$} & X_0 \ar[ld]^{x_0}\\
& \X&}$$
Il faut donc construire un 2-isomorphisme de $F_{X_1, x_1}$ vers $F_0\times_{X_0,p_2} X_1$. Or ces espaces algébriques
représentent tous les deux la restriction de $\Fc$ au site lisse-étale de $X_1$ (car $X_1$ est lisse sur $\X$),
l'un par $x_1$ et l'autre par $x_0\circ p_2$. Le 2-isomorphisme canonique $\gamma$ de $x_1$ vers $x_0\circ p_2$
induit un isomorphisme de $x_1^{-1}\Fc$ vers $(x_0\circ p_2)^{-1}\Fc$. De plus il est clair que le foncteur naturel
de la catégorie des espaces algébriques \emph{lisses} sur $X_1$ vers la catégorie des faisceaux lisses-étales sur $X_1$
est pleinement fidèle. L'isomorphisme évoqué précédemment provient donc d'un unique isomorphisme
$$\xymatrix{\theta_{p_2, \gamma} : F_{X_1,x_1} \ar[r]^{\sim} & F_0\times_{X_0,p_2} X_1.}$$
\end{etape}

\begin{etape}{Cas des objets qui se factorisent par $X_0$.}
Soit maintenant $u : U \fleche \X$ tel qu'il existe une factorisation $(\varphi, \alpha)$ de $u$ par $X_0$ :
$$\shorthandoff{!;:?}
\xymatrix@R=0.9pc@C=0.9pc{U \ar[rr]^{\varphi} \ar[rd]_u & \raisebox{-3ex}{$^{\alpha} \FlecheNE$} & X_0 \ar[ld]^{x_0}\\
& \X&}$$
On pose
$$F_{U,u,\varphi,\alpha}= F_{X_0,x_0}\times_{X_0,\varphi} U.$$
Si $(\varphi_1, \alpha_1)$ et $(\varphi_2, \alpha_2)$ sont deux factorisations différentes de $u$ par $X_0$,
elles induisent une factorisation de $u$ par $X_1$ et on en déduit des isomorphismes canoniques entre
$F_{U,u,\varphi_1, \alpha_1}$ et $F_{X_1,x_1}$ d'une part et entre $F_{U,u,\varphi_2, \alpha_2}$ et $F_{X_0,x_0}\times_{X_0,p_2} X_1$
d'autre part. \`A un isomorphisme canonique près, $F_{U,u,\varphi,\alpha}$ ne dépend donc pas de la factorisation
$(\varphi, \alpha)$ choisie et nous le noterons $F_{U,u}$. C'est un espace algébrique lisse sur $U$.

De plus, si $(U_1, u_1)$ et $(U_2, u_2)$ se factorisent par $X_0$, alors tout diagramme 2-commu\-tatif
$$\shorthandoff{!;:?}
\xymatrix@R=0.9pc@C=0.9pc{U_1 \ar[rr]^{\varphi} \ar[rd]_{u_1} & \raisebox{-3ex}{$^{\alpha} \FlecheNE$} & U_2 \ar[ld]^{u_2}\\
& \X&}$$
induit un isomorphisme canonique
$$\theta_{\varphi,\alpha} : F_{U_1,u_1} \flechelongue F_{U_2,u_2}\times_{U_2,\varphi} U_1.$$
\end{etape}

\begin{etape}{Construction d'un ensemble de sections $S(U,u)$ pour $(U,u)$ quelconque.}
 Pour tout $U\in\ob\aff$ et tout $u\in\ob\X_U$ on pose
$$S(U,u) =\lpro F_{U',u'}(U')$$
où la limite projective est prise sur les diagrammes 2-commutatifs
$$\xymatrix{U' \ar[d]_{\varphi} \ar[rd]^{u'} &\\
U \ar[r]_u & \X}$$
dans lesquels $U'$ est un espace algébrique, $\varphi$ est lisse et $u'$ se factorise par $X_0$.

\`A ce stade, nous pouvons faire deux observations importantes pour la suite.
\begin{description}
\item[Observation 1 :] Si $(U,u)$ lui-même se factorise par $X_0$, alors $S(U,u)$ est canoniquement isomorphe à $F_{U,u}(U)$.
\item[Observation 2 :] Si $u$ est lisse alors $S(U,u)$ s'identifie canoniquement à $\Fc(U,u)$. Il est en effet clair que $S(U,u)$ s'identifie au noyau des deux flèches :
$$\xymatrix{F_{U_0,u_0}(U_0)\ar@<0.5ex>[r] \ar@<-0.5ex>[r] &F_{U_1,u_1}(U_1)}$$
où $U_0$ désigne le produit fibré $U \times_{\X} X_0$ et $U_1$ le produit $U_0\times_U U_0$. Comme $U_0$ et $U_1$ sont lisses sur $X_0$ il résulte des définitions que $F_{U_0,u_0}(U_0)$ et $F_{U_1,u_1}(U_1)$ s'identifient respectivement à $\Fc(U_0,u_0)$ et $\Fc(U_1,u_1)$, si bien que $S(U,u)$ est isomorphe au noyau du couple de flèches :
$$\xymatrix{\Fc(U_0,u_0)\ar@<0.5ex>[r] \ar@<-0.5ex>[r]& \Fc(U_1,u_1)}$$
qui lui-même est isomorphe à $\Fc(U,u)$ puisque la famille $U_0 \fleche U$ est couvrante.
\end{description}
\end{etape}

\begin{etape}{Fonctorialité.}
 Soient $U,V\in\ob\aff$, $u$ un objet de $\X_U$ et $\psi$ un morphisme de $V$ dans $U$. On pose 
$v=u\circ \psi$. On va construire un morphisme de transition 
$$S(\psi) : S(U,u) \flechelongue S(V,v).$$
Soit $(s_{U',u'})$ un élément de $S(U,u)$. On veut lui associer un élément de $S(V,v)$.
Soit
$$\xymatrix{V' \ar[rrd]^{v'} \ar[d]_{\varphi} && X_0 \ar[d]\\ V\ar[r]_{\psi} &U \ar[r]_u &\X}$$
un diagramme 2-commutatif dans lequel $\varphi$ est lisse et $v'$ se factorise par $X_0$.
En notant $U_0$ le produit fibré de $U$ par $X_0$ on voit que se donner une factorisation
de $v'$ par $X_0$ revient à se donner une factorisation de $v'$ par $U_0$. On choisit un morphisme $f$ de $V'$ dans $U_0$ qui
factorise $v'$. La section $s_{U_0,u_0}$ de $F_{U_0,u_0}(U_0)$ induit alors via $f$ un élément $f^*s_{U_0,u_0}$ de $F_{V',v'}(V')$
que l'on note $t_{V',v'}$ (on rappelle que $F_{V',v'}=F_{U_0,u_0}\times_{U_0}V'$).
Il faut juste vérifier que $t_{V',v'}$ ne dépend pas de la factorisation $f$ choisie.
C'est évident car si $f_1$ et $f_2$ sont deux telles factorisations, alors vu que la famille $(s_{U',u'})$
est compatible on a nécessairement $f_1^*s_{U_0,u_0}=(f_1\times f_2)^*s_{U_1,u_1}=f_2^*s_{U_0,u_0}$ où $U_1$
est le produit fibré $U_0 \times_U U_0$.

Il est clair que les morphismes ainsi construits vérifient
$S(\psi_2\circ \psi_1)= S(\psi_1)\circ S(\psi_2)$ dès que $\psi_1$ et $\psi_2$ sont deux morphismes composables.
\end{etape}

\begin{etapefinale}{Construction d'un faisceau $F_{U,u}$ pour $(U,u)$ quelconque.}
 Soient $U\in\ob\aff$ et $u\in\ob\X_U$. On définit un préfaisceau $F_{U,u}$ sur $(\textrm{Aff}/U)$
en posant pour tout objet $V$ de $(\textrm{Aff}/U)$ :
$$F_{U,u}(V)=S(V,v)$$ où $v$ est le morphisme composé $V\fleche U \fleche \X$. Les morphismes de transition sont
donnés par les $S(\psi)$. L'observation~1 ci-dessus montre que dans le cas où $u$ se factorise par $X_0$,
le préfaisceau que nous venons de définir est canoniquement isomorphe au préfaisceau $F_{U,u}$ défini précédemment. L'observation~2 permet quant à elle de s'assurer que lorsque $u$
est lisse l'ensemble $F_{U,u}(U)$ des sections de $F_{U,u}$ au-dessus de $U$ s'identifie canoniquement à l'ensemble
$\Fc(U,u)$. Par ailleurs il est évident qu'un diagramme 2-commutatif de la forme
$$\shorthandoff{!;:?}
\xymatrix@R=0.9pc@C=0.9pc{V \ar[rr]^{\varphi} \ar[rd]_v & \raisebox{-3ex}{$^{\alpha} \FlecheNE$} & U \ar[ld]^u\\
& \X&}$$
induit un isomorphisme canonique $\theta_{\varphi,\alpha}$ de $F_{V,v}$ vers $F_{U,u}\times_U V$ et que ces
isomorphismes vérifient une condition de cocycle que seule la crainte de voir le lecteur céder à un agacement
bien compréhensible face à ce fleuve de trivialités nous dissuade d'expliciter. 

Pour en finir avec la démonstration de notre lemme, il reste juste à montrer que les préfaisceaux
$F_{U,u}$ que nous venons de construire sont des faisceaux pour la topologie étale. Ceci revient à montrer que
si $(U_i \fleche U)_i$ est une famille couvrante dans $\aff$ et si $u$ est un objet de $\X_U$ alors le diagramme
(avec les notations évidentes)
$$\xymatrix{S(U,u) \ar[r] & \prod_i S(U_i,u_i) \ar@<0.5ex>[r] \ar@<-0.5ex>[r]&\prod_{i,j} S(U_{i,j},u_{i,j})}$$
est exact.
On note $U_0$ (resp. $U_0^i$, $U_1$, $U_1^i$, $U_0^{i,j}$) le produit fibré $U \times_{\X} X_0$ (resp.
$U_i   \times_{\X} X_0$, $U_0\times_U U_0$, $U_1\times_U U_i$, $U_0 \times_U U_{i,j}$). On sait par ailleurs que $S(U,u)$ s'identifie au noyau des deux flèches :
$$\xymatrix{F_{U_0,u_0}(U_0) \ar@<0.5ex>[r] \ar@<-0.5ex>[r]&F_{U_1,u_1}(U_1)}.$$
Le résultat découle alors immédiatement du fait que $F_{U_0, u_0}$ est un faisceau étale et que 
les familles $(U_0^i \fleche U_0)$ et $(U_1^i \fleche U_1)$ sont couvrantes.
\end{etapefinale}
\end{demo}

Reprenons la démonstration du théorème \ref{descente_pour_les_champs}.

\begin{etape}{Construction d'un champ $\Pc$ candidat.}
On définit pour tout $U\in\ob\aff$ une catégorie $\Pc_U$ de la manière suivante.
Un objet est un couple $(x,s)$ où $x$ est un objet de la catégorie $\X_U$ et
$s$ est un élément de $F_{U,x}(U)$. Un morphisme entre deux tels objets
$(x,s)$ et $(x',s')$ est un 2-isomorphisme
$$\alpha : x \flechelongue x'$$
tel que l'isomorphisme $\theta_{\id_U,\alpha}$ de $F_{U,x}$ dans $F_{U,x'}$
envoie $s$ sur $s'$.

Ainsi défini, il est clair que $\Pc$ (muni des flèches de transition évidentes)
est un $S$-espace en groupoïdes. Le fait que $\Pc$ soit un $S$-champ résulte
alors facilement du fait que $\X$ lui-même en est un et que les
foncteurs $F_{U,x}$ sont des faisceaux étales.

Enfin pour tout $U\in\ob\aff$ on a un foncteur d'oubli évident de $\Pc_U$
vers $\X_U$. Ces foncteurs induisent un morphisme de $S$-champs :
$$p : \Pc \flechelongue \X$$.
\end{etape}

\begin{etape}{Le $S$-champ $\Pc$ est algébrique.}
Commençons par remarquer que pour tout $U\in\ob\aff$ et tout objet $u$ de
$\X_U$ le produit fibré $\Pc\times_{\X,u} U$ s'identifie canoniquement
au $U$-champ associé au $U$-espace $F_{U,u}$. En d'autres termes on a un carré
2-cartésien
$$\xymatrix{F_{U,u} \ar[r] \ar[d] \cartesien & U\ar[d]^u \\
\Pc \ar[r]_p &\X. }$$
La vérification de ce fait est aisée, quoique fastidieuse.
En particulier en prenant $U=X_0$ et $u=x_0$ on obtient un carré 2-cartesien
$$\xymatrix{F_{X_0,x_0} \ar[r] \ar[d] \cartesien & X_0\ar[d]^{x_0} \\
\Pc \ar[r]_p &\X. }$$

Il résulte alors de \cite{LMB}~(4.3.3) et~(4.5) que le morphisme $p$
est représentable et lisse et que le champ $\Pc$ est algébrique. Ceci montre au passage, \emph{a posteriori},
que chacun des $S$-espaces $F_{U,u}$ construits
ci-dessus est algébrique et lisse sur $U$.
\end{etape}

\begin{etapefinale}{L'objet $(\Pc,p)$ de $\gll(\X)$ représente le faisceau $\Fc$.}
Il faut construire, pour tout objet $(U,u)$ du site lisse-étale de $\X$, un isomorphisme canonique
$$\xymatrix{\Fc(U,u) \ar[r]^-{\sim} &\Hom_{\gll(\X)}((U,u),(\Pc,p)).}$$
Un élément du membre de droite est représenté par un triplet $(x,s,\alpha)$ où $x$ est un objet de $\X_U$, $s$ est un élément de $F_{U,x}(U)$ et $\alpha$ est un 2-isomorphisme
entre $u$ et $x$. Deux triplets $(x,s,\alpha)$ et $(x',s',\alpha')$ représentent le même élément
s'il existe un 2-isomorphisme $\beta$ de $x$ vers $x'$ tel que $\theta_{\id_U, \beta}$ envoie $s$ sur $s'$
et tel que $\beta \circ \alpha=\alpha'$.
On construit un morphisme de $F_{U,u}(U)$ dans cet ensemble en associant à une section $s$ de $F_{U,u}(U)$ l'élément représenté par
le triplet $(u, s,\id_U)$. Il est immédiat qu'il s'agit là d'un  isomorphisme, fonctoriel en $(U,u)$, si bien qu'en le composant
avec l'isomorphisme~(ii) du lemme~(\ref{lemme_descente_pour_les_champs}) nous obtenons le résultat voulu.
\end{etapefinale}
\end{demo}

\subsection{Cohomologie et torseurs}
\label{Cohomologie_et_torseurs}

Soit $\X$ un $S$-champ algébrique et soit $\Fc$ un faisceau lisse-étale abélien sur $\X$. On aimerait calculer
le groupe $H^1(\X,\Fc)$ en termes de torseurs. Bien évidemment, le résultat très général de Giraud
(\cite{Giraud}~III~3.5.4) s'applique, et l'on sait déjà que le groupe $H^1(\X,\Fc)$ s'identifie au groupe des classes
à isomorphisme près de $\Fc$-torseurs du topos $\widetilde{\liset(\X)}$ des faisceaux sur $\liset(\X)$. (Pour la définition
d'un torseur d'un topos, nous renvoyons à \emph{loc. cit.}~III~1.4.1.) Nous allons dans les lignes qui suivent essayer de
remplacer les $\Fc$-torseurs du topos $\widetilde{\liset(\X)}$ par des objets offrant plus de prise à l'intuition
géométrique.

On sait (cf. \emph{loc. cit.}~III remarque~1.7.2) que si $E$ est un site standard, c'est-à-dire un site dont
la topologie est moins fine que la topologie canonique, et qui admet des produits fibrés finis,
alors la catégorie des torseurs du site $E$ (c'est-à-dire des objets de $E$ munis d'une action d'un groupe $G$
de $E$ vérifiant les conditions que l'on imagine) est équivalente à la sous-catégorie pleine de la catégorie
des torseurs du topos $\widetilde{E}$ formée des torseurs représentables. Il en résulte que si $G$ est un objet en groupes
de $E$ tel que tout $G$-torseur de $\widetilde{E}$ soit représentable, la catégorie
des $G$-torseurs de $E$ est équivalente à la catégorie des $G$-torseurs de $\widetilde{E}$.

Malheureusement, le site $\liset(\X)$ n'est pas un site standard : il n'a pas d'objet final, donc pas de
produits fibrés finis. On s'en sort encore une fois en utilisant le site lisse-lisse champêtre $\gll(\X)$,
qui, lui, est bien un site standard, comme le lecteur le vérifiera facilement.
Le théorème~(\ref{descente_pour_les_champs}) montre que tout faisceau localement représentable sur $\gll(\X)$ est représentable.
En particulier si $G$ est un espace algébrique en groupes lisse sur $S$, alors tout $G$-torseur du topos
$\widetilde{\gll(\X)}$ est représentable par un objet du site $\gll(\X)$.

La définition suivante
ne fait qu'expliciter ce qu'est un $G$-torseur de $\gll(\X)$, dans le cas où le groupe $G$ provient de la base $S$.

\begin{sousdefi}
\label{def_torseur}
Soient $S$ un schéma, $G$ un $S$-espace algébrique en groupes lisse sur $S$ et $\X$ un $S$-champ algébrique.
Soit $p : \Pc \fleche \X$ un 1-morphisme représentable et lisse de $S$-champs algébriques. Une action de $G$
sur $(\Pc,p)$ est un quadruplet $(\mu,\varphi_{\mu}, \varphi_e, \varphi_{\text{ass}})$, où
$\mu$ est un 1-morphisme de $G\times_S \Pc$ vers $\Pc$, et où $\varphi_{\mu}, \varphi_e$ et $\varphi_{\text{ass}}$
sont des 2-isomorphismes faisant 2-commuter les diagrammes suivants :
$$\xymatrix{G\times_S \Pc \ar[r]^-{\mu} \ar[d]_{\pr_2} \ar@{}[dr]|{\varphi_{\mu} \FlecheNE}& \Pc \ar[d]^p \\
\Pc\ar[r]_p & \X}$$
$$\shorthandoff{!;:?}
\xymatrix@R=0.9pc@C=0.9pc{\Pc \ar[rr]^-{e\times\Id_{\Pc}} \ar[rd]_{\Id_{\Pc}}
& \raisebox{-3ex}{$^{ \varphi_e} \FlecheNE$}
&G\times_S\Pc \ar[ld]^{\mu}\\ & \Pc&}$$
$$\xymatrix{G\times_S G \times_S \Pc \ar[r]^-{\Id_G\times \mu} \ar[d]_{m_G\times \Id_{\Pc}}
 \ar@{}[dr]|{\varphi_{\text{ass}} \FlecheNE}&
G\times_S \Pc  \ar[d]^{\mu} \\ G\times_S \Pc \ar[r]_{\mu} & \Pc}$$
les 2-isomorphismes $\varphi_e$ et $\varphi_{\text{ass}}$ étant assujettis aux conditions de compatibilité, que nous nous dispenserons d'écrire, assurant
que ce sont des 2-morphismes dans la 2-catégorie des ouverts lisses champêtres de $\X$ (cf. début du paragraphe~(\ref{site_llc})).

On dit qu'un tel couple $(\Pc, p)$ muni d'une action de $G$ est un $G$-torseur sur $\X$ si les deux conditions supplémentaires
suivantes sont vérifiées :
\begin{itemize}
\item[(a)] $p$ est surjectif ;
\item[(b)] le morphisme naturel
$$G\times_S \Pc \flechelongue \Pc \times_{\X}\Pc$$
induit par le triplet $(\mu, \pr_2, \varphi_{\mu})$ est un isomorphisme.
\end{itemize}
\end{sousdefi}

La discussion précédente permet alors d'affirmer :

\begin{sousprop}[\cite{Giraud}]
\label{thm_coh_torseur}
Soient $S$ un schéma, $G$ un $S$-espace algébrique en groupes commutatifs lisse sur $S$ et $\X$ un $S$-champ algébrique.
Alors le groupe $H^1(\X,G)$ est canoniquement isomorphe au groupe des classes d'isomorphie de $G$-torseurs sur $\X$
(muni de la loi de groupe induite par le produit contracté de torseurs, cf.~\cite{Giraud}~III~2.4.5).
\end{sousprop}

\begin{sousremarque}\rm
Comme d'habitude, un $G$-torseur $p : \Pc \fleche \X$ est trivial si et seulement si le morphisme structural $p$ a une section.
\end{sousremarque}

\section{Cohomologie fppf}
\label{Cohomologie_plate}
\subsection{Sorites sur la cohomologie plate}

Soit $\X$ un $S$-champ algébrique. On définit le gros site \emph{fppf} de $\X$,
noté $\fppf(\X)$, de la manière suivante. Les ouverts sont les couples $(U,u)$,
où $U$ est un espace algébrique et $u : U \fleche \X$ est un morphisme
localement de présentation finie. Un morphisme entre deux tels ouverts $(U,u)$
et $(V,v)$ est un couple $(\varphi, \alpha)$ où $\varphi : U \fleche V$ est un
morphisme d'espaces algébriques et $\alpha$ est un 2-isomorphisme de $u$ vers
$v\circ \varphi$. Une famille couvrante est une collection de morphismes
$((\varphi_i, \alpha_i) : (U_i,u_i)\fleche (U,u))_{i\in I}$ telle que le
morphisme
$$\coprod_{i\in I} \varphi_i : \coprod_{i\in I} U_i \flechelongue U$$
soit fidèlement plat et localement de présentation finie.

\begin{sousremarque}\rm
Il est évident que l'on aurait obtenu un topos équivalent en ne prenant pour
ouverts que les couples $(U,u)$ où $U$ est un schéma (appliquer le lemme de
comparaison de SGA~4, \cite{SGA4_1}~III~4.1). En particulier si $\X$ est un
schéma on retrouve le topos des faisceaux sur le gros site \emph{fppf}
usuel (considéré par exemple dans \cite{Dix}, exposé~VI, paragraphe~5, p.124).
\end{sousremarque}

On définit de manière évidente les faisceaux d'anneaux $\O_{\X}$ et $\Z$ sur
$\fppf(\X)$. On note $\X_{\pl}$ le topos des faisceaux sur $\fppf(\X)$. Si
$\Ac$ est un anneau du topos $\X_{\pl}$, on note $\Mod_{\Ac}^{\pl}(\X)$ la
catégorie des faisceaux de modules sur le site annelé $(\fppf(\X),\Ac)$. Elle
sera notée $\Mod^{\pl}(\X)$ lorsque $\Ac=\O_{\X}$, et $\Ab^{\pl}(\X)$
lorsque $\Ac=\Z$.

On a pour tout $\Ac$ un foncteur d'oubli évident
$$\Mod_{\Ac}^{\pl}(\X)\flechelongue \Ab^{\pl}(\X).$$
On définit le foncteur \og sections globales\fg\ en posant pour tout faisceau
$\Fc$ sur $\fppf(\X)$ :
$$\Gamma_{\pl}(\X,\Fc)=\lpro \Fc(U,u)$$
où la limite projective est prise sur l'ensemble des couples $(U,u)$ de
$\fppf(\X)$. Il est clair que ce foncteur commute aux limites projectives
quelconques. En particulier il est exact à gauche. De même que dans le cas des
faisceaux lisses-étales, il résulte de SGA~4 (\cite{SGA4_1}~II~6.7) que la
catégorie $\Mod_{\Ac}^{\pl}(\X)$ est une catégorie abélienne avec suffisament
d'objets injectifs. On définit alors $H^i_{\pl}(\X,.)$ comme étant le
$i^{\textrm{ième}}$ foncteur dérivé à droite de $\Gamma_{\pl}(\X,.) : 
\Ab^{\pl}(\X)\fleche \Ab$. Il coïncide sur $\Mod_{\Ac}^{\pl}(\X)$ avec le
$i^{\textrm{ième}}$ foncteur dérivé à droite de $\Gamma_{\pl}(\X,.) : 
\Mod_{\Ac}^{\pl}(\X)\fleche (\Gamma_{\pl}(\X,\Ac)$-$\Mod)$
(\cite{SGA4_2} V 3.5).

\medskip
\noindent
{\sc Fonctorialité.}

Soit $f: \X\fleche \Y$ un 1-morphisme de $S$-champs algébriques. On lui associe
un couple de foncteurs adjoints $(f_{\pl}^{-1},f_*^{\pl})$ d'une manière tout à
fait analogue à ce qui a été fait au paragraphe
(\ref{fonctorialite_lisse_etale}) pour les faisceaux lisses-étales. Le foncteur
$f_*^{\pl}$ est alors exact à gauche, et on note $R^if_*^{\pl}$ ses foncteurs
dérivés à droite. 

\begin{sousremarque}\rm
Contrairement à ce qu'il se passe dans le cas des faisceaux lisses-étales, la limite inductive qui définit le
foncteur image inverse $f_{\pl}^{-1}$ \emph{est} filtrante. Ceci est essentiellement dû au fait qu'un morphisme
entre deux objets localement de présentation finie est lui-même localement de présentation finie (ce qui se
produisait aussi pour le site étale, mais n'était plus vrai en remplaçant étale par lisse). On en déduit
(cf. par exemple \cite{Milne_Etale_coh} annexe~A) que le foncteur $f_{\pl}^{-1}$ est exact, et donc que le
couple $(f_{\pl}^{-1},f_*^{\pl})$ est un morphisme de topos.
\end{sousremarque}

Notons ici aussi deux cas particuliers dans lesquels les foncteurs image directe
ou image inverse ont une expression plus simple. Lorsque $f$ est localement de
présentation finie, le foncteur $f_{\pl}^{-1}$ est simplement le foncteur de
restriction au site \emph{fppf} de $\X$ via le foncteur $(U,u)\mapsto (U,f\circ
u)$. Si $f$ est représentable, le foncteur $f_*^{\pl}$ provient d'une
application continue $\fppf(\X) \fleche \fppf(\Y)$, qui à un ouvert \emph{fppf}
$u : U\fleche \Y$ associe l'ouvert formé de l'espace algébrique $U\times_{\Y}\X$
muni de la projection sur $\X$. Naturellement ceci n'est pas vrai si $f$ n'est
pas représentable, puisqu'alors $U\times_{\Y}\X$ n'est pas nécessairement un
espace algébrique. Ce défaut s'avère gênant dans le calcul des images directes
supérieures, et motive l'introduction du site \emph{fppf} champêtre ci-dessous.

\medskip
\noindent
{\sc Le gros site \emph{fppf} champêtre d'un $S$-champ algébrique.}

Soit $\X$ un $S$-champ algébrique. On définit le gros site \emph{fppf} champêtre
de $\X$ de la manière suivante. Un ouvert de $\X$ est un couple  $(\Uc,u)$
o\`u $\Uc$ est un $S$-champ alg\'ebrique et $u : \Uc \fleche \X$ est
un morphisme repr\'esentable et localement de présentation finie. Un 1-morphisme
de $(\Uc,u)$ vers $(\Vc,v)$ est un couple $(\varphi,\alpha)$ o\`u
$\varphi : \Uc \fleche \Vc$ est un 1-morphisme de $S$-champs alg\'ebriques et
$\alpha$ est un 2-isomorphisme de $u$ dans $v\circ \varphi$.
Si $(\varphi,\alpha)$ et $(\psi,\beta)$ sont deux 1-morphismes de $(\U,u)$ dans
$(\Vc,v)$, un 2-morphisme entre $(\varphi,\alpha)$ et $(\psi,\beta)$ est un
2-isomorphisme $\gamma : \varphi \Rightarrow \psi$ tel que
$\beta=(v_*\gamma)\circ \alpha$.

L'analogue \emph{fppf} du lemme (\ref{lemme_site_champetre}) est alors valable,
de sorte que la 2-catégorie que nous venons de décrire est en fait équivalente à
une catégorie, que nous appellerons la catégorie des ouverts \emph{fppf}
champêtres de $\X$. On définit maintenant le gros site \emph{fppf}
champêtre de $\X$ de manière évidente, et on le note $\fppfc(\X)$.

\begin{sousprop}
Le foncteur d'inclusion $\fppf(\X) \fleche \fppfc(\X)$ induit une équivalence de
topos de la catégorie des faisceaux sur $\fppfc(\X)$ vers la catégorie des
faisceaux sur $\fppf(\X)$.
\end{sousprop}
\begin{demo}
Il suffit encore une fois d'appliquer \cite{SGA4_1}~III~4.1.
\end{demo}

On voit alors facilement que via cette équivalence de catégories, le foncteur
image directe $f_*^{\pl}$ provient d'une application continue
$$\fppfc(\X) \flechelongue \fppfc(\Y)$$
qui envoie un ouvert $u : \Uc\fleche \Y$ de $\Y$ sur l'ouvert $\X\times_{\Y}
\Uc\fleche \X$ de $\X$. Maintenant les mêmes démonstrations que celles des
propositions~(\ref{prop_image_directe_sup})~1)
et~(\ref{images_directes_sup_compatibles}) permettent d'obtenir le résultat
suivant.

\begin{sousprop} 
\label{images_directes_superieures_fppf}
Soit $f : \X \fleche \Y$ un morphisme de $S$-champs
algébriques.
\begin{itemize}
\item[1)] Soit $\Fc$ un faisceau abélien sur $\fppf(\X)$. Alors le faisceau
$R^qf_*^{\pl} \Fc$ est le faisceau associé au préfaisceau qui à tout ouvert
\emph{fppf} $u : U\fleche \Y$ de $\Y$ associe $H^q_{\pl}(\X\times_{\Y} U,
\Fc_{(\X\times_{\Y} U,\pr_{\X})})$.
\item[2)] La restriction du foncteur $R^qf_*^{\pl}$ à la catégorie
$\Mod^{\pl}(\X)$ coïncide avec le $q^{\textrm{ième}}$ foncteur dérivé à droite
du foncteur $f_*^{\pl} : \Mod(\X) \fleche \Mod(\Y)$.\ $\square$
\end{itemize}
\end{sousprop}

\subsection{Comparaison avec la cohomologie lisse-étale}
Soit $\X$ un $S$-champ algébrique. On a une application continue évidente :
$$p : \fppf(\X) \flechelongue \liset(\X)$$
induite par le foncteur d'inclusion de $\liset(\X)$ dans $\fppf(\X)$.
En particulier $p$ induit un couple de foncteurs adjoints :
$$\begin{array}{c} p_* : \X_{\pl} \flechelongue \X_{\text{lis-ét}}\\
p^{-1} : \X_{\text{lis-ét}}\flechelongue \X_{\pl}.
\end{array}$$
Le foncteur $p^{-1}$ peut être décrit de la manière suivante. On définit d'abord un adjoint à gauche au foncteur
$p_*$ pour les préfaisceaux, que l'on note $\widehat{p^{-1}}$, en associant à tout préfaisceau $\Fc$ sur $\liset(\X)$
le préfaisceau sur $\fppf(\X)$ qui à tout ouvert \emph{fppf} $u : U \fleche \X$ associe
$$(\widehat{p^{-1}}\Fc)(U,u) = \lind \Fc(V,v)$$
où la limite inductive est prise sur l'ensemble des diagrammes 2-commutatifs :
$$\shorthandoff{!;:?}
\xymatrix@R=0.9pc@C=0.9pc{U \ar[rr]^{\varphi} \ar[rd]_u  & \raisebox{-3ex}{$^{\alpha} \FlecheNE$} & V \ar[ld]^v\\
& \X&}$$
avec $(V,v)\in \ob\liset(\X)$. On définit alors $p^{-1}$ en posant $p^{-1}=\underline{a}\widehat{p^{-1}}$ où
$\underline{a}$ est le foncteur \og faisceau associé\fg. Encore une fois, on voit que $\widehat{p^{-1}}$ est défini par
une limite inductive sur $\liset(\X)$ qui n'est pas filtrante. Le couple $(p^{-1},p_*)$ n'est donc pas un morphisme de
topos a priori, et il n'y a aucune raison formelle pour que $p_*$ transforme les objets injectifs de $\Ab^{\pl}(\X)$ en
injectifs de $\Ab(\X)$, ce qui nous oblige à travailler un peu plus pour obtenir la suite spectrale (\ref{ss_comparaison_fppf_vs_liset})
ci-dessous.

\begin{souslem}
Si $\Fc$ est un objet injectif de $\Ab^{\pl}(\X)$, alors $p_*\Fc$ est un faisceau acyclique pour le foncteur
$\Gamma_{\text{\rm lis-ét}}(\X,.)$.
\end{souslem}
\begin{demo}
Soit $\Fc$ un objet injectif de $\Ab^{\pl}(\X)$. En particulier pour tout ouvert plat
champêtre $(\Uc,u)$ et pour tout $q>0$ on a $H^q_{\pl}(\Uc,\Fc_{(\Uc,u)})=0$. Remarquons que
la proposition~(\ref{prop_acyclicite}) ne faisait intervenir que des propriétés générales 
des faisceaux sur un site, et elle admet un analogue évident pour les faisceaux sur le
site $\fppfc(\X)$. En particulier pour toute famille couvrante $\Uc'\fleche \Uc$
de $\fppfc(\X)$ et pour tout $q>0$, le groupe $H^q(\Uc'/\Uc,\Fc)$ est nul.
Maintenant, toute famille couvrante  $\Uc'\fleche \Uc$ de $\gll(\X)$ est une famille
couvrante de $\fppfc(\X)$, et on a pour tout $q>0$
$$H^q(\Uc'/\Uc,p_*\Fc)=H^q(\Uc'/\Uc,\Fc)=0.$$
Donc d'après la proposition~(\ref{prop_acyclicite}) le faisceau $p_*\Fc$ est $\gll$-acyclique.
En particulier il est acyclique pour le foncteur $\Gamma_{\text{\rm lis-ét}}(\X,.)$.
\end{demo}

\begin{souslem}
Le morphisme naturel
$$\Gamma_{\pl}(\X,\Fc)\flechelongue\Gamma_{\text{\rm lis-ét}}(\X,p_*\Fc)$$
est un isomorphisme.
\end{souslem}
\begin{demo}
C'est complètement évident quand on regarde $\Fc$ et $p_*\Fc$ comme des faisceaux sur les sites champêtres $\fppfc(\X)$ et $\gll(\X)$.
\end{demo}

On déduit des deux résultats précédents une suite spectrale de Leray, qui résume les relations générales
entre cohomologie lisse-étale et cohomologie \emph{fppf}.

\begin{sousthm}
\label{ss_comparaison_fppf_vs_liset}
Soit $\X$ un $S$-champ algébrique, et soit $\Fc$ un faisceau abélien sur $\fppf(\X)$. On a alors une suite
spectrale (fonctorielle en $\Fc$) :
$$H^p_{\text{\rm lis-ét}}(\X,R^qp_*\Fc) \Rightarrow H^{p+q}_{\pl}(\X,\Fc). $$
\vskip-8mm\hfill $\square$
\end{sousthm}

\begin{sousremarque}\rm
\label{rmq_comparaison_coh_fppf_vs_liset}
On en déduit en particulier des morphismes canoniques :
$$H^p_{\text{lis-ét}}(\X,p_*\Fc) \flechelongue H^{p}_{\pl}(\X,\Fc).$$
De plus ces morphismes sont des isomorphismes si pour tout $q>0$, le faisceau lisse-étale $R^qp_*\Fc$ est nul.
\end{sousremarque}

Notre objectif à présent est de généraliser aux champs algébriques le résultat de Grothendieck (\cite{Dix}, exposé~VI,
paragraphe~11) selon lequel si $G$ est un groupe lisse sur un schéma $X$, le morphisme canonique
$$H^q_{\text{ét}}(X,p_*G) \flechelongue H^q_{\pl}(X,G)$$
est un isomorphisme pour tout $q$. Dans la mesure où nous n'utiliserons ce résultat que pour le groupe $\gm$, nous
n'avons pas cherché à le démontrer pour les \og champs en groupes\fg\ lisses sur $\X$ (notion qui resterait d'ailleurs
à définir), mais nous nous sommes contenté de considérer un groupe lisse sur la base $S$. Nous allons voir dans les
lignes qui suivent que dans ce cadre élémentaire le résultat se déduit assez facilement du cas des schémas.

\begin{sousthm}
\label{coh_fppf_groupe_lisse}
Soient $S$ un schéma, $G$ un schéma en groupes lisse sur $S$, et $\X$ un $S$-champ algébrique.
On note encore $G$ le faisceau défini sur le site $\fppf(\X)$ par :
$$\forall (U,u)\in\ob\fppf(\X)\quad G(U,u):=G(U,f\circ u)=\Hom_S(U,G).$$
Alors pour tout $q$ le morphisme canonique
$$H^q_{\text{\rm lis-ét}}(\X, p_* G) \flechelongue H^q_{\pl}(\X,G)$$
est un isomorphisme.
\end{sousthm}
\begin{demo}
Vu la remarque (\ref{rmq_comparaison_coh_fppf_vs_liset}), il suffit de montrer que pour tout $q>0$,
le faisceau lisse-étale $R^qp_*G$ est nul. Or par des arguments tout à fait analogues à ceux de
(\ref{prop_image_directe_sup}) on voit que $R^qp_*G$ est le faisceau lisse-étale associé au préfaisceau
qui à tout ouvert $u : U \fleche \X$ associe $H^q_{\pl}(U,G\times_S U)$ (on notera encore $H^q_{\pl}(U,G)$
ce dernier groupe). Alors d'après le lemme (\ref{lemme_nullite_faisceau_associe}) il suffit de montrer que pour tout morphisme lisse
$u : U\fleche \X$ où $U$ est un schéma et pour tout $\xi\in H^q_{\pl}(U,G)$ il existe un schéma $U'$
et un morphisme étale et surjectif $U'\fleche U$ tel que $\xi_{|_{U'}}$ soit nul.
Fixons un tel couple $(U,u)$. On note $p_U$ l'application continue :
$$p_U : \fppf(U) \flechelongue \et(U)$$
induite par le foncteur d'inclusion de $\et(U)$ dans $\fppf(U)$.
D'après le théorème~(11.7) de \cite{Dix}, pour tout $U'$ sur $U$ les morphismes canoniques
$$H^q_{\text{ét}}(U',G) \flechelongue H^q_{\pl}(U',G)$$
sont des isomorphismes. Puis d'après le lemme~(11.1) de \cite{Dix}, on en déduit que pour tout $q>0$,
le faisceau $R^qp_{U*}G$ est nul. Or ce dernier n'est autre que le faisceau étale associé au préfaisceau
$$\xymatrix{U' \ar@{|->}[r]& H^q_{\pl}(U',G)}.$$
En particulier (lemme~(\ref{lemme_nullite_faisceau_associe})) pour tout $\xi\in H^q_{\pl}(U,G)$ il existe
une famille couvrante étale $U'\fleche U$ telle que $\xi_{|_{U'}}$ soit nul.
\end{demo}

\bibliographystyle{latex/plain-fr}
\addcontentsline{toc}{chapter}{Bibliographie}
\bibliography{latex/mabiblio}

\begin{thebibliography}{10}
\expandafter\ifx\csname fonteauteurs\endcsname\relax
\def\fonteauteurs{\scshape}\fi

\bibitem{SGA3}
{\em Sch\'emas en groupes. {II}: {G}roupes de type multiplicatif, et structure
  des sch\'emas en groupes g\'en\'eraux}.
\newblock S\'eminaire de G\'eom\'etrie Alg\'ebrique du Bois Marie 1962/64 (SGA
  3). Dirig\'e par M. Demazure et A. Grothendieck. Lecture Notes in
  Mathematics, Vol. 152. Springer-Verlag, Berlin, 1962/1964.

\bibitem{Dix}
{\em Dix expos\'es sur la cohomologie des sch\'emas}.
\newblock Advanced Studies in Pure Mathematics, Vol. 3. North-Holland
  Publishing Co., Amsterdam, 1968.

\bibitem{SGA4_1}
{\em Th\'eorie des topos et cohomologie \'etale des sch\'emas. {T}ome 1:
  {T}h\'eorie des topos}.
\newblock Springer-Verlag, Berlin, 1972.
\newblock S\'eminaire de G\'eom\'etrie Alg\'ebrique du Bois-Marie 1963--1964
  (SGA 4), Dirig\'e par M. Artin, A. Grothendieck, et J. L. Verdier. Avec la
  collaboration de N. Bourbaki, P. Deligne et B. Saint-Donat, Lecture Notes in
  Mathematics, Vol. 269.

\bibitem{SGA4_2}
{\em Th\'eorie des topos et cohomologie \'etale des sch\'emas. {T}ome 2}.
\newblock Springer-Verlag, Berlin, 1972.
\newblock S\'eminaire de G\'eom\'etrie Alg\'ebrique du Bois-Marie 1963--1964
  (SGA 4), Dirig\'e par M. Artin, A. Grothendieck et J. L. Verdier. Avec la
  collaboration de N. Bourbaki, P. Deligne et B. Saint-Donat, Lecture Notes in
  Mathematics, Vol. 270.

\bibitem{SGA1}
{\em Rev\^etements \'etales et groupe fondamental ({SGA} 1)}.
\newblock Documents Math\'ematiques (Paris) , 3. Soci\'et\'e Math\'ematique de
  France, Paris, 2003.
\newblock S\'eminaire de g\'eom\'etrie alg\'ebrique du Bois Marie 1960--61.
  Dirigé par A. Grothendieck, augmenté de deux exposés de Mme M. Raynaud.
  \'Edition recomposée et annotée du volume 224 des Lectures Notes in
  Mathematics publié en 1971 par Springer-Verlag.

\bibitem{Abramovich_Vistoli_note}
Dan \bgroup\fonteauteurs\bgroup Abramovich\egroup\egroup{} et Angelo
  \bgroup\fonteauteurs\bgroup Vistoli\egroup\egroup{} :
\newblock Complete moduli for families over semistable curves.
\newblock http://fr.arxiv.org/abs/math/9811059, 1998.

\bibitem{Abramovich_Vistoli_CMFFS}
Dan \bgroup\fonteauteurs\bgroup Abramovich\egroup\egroup{} et Angelo
  \bgroup\fonteauteurs\bgroup Vistoli\egroup\egroup{} :
\newblock Complete moduli for fibered surfaces.
\newblock \emph{In} {\em Recent progress in intersection theory (Bologna,
  1997)}, Trends Math., pages 1--31. Birkh\"auser Boston, Boston, MA, 2000.

\bibitem{Abramovich_Vistoli_CSSM}
Dan \bgroup\fonteauteurs\bgroup Abramovich\egroup\egroup{} et Angelo
  \bgroup\fonteauteurs\bgroup Vistoli\egroup\egroup{} :
\newblock Compactifying the space of stable maps.
\newblock {\em J. Amer. Math. Soc.}, 15(1)\string:\penalty500\relax 27--75
  (electronic), 2002.

\bibitem{Aoki_defm}
Masao \bgroup\fonteauteurs\bgroup Aoki\egroup\egroup{} :
\newblock Deformation theory of algebraic stacks.
\newblock {\em Compos. Math.}, 141(1)\string:\penalty500\relax 19--34, 2005.

\bibitem{Aoki_erratum}
Masao \bgroup\fonteauteurs\bgroup Aoki\egroup\egroup{} :
\newblock Erratum: ``{H}om stacks''.
\newblock {\em Manuscripta Math.}, 121(1)\string:\penalty500\relax 135, 2006.

\bibitem{Aoki_Hom}
Masao \bgroup\fonteauteurs\bgroup Aoki\egroup\egroup{} :
\newblock Hom stacks.
\newblock {\em Manuscripta Math.}, 119(1)\string:\penalty500\relax 37--56,
  2006.

\bibitem{Global_Analysis_1}
Michael \bgroup\fonteauteurs\bgroup Artin\egroup\egroup{} :
\newblock Algebraization of formal moduli. {I}.
\newblock \emph{In} {\em Global Analysis (Papers in Honor of K. Kodaira)},
  pages 21--71. Univ. Tokyo Press, Tokyo, 1969.

\bibitem{Artin_implicit_fct_thm}
Michael \bgroup\fonteauteurs\bgroup Artin\egroup\egroup{} :
\newblock The implicit function theorem in algebraic geometry.
\newblock \emph{In} {\em Algebraic Geometry (Internat. Colloq., Tata Inst.
  Fund. Res., Bombay, 1968)}, pages 13--34. Oxford Univ. Press, London, 1969.

\bibitem{Artin_Montreal}
Michael \bgroup\fonteauteurs\bgroup Artin\egroup\egroup{} :
\newblock {\em Th\'eor\`emes de repr\'esentabilit\'e pour les espaces
  alg\'ebriques}.
\newblock Les Presses de l'Universit\'e de Montr\'eal, Montreal, Que., 1973.
\newblock En collaboration avec Alexandru Lascu et Jean-Fran\c cois Boutot,
  S\'eminaire de Math\'ematiques Sup\'erieures, No. 44 (\'Et\'e, 1970).

\bibitem{Artin_Versal_defm}
Michael \bgroup\fonteauteurs\bgroup Artin\egroup\egroup{} :
\newblock Versal deformations and algebraic stacks.
\newblock {\em Invent. Math.}, 27\string:\penalty500\relax 165--189, 1974.

\bibitem{Cadman_USTITCOC}
Charles \bgroup\fonteauteurs\bgroup Cadman\egroup\egroup{} :
\newblock Using stacks to impose tangency conditions on curves.
\newblock \`A paraître dans American Journal of Mathematics.
  http\!\!://www.math.lsa.umich.edu/$\sim$ cdcadman/research/stacks.pdf.

\bibitem{Cartan_Eilenberg}
Henri \bgroup\fonteauteurs\bgroup Cartan\egroup\egroup{} et Samuel
  \bgroup\fonteauteurs\bgroup Eilenberg\egroup\egroup{} :
\newblock {\em Homological algebra}.
\newblock Princeton University Press, Princeton, N. J., 1956.

\bibitem{Chiodo_twisted_curves_spin_structures}
Alessandro \bgroup\fonteauteurs\bgroup Chiodo\egroup\egroup{} :
\newblock Stable twisted curves and their $r$-spin structures.
\newblock http\!\!://front.math.ucdavis.edu/math.AG/0603687.

\bibitem{Conrad_thm_Chevalley}
Brian \bgroup\fonteauteurs\bgroup Conrad\egroup\egroup{} :
\newblock A modern proof of {C}hevalley's theorem on algebraic groups.
\newblock {\em J. Ramanujan Math. Soc.}, 17(1)\string:\penalty500\relax 1--18,
  2002.

\bibitem{Deligne_Mumford}
Pierre \bgroup\fonteauteurs\bgroup Deligne\egroup\egroup{} et David
  \bgroup\fonteauteurs\bgroup Mumford\egroup\egroup{} :
\newblock The irreducibility of the space of curves of given genus.
\newblock {\em Inst. Hautes \'Etudes Sci. Publ. Math.},
  (36)\string:\penalty500\relax 75--109, 1969.

\bibitem{Duskin}
John~W. \bgroup\fonteauteurs\bgroup Duskin\egroup\egroup{} :
\newblock Simplicial methods and the interpretation of ``triple''\ cohomology.
\newblock {\em Mem. Amer. Math. Soc.}, 3(issue 2, 163)\string:\penalty500\relax
  v+135, 1975.

\bibitem{Faltings_finiteness_coco}
Gerd \bgroup\fonteauteurs\bgroup Faltings\egroup\egroup{} :
\newblock Finiteness of coherent cohomology for proper fppf stacks.
\newblock {\em J. Algebraic Geom.}, 12(2)\string:\penalty500\relax 357--366,
  2003.

\bibitem{Ferrand}
Daniel \bgroup\fonteauteurs\bgroup Ferrand\egroup\egroup{} :
\newblock Conducteur, descente et pincement.
\newblock {\em Bull. Soc. Math. France}, 131(4)\string:\penalty500\relax
  553--585, 2003.

\bibitem{Giraud}
Jean \bgroup\fonteauteurs\bgroup Giraud\egroup\egroup{} :
\newblock {\em Cohomologie non ab\'elienne}.
\newblock Springer-Verlag, Berlin, 1971.
\newblock Die Grundlehren der mathematischen Wissenschaften, Band 179.

\bibitem{Godement_faisceaux}
Roger \bgroup\fonteauteurs\bgroup Godement\egroup\egroup{} :
\newblock {\em Topologie alg\'ebrique et th\'eorie des faisceaux}.
\newblock Hermann, Paris, 1973.
\newblock Troisi\`eme \'edition revue et corrig\'ee, Publications de l'Institut
  de Math\'ematique de l'Universit\'e de Strasbourg, XIII, Actualit\'es
  Scientifiques et Industrielles, No. 1252.

\bibitem{Tohoku}
Alexander \bgroup\fonteauteurs\bgroup Grothendieck\egroup\egroup{} :
\newblock Sur quelques points d'alg\`ebre homologique.
\newblock {\em T\^ohoku Math. J. (2)}, 9\string:\penalty500\relax 119--221,
  1957.

\bibitem{EGA2}
Alexander \bgroup\fonteauteurs\bgroup Grothendieck\egroup\egroup{} :
\newblock \'{E}l\'ements de g\'eom\'etrie alg\'ebrique. {II}. \'{E}tude globale
  \'el\'ementaire de quelques classes de morphismes.
\newblock {\em Inst. Hautes \'Etudes Sci. Publ. Math.},
  (8)\string:\penalty500\relax 222, 1961.

\bibitem{FGA}
Alexander \bgroup\fonteauteurs\bgroup Grothendieck\egroup\egroup{} :
\newblock {\em Fondements de la g\'eom\'etrie alg\'ebrique. [{E}xtraits du
  {S}\'eminaire {B}ourbaki, 1957--1962.]}.
\newblock Secr\'etariat math\'ematique, Paris, 1962.

\bibitem{EGA4_2}
Alexander \bgroup\fonteauteurs\bgroup Grothendieck\egroup\egroup{} :
\newblock \'{E}l\'ements de g\'eom\'etrie alg\'ebrique. {IV}. \'{E}tude locale
  des sch\'emas et des morphismes de sch\'emas. {II}.
\newblock {\em Inst. Hautes \'Etudes Sci. Publ. Math.},
  (24)\string:\penalty500\relax 231, 1965.

\bibitem{EGA4_3}
Alexander \bgroup\fonteauteurs\bgroup Grothendieck\egroup\egroup{} :
\newblock \'{E}l\'ements de g\'eom\'etrie alg\'ebrique. {IV}. \'{E}tude locale
  des sch\'emas et des morphismes de sch\'emas. {III}.
\newblock {\em Inst. Hautes \'Etudes Sci. Publ. Math.},
  (28)\string:\penalty500\relax 255, 1966.

\bibitem{EGA4_4}
Alexander \bgroup\fonteauteurs\bgroup Grothendieck\egroup\egroup{} :
\newblock \'{E}l\'ements de g\'eom\'etrie alg\'ebrique. {IV}. \'{E}tude locale
  des sch\'emas et des morphismes de sch\'emas {IV}.
\newblock {\em Inst. Hautes \'Etudes Sci. Publ. Math.},
  (32)\string:\penalty500\relax 361, 1967.

\bibitem{EGA1}
Alexander \bgroup\fonteauteurs\bgroup Grothendieck\egroup\egroup{} :
\newblock {\em \'{E}l\'ements de g\'eom\'etrie alg\'ebrique. {I}. {L}e langage
  des sch\'emas}, volume 166 de {\em Grundlehren Math. Wiss.}
\newblock Springer-Verlag, 1971.

\bibitem{Illusie_CCD}
Luc \bgroup\fonteauteurs\bgroup Illusie\egroup\egroup{} :
\newblock {\em Complexe cotangent et d\'eformations. {I}}.
\newblock Springer-Verlag, Berlin, 1971.
\newblock Lecture Notes in Mathematics, Vol. 239.

\bibitem{poly_Kleiman}
Steven \bgroup\fonteauteurs\bgroup Kleiman\egroup\egroup{} :
\newblock The picard scheme.
\newblock http\!\!://front.math.ucdavis.edu/math.AG /0504020.

\bibitem{Knutson}
Donald \bgroup\fonteauteurs\bgroup Knutson\egroup\egroup{} :
\newblock {\em Algebraic spaces}.
\newblock Springer-Verlag, Berlin, 1971.
\newblock Lecture Notes in Mathematics, Vol. 203.

\bibitem{Lang_Algebra}
Serge \bgroup\fonteauteurs\bgroup Lang\egroup\egroup{} :
\newblock {\em Algebra}, volume 211 de {\em Graduate Texts in Mathematics}.
\newblock Springer-Verlag, New York, third \'edition, 2002.

\bibitem{Laszlo_Olsson_Six_operationsI}
Yves \bgroup\fonteauteurs\bgroup Laszlo\egroup\egroup{} et Martin~C.
  \bgroup\fonteauteurs\bgroup Olsson\egroup\egroup{} :
\newblock The six operations for sheaves on artin stacks i: Finite
  coefficients.
\newblock http\!\!://www.math.polytechnique.fr/$\sim$laszlo/.

\bibitem{LMB}
G{\'e}rard \bgroup\fonteauteurs\bgroup Laumon\egroup\egroup{} et Laurent
  \bgroup\fonteauteurs\bgroup Moret-Bailly\egroup\egroup{} :
\newblock {\em Champs alg\'ebriques}, volume~39 de {\em Ergebnisse der
  Mathematik und ihrer Grenzgebiete. 3. Folge. A Series of Modern Surveys in
  Mathematics [Results in Mathematics and Related Areas. 3rd Series. A Series
  of Modern Surveys in Mathematics]}.
\newblock Springer-Verlag, Berlin, 2000.

\bibitem{Lieblich}
Max \bgroup\fonteauteurs\bgroup Lieblich\egroup\egroup{} :
\newblock Moduli of twisted sheaves.
\newblock \`A paraître dans Duke Math. Journal.
  http\!\!://front.math.ucdavis.edu/math.AG/0411337.

\bibitem{Lieblich_coherent_algebras}
Max \bgroup\fonteauteurs\bgroup Lieblich\egroup\egroup{} :
\newblock Remarks on the stack of coherent algebras.
\newblock {\em Int. Math. Res. Not.}, pages Art. ID 75273, 12, 2006.

\bibitem{Milne_Etale_coh}
James~S. \bgroup\fonteauteurs\bgroup Milne\egroup\egroup{} :
\newblock {\em \'{E}tale cohomology}, volume~33 de {\em Princeton Mathematical
  Series}.
\newblock Princeton University Press, Princeton, N.J., 1980.

\bibitem{Mumford_Picard_groups}
David \bgroup\fonteauteurs\bgroup Mumford\egroup\egroup{} :
\newblock Picard groups of moduli problems.
\newblock \emph{In} {\em Arithmetical Algebraic Geometry (Proc. Conf. Purdue
  Univ., 1963)}, pages 33--81. Harper \& Row, New York, 1965.

\bibitem{Varietes_abeliennes}
David \bgroup\fonteauteurs\bgroup Mumford\egroup\egroup{} :
\newblock {\em Abelian varieties}.
\newblock Tata Institute of Fundamental Research Studies in Mathematics, No. 5.
  Published for the Tata Institute of Fundamental Research, Bombay, 1970.

\bibitem{Olsson_log_twisted_curves}
Martin~C. \bgroup\fonteauteurs\bgroup Olsson\egroup\egroup{} :
\newblock On (log) twisted curves.

\bibitem{Olsson_Sheaves_on_Artin_stacks}
Martin~C. \bgroup\fonteauteurs\bgroup Olsson\egroup\egroup{} :
\newblock Sheaves on artin stacks.
\newblock \`A paraître dans J. Reine Angew. Math. (Journal de Crelle).
  http\!\!://www.ma.utexas.edu/~molsson/qcohrevised.pdf.

\bibitem{Olsson_lemme_chow}
Martin~C. \bgroup\fonteauteurs\bgroup Olsson\egroup\egroup{} :
\newblock On proper coverings of {A}rtin stacks.
\newblock {\em Adv. Math.}, 198(1)\string:\penalty500\relax 93--106, 2005.

\bibitem{Olsson_defm}
Martin~C. \bgroup\fonteauteurs\bgroup Olsson\egroup\egroup{} :
\newblock Deformation theory of representable morphisms of algebraic stacks.
\newblock {\em Math. Z.}, 253(1)\string:\penalty500\relax 25--62, 2006.

\bibitem{Springer_Linear_Algebraic_Groups}
Tonny~A. \bgroup\fonteauteurs\bgroup Springer\egroup\egroup{} :
\newblock {\em Linear algebraic groups}, volume~9 de {\em Progress in
  Mathematics}.
\newblock Birkh\"auser Boston Inc., Boston, MA, second \'edition, 1998.

\end{thebibliography}
\end{document}